\documentclass{amsbook}
\usepackage{nohyperref}
\usepackage{array}
\usepackage{amssymb}
\usepackage{float}
\usepackage{url}
\newtheorem*{question}{Question}
\newtheorem*{theorem}{Theorem}
\newcommand{\myaxiom}[1]{{\sffamily\upshape\uppercase{#1}}}
\newcommand{\manotch}{$\text{\myaxiom{MA}}+\neg\text{\myaxiom{CH}}$}
\newcommand{\visl}{\myaxiom{V}~$=$~\myaxiom{L}}
\newcommand{\bylabel}[1]{%
\def\a{#1}%
\def\b{}%
\ifx\a\b{}\else~(#1)\fi}

\newcommand{\myproblem}[2]{%
\noindent
\textbf{#1.}\bylabel{#2}}

\newcommand{\problembreak}{\medbreak}
\newcommand{\startproblem}{\problembreak}
\newenvironment{myprob}{}{\problembreak}

\newcommand{\mynote}[1]{\emph{#1.}}

\newcommand{\mypreface}{%
\noindent
\emph{Editor's notes.} }

\newfloat{myfoot}{pb}{lof}
\floatstyle{plain}
\newenvironment{myfooter}{%
\hrule
\vspace{6pt}
\raggedright
\small}{\vfill}

\newcounter{myenum}
\newenvironment{myenumerate}%
{%
 \begin{list}{(\arabic{myenum})}{%
 \usecounter{myenum}%
 \setlength{\topsep}{0in}%
 \setlength{\parsep}{0in}%
 \setlength{\partopsep}{0in}%
 }%
}%
{\end{list}}

\newenvironment{myitemize}%
{%
 \begin{list}{$\bullet$}{%
 \setlength{\topsep}{0in}%
 \setlength{\parsep}{0in}%
 \setlength{\partopsep}{0in}%
 }%
}%
{\end{list}}

\renewenvironment{thebibliography}[1]{%
 \expandafter\section\expandafter*\expandafter{\bibname}%
 \normalfont\footnotesize\labelsep .5em\relax
 \renewcommand\theenumiv{\arabic{enumiv}}%
 \list{\@biblabel{\theenumiv}}{\settowidth\labelwidth{\@biblabel{#1}}%
 \leftmargin\labelwidth \advance\leftmargin\labelsep
 \usecounter{enumiv}}%
}{%
 \endlist 
}
\def\@biblabel#1{[#1]}

\begin{document}
\frontmatter
\begin{titlepage}
\mbox{}

\vspace{4in}
\noindent
{\Huge\bfseries
Problems from\\[0.5\baselineskip]
Topology Proceedings}

\vspace{0.5in}
\noindent
{\itshape Edited by Elliott Pearl}

\vfill
\begin{flushright}
{\itshape Topology Atlas, Toronto, 2003}
\end{flushright}
\end{titlepage}

\begin{titlepage}
\noindent
Topology Atlas\\
Toronto, Ontario, Canada\\
http://at.yorku.ca/topology/\\
atlas@at.yorku.ca

\vspace{1in}
\noindent
\textbf{Cataloguing in Publication Data}\\
Problems from topology proceedings / 
edited by Elliott~Pearl.

vi, 216 p.

Includes bibliographical references.

ISBN 0-9730867-1-8

1. Topology---Problems, exercises, etc. 
I. Pearl, Elliott.
II. Title. 

\noindent
Dewey 514 20\\
LC QA611\\
MSC (2000) 54-06

\vspace{1in}
\noindent
Copyright \copyright\ 2003 Topology Atlas. All rights reserved.
Users of this publication are permitted to make fair use of the material 
in teaching, research and reviewing.
No part of this publication may be distributed for commercial purposes
without the prior permission of the publisher.

\vspace{1in}
\noindent
ISBN 0-9730867-1-8\\
Produced November 2003.
Preliminary versions of this publication were distributed on the 
Topology Atlas website. 
This publication is available in several electronic formats on the 
Topology Atlas website.

\vfill
\noindent
Produced in Canada
\end{titlepage}

\setcounter{page}{3}
\chapter*{Contents}

\noindent
Preface
\dotfill
\pageref{tppreface}

\vspace{\baselineskip}
\noindent
Contributed Problems in \emph{Topology Proceedings}
\dotfill
\pageref{tpcontributed}

\emph{Edited by Peter J. Nyikos and Elliott Pearl.}

\vspace{\baselineskip}
\noindent
Classic Problems
\dotfill
\pageref{tpclassic}

\emph{By Peter J. Nyikos.}

\vspace{\baselineskip}
\noindent
New Classic Problems
\dotfill
\pageref{tpnewclassic}

\emph{Contributions by 
Z.T.~Balogh, S.W.~Davis, A.~Dow, G.~Gruenhage, P.J.~Nyikos,\\ 
M.E.~Rudin, F.D.~Tall, S.~Watson.}

\vspace{\baselineskip}
\noindent
Problems from M.E.~Rudin's \emph{Lecture notes in set-theoretic topology}
\dotfill
\pageref{tprudin}

\emph{By Elliott Pearl.}

\vspace{\baselineskip}
\noindent
Problems from A.V.~Arhangel$'$\kern-.1667em ski\u\i's 
\emph{Structure and classification of topological spaces and cardinal 
invariants}
\dotfill
\pageref{tpshura}

\emph{By A.V. Arhangel$'$\kern-.1667em ski\u{\i} and Elliott Pearl.}

\vspace{\baselineskip}
\noindent
A note on P.~Nyikos's \emph{A survey of two problems in topology}
\dotfill
\pageref{tptwonyikos}

\emph{By Elliott Pearl.}

\vspace{\baselineskip}
\noindent
A note on \emph{Open problems in infinite-dimensional topology}
\dotfill
\pageref{tpinfdim}

\emph{By Elliott Pearl.}

\vspace{\baselineskip}
\noindent
Non-uniformly continuous homeomorphisms with uniformly continuous iterates
\dotfill
\pageref{tputz}

\emph{By W.R. Utz.}

\vspace{\baselineskip}
\noindent
Questions on homeomorphism groups of chainable and homogeneous continua
\dotfill
\pageref{tpbrechner}

\emph{By Beverly L. Brechner.}

\vspace{\baselineskip}
\noindent
Some problems in applied knot theory and geometric topology
\dotfill
\pageref{tpsumners}

\emph{Contibutions by 
D.W.~Sumners, J.L.~Bryant, R.C.~Lacher, R.F.~Williams,\\
J.~Vieitez.}

\vspace{\baselineskip}
\noindent
Problems from Chattanooga, 1996
\dotfill
\pageref{tpchattanooga}

\emph{Contributions by
W.W.~Comfort, F.D.~Tall, D.J.~Lutzer, C.~Pan,\\
G.~Gruenhage, S.~Purisch, P.J.~Nyikos.}

\newpage
\vspace{\baselineskip}
\noindent
Problems from Oxford, 2000
\dotfill
\pageref{tpoxford}

\emph{Contributions by 
A.V.~Arhangel$'$\kern-.1667em ski\u\i, S.~Antonyan, K.P.~Hart, L.~Ludwig,\\
M.~Matveev, J.T.~Moore, P.J.~Nyikos, S.A.~Peregudov, R.~Pol, J.T.~Rogers,\\
M.E.~Rudin, K.~Shankar.}

\vspace{\baselineskip}
\noindent
Continuum theory problems
\dotfill
\pageref{tplewis}

\emph{By Wayne Lewis.}

\vspace{\baselineskip}
\noindent
Problems in continuum theory
\dotfill
\pageref{tpprajs}

\emph{By Janusz R. Prajs.}

\vspace{\baselineskip}
\noindent
The plane fixed-point problem
\dotfill
\pageref{tphagopian}

\emph{By Charles L. Hagopian.}

\vspace{\baselineskip}
\noindent
On an old problem of Knaster
\dotfill
\pageref{tpcharatonik1}

\emph{By Janusz J. Charatonik.}

\vspace{\baselineskip}
\noindent
Means on arc-like continua
\dotfill
\pageref{tpcharatonik2}

\emph{By Janusz J.~Charatonik.}

\vspace{\baselineskip}
\noindent
Classification of homogeneous continua
\dotfill
\pageref{tprogers}

\emph{By James T. Rogers, Jr.}

\chapter*{Preface}
\label{tppreface}

I hope that this collection of problems will be an interesting and useful
resource for researchers. 

This volume consists of material from the \emph{Problem Section} of the
journal \emph{Topology Proceedings} originally collected and edited by
Peter Nyikos and subsequently edited by Elliott Pearl for this
publication. This volume also contains some other well-known problems
lists that have appeared in \emph{Topology Proceedings}.

Some warnings and acknowledgments are in order.

I have made some changes to the original source material. The original
wording of the problems is mostly intact. I have rewritten many of the
solutions, originally contributed by Peter Nyikos, in order to give a more
uniform current presentation. I have contributed some new reports of
solutions. I have often taken wording from abstracts of articles and from
reviews (\emph{Mathematical Reviews} and \emph{Zentralblatt MATH}) without
specific attribution. In cases where the person submitting a problem was
not responsible for first asking the problem, I have tried to provide a
reference to the original source of the problem.

Regrettably, I cannot guarantee that all assumptions regarding lower
separation axioms have been reported accurately from the original sources.

Some portions of this volume have been checked by experts for accuracy of
updates and transcription.

I have corrected some typographical errors from the original source
material. I have surely introduced new typographical errors during the
process of typesetting the original documents.

The large bibliography sections were prepared using some of the features
of \emph{MathSciNet} and \emph{Zentralblatt MATH}.

No index has been prepared for this volume. This volume is distributed in
several electronic formats some of which are searchable with viewing
applications.

I thank York University for access to online resources. I thank York
University and the University of Toronto for access to their libraries.

I thank Dmitri Shakhmatov and Stephen Watson for developing Topology
Atlas as a research tool for the community of topologists.

I thank Gary Gruenhage, John C.\ Mayer, Peter Nyikos, Murat Tuncali and 
the editorial board of \emph{Topology Proceedings} for permission to 
reprint this material from \emph{Topology Proceedings} and to proceed 
with this publishing project.

I thank Peter Nyikos for maintaining the problem section for twenty years.

The material from Mary Ellen Rudin's \emph{Lecture notes in set-theoretic
topology} are distributed with the permission of the American Mathematical
Society.

A.V.\ Arhangel$'$\kern-.1667em ski\u{\i} has given his permission to include in this
volume the material from his survey article \emph{Structure and
classification of topological spaces and cardinal invariants}.

The chapter \emph{Problems in continuum theory} consists of material from
the article \emph{Several old and new problems in continuum theory} by
Janusz J.\ Charatonik and Janusz R.\ Prajs and from the website that they
maintain. This material is distributed with the permission of the authors.

The original essay \emph{The plane fixed-point problem} by Charles
Hagopian is distributed with the permission of the author.

The original essays \emph{On an old problem of Knaster} and \emph{Means on
arc-like continua} by Janusz J.\ Charatonik are distributed with the
permission of the author.

The original essay \emph{Expansive diffeomorphisms on $3$-manifolds} by
Jos\'e Vieitez is distributed with the permission of the author.

I thank many people for contributing solutions and checking portions 
(small and large) of this edition:
A.V.~Arhangel$'$\kern-.1667em ski\u\i,
Christoph Bandt,
Paul Bankston,
Carlos Borges,
Raushan Buzyakova,
Dennis Burke,
Max Burke,
Janusz Charatonik,
Chris Ciesielski,
Sheldon Davis,
Alan Dow,
Alexander Dranishnikov,
Todd Eisworth,
Gary Gruenhage, 
Charles Hagopian,
K.P.~Hart,
Oleg Okunev,
Piotr Koszmider,
Paul Latiolais,
Arkady Leiderman,
Ronnie Levy,
Wayne Lewis,
Lew Ludwig,
David Lutzer,
Mikhail Matveev,
Justin Moore,
Grzegorz Plebanek,
Janusz Prajs,
Jim Rogers,
Andrzej Roslanowski,
Mary Ellen Rudin,
Masami Sakai,
John Schommer,
Dmitri Shakhmatov,
Weixiao Shen,
Alex Shibakov,
Petr Simon,
Greg Swiatek,
Paul Szeptycki,
Frank Tall,
Gino Tironi,
Artur Tomita,
Vassilis Tzannes,
Vladimir Uspenskij, 
W.R.~Utz,
Stephen Watson,
Bob Williams,
Scott Williams.

I welcome any corrections or new information on solutions.
Indeed, I hope to use your contributions to prepare a revised edition of 
this volume.

\bigskip
\begin{flushright}
\noindent
\emph{Elliott Pearl\\
November, 2003\\
Toronto, ON, Canada\\
elliott@at.yorku.ca}
\end{flushright}

\mainmatter
\chapter*{Contributed Problems in \emph{Topology Proceedings}}
\label{tpcontributed}
\markboth{\normalsize\textsc{\lowercase{Contributed Problems}}}{}
\begin{myfoot}
\begin{myfooter}
Peter J. Nyikos and Elliott Pearl, 
\emph{Contributed Problems in \emph{Topology Proceedings}},\\
Problems from Topology Proceedings, Topology Atlas, 2003, 
pp.\ 1--68.
\end{myfooter}
\end{myfoot}

\mypreface
This is a collection of problems and solutions that appeared in the 
problem section of the journal \emph{Topology Proceedings}.
The problem section was edited by Peter~J.~Nyikos for twenty years from 
the journal's founding in 1976.
John~C.~Mayer began editing the problem section with volume 21 in 1996.
In this version, the notes and solutions collected throughout the 
twenty-seven year history of the problem section have been updated with 
current information.

\mynote{Conventions and notation}
The person who contributed each problem is mentioned in parentheses after 
the respective problem number. 
This is not necessarily the person who first asked the problem.
Usually there is a reference to a relevant article in 
\emph{Topology Proceedings}.
Sometimes there is a reference to other relevant articles.
There are a few discontinuities in the numbering of the problems.
Some problems have been omitted.

\section*{A. Cardinal invariants}

\begin{myprob}
\myproblem{A1}{K.~Kunen \cite{MR56:8362}}
Does \manotch\ imply that there are no $L$-spaces? 

\mynote{Notes}
Kunen \cite{MR56:8362} showed that \manotch\ implies that there are no 
Luzin spaces (hence there are no Souslin lines either).
A \emph{Luzin space} is an uncountable Hausdorff space in which 
every nowhere dense subset is countable and which has at most countably 
many isolated points.

\mynote{Solution}
U.~Abraham and S.~Todor\v{c}evi\'c \cite{MR86h:03092} showed that the 
existence of an $L$-space is consistent with \manotch.
\end{myprob}

\begin{myprob}
\myproblem{A3}{E. van~Douwen \cite{MR56:1262}}
Is every point-finite open family in a c.c.c.\ space $\sigma$-centered
(i.e., the union of countably many centered families)?

\mynote{Solution}
No (Ortwin F\"orster).
J. Stepr\={a}ns and S. Watson \cite{MR89e:54061} described a subspace of 
the Pixley-Roy space on the irrationals that is a first countable c.c.c.\ 
space which does not have a $\sigma$-linked base.
\end{myprob}

\begin{myprob}
\myproblem{A4}{E. van~Douwen \cite[Problem~391]{MR1078652}}
For which $\kappa > \omega$ is there a compact homogeneous Hausdorff
space $X$ with $c(X) = \kappa$?

\mynote{Notes}
This is known as van~Douwen's problem.
Here $c(X)$ denotes cellularity, i.e., the supremum of all possible
cardinalities of collections of disjoint open sets.
There is an example with $c(X) = 2^{\aleph_0}$.
\end{myprob}

\begin{myprob}
\myproblem{A5}{A.V.~Arhangel$'$\kern-.1667em ski\u{\i}}
Let $c(X)$ denote the cellularity of $X$.
Does there exist a space $X$ such that $c(X^2) > c(X)$?

\mynote{Solution}
Yes (S. Todor\v{c}evi\'c \cite{MR87b:04003}).
\end{myprob}

\begin{myprob}
\myproblem{A6}{A.V.~Arhangel$'$\kern-.1667em ski\u{\i}}
Let $d(X)$ denote the density of $X$ and let $t(X)$ denote the tightness
of $X$, $\omega \cdot \min\{\kappa : 
(\forall A \subset X)
(\forall x \in \operatorname{cl}A)
(\exists B \subset A) {x \in \operatorname{cl}B}, {|B| \leq \kappa}\}$.
Does there exist a compact space $X$ such that $c(X) = t(X) < d(X)$.
Yes, if \myaxiom{CH} or there exists a Souslin line.
\end{myprob}

\begin{myprob}
\myproblem{A7}{T. Przymusi\'nski}
Does there exist for every cardinal $\lambda$ an isometrically universal
metric space of weight $\lambda$?
Yes, if \myaxiom{GCH}.
\end{myprob}

\begin{myprob}
\myproblem{A8}{V. Saks \cite{MR82a:54013}}
A set $C \subset \beta\omega \setminus \omega$ is a \emph{cluster set} 
if there exist $x \in \beta\omega \setminus \omega$ and a sequence 
$\{x_n : n \in \omega\}$ in $\beta\omega$ such that
$C = \{ \mathcal{D} \in \beta\omega \setminus \omega : 
x = \mathcal{D} \setminus \lim x_n, 
\{n : x_n \neq x\} \in \mathcal{D} \}$.
Here a point of $\beta\omega$ is identified with the ultrafilter on
$\omega$ that converges to it.
Is it a theorem of \myaxiom{zfc} that $\beta\omega \setminus \omega$ 
is not the union of fewer than $2^\mathfrak{c}$ cluster sets?

\mynote{Notes}
See especially \cite[Theorem 3.1]{MR82a:54013}.
\end{myprob}

\begin{myprob}
\myproblem{A9}{E. van~Douwen \cite{MR82a:54009}}
If $G$ is an infinite countably compact group, is $|G|^\omega = |G|$?
Yes, if \myaxiom{GCH}.

\mynote{Solution}
No is consistent.
A. Tomita \cite{MR1974663} showed that there is a model of \myaxiom{CH} 
in which there is a countably compact group of cardinality $\aleph_\omega$.
\end{myprob}

\begin{myprob}
\myproblem{A10}{E. van~Douwen \cite{MR83h:54004}}
Is the character, or hereditary Lindel\"of degree, or spread, equal to the
weight for a compact $F$-space?
for a compact basically disconnected space?

\mynote{Notes}
Yes, for compact extremally disconnected spaces by a result of B. Balcar
and F. Fran\v{e}k \cite{MR83m:06020}.
\end{myprob}

\begin{myprob}
\myproblem{A11}{G. Grabner \cite{MR83i:54005}}
Suppose that $X$ is a wrb space.
Does $\chi(X) = t(X)$?

\mynote{Notes}
A space is \emph{wrb} if each point has a local base which is the 
countable union of Noetherian collections of subinfinite rank.
\end{myprob}

\begin{myprob}
\myproblem{A12}{P. Nyikos \cite{MR89i:54034}}
Does there exist, for each cardinal $\kappa$, a first countable, locally
compact, countably compact space of cardinality $\geq \kappa$?

\mynote{Notes}
Yes if $\Box_\kappa$ and $\operatorname{cf} [\kappa]^\omega = \kappa^+$
for all singular cardinals of countable cofinality (P. Nyikos), hence 
yes if the Covering Lemma holds over the Core Model.
A negative answer in some model would thus imply the presence of inner 
models with a proper class of measurable cardinals.
An affirmative answer is compatible with any possible cardinal arithmetic 
(S. Shelah).
\end{myprob}

\begin{myprob}
\myproblem{A13}{E. van~Douwen}
Let $\exp_Y X$ stand for the least cardinal $\kappa$ (if it exists) such
that $X$ can be embedded as a closed subspace in a product of $\kappa$
copies of $Y$.
Does there exist an $N$-compact space $X$ such that
$\exp_\mathbb{N} X \neq \exp_\mathbb{R} X$?

\mynote{Notes}
Such a space cannot be strongly zero-dimensional.
\end{myprob}

\begin{myprob}
\myproblem{A14}{E. van~Douwen}
Is every compact Hausdorff space a continuous image of some 
zero-dimensional compact space of the same cardinality? 
of the same character?
The answer is well-known to be yes for weight.
\end{myprob}

\begin{myprob}
\myproblem{A15}{E. van~Douwen}
Is there for each $\kappa \geq \omega$ a (preferably homogeneous, or even
groupable) hereditarily paracompact (or hereditarily normal) space $X$
with $w(X) = \kappa$ and $|X| = 2^\kappa$?

\mynote{Notes}
Yes to all questions if $2^\kappa = \kappa^+$. 
Also, $w(X) \leq \kappa < |X|$ is always possible.
\end{myprob}

\begin{myprob}
\myproblem{A16}{E. van~Douwen}
Is there for each $\kappa \geq \omega$ a homogeneous compact Hausdorff 
space $X$ with $\chi(X) = \kappa$ and $w(X) = 2^\kappa$?
Or is $\omega$ the only value of $\kappa$ for which this is true?
\end{myprob}

\begin{myprob}
\myproblem{A17}{E. van~Douwen \cite{vanDouwen:Noetherian}}
Is there always a regular space without a Noetherian base?
(Noetherian: no infinite ascending chains.)

\mynote{Notes}
For any ordinal $\alpha$, the space $\alpha$ has a Noetherian base if and
only if $\alpha+1$ does not contain a strongly inaccessible cardinal.
A. Tamariz-Mascar\'ua and R.G.~Wilson \cite{MR89k:54091} showed that 
there is a $T_1$ space without a Noetherian base.
\end{myprob}

\begin{myprob}
\myproblem{A19}{E. van~Douwen}
Is a first countable $T_1$ space normal if every two disjoint closed sets
of size $\leq \mathfrak{c}$ can be put into disjoint open sets?

\mynote{Solution}
If there is no counterexample then there is an inner model with a proper
class of measurable cardinals.
But if the consistency of a supercompact cardinal is assumed, then an
affirmative answer is consistent (I. Juh\'asz).
\end{myprob}

\begin{myprob}
\myproblem{A20}{A.~Garc\'{\i}a-Maynez \cite{MR87i:54004}}
Let $X$ be a $T_3$-space and let $X$ be an infinite cardinal.
Assume the pluming degree of $X$ is $\leq \lambda$.
Is it true that every compact subset of $X$ lies in a compact set which
has a local basis for its neighborhood system consisting of at most
$\lambda$ elements?
\end{myprob}

\begin{myprob}
\myproblem{A21}{B. Shapirovski\u{\i} \cite{MR55:9025}}
Let $A$ be a subset of a space $X$ and let $x \in A^\prime$.
$A^\prime$ denotes the derived set.
Define the accessibility number $a(x,A)$ to be
$\min\{|B| : B \subset A, x \in B^\prime\}$.
Define 
$t_c(x,X)$ to be 
$\sup\{a(x,F) : F\ \text{is closed}, x \in F^\prime\}$.
As usual, define $t(x,X)$ as $\sup\{a(x,A) : x \in A^\prime\}$.
Can we ever have $t_c(x,X) < t(x,X)$ in a compact Hausdorff space?

\mynote{Notes}
No, for c.c.c.\ compact spaces under \myaxiom{GCH} \cite{MR82b:54014}.
\end{myprob}

\begin{myprob}
\myproblem{A22}{D. Shakhmatov \cite{MR94d:54010}}
Assume that $\tau$ is a Tychonoff [resp.\ Hausdorff, regular, $T_1$ etc.]\ 
homogeneous topology on a set $X$.
Are there Tychonoff [resp.\ Hausdorff, regular, $T_1$ etc.]\
homogeneous topologies $\tau_*$ and $\tau^*$ on $X$ such that
$\tau_* \subset \tau \subset \tau^*$,
$w(X, \tau_*) \leq nw(X, \tau)$ and
$w(X, \tau^*) \leq nw(X, \tau)$?

\mynote{Notes}
For background on this problem for the case of topological groups and
other topological algebras, see papers by A.V.~Arhangel$'$\kern-.1667em ski\u{\i} in
\cite{MR80m:54005}
where the ``left half'' is achieved in the category of topological groups
and continuous homeomorphisms.
In (D. Shakhmatov \cite{MR85h:22003}) 
this is extended to many other categories.
In (V.~Pestov and D. Shakhmatov \cite{MR87j:22003}),
the right half is shown to fail in the categories of topological
groups and topological vector spaces, for countable net weight;
in the latter case, $\mathbb{R}^\infty$ provides a counterexample.
\end{myprob}

\section*{B. Generalized metric spaces and metrization}

\begin{myprob}
\myproblem{B1}{T. Przymusi\'nski \cite{MR56:6617}}
Can each normal (or metacompact) Moore space of weight $\leq \mathfrak{c}$
be embedded into a separable Moore space? 

\mynote{Notes}
Under \myaxiom{CH}, the answer is yes even if ``normal'' and 
``metacompact'' are completely dropped (E. van~Douwen and T. Przymusi\'nski). 

\mynote{Solution}
B. Fitzpatrick, J.W.~Ott and G.M.~Reed asked
``Can each Moore space with weight at most $\mathfrak{c}$ be embedded
in a separable Moore space?''
The answer to this question is independent of \myaxiom{zfc} 
(E. van~Douwen and T. Przymusi\'nski \cite{MR82j:54051}).
\end{myprob}

\begin{myprob}
\myproblem{B2}{D. Burke \cite{MR56:6621}}
Is the perfect image of a quasi-developable space also quasi-developable? 

\mynote{Solution}
Yes (D. Burke \cite{MR85d:54009}).
\end{myprob}

\begin{myprob}
\myproblem{B3}{K. Alster and P. Zenor \cite{MR56:9486}}
Is every locally connected and locally rim-compact normal Moore space 
metrizable? 

\mynote{Solution}
Yes (P. Zenor \cite{MR82i:54040}).
\end{myprob}

\begin{myprob}
\myproblem{B4}{D.~Burke and D.~Lutzer \cite{MR55:1318}}
Must a strict $p$-space with a $G_\delta$-diagonal be developable
(equivalently, $\theta$-refinable (=submetacompact))?

\mynote{Notes}
It was erroneously announced in \cite{MR55:1318} that J.~Chaber had given
an affirmative answer; however, Chaber did not claim to settle the question
except in the cases where the space is locally compact or locally second
countable \cite{MR56:9493, MR82j:54046}.

\mynote{Solution}
Yes, because every strict $p$-space is submetacompact
(S.L.~Jiang~\cite{MR89j:54030}). 
\end{myprob}

\begin{myprob}
\myproblem{B5}{H. Wicke \cite{MR80k:54053}}
Is every monotonically semi-stratifiable hereditarily submetacompact
space semi-stratifiable?
\end{myprob}

\begin{myprob}
\myproblem{B6}{H. Wicke \cite{MR80k:54053}}
Is every monotonic $\beta$-space which is hereditarily submetacompact a
$\beta$-space?
\end{myprob}

\begin{myprob}
\myproblem{B7}{H. Wicke \cite{MR80k:54053}}
Does every primitive $q$-space with a $\theta$-diagonal have a primitive
base?

\mynote{Notes}
R.~Ruth \cite{MR86c:54026} proved that a space has a primitive base if 
and only if it is both a $\theta$-space and a primitive $\sigma$-space. 
Also, a primitive $\sigma$-space with a $\theta$-diagonal has a primitive 
diagonal.
\end{myprob}

\begin{myprob}
\myproblem{B8}{C.E.~Aull \cite{MR80m:54044}}
For all base axioms such that countably compact regular $+$ base axiom
$\Rightarrow$ metrizable,
is it true that regular $+$ $\beta$ $+$ collectionwise normal $+$ base
axiom $\Rightarrow$ metrizable?
In particular, what about quasi-developable spaces, or those with
$\delta\theta$-bases or point-countable bases?
\end{myprob}

\begin{myprob}
\myproblem{B9}{C.E.~Aull \cite{MR42:2425, MR80m:54044}}
Is every space in the class MOBI quasi-developable?
\end{myprob}

\begin{myprob}
\myproblem{B10}{C.E.~Aull \cite{MR80m:54044}}
Is every space with a $\sigma$-locally countable base quasi-develop\-able?

\mynote{Notes}
D. Burke \cite[p.~25]{MR81h:54020} showed that a submetacompact 
(=$\theta$-refinable) regular space with a $\sigma$-locally countable 
base is developable. 
Thus Problems B10 and B11 have affirmative answers where submetacompact
regular spaces are concerned.
D. Burke \cite{MR85d:54009} showed that the class of spaces with 
primitive bases is closed under perfect maps. 
J. Kofner \cite{MR81m:54060} showed that the class of quasi-metrizable spaces 
is also closed under perfect maps. 
Also, H.R.~Bennett's example of a paracompact, nonmetrizable space in
MOBI \cite{MR42:2425} shows that the class of spaces with $\sigma$-locally 
countable bases is not preserved under compact open mappings.
\end{myprob}

\begin{myprob}
\myproblem{B11}{C.E.~Aull \cite{MR80m:54044}}
Is every collectionwise normal space with a $\sigma$-locally countable
base metrizable (equivalently, paracompact)?
\end{myprob}

\begin{myprob}
\myproblem{B12}{C.E.~Aull \cite{MR80m:54044}}
Is every first countable space with a weak uniform base (WUB) 
quasi-developable?

\mynote{Notes}
A base $\mathcal{B}$ for a space $X$ is a \emph{(weakly) uniform base} if 
for each $x \in X$ and each infinite subcollection $\mathcal{H}$ of 
$\mathcal{B}$, each member of which contains $X$, $\mathcal{H}$ is a 
local base for $X$ (resp.\ $\bigcap \mathcal{H} = \{x\}$).
A $T_3$ space has a uniform base iff it is a metacompact Moore space 
(P.S.~Alexandroff, R.W.~Heath). 
\end{myprob}

\begin{myprob}
\myproblem{B13}{C.E.~Aull \cite{MR80m:54044}}
Does every developable space with a WUB and without isolated points have a
uniform base? 
Equivalently, is it metacompact?
\end{myprob}

\begin{myprob}
\myproblem{B14}{A.V.~Arhangel$'$\kern-.1667em ski\u{\i}}
Let $X$ be regular, Lindel\"of, and symmetrizable.
Is $X$ separable?
Does $X$ have a $G_\delta$-diagonal?

\mynote{Notes}
It is consistent that the answer to the first is negative, but the 
construction does have a $G_\delta$-diagonal 
(D. Shakhmatov \cite{MR89c:54009}). 
There is a Hausdorff Lindel\"of, and symmetrizable space that is not
separable (Z. Balogh, D. Burke and S. Davis \cite{MR91e:54076}).
\end{myprob}

\begin{myprob}
\myproblem{B15}{H. Junnila \cite{MR80j:54015}}
Is every strict $p$-space submetacompact? 

\mynote{Notes}
This problem is a generalization of Problem B4.

\mynote{Solution}
Yes (S.L.~Jiang \cite{MR89j:54030}).
\end{myprob}

\begin{myprob}
\myproblem{B16}{H. Junnila \cite{MR80j:54015}}
Does there exist, in \myaxiom{zfc}, a set $X$ and two topologies $\tau$ 
and $\pi$ on $X$ such that $\tau \subset \pi$, every $\pi$-open set is an 
$F_\sigma$ set with respect to $\tau$, the space $(X,\pi)$ is metrizable 
but the space $(X,\tau)$ is not a $\sigma$-space?
\end{myprob}

\begin{myprob}
\myproblem{B17}{J.M.~van Wouwe \cite{MR82a:54058}}
Is each GO-space $X$, that is hereditarily a $\Sigma$-space, metrizable?
What if $X$ is compact?

\mynote{Solution}
Yes (Z. Balogh \cite{MR85m:54025}).
\end{myprob}

\begin{myprob}
\myproblem{B18}{D. Burke \cite{MR81h:54020}}
Does every regular space $X$ with a $\sigma$-locally countable base have
a $\sigma$-disjoint base?

\mynote{Notes}
No, if there is a $Q$-set, because then there is a para-Lindel\"of
nonmetrizable normal Moore space (C. Navy \cite{Navy}) and no
nonmetrizable normal Moore space can have a $\sigma$-disjoint base or 
even be screenable.
\end{myprob}

\begin{myprob}
\myproblem{B19}{H.-X. Zhou \cite[M.~Hu\v{s}ek]{MR58:24198}}
A space $X$ is said to have a small diagonal if any uncountable subset of 
$X^2 \setminus \Delta$ has an uncountable subset with closure disjoint
from the diagonal. (This definition is due to M. Hu\v{s}ek.)
Must a compact $T_2$-space with a small diagonal be metrizable?

\mynote{Notes}
Yes (H.-X. Zhou \cite[Theorem~7.5]{MR86j:54008}), if \myaxiom{CH} and the 
character of $\omega_1$ is at most $\omega_1$ in every first countable 
space (as in a model involving inaccessible cardinals due to W. Fleissner 
\cite{MR55:11190}).
Yes if \myaxiom{CH} (I. Juh\'asz and Z. Szentmikl\'ossy \cite{MR93b:54024}).
\end{myprob}

\begin{myprob}
\myproblem{B20}{H. -X. Zhou}
Is a strongly $\omega_1$-compact, locally compact space with a
$G_\delta$-diagonal metrizable?
\end{myprob}

\begin{myprob}
\myproblem{B21}{R.M.~Stephenson \cite{MR81m:54056}}
Is every regular, feebly compact, symmetrizable space first countable
(equivalently, developable)?
\end{myprob}

\begin{myprob}
\myproblem{B22}{P. Nyikos}
Is every weakly $\theta$-refinable space (=weakly submetacompact space) 
with a base of countable order quasi-developable (equivalently, 
by an old theorem of Bennett and Berney \cite{MR49:6177}, 
hereditarily weakly submetacompact)?
\end{myprob}

\begin{myprob}
\myproblem{B23}{P. Nyikos}
Is every collectionwise normal, countably paracompact space with a
$\sigma$-locally countable base metrizable (equivalently, by an old
theorem of Fedor\v{c}uk, paracompact)?
\end{myprob}

\begin{myprob}
\myproblem{B24}{T.J.~Peters \cite{MR87h:54050}}
Do there exist spaces $X$ and $Y$ such that neither $X$ nor $Y$ has a
$\sigma$-discrete $\pi$-base (equivalently, a $\sigma$-locally finite 
$\pi$-base) but $X \times Y$ has one?

\mynote{Solution}
Yes (A. Dow and T.J.~Peters \cite{MR90a:54069}).
\end{myprob}

\begin{myprob}
\myproblem{B25}{P. Nyikos}
A space is paranormal if every countable discrete collection of closed
sets
$\{F_n : n \in \omega\}$ can be expanded to a locally finite 
collection of open sets
$\{G_n : n \in \omega\}$, i.e., $F_n \subset G_n$ and 
$G_n \cap F_m \neq 0$ iff $F_m = F_n$.
Is there a real example of a nonmetrizable paranormal Moore space?
\end{myprob}

\begin{myprob}
\myproblem{B26}{J. Porter and G. Woods \cite{MR86h:54002}}
A space is \emph{$RC$-perfect} if each of its open sets is a union of 
countably many regular closed subsets of the space.
Is there a \myaxiom{ZFC} example of a feebly compact, $RC$-perfect, 
regular space that is not separable? 
A compact $L$-space is a consistent example.

\mynote{Notes}
Does \manotch\ imply that any $RC$-perfect, feebly compact space is compact 
(or separable) (J. Porter and G. Woods \cite{MR86h:54002})?
Is there a \myaxiom{ZFC} example of a feebly compact, $RC$-perfect, 
regular space that is not normal?
\end{myprob}

\begin{myprob}
\myproblem{B27}{K. Tamano \cite{MR87j:54011}}
Find an internal characterization of subspaces of the product of countably
many La\v{s}nev spaces.
\end{myprob}

\begin{myprob}
\myproblem{B28}{K. Tamano \cite{MR87j:54011}}
Does the product of countably many La\v{s}nev spaces have a
$\sigma$-hereditarily closure-preserving $k$-network?

\mynote{Solution}
No.
S. Lin \cite{MR89e:54025} proved that for any La\v{s}nev space $X$ the 
product $X \times \mathbb{I}$ has a $\sigma$-hereditarily closure 
preserving $k$-network if and only if $X$ has a $\sigma$-locally finite 
$k$-network.
There are La\v{s}nev spaces that do not have a $\sigma$-locally finite 
$k$-network.
\end{myprob}

\begin{myprob}
\myproblem{B29}{P. Nyikos}
Is every locally compact, locally connected, countably paracompact Moore
space metrizable?
Yes is consistent.
\end{myprob}

\begin{myprob}
\myproblem{B30}{M.E.~Rudin \cite[The Point-Countable Base Problem]{MR87b:54018}}
A \emph{Collins space} is one in which each point $x$ has a special countable 
open base $W_x$ with the property that, if $U$ is a neighborhood of a 
point $y$, there is a neighborhood $V$ of $y$ such that, for all $x \in V$ 
there is a $W \in W_x$ with $y \in V \subset U$. 
Recall that a Collins space is metrizable precisely if $W_x$ can be made
a nested decreasing sequence for each $x$.
It is easy to see that every space with a point-countable base is a
Collins space.
Is the converse true?

\mynote{Notes}
This problem is in \cite[Problem 378]{MR1078650}.
M.E.~Rudin wrote:
``The conjecture [that the converse is true] has become doubly interesting
to me since I now know that I do not know how to construct a 
counterexample.''
\end{myprob}

\begin{myprob}
\myproblem{B31}{C.R.~Borges \cite{MR89f:54059}}
If $(X, \tau)$ is a topologically complete submetrizable topological
space, is there a complete metric for $X$ whose topology is coarser than
$\tau$?
\end{myprob}

\begin{myprob}
\myproblem{B32}{P. Nyikos}
Is it consistent that every compact space with hereditarily collectionwise
Hausdorff square is metrizable?

\mynote{Notes}
If \manotch, then every compact space with hereditarily strongly 
collectionwise Hausdorff square is metrizable, but this is false 
under \myaxiom{CH}.
\end{myprob}

\begin{myprob}
\myproblem{B33}{P. Nyikos}
Can the consistency of ``all normal Moore spaces of cardinality $\le
\kappa$ are metrizable'' be established without using large cardinals if
$\kappa = \mathfrak{c}$?
$\kappa = 2^\mathfrak{c}$?
$\kappa = \beth_\omega$?
\end{myprob}

\begin{myprob}
\myproblem{B34}{T. Hoshina, communicated by T. Goto \cite{MR93f:54013}}
Can every La\v{s}nev space be embedded in a La\v{s}nev space that is the
closed continuous image of a complete metric space?
\end{myprob}

\begin{myprob}
\myproblem{B35}{A. Okuyama \cite{MR94d:54033}}
Is every Lindel\"of Hausdorff space a weak P($\aleph_0$)-space?
No, if \myaxiom{MA} or $\mathfrak{b} = \omega_1$.

\mynote{Notes}
For a paracompact [resp.\ Lindel\"of] regular space $X$, the product 
$X \times \mathbb{P}$ is paracompact [resp.\ Lindel\"of] 
iff $X$ is a weak P($\aleph_0$)-space.
For further information see
B. Lawrence \cite{MR90m:54014} and 
K. Alster \cite{MR90m:54012,MR93g:54013}.
\end{myprob}

\begin{myprob}
\myproblem{B36}{W. Just and H. Wicke \cite{MR95c:54021}}
Is every bisequential space the continuous image of a metrizable
space under a map with completely metrizable (or even discrete) fibers?
\end{myprob}

\begin{myprob}
\myproblem{B37}{S. Lin \cite{MR1429180}}
Suppose $X$ is a space with a point-countable closed $k$-network.
Does $X$ have a point-countable compact $k$-network if every first
countable closed subspace of $X$ is locally compact?

\mynote{Solution}
No.
M. Sakai \cite{MR2001h:54058} showed that there is a space $X$ satisfying 
the following conditions: 
$X$ has a point-countable closed $k$-network,
every first countable closed subspace of $X$ is compact, and 
$X$ does not have any point-countable compact $k$-network. 
H. Chen \cite{MR2002i:54012} also gave a negative answer.
\end{myprob}

\begin{myprob}
\myproblem{B38}{S. Lin \cite{MR1429180}}
Suppose $X$ is a quotient $s$-image of a metric space.
Does $X$ have a point-countable closed $k$-network if every first
countable closed subspace of $X$ is locally compact?

\mynote{Solution}
H. Chen \cite{chen} showed that a negative answer is consistent.
\end{myprob}

\begin{myprob}
\myproblem{B39}{S. Lin \cite{MR1429180}}
Suppose $X$ has a $\sigma$-closure-preserving compact $k$-network.
Is $X$ a $k$-space if $X$ is a $k_R$-space?
\end{myprob}

\begin{myprob}
\myproblem{B40}{H. Hung \cite{MR98j:54049}}
Is there a metrization theorem in terms of weak, non-uniform factors?

\mynote{Notes}
This paper \cite{MR98j:54049} underlines once again the desirability 
of a non-uniform metrization theorem; Theorem 1.1 being uniform, following 
immediately from \cite[Corollary~2.3]{MR56:13172}, and Theorem 0.2 being 
non-uniform. 
See also \cite{MR96k:54048}.
\end{myprob}

\begin{myprob}
\myproblem{B41}{H. Bennett and D. Lutzer \cite{MR99i:54043}}
Is it consistently true that if $X$ is a Lindel\"of LOTS that is
paracompact off of the diagonal, then $X$ has a $\sigma$-point finite base?
\end{myprob}

\begin{myprob}
\myproblem{B42}{H. Bennett and D. Lutzer \cite{MR99i:54043}}
Can there be a Souslin space (i.e., a nonseparable LOTS with countable
cellularity, no completeness or connectedness assumed) such that
$X^2 \setminus \Delta$ is paracompact? hereditarily paracompact?

\mynote{Solution}
Yes, consistently.
G. Gruenhage showed that if there is a Souslin space, then there is a 
Souslin space $X$ such that $X^2 \setminus \Delta$ is hereditarily 
paracompact. 
The proof appeared in a paper by H. Bennett, D. Lutzer, and
M.E.~Rudin \cite{MR1964446}.
\end{myprob}

\begin{myprob}
\myproblem{B43}{H. Bennett and D. Lutzer \cite{MR99i:54043}}
Suppose $X$ is a LOTS that is first countable and hereditarily paracompact 
off of the diagonal (i.e., $X^2 \setminus \Delta$ is hereditarily
paracompact).
Must $X$ have a point-countable base?
\end{myprob}

\section*{BB. Metric spaces} 

\begin{myprob}
\myproblem{BB1}{Y. Hattori and H. Ohta \cite{MR95k:54052}}
A metric space is said to have UMP (resp.\ WUMP) if for every pair of 
distinct points $x$, $y$ there exists exactly (resp.\ at most) one point
$p$ such that $d(x, p) = d(y, p)$.
Is a separable metric space having UMP homeomorphic to a subspace of the
real line?
\end{myprob}

\begin{myprob}
\myproblem{BB2}{Y. Hattori and H. Ohta \cite{MR95k:54052}} 
Is a rim-compact (i.e., each point has a neighborhood base consisting of
sets with compact boundary) and separable metric space having WUMP
homeomorphic to a subspace of the real line?
\end{myprob}

\section*{C. Compactness and generalizations}

\begin{myprob}
\myproblem{C1}{T. Przymusi\'nski \cite{TP1.321}}
Can each first countable compact space be embedded into a separable first
countable space? 
A separable first countable compact space? 

\mynote{Notes}
Yes to both questions, if \myaxiom{CH} is assumed.
The first answer can be found in the research announcement by
Przymusi\'nski \cite{TP1.321}.
The second answer can be found in R. Walker's book 
\cite[p.~143]{MR52:1595}.
However, the proof of Parovi\v{c}enko's result on which this relies 
\cite[p.~82]{MR52:1595} has a gap in it; but this gap can be filled. 
\end{myprob}

\begin{myprob}
\myproblem{C2}{G. Woods \cite{MR56:13164}}
Is it consistent that there exists a normal countable compact Hausdorff
$F$-space $X$ such that $|C^*(X)| = 2^{\aleph_0}$ and $X$ is not compact? 

\mynote{Solution}
(E. van~Douwen) 
Yes, in fact the assertion is equivalent to $\neg\text{\myaxiom{CH}}$ 
\cite{MR81j:54057}.
There is an absolute example of a countably compact normal basically
disconnected space which is not compact and satisfies 
$|C^*(X)| = \aleph_2 \cdot 2^{\aleph_0}$.
This example may shed some light on D1. 
\end{myprob}

\begin{myprob}
\myproblem{C3}{E. van~Douwen \cite{MR58:30998}}
Is a compact Hausdorff space nonhomogeneous if it can be mapped
continuously onto $\beta\mathbb{N}$?

Yes, if $w(X) \leq \mathfrak{c}$.
This is Problem 247 from \emph{Open Problems in Topology} 
\cite{MR1078643}.
\end{myprob}

\begin{myprob}
\myproblem{C4}{W.W.~Comfort \cite{MR80m:54035}}
Let $\beta\kappa$ denote the Stone-\v{C}ech compactification of the
discrete space of cardinal $\kappa$.
Let $U_\lambda(\kappa) = \{ p \in \beta\kappa : 
(\forall A \in p)\ |A| \geq \lambda \}$, 
let $U(\kappa) = U_\kappa(\kappa)$ and 
$\kappa^* = \beta\kappa \setminus \kappa$. 
Is it a theorem in \myaxiom{ZFC} that if $\lambda \neq \kappa$ then
$U(\lambda) \not\cong U(\kappa)$? 

\mynote{Notes}
The symbol $\cong$ denotes homeomorphism.
This is true if $\operatorname{cf}(\lambda) \neq \operatorname{cf}(\kappa)$.
van~Douwen \cite{MR92d:06032} showed that there is at most one $n \in 
\omega$ for which there is a $\kappa > \omega_n$ with 
$U(\omega_n) \cong U(\kappa)$.
\end{myprob}

\begin{myprob}
\myproblem{C5}{W.W.~Comfort \cite{MR80m:54035}}
With notation as in C4, is it a theorem in \myaxiom{ZFC} that 
$\omega_1^* \not\simeq \omega_0^*$?

\mynote{Notes}
This is an old problem.
See Problem 242 from \emph{Open Problems in Topology} \cite{MR1078643}.
Equivalently, are the Boolean algebras $\mathcal{P}(\omega)/\text{fin}$ 
and $\mathcal{P}(\omega_1)/{[\omega_1]^{<\omega}}$.
It is known \cite{MR80b:54026}
that if $\kappa > \lambda \geq \omega_0$, and 
$\kappa^* \simeq \lambda^*$ then 
$\lambda = \omega_0$ and $\kappa = \omega_1$.
\end{myprob}

\begin{myprob}
\myproblem{C6}{W.W.~Comfort \cite{MR80m:54035}}
More generally, is it a theorem in \myaxiom{ZFC} that if 
$\kappa > \alpha \geq \omega_0$, $\lambda > \beta \geq \omega_0$, and 
$U_\alpha(\kappa) \cong U_\beta(\lambda)$, then 
$\lambda = \kappa$ and $\alpha = \beta$?
\end{myprob}

\begin{myprob}
\myproblem{C7}{W.W.~Comfort \cite{MR80m:54035, MR56:13136}}
It is known that if $\{X_i : i\in I\}$ is a family of Tychonoff 
spaces such that $X_J = \prod_{i \in J} X_i$ is countably compact for 
all $J\subseteq I$ with $|J|\leq 2^\mathfrak{c}$, then 
$X_I = \prod_{i \in I} X_i$ is countably compact. 
See J. Ginsburg and V. Saks \cite{MR52:1633}.
Is $2^\mathfrak{c}$ the optimal test cardinal in this respect?
Is there $\{X_i : i \in I\}$ with $|I| = 2^\mathfrak{c}$, $X_J$ is 
countably compact whenever $J \subsetneq I$, and $X_I$ not countably
compact?
Is there $X$ such that $X^\alpha$ is countably compact iff 
$\alpha < 2^\mathfrak{c}$?
\end{myprob}

\begin{myprob}
\myproblem{C8}{W.W.~Comfort \cite{MR56:13136}, \cite[communicated 
independently by N.~Hindman and S.~Glazer]{MR80m:54035}}
For $p, q \in \beta\mathbb{N}$, define $p + q$ by $A \in p + q$ if 
$\{n : A - n \in p\} \in q$.
Then $p + q \in \beta\mathbb{N}$, and it is known that there exists
$\bar{p} \in \beta\mathbb{N}$ such that $\bar{p} + \bar{p} = \bar{p}$.
Similarly (with $\cdot$ defined analogously) there is 
$\bar{q} \in \beta\mathbb{N}$ such that $\bar{q} \cdot \bar{q} = \bar{q}$. 
Is there $p \in \beta\mathbb{N}$ such that 
$\bar{p} + \bar{p} = \bar{p} \cdot \bar{p} = \bar{p}$?

\mynote{Solution}
No, N. Hindman \cite{MR80f:05005, MR81m:54040} proved there do not exist 
points $p,q\in\beta\mathbb{N}\setminus\mathbb{N}$ such that 
$p + q = p \cdot q$.
\end{myprob}

\begin{myprob}
\myproblem{C9}{D. Cameron \cite{MR80k:54034}}
Under what conditions is $\beta X$ maximal countably compact?
\end{myprob}

\begin{myprob}
\myproblem{C10}{D. Cameron \cite{MR80k:54034}}
Are all compact spaces strongly compact?
\end{myprob}

\begin{myprob}
\myproblem{C11}{D. Cameron \cite{MR80k:54034}}
Are all countably compact spaces strongly countably compact?
\end{myprob}

\begin{myprob}
\myproblem{C12}{D. Cameron \cite{MR80k:54034}}
Are all sequentially compact spaces strongly sequentially compact? 
\end{myprob}

\begin{myprob}
\myproblem{C13}{D. Cameron \cite{MR80k:54034}}
Are there maximal countably compact spaces which are not sequentially 
compact?
\end{myprob}

\begin{myprob}
\myproblem{C14}{D. Cameron \cite{MR80k:54034}}
What are intrinsic necessary and sufficient conditions for a space to
be maximal pseudocompact?
\end{myprob}

\begin{myprob}
\myproblem{C15}{P. Nyikos}
Does there exist a first countable compact $T_1$ space of cardinality 
$> \mathfrak{c}$? 
a compact $T_1$ space with points $G_\delta$ and cardinality 
$> \mathfrak{c}$?
How large can the cardinality be in either case?

\mynote{Solution}
No,
A.A.~Gryzlov \cite{MR81j:54008} proved that for every compact $T_1$-space
$X$, $|X| \leq 2^{\psi(X)}$.
\end{myprob}

\begin{myprob}
\myproblem{C16}{J. Hagler}
Does there exist a compact space $K$ with countable dense subset $D$ such
that every sequence in $D$ has a convergent subsequence, but $K$ is not
sequentially compact?
We may assume without loss of generality that $K$ is a compactification of
$\omega$, i.e., that the points of $D$ are isolated.

\mynote{Notes}
Yes if $\mathfrak{s} = \mathfrak{c}$.
In fact, $\mathfrak{s} = \mathfrak{c}$ implies that $2^\mathfrak{c}$ 
itself is an example of such a $K$ (P. Nyikos).

\mynote{Solution}
Yes, (A. Dow).
\end{myprob}

\begin{myprob}
\myproblem{C17}{P. Nyikos \cite[Problem 356]{MR1078647}}
If a compact space has the property that all countably compact subsets
are compact, is the space sequentially compact?
Yes, if $\mathfrak{c} < 2^\mathfrak{t}$.
\end{myprob}

\begin{myprob}
\myproblem{C18}{P. Nyikos \cite{MR10:315h}, 
Kat\v{e}tov's Problem \cite{MR10:315h}}
Is there a compact nonmetrizable space $X$ such that $X^2$ is 
hereditarily normal?

\mynote{Notes}
Yes, if \manotch\ \cite{TP2.1.359}.
See \cite{MR94b:54009} for a complete proof.
Yes, if there is an uncountable $Q$-set, or assuming \myaxiom{CH}
\cite{MR94b:54009}.

\mynote{Solution}
P. Larson and S. Todor\v{c}evi\'c \cite{MR2003b:54033}
proved that it is consistent that the answer is negative.
\end{myprob}

\begin{myprob}
\myproblem{C19}{E. van~Douwen}
Is a compact space metrizable if its square is:
(1) hereditarily collectionwise normal?
(2) hereditarily collectionwise Hausdorff?

\mynote{Notes}
(1) Yes, if \manotch\ \cite{MR82g:54010} (P. Nyikos).
(2) No, if \myaxiom{CH} (K. Kunen).
\end{myprob}

\begin{myprob}
\myproblem{C20}{E. van~Douwen}
Consider the following statements about an infinite compact space $X$:
\begin{myenumerate}
\item
there are $Y \subset X$ and $y \in Y$ such that 
$\chi(y, Y) \in \{\omega, \omega_1\}$;
\item
there is a decreasing family $\mathcal{F}$ of closed sets with
$|\mathcal{F}| \in \{\omega, \omega_1\}$ and $|\bigcap \mathcal{F}| = 1$.
\end{myenumerate}
Without loss of generality, $X$ is separable, hence \myaxiom{CH} implies (1). 
Clearly (1) implies (2).
What happens under $\neg\text{\myaxiom{CH}}$?

\mynote{Notes}
I. Juh\'asz and Z. Szentmikl\'ossy have shown that if $X$ is of
uncountable tightness, then $X$ has a convergent free
$\omega_1$-sequence, providing a closed $Y$ as the (1) \cite{MR93b:54024}.
Hence \myaxiom{PFA} implies (1), hence (2), by Balogh's theorem that 
\myaxiom{PFA} implies every compact Hausdorff space of countable tightness 
is sequential.
Also, (1) holds in a model obtained by adding uncountably many Cohen
reals in any model of set theory since Juh\'asz showed that every compact 
Hausdorff space of countable tightness has a point of character
$\leq \omega_1$ in that model \cite{MR94k:54004}.
Juh\'asz and Szentmikl\'ossy have also shown that (1) has an affirmative 
solution under $\clubsuit$.
\end{myprob}

\begin{myprob}
\myproblem{C21}{E. van~Douwen}
Is it true that for all infinite cardinals $\kappa$ we have:
$\kappa$ is singular iff
initial $\kappa$-compactness is productive iff
initial $\kappa$-compactness is finitely productive?

\mynote{Solution}
(E. van~Douwen)
Yes if \myaxiom{GCH} but no if 
\myaxiom{MA} + $\mathfrak{c} > \aleph_\omega$ \cite{MR93f:54012}.
Moreover, there is no known model in which initial $\kappa$-compactness
is finitely productive for any cardinals other than singular strong limit
cardinals. 
Compare Problem C37.
\end{myprob}

\begin{myprob}
\myproblem{C22}{E. van~Douwen}
Is initial $\kappa$-compactness productive if $\kappa$ is singular?

\mynote{Notes}
Yes if for all $\mu < \kappa$, $2^\mu < \kappa$ hence yes if \myaxiom{GCH}
(V. Saks and R.M.~Stephenson, \cite{MR42:8448}).

\mynote{Solution}
The statement in the problem is independent of \myaxiom{ZFC}.
Assuming $\text{\myaxiom{MA}} + \mathfrak{c} > \aleph_\omega$, there are 
two initially $\aleph_\omega$-compact normal spaces whose product is not 
initially $\aleph_\omega$-compact 
(E. van~Douwen \cite{MR93f:54012}).
\end{myprob}

\begin{myprob}
\myproblem{C23}{E. van~Douwen}
Does there exist a normal space which is not initially $\kappa$-compact 
but which has a dense initially $\kappa$-compact subspace, for some (each) 
$\kappa > \omega$?
This cannot happen if $\kappa = \omega$ of course.
\end{myprob}

\begin{myprob}
\myproblem{C24}{M. Pouzet \cite{MR83j:04002}}
A space $X$ is called \emph{impartible} if for every partition $\{A,B\}$ 
of $X$, there is a homeomorphism from $X$ into $A$ or into $B$.
Is there a compact impartible space?

\mynote{Notes}
No is consistent (G. Balasubramanian \cite{MR89f:54040}).
\end{myprob}

\begin{myprob}
\myproblem{C24}{V. Saks \cite[attributed to W.W.~Comfort]{MR82a:54013}}
Does there exist a family of spaces $\{X_i : i \in I\}$ with 
$|I| = 2^\mathfrak{c}$, $\prod_{i \in I} X_i$ is not countably compact, 
and $\prod_{i \in J} X_i$ is countably compact, whenever $J \subset I$ 
and $|J| < 2^\mathfrak{c}$?

\mynote{Notes}
This is a special case of Problem C7. 
An affirmative answer to any of A8, P10, or P11 would be sufficient to
construct such a family.

Yes if $2^\mathfrak{c} = \aleph_2$:
The product of $\aleph_1$ sequentially compact spaces is countably compact 
(C.T.~Scarborough and A.H.~Stone \cite{MR34:3528})
and if \myaxiom{CH} then there is a family of $2^\mathfrak{c}$ 
sequentially compact spaces whose product is not countably compact 
(M.~Rajagopalan \cite{MR56:9506}).
The proofs and constructions generalize to models of
$\text{\myaxiom{MA}} + 2^\mathfrak{c} = \mathfrak{c}^+$.
\end{myprob}

\begin{myprob}
\myproblem{C25}{V. Saks \cite{MR82a:54013}}
Do there exist spaces $X$ and $Y$ such that $X^\kappa$ and $Y^\kappa$ are
countably compact for all cardinals $\kappa$, but $X \times Y$ is not
countably compact?
\end{myprob}

\begin{myprob}
\myproblem{C26}{W.W.~Comfort \cite{MR81j:54039}}
Let $\alpha \geq \beta \geq \omega$.
An infinite space $X$ is called pseudo-$(\alpha, \beta)$-compact if for 
every family $\{U_\xi : \xi < \alpha\}$ of nonempty open subsets of 
$X$, there exists $x \in X$ such that
$|\{\xi < \alpha : W \cap U_\xi \neq \emptyset\}| \geq \beta$ 
for every neighborhood $W$ of $X$.
If $\beta$ is singular and $1 < m < \omega$, does there exist a Tychonoff
space $X$ such that $X^{m-1}$ is pseudo-$(\beta, \beta)$-compact and $X^m$
is not pseudo-$(\alpha,\omega)$-compact?
\end{myprob}

\begin{myprob}
\myproblem{C27}{W.W.~Comfort \cite{MR81j:54039}}
Let $\alpha > \beta \geq \omega$ with 
$\operatorname{cf}(\alpha) = \omega$. 
Is there a Tychonoff space $X$ such that $X^m$ is 
pseudo-$(\alpha, \beta)$-compact for all $m < \omega$ and $X^\omega$ is
not pseudo-$(\alpha, \alpha)$-compact?
\end{myprob}

\begin{myprob}
\myproblem{C28}{P. Nyikos}
Does there exist a separable, first countable, countably compact, $T_2$
(hence regular) space which is not compact?

\mynote{Notes}
Yes, if $\mathfrak{b} = \mathfrak{c}$ and other models of set theory.
See the series of articles \emph{On first countable, countably compact 
spaces} by P. Nyikos \cite{MR82k:54028,MR85c:54009,MR91j:54066,MR1078644}.
\end{myprob}

\begin{myprob}
\myproblem{C29}{P. Nyikos}
Does there exist a first countable, countably compact, noncompact regular
space which does not contain a copy of $\omega_1$?

\mynote{Notes}
Yes, if $\clubsuit$; also yes in any model which is obtained from a model
of $\clubsuit$ by iterated c.c.c.\ forcing, so that yes is compatible
with \manotch.

\mynote{Solution}
No is also consistent.
It follows from \myaxiom{PFA} that no such space exists (Z. Balogh) and a 
negative answer is also equiconsistent with \myaxiom{ZFC} (A. Dow).
\end{myprob}

\begin{myprob}
\myproblem{C30}{S. Watson}
Is there a pseudocompact, meta-Lindel\"of space which is not compact?

\mynote{Notes}
Yes, if \myaxiom{CH} (B. Scott \cite{MR81m:54034}).
\end{myprob}

\begin{myprob}
\myproblem{C31}{S. Watson}
Is there a pseudocompact, para-Lindel\"of space which is not compact?

\mynote{Solution}
No, (D.~Burke and S.~Davis \cite{BurkeDavis},
\cite[Theorem~9.7]{MR86e:54030}).
\end{myprob}

\begin{myprob}
\myproblem{C32}{P. Nyikos}
Is every separable, first countable, normal, countably compact space compact?

\mynote{Notes}
No, if $\mathfrak{p} = \omega_1$.
No, if $\mathfrak{p} = \omega_1$
(S.P.~Franklin and M. Rajagopalan \cite{MR44:972}).

\mynote{Solution}
Yes if \myaxiom{PFA} (D. Fremlin) and an affirmative answer is 
equiconsistent with \myaxiom{ZFC} (A. Dow).
But also, a negative answer is consistent with \manotch\ and with 
$\text{\myaxiom{PFA}}^-$ (P. Nyikos).
\end{myprob}

\begin{myprob}
\myproblem{C33}{J. Vaughan}
Is there a separable, first countable, countably compact, non-normal space?

\mynote{Notes}
Yes if $\mathfrak{p} = \omega_1$ or $\mathfrak{b} = \mathfrak{c}$,
hence yes if $\mathfrak{c} \leq \omega_2$
\end{myprob}

\begin{myprob}
\myproblem{C34}{T. Przymusi\'nski}
A space is sequentially separable if it has a countable subset $D$ such
that every point is the limit of a sequence from $D$.
Can every first countable compact space be embedded in a sequentially
separable space?
Yes, if \myaxiom{CH}.
\end{myprob}

\begin{myprob}
\myproblem{C35}{P. Nyikos}
Is \myaxiom{CH} alone enough to imply the existence of a locally compact, 
countably compact, hereditarily separable space which is not compact?
a perfectly normal, countably compact space which is not compact? 

\mynote{Notes}
Under ``\myaxiom{CH} + there exists a Souslin tree'' there is a single 
example with all these properties, and various non-Lindel\"of spaces have 
been constructed under \myaxiom{CH} that are countably compact and 
hereditarily separable, or perfectly normal, locally compact and 
hereditarily separable.
\end{myprob}

\begin{myprob}
\myproblem{C36}{E. van~Douwen}
Does there exist in \myaxiom{ZFC} a separable normal countably compact 
noncompact space?
Examples exist if \myaxiom{MA} or if $\mathfrak{p} = \omega_1$.

\mynote{Solution}
Yes, (S.P.~Franklin and M.~Rajagopalan \cite[Ex.~1.5]{MR44:972}).
Their example is also locally compact and scattered, hence sequentially
compact.
van~Douwen probably wanted a first countable example.
\end{myprob}

\begin{myprob}
\myproblem{C37}{P. Nyikos}
Is initial $\kappa$-compactness productive if and only if $\kappa$ is a
singular strong limit cardinal?

\mynote{Notes}
For `if', the answer is affirmative in \myaxiom{ZFC} 
(V. Saks and R.M.~Stephenson \cite{MR42:8448}). 
See also \cite{MR54:13848}.
For `only if', there is an affirmative answer under \myaxiom{GCH} 
(E.~van~Douwen) and in numerous other models of set theory 
(E.~van~Douwen, P.~Nyikos).
\end{myprob}

\begin{myprob}
\myproblem{C38}{E. van~Douwen}
Is there a (preferably separable locally compact) first countable
pseudocompact space that is $\aleph_1$-compact (i.e., has no uncountable
closed discrete subset) but is not countably compact?

\mynote{Notes}
Yes if $\mathfrak{b} = \omega_1$, see the example by
P. Nyikos described in \cite[Notes to \S~13]{MR87f:54008} 
or $\mathfrak{b} = \mathfrak{c}$ (E. van~Douwen).
\end{myprob}

\begin{myprob}
\myproblem{C39}{E. van~Douwen \cite{MR87f:54008}}
Let $\mu$ be the least cardinality of a compact space that is not
sequentially compact.
It is known that $2^\mathfrak{t} \leq \mu \leq 2^\mathfrak{s}$.
What else can be said about $\mu$?

\mynote{Notes}
Here $\mathfrak{t}$ denotes the least cardinality of a tower: a chain of
subsets of $\omega$ with respect to almost-containment
($A \subset^* B$ iff $A \setminus B$ is finite) such that no infinite
subset of $\omega$ is almost contained in every one. 
$\mathfrak{s}$ is the least cardinality of a splitting family 
$\mathcal{S}$ of subsets of $\omega$: a family such that for each infinite 
$A \subset \omega$, there exists $S \in \mathcal{S}$ such that $A \cap S$
and $A \setminus S$ are both infinite.

Let $\mathfrak{h}$ denote the least height of a tree $\pi$-base for 
$\omega^*$.
Then $\mathfrak{h} \leq \mathfrak{s}$, and there is a family of
$\mathfrak{h}$ compact sequential spaces of cardinality 
$\leq \mathfrak{c}$ whose product is not sequentially compact.
Thus $\mu \leq 2^\mathfrak{h}$.
Also, $\mathfrak{h}$ is equal to the least cardinality of a family of 
sequentially compact spaces whose product is not sequentially compact, 
as well as the least cardinality of a family of nowhere dense subsets of 
$\omega^*$ whose union is dense (i.e., the weak Nov\'ak number).
For additional information on $\mathfrak{h}$, see the paper by B. Balcar,
J. Pelant and P. Simon \cite{MR82c:54003},
where it is denoted by $\kappa(\mathbb{N}^*)$, 
and Peter Dordal's thesis \cite{Dordal}, 
where it is denoted $d$, and where it is shown that 
$\mathfrak{t} < \mathfrak{h}$ is consistent. 
S. Shelah's model of $\mathfrak{b} < \mathfrak{s}$ \cite{MR86b:03064} has 
$\mathfrak{h} < \mathfrak{s}$ because $\mathfrak{h} \leq \mathfrak{b}$.

A further improvement is that $\mu \leq \beta$, where 
$\beta$ $=$ $\min\{$ $|B|$ : $B$ is the set of branches in some tree
$\pi$-base for $\omega^*\}$.
It is easy to see that $\beta \leq 2^\mathfrak{h}$.
Moreover, it is consistent to have $\beta < 2^\mathfrak{h}$
(P. Nyikos and S. Shelah).

It is possible to have $\mu = \mathfrak{s} = \mathfrak{c}$, hence
$\mu < 2^\mathfrak{s}$ (S. Shelah \cite{MR86b:03064}).

$\mathfrak{n} \leq \mu$ and there is a model where 
$2^\mathfrak{t} < \mathfrak{n}$ (A. Dow \cite{MR91b:54054}).
$\mathfrak{n}$ is the Nov\'ak number of $\omega^*$, i.e., the minimum 
cardinality of a family of nowhere dense sets covering $\omega^*$.
\end{myprob}

\begin{myprob}
\myproblem{C40}{P. Nyikos}
Is there a first countable, $H$-closed space of cardinality $\aleph_1$?
Equivalently: is there a compact Hausdorff space that can be partitioned
in $\aleph_1$ nonempty zero-sets?
Yes, if \myaxiom{CH}.

\mynote{Solution}
No if under \manotch.
G. Gruenhage \cite{MR94a:54015} showed that if the real line is not the 
union of $\kappa$ many nowhere dense sets, then no compact Hausdorff space 
can be partitioned into $\kappa$ many disjoint $G_\delta$ sets 
(equivalently, zero sets).
\end{myprob}

\begin{myprob}
\myproblem{C41}{E. van~Douwen}
Is there a regular (noncompact, countably compact) space which is 
homeomorphic to each of its closed noncompact subspaces, and is not 
orderable?

\mynote{Notes}
The orderable such spaces are regular cardinals.
\end{myprob}

\begin{myprob}
\myproblem{C42}{T.J.~Peters \cite{MR86c:54024}}
Is the class of $G$-spaces finitely productive?
\end{myprob}

\begin{myprob}
\myproblem{C43}{T.J.~Peters \cite{MR86c:54024}}
Determine conditions on an infinite family of $G$-spaces which will ensure
that their product is $G$.
Specifically, if every countable partial product of some family
$\{X_\xi : \xi < \alpha\}$ of spaces is also a $G$-space, then must 
their full product be one also.
\end{myprob}

\begin{myprob}
\myproblem{C44}{T.J.~Peters \cite{MR86c:54024}}
Do there exist non-$G$-spaces $X$ and $Y$ such that $X \times Y$ is a
$G$-space?
\end{myprob}

\begin{myprob}
\myproblem{C45}{E. van~Douwen}
Is there a compact Fr\'echet-Urysohn space with a pseudocompact noncompact
subspace?
Yes, if $\mathfrak{b} = \mathfrak{c}$.

\mynote{Solution}
Yes, there is even a Talagrand compact space $X$ with a point $p$ such
that $X = \beta(X \setminus \{p\})$ 
(E. Reznichenko).
\end{myprob}

\begin{myprob}
\myproblem{C46}{E. van~Douwen}
Suppose every pseudocompact subspace of a compact space $X$ is compact. 
Must $X$ be hereditarily realcompact?
No if $\clubsuit$.

\mynote{Solution}
No (P. Nyikos \cite{MR81j:58012}).
The subspace $T^+$ of the tangent bundle on the long line is a Moore 
manifold in which every separable subspace is metrizable and so every 
pseudocompact subspace is compact, yet it is not realcompact.
Its one-point compactification is the counterexample.
\end{myprob}

\begin{myprob}
\myproblem{C47}{E. van~Douwen}
Is there a regular Baire space $X$ which has a $1$-$1$ regular continuous 
image $Y$ of smaller weight but no such image that is Baire?
\end{myprob}

\begin{myprob}
\myproblem{C48}{P. Nyikos}
Is there a compact non-scattered space that is the union of a chain of
compact scattered subspaces?

\mynote{Solution}
No, I. Juh\'asz and E. van~Douwen have pointed out that a compact
nonscattered space $X$ has a separable nonscattered subspace, because $X$
admits a continuous map onto $[0,1]$ and any closed subspace $Y$ to which
the restriction is irreducible must be separable.
\end{myprob}

\begin{myprob}
\myproblem{C49}{J. Porter \cite{MR87f:54034}}
Can each Hausdorff space be embedded in some CFC space?

\mynote{Notes}
A space $X$ is \emph{compactly functionally compact} (CFC) if continuous 
function $f\colon X \to Y$ with compact fibers is a closed function.
\end{myprob}

\begin{myprob}
\myproblem{C50}{J. Porter \cite{MR87f:54034}}
Is the product of CFC spaces a CFC space?
\end{myprob}

\begin{myprob}
\myproblem{C51}{V. Malykhin}
Recall that a space is \emph{weakly first countable} if to each point 
$x$ one can assign a countable filterbase $F_x$ of sets containing $x$ 
such that a set $U$ is open iff for each $x \in U$ there is $P \in P_x$ 
such that $P \subset U$ \cite{MR37:3534}.
Is there a weakly first countable compact space which is not first
countable?
One that is of cardinality $> \mathfrak{c}$?
Yes, if \myaxiom{CH} \cite{MR84a:54008}.

\mynote{Notes}
Yes to the first question if $\mathfrak{b} = \mathfrak{c}$ (H.-X. Zhou).
If $\aleph_1$ dominating reals are iteratively added and every countable 
subset of $\omega$ appears at some initial stage, then arbitrarily
large weakly first countable compact Hausdorff spaces exist (P.~Nyikos).
\end{myprob}

\begin{myprob}
\myproblem{C52}{B. Shapirovski\u{\i}}
Is it true that every infinite compact Hausdorff space contains either
$\beta\omega$, or a point with countable $\pi$-character, or a nontrivial
convergent sequence?
\end{myprob}

\begin{myprob}
\myproblem{C53}{V. Uspenskij}
Is every Eberlein compact space of nonmeasurable cardinal bisequential?

\mynote{Solution}
No (P. Nyikos).
The result does hold, however, for uniform Eberlein compacta.
\end{myprob}

\begin{myprob}
\myproblem{C54}{P. Nyikos}
A space is called \emph{$\alpha$-realcompact} if every maximal family of 
closed sets with the c.i.p.\ has nonempty intersection.
Is there a compact sequential space of nonmeasurable cardinal that is not
hereditarily $\alpha$-realcompact?
Yes, if $\clubsuit$.

\mynote{Solution}
Yes, (A. Dow \cite{MR92f:54027}).
\end{myprob}

\begin{myprob}
\myproblem{C55}{P. Nyikos}
Is $2^\mathfrak{s}$ always the smallest cardinality of an infinite
compact Hausdorff space with no nontrivial convergent sequences?

\mynote{Notes}
Here $\mathfrak{s}$ denotes the \emph{splitting number}, which can be 
characterized as the least cardinal $\kappa$ such that $2^\kappa$ is not
sequentially compact.
Fedor\v{c}uk showed, in effect, that if $\mathfrak{s} = \aleph_1$ then
there is a compact Hausdorff space of cardinality $2^\mathfrak{s}$ with
no nontrivial convergent sequences.
\end{myprob}

\begin{myprob}
\myproblem{C56}{P. Nyikos}
Is it consistent that every separable, hereditarily normal,
countably compact space is compact?

\mynote{Solution}
Yes (P. Nyikos, B. Shapirovski\u{\i}, Z. Szentmikl\'ossy, and
B. Veli\v{c}kovi\'c \cite{MR1255822}).
\end{myprob}

\begin{myprob}
\myproblem{C57}{P. Nyikos}
Is there an internal characterization of Rosenthal compacta?

\mynote{Notes}
A \emph{Rosenthal compact space} if it homeomorphic to a compact subset, 
in the topology of pointwise convergence, of the set of Baire class $1$ 
functions on a Polish space. 
See \cite{MR82f:54030,MR84g:46019}.
\end{myprob}

\begin{myprob}
\myproblem{C58}{P. Nyikos}
Is it consistent that every separable, hereditarily normal, 
countably compact space is compact?

\mynote{Notes}
Yes, this is C56.
\end{myprob}

\begin{myprob}
\myproblem{C59}{P. Nyikos}
Is it consistent that every hereditarily normal, countably compact space
is either compact or contains a copy of $\omega_1$?
\end{myprob}

\begin{myprob}
\myproblem{C60}{L. Friedler, M. Girou, D. Pettey, and J. Porter 
\cite{MR94m:54058}}
A regular $T_1$ [resp.\ Urysohn] space $X$ is \emph{$R$-closed}
[resp.\ \emph{$U$-closed}] if $X$ is a closed subspace of every regular
$T_1$ [resp.\ Urysohn] space containing $X$ as a subspace.
Is a space in which each closed set is $R$-closed [resp.\ $U$-closed]
necessarily compact?
\end{myprob}

\begin{myprob}
\myproblem{C61}{L. Friedler, M. Girou, D. Pettey, and J. Porter
\cite{MR94m:54058}}
A regular $T_1$ space is \emph{$RC$-regular} if it can be embedded in an 
$R$-closed space.
Find an internal characterization of $RC$-regular spaces.
\end{myprob}

\begin{myprob}
\myproblem{C62}{L. Friedler, M. Girou, D. Pettey, and J. Porter
\cite{MR94m:54058}}
Is the product of two $R$-closed spaces necessarily $RC$-regular?
\end{myprob}

\begin{myprob}
\myproblem{C63}{L. Friedler, M. Girou, D. Pettey, and J. Porter
\cite{MR94m:54058}}
Is there only one minimal regular topology coarser than an $R$-closed
topology that has a proper regular subtopology?
\end{myprob}

\begin{myprob}
\myproblem{C64}{L. Friedler, M. Girou, D. Pettey, and J. Porter
\cite{MR94m:54058}}
If the product of spaces $X$ and $Y$ is strongly minimal regular 
[resp.\ $RC$-regular] then must each of $X$ and $Y$ be strongly minimal
regular [resp.\ $RC$-regular]?
\end{myprob}

\begin{myprob}
\myproblem{C65}{V. Tzannes \cite{MR98j:54036}}
Does there exist a regular (first countable, separable) countably compact
space on which every continuous real-valued function is constant?
\end{myprob}

\begin{myprob}
\myproblem{C66}{V. Tzannes \cite{MR98j:54036}}
Does there exist, for every Hausdorff space $R$, a regular (first 
countable, separable) countably compact space on which every continuous 
function into $R$ is constant?
\end{myprob}

\begin{myprob}
\myproblem{C67}{V. Tzannes \cite{MR98j:54036}}
Characterize the Hausdorff (regular, normal) spaces which can be 
represented as closed subspaces of Hausdorff (regular, normal) 
star-Lindel\"of spaces.
\end{myprob}

\begin{myprob}
\myproblem{C68}{V. Tzannes \cite{MR98j:54036}}
How big can be the extent of a Hausdorff (regular, normal) star-Lindel\"of 
space?

\mynote{Notes} 
We say that a space is \emph{star-Lindel\"of} if for every open cover
$\mathcal{U}$ of $X$ there exists a countable subset $F \subset X$ such
that $\operatorname{St}^1(F, \mathcal{U}) = X$.
Star-Lindel\"ofness is a joint generalization of Lindel\"ofness, countable
compactness and separability.
Partial answers to C67 and C68 were obtained by M. Bonanzinga in 
\cite{MR99e:54015}.
\end{myprob}

\begin{myprob}
\myproblem{C69}{L. Feng and S. Garcia-Ferreira \cite{MR2002a:54002}}
What kind of spaces can be extended to maximal Tychonoff $MI$ spaces? 

\mynote{Notes}
An \emph{$MI$ space} (E. Hewitt \cite{MR5:46e}) is a crowded space in 
which every dense subset is open. 
If a space is $MI$ then every Tychonoff crowded extension of it is $MI$.
A Tychonoff space is Hausdorff maximal iff it is a maximal Tychonoff $MI$ 
space.
\end{myprob}

\section*{D. Paracompactness and generalizations}

\begin{myprob}
\myproblem{D1}{G. Woods \cite{MR56:13164}}
Is there a real (i.e., not using any set-theoretic hypotheses other than
\myaxiom{ZFC}) example of an extremally disconnected locally compact 
normal nonparacompact Hausdorff space? 
\end{myprob}

\begin{myprob}
\myproblem{D2}{J.C.~Smith \cite{MR56:13160}}
Let $X$ be a regular $q$-space. 
If $X$ is $\aleph$-preparacompact and weakly $\theta$-refinable 
(=weakly submetacompact), then is $X$ paracompact? 

\mynote{Notes}
A $T_2$ space $X$ is said to be \emph{preparacompact}
(\emph{$\aleph$-preparacompact}) if each open cover of $X$ has an open 
refinement $\mathcal{H} = \{H_\alpha : \alpha\in A\}$ such that, if 
$B$ is any infinite (uncountable) subset of $A$ and if $p_\beta$ and 
$q_\beta \in H_\beta$ for each $\beta \in B$ with $p_\alpha \neq p_\beta$ 
and $q_\alpha \neq q_\beta$ for $\alpha \neq \beta$, then the set
$Q = \{q_\beta : \beta\in B\}$ has a limit point whenever the set 
$P = \{p_\beta : \beta\in B\}$ has a limit point.
A space $X$ is called a \emph{$q$-space} if each point $x \in X$ has a 
sequence of neighborhoods $\{N_i\}_{i \in \omega}$ such that, if 
$y_i \in N_i$ for each $i$ with $y_i \neq y_j$ for $i \neq j$, then the 
set $\{y_i\}_{i \in \omega}$ has a limit point.
If $X$ is a regular $q$-space then the following
statements are equivalent \cite{MR52:6662}:
$X$ is paracompact;
$X$ is $\aleph$-preparacompact and subparacompact; 
$X$ is $\aleph$-preparacompact and metacompact. 
\end{myprob}

\begin{myprob}
\myproblem{D3}{J.C.~Smith \cite{MR56:13160}}
Are $\aleph$-preparacompact or preparacompact spaces countably paracompact 
or expandable? 

\mynote{Notes}
A space is \emph{expandable} if every locally finite collection of 
subsets can be expanded to locally finite collection of open sets.
\end{myprob}

\begin{myprob}
\myproblem{D4}{J.C.~Smith \cite{MR56:13160}}
What class of spaces, weaker than irreducible spaces, imply
paracompactness in the presence of $\aleph$-paracompactness?

\mynote{Notes}
$X$ is called \emph{irreducible} if every open cover has an irreducible 
refinement (a cover is irreducible if no proper subfamily is a cover).
\end{myprob}

\begin{myprob}
\myproblem{D5}{P. Bankston \cite{MR56:13141}}
Can an ultrapower of a (paracompact) space be normal without being 
paracompact? 
\end{myprob}

\begin{myprob}
\myproblem{D6}{K. Alster and P. Zenor \cite{MR56:9486}}
Is every perfectly normal manifold collectionwise normal? 

\mynote{Notes}
If \manotch, perfectly normal manifolds are metrizable 
(M.E.~Rudin \cite{MR80j:54014}).
\end{myprob}

\begin{myprob}
\myproblem{D7}{K. Alster and P. Zenor \cite{MR56:9486}}
Is every locally compact and locally connected normal $T_2$-space
collectionwise normal with respect to compact sets?
\end{myprob}

\begin{myprob}
\myproblem{D8}{E. van~Douwen}
Is there a paracompact (metacompact or subparacompact or hereditarily
Lindel\"of) space that is not a $D$-space?

\mynote{Notes}
$X$ is a \emph{$D$-space} if for every $V\colon X \to \tau X$ with 
$x \in V(x)$ for all $x$, there is a closed discrete $D \subset X$ such 
that $\bigcup\{V(x) : x \in D\} = X$.
A generalized ordered space is paracompact iff it is a $D$-space.
(E. van~Douwen and D. Lutzer \cite{MR97f:54039}).
\end{myprob}

\begin{myprob}
\myproblem{D9}{P. Nyikos}
Is the finite product of metacompact $\sigma$-scattered spaces likewise
metacompact? 
What if (weakly) submetacompact, or screenable, or
$\sigma$-metacompact, or meta-Lindel\"of is substituted for metacompact?
\end{myprob}

\begin{myprob}
\myproblem{D10}{P. Nyikos}
Is the product of a metacompact space and a metacompact scattered space
likewise metacompact?
(What about the other covering properties mentioned in D9?)

\mynote{Notes}
A space is called \emph{$C$-scattered} if each closed subspace has a 
point with a neighborhood in the relative topology which is locally 
compact.
A subspace $A$ of a space $X$ is \emph{metacompact relative to $X$} if 
for each open (in $X$) cover of $A$ there is a point-finite (in $X$) open 
(in $X$) refinement which covers $A$.

\mynote{Solution}
Yes for regular spaces (H. Hdeib):
If $A$ is a closed $C$-scattered subset of a regular metacompact space 
$X$, then $A \times Y$ is metacompact relative to $X \times Y$ for any 
regular metacompact space $Y$.
As a corollary, the product of a regular metacompact $C$-scattered 
(in particular, scattered) space with a regular metacompact space is 
metacompact.
\end{myprob}

\begin{myprob}
\myproblem{D11}{P. Nyikos}
Is the finite product of hereditarily (weakly) $\delta\theta$-refinable
$\sigma$-scattered spaces likewise hereditarily (weakly) 
$\delta\theta$-refinable?
What about (weakly) $\delta\theta$-refinable spaces?
\end{myprob}

\begin{myprob}
\myproblem{D12}{S. Williams \cite{MR81f:54002}}
Is $\Box^\omega (\omega+1)$ always paracompact or normal?
\end{myprob}

\begin{myprob}
\myproblem{D13}{S. Williams \cite{MR81f:54002}}
Is $\Box^{\omega_1} (\omega+1)$ normal in any model of \myaxiom{ZFC}?

\mynote{Solution}
No.
L.B.~Lawrence \cite{MR96f:54027} proved that 
$\Box^{\omega_1} (\omega+1)$ is neither normal nor collectionwise Hausdorff.
\end{myprob}

\begin{myprob}
\myproblem{D14}{S. Williams \cite{MR81f:54002}}
Can there be a normal nonparacompact box product of compact spaces?
\end{myprob}

\begin{myprob}
\myproblem{D15}{S. Williams \cite{MR81f:54002}}
Is the box product of countably many compact linearly ordered topological
spaces paracompact?
\end{myprob}

\begin{myprob}
\myproblem{D16}{S. Williams \cite{MR81f:54002}}
For directed sets $D$ and $E$, define $D \leq E$ if there exists a 
function $T\colon D \to E$ preserving bounded sets; allow $D \equiv E$ if 
$D \leq E$ and $E \leq D$.
For which directed sets $D$ does $D \equiv {}^\omega \omega$ imply
$\Box^\omega (\omega+1)$ is paracompact?
Does $\omega \times \omega_2 \equiv {}^\omega \omega$ imply
$\Box^\omega (\omega + 1)$ is paracompact?

\mynote{Notes}
If $\kappa \leq \mathfrak{c}$ is an ordinal of uncountable cofinality,
then each of $\kappa \equiv {}^\omega \omega$ and 
$\kappa \times \mathfrak{c} \equiv {}^\omega \omega$
imply $\Box^\omega (\omega+1)$ is paracompact.
\end{myprob}

\begin{myprob}
\myproblem{D17}{C.E.~Aull \cite{MR80m:54044}}
Is every collectionwise normal space with an orthobase paracompact?
Is it consistent that every normal space with an orthobase is paracompact?

\mynote{Notes}
A base $\mathcal{B}$ for a topological space $X$ is an orthobase if for 
each $\mathcal{B}' \subseteq \mathcal{B}$, either 
$\bigcap\mathcal{B}'$ is an open set of $X$, or 
$\bigcap\mathcal{B}'=\{p\}$ and $\mathcal{B}'$ is a local base at $p$. 
G.~Gruenhage \cite{MR80k:54056} proved that monotonically normal spaces 
with an orthobase are paracompact.
\end{myprob}

\begin{myprob}
\myproblem{D18}{H.~Junnila \cite{MR80j:54015}}
Is a space submetacompact if every directed open cover has a 
$\sigma$-cushioned refinement? 

\mynote{Notes}
See the surveys by H.~Junnila \cite{MR81e:54019} and
S.~Jiang \cite{MR2001j:54023}.
\end{myprob}

\begin{myprob}
\myproblem{D19}{C.E.~Aull \cite{MR80m:54044}}
For Tychonoff spaces, does pseudocompact plus metacompact equal compact?
In a pseudocompact Tychonoff space, does every point-finite collection
$\mathcal{U}$ of open sets have a finite subcollection $\mathcal{V}$ such 
that $\bigcup \mathcal{U}$ is dense in $\bigcup \mathcal{V}$?

\mynote{Solution}
Yes to the first; independently answered by B.~Scott, 
O.~F\"orster and S.~Watson 
\cite{MR81m:54034,Forster,MR81j:54042}.
No to the second problem (B. Scott).
\end{myprob}

\begin{myprob}
\myproblem{D20}{G.M.~Reed}
Does \manotch\ imply either
perfect (normal), locally compact spaces are subparacompact or
that there is no Dowker manifold? 
\end{myprob}

\begin{myprob}
\myproblem{D21}{G.M.~Reed \cite{MR80j:54020}}
Does there exist in \myaxiom{ZFC} a normal space of cardinality $\aleph_1$ 
with a point-countable base which is not perfect?

\mynote{Notes}
With $\mathfrak{c}$ in place of $\omega_1$, there are many examples, such
as the Michael line.
P. Davies \cite{MR80j:54024} constructed a completely regular space of 
cardinality $\aleph_1$ with a point-countable base which is not perfect.
If there is a normal counterexample, then the closed set which is not a 
$G_\delta$ cannot be discrete.
\end{myprob}

\begin{myprob}
\myproblem{D22}{G.M.~Reed}
Does there exist a strongly collectionwise Hausdorff Moore space which is
not normal?

\mynote{Solution}
Yes, if there is a $Q$-set;
C. Navy \cite{Navy} proved that every para-Lindel\"of Moore space is
strongly collectionwise Hausdorff.
\end{myprob}

\begin{myprob}
\myproblem{D23}{M.A.~Swardson, attributed to R. Blair \cite{MR82i:54015}}
Does \manotch\ imply that every perfectly normal space of nonmeasurable 
cardinality is realcompact?

\mynote{Solution}
F. Hern\'andez-Hern\'andez and T. Ishiu \cite{hhi} showed that is it 
consistent with \manotch\ that there is a perfectly normal 
non-realcompact space of cardinality $\aleph_1$. 
The example is obtained by refining the order topology on $\omega_1$ in a 
forcing extension. 
\end{myprob}

\begin{myprob}
\myproblem{D24}{S. Watson, The Arhangel$'$\kern-.1667em ski\u{\i}-Tall Problem}
Is every normal, locally compact, metacompact space paracompact?

\mynote{Solution}
The answer is independent.
See Watson's contribution to \emph{New Classic Problems}.
\end{myprob}

\begin{myprob}
\myproblem{D25}{S. Watson \cite[Problem 88]{MR1078640}}
Is there a locally compact, perfectly normal space which is not
paracompact?

Yes if \myaxiom{MA} or if there exists a Souslin tree. 
Yes, if $\diamondsuit^*$ (G. Gruenhage and P.~Daniels \cite{MR86m:54021}).
A real example must be collectionwise Hausdorff under \visl\ but must not
be under \manotch; if one adds $\aleph_2$ random reals to 
a model of \visl\ the example must be collectionwise normal in the model.
\end{myprob}

\begin{myprob}
\myproblem{D26}{W. Fleissner and G.M.~Reed \cite{MR80j:54020}}
Is every collectionwise normal para-Lindel\"of space paracompact?

\mynote{Notes}
C. Navy \cite{Navy} gave an example of a normal para-Lindel\"of 
nonparacompact space.
This problem is in \cite[Problem 109]{MR1078640}.
\end{myprob}

\begin{myprob}
\myproblem{D27}{H. Wicke \cite[R. Hodel]{MR50:8400}}
Is every collectionwise normal meta-Lindel\"of space paracompact?
What if it is first countable?

\mynote{Notes}
This problem is in \cite[Problem 110]{MR1078640}.

\mynote{Solution}
No, Z. Balogh gave a \myaxiom{ZFC} example \cite{MR2003c:54047}.
See Watson's contribution to \emph{New Classic Problems}.
\end{myprob}

\begin{myprob}
\myproblem{D28}{H. Wicke}
Is there a meta-Lindel\"of space which is not weakly $\theta$-refinable
(=not weakly submetacompact)?

\mynote{Notes}
Yes, if \myaxiom{CH} (R.J.~Gardner and G. Gruenhage \cite{MR80d:54021}).

\mynote{Solution}
Yes, G. Gruenhage \cite{MR88j:54039} showed that for the Corson compact 
space $X$ constructed by S. Todor\v{c}evi\'c \cite{MR84g:03078},
$X^2 \setminus \Delta$ is a meta-Lindel\"of space which is not 
weakly submetacompact.
\end{myprob}

\begin{myprob}
\myproblem{D28}{P. de Caux \cite{MR83h:54032}}
Is every Lindel\"of space a $D$-space?

\mynote{Notes}
Compare D8.
\end{myprob}

\begin{myprob}
\myproblem{D29}{G. Grabner}
Suppose that $X$ is a regular wrb space.
Are the following equivalent?
$X$ is paracompact;
$X$ is irreducible and $\aleph$-preparacompact;
$X$ is submetacompact and $\aleph$-preparacompact.

\mynote{Notes}
See problems A11, D2 and D4 for definitions.
\end{myprob}

\begin{myprob}
\myproblem{D30}{P. Nyikos}
Is there a first countable space (or even a space of countable 
pseudocharacter) that is weakly $\theta$-refinable (weakly 
submetacompact) and countably metacompact, but not subparacompact?

\mynote{Notes}
Yes to the \emph{countable pseudocharacter} version (P. Nyikos).

\mynote{Solution}
Yes if $\neg\text{\myaxiom{CH}}$. 
In fact, G. Gruenhage and Z. Balogh have shown that \myaxiom{CH} is 
equivalent to the statement that every locally compact, first countable, 
$\theta$-refinable (=submetacompact) space is subparacompact.
Gruenhage's $\neg\text{\myaxiom{CH}}$ example is, in addition, 
metacompact.
\end{myprob}

\begin{myprob}
\myproblem{D31}{P. Nyikos}
Is there a quasi-developable countably metacompact space which is not
subparacompact?

\mynote{Solution}
P. Gartside, C. Good, R. Knight and A. Mohamad \cite{MR2002c:54002} 
constructed a quasi-developable manifold which is not developable 
(hence not subparacompact).
Furthermore, it is consistent that the example can be made countably 
metacompact.
\end{myprob}

\begin{myprob}
\myproblem{D32}{P. Nyikos}
Is every quasi-developable collectionwise normal countably paracompact
space paracompact?
\end{myprob}

\begin{myprob}
\myproblem{D33}{P. Nyikos}
Does \myaxiom{MA} imply every locally compact Hausdorff space of weight 
$< \mathfrak{c}$ is either subparacompact or contains a countably compact
noncompact subspace?
If one substitutes ``cardinality'' or ``weight'' the answer is affirmative
(Z. Balogh).

\mynote{Solution}
No (P. Nyikos): there is a \myaxiom{ZFC} example of a manifold of weight 
$\aleph_1$ which is quasi-developable but not even countably metacompact.
\end{myprob}

\begin{myprob}
\myproblem{D34}{E. van~Douwen}
Is there a nonparacompact, collectionwise normal space that is \emph{not
trivially so}?
Such a space would be realcompact and countably paracompact, and each
closed subspace $F$ would be irreducible (i.e., every open cover has an
open refinement with no proper subcover) or at least satisfy 
$L(F) = \hat{e}(F)$ where
$L(F)
= \min\{\kappa : \text{each open cover of $F$ has a subcover of 
cardinality $\leq \kappa$}\}$
and
$\hat{e}(F)
= \min\{\kappa : \text{no closed discrete subspace of $F$ has 
cardinality $\kappa$}\}$.
It would be even better if the space is a $D$-space, i.e., for every
neighbornet there is a closed discrete subspace $D$ such that the
restriction of the neighbornet to $D$ covers the space.

\mynote{Solution}
G. Gruenhage \cite{gruenhaged} proved that R. Pol's 1977 example 
\cite{MR57:4113.1} of a perfectly normal, collectionwise normal, 
nonparacompact space is a $D$-space.
\end{myprob}

\begin{myprob}
\myproblem{D35}{P. Nyikos}
Does there exist a screenable anti-Dowker space? 
That is, does there exist a screenable space that is countably
paracompact but not normal? 
If \myaxiom{PMEA}, any example must be of character $\geq \mathfrak{c}$.

\mynote{Solution}
Yes, applying the Wage machine \cite{MR53:9158}
to Bing's example $G$ gives a screenable space.
\end{myprob}

\begin{myprob}
\myproblem{D36}{D. Burke and P. Nyikos}
In a regular, first countable, countably metacompact space, must every
closed discrete subspace be a $G_\delta$?
What if the space is countably paracompact? normal?

\mynote{Notes}
Yes to each question if \myaxiom{PMEA}.
Yes to the second (S. Watson) and third (W. Fleissner) if \visl:
countably paracompact (resp.\ normal) first countable, Hausdorff spaces
are collectionwise Hausdorff.

\mynote{Solution}
No to the first question, if \visl.
P. Szeptycki \cite{MR95h:54016} constructed from $\diamondsuit^*$ a
first countable, regular, countably metacompact space with a closed
discrete set that is not a $G_\delta$-set.
\end{myprob}

\begin{myprob}
\myproblem{D37}{P. Nyikos}
Is there a real example of a locally compact, realcompact, first
countable space of cardinality $\aleph_1$ that is not normal?

\mynote{Solution}
Yes, there is a Moore space obtained by splitting nonisolated points of
the Cantor tree, which has all the desired properties.
See S. Shelah's \cite[Theorem~11.4.2]{MR84h:03002}.
\end{myprob}

\begin{myprob}
\myproblem{D38}{C.R.~Borges and A. Wehrly \cite{MR94a:54059}}
Are subparacompact spaces $D$-spaces?
\end{myprob}

\begin{myprob}
\myproblem{D39}{C.R.~Borges and A. Wehrly \cite{MR94a:54059}}
Are monotonically normal paracompact spaces $D$-spaces?
\end{myprob}

\begin{myprob}
\myproblem{D40}{C.R.~Borges and A. Wehrly \cite{MR94a:54059}}
Is the countable product of Sorgenfrey lines a $D$-space?

\mynote{Notes}
In the article by Borges and Wehrly, it was also asked whether the finite
product of irrational Sorgenfrey lines is a $D$-space, but this was
answered affirmatively by P. de Caux \cite{MR83h:54032}
where he showed that each subspace of each finite power of the Sorgenfrey
line is a $D$-space.
\end{myprob}

\begin{myprob}
\myproblem{D41}{K. Tamano \cite{MR96d:54020}}
Is the space $\omega^{\omega_1}$ weakly $\delta\theta$-refinable?

\mynote{Solution}
No (J. Chaber, G. Gruenhage, R. Pol \cite{MR96k:54061}).
\end{myprob}

\begin{myprob}
\myproblem{D42}{P. Szeptycki \cite[P.~Nyikos]{MR97a:54022}}
Does \visl\ imply that first countable, countably paracompact spaces are
strongly collectionwise Hausdorff?
\end{myprob}

\begin{myprob}
\myproblem{D43}{P. Szeptycki \cite{MR97a:54022}}
Are first countable, countably paracompact, collectionwise Hausdorff 
spaces strongly collectionwise Hausdorff?

\mynote{Notes}
(P. Szeptycki)
A space is strongly collectionwise Hausdorff if closed discrete sets can 
be separated by a discrete family of open sets.
The structure of closed discrete sets in first countable spaces has a long
and interesting history beginning with the normal Moore space conjecture.
The question whether normal, first countable spaces are collectionwise 
Hausdorff and whether countably paracompact, first countable spaces are 
collectionwise Hausdorff is particularly interesting.
A series of results by D.~Burke, W.~Fleissner, P.~Nyikos, 
F.~Tall, and S.~Watson address these questions under \visl, 
\myaxiom{PMEA}, and other assumptions.
D42 of Nyikos appears to be one of the last important questions concerning
the effect of \visl\ on the separation of closed discrete sets in first
countable spaces.

While Burke has shown that \myaxiom{PMEA} provides a consistent positive 
answer (even without the assumption of collectionwise Hausdorff) 
\cite{MR85h:54032}, a positive answer to D43 assuming \visl\ would yield 
a positive answer to Nyikos's question.
However, any consistent counterexample would go a long way toward clarifying 
the distinction between normality and countable paracompactness.
Note that the assumption of first countability is essential as a 
\myaxiom{ZFC} example with uncountable character has been constructed by 
Watson \cite{MR92d:54024}.
Also, if we weaken countable paracompactness to paranormality in D42 or D43, 
then \cite{MR98j:54038} gives consistent negative answers, respectively.
\end{myprob}

\begin{myprob}
\myproblem{D44}{K. Yamazaki \cite{MR2000a:54026}}
Let $X$ be a collectionwise normal space and $Y$ a paracompact
$\Sigma$-space (or a paracompact $\sigma$-space, or a $M_3$-space).
Suppose $X \times Y$ is normal and countably paracompact. 
Then is $X \times Y$ collectionwise normal?

\mynote{Notes}
Yes, if $X$ is also a $P$-space.
See also the author's second article \cite{MR2000a:54027}.
\end{myprob}

\section*{E. Separation and disconnectedness}

\begin{myprob}
\myproblem{E1}{M. Wage \cite{MR81d:54019,MR55:6358,MR56:16595}}
Is there an extremally disconnected Dowker space?

\mynote{Solution}
Yes (A. Dow and J. van~Mill \cite{MR84a:54028}).
\end{myprob}

\begin{myprob}
\myproblem{E2}{M. Wage \cite{MR56:16595}}
Is there a strong $S$-space that is extremally disconnected?

\mynote{Notes}
If \manotch, there are no strong $S$-spaces 
(K. Kunen \cite{MR55:13362}).
\end{myprob}

\begin{myprob}
\myproblem{E3}{P. Bankston \cite{MR56:13141,MR80h:54007}}
Are ultraproducts of scattered Hausdorff spaces scattered?
Non-Hausdorff counterexamples are known.

\mynote{Notes}
No. E. van~Douwen showed that an ultrapower of a scattered space $X$ is 
scattered if and only if the Cantor-Bendixson rank of $X$ is finite.
Bankston had translated to ultraproducts a question of R.W.~Button 
\cite{MR80a:54096} in nonstandard topology: if $X$ is scattered, is then
${}^*\!X$, endowed with the $Q$-topology, scattered?
R. \v{Z}ivaljevi\'c \cite{MR86c:54053} showed that ${}^*\!X$ is scattered 
iff $X$ has finite Cantor-Bendixson rank.
\end{myprob}

\begin{myprob}
\myproblem{E4}{B. Smith-Thomas \cite{MR80k:54044}}
If $X$ is a $k_W$-space, is $\beta X \setminus X$ necessarily an 
$F$-space?

\mynote{Solution}
A $k_\omega$-space has the weak topology determined by an increasing 
sequence of compact, $T_2$ subspaces of which it is the union. 
E.~van~Douwen showed that the answer to E4 is negative.
A proof similar to that for the rationals shows that no countable
dense-in-itself $k_\omega$-space has an $F$-space for its growth.
\end{myprob}

\begin{myprob}
\myproblem{E5}{A.V.~Arhangel$'$\kern-.1667em ski\u{\i}}
Does every zero-dimensional space have a strongly zero-dimensional
subtopology?

\mynote{Notes}
(P. Nyikos)
All examples of zero-dimensional spaces which are known to the Problems
Editor have strongly zero-dimensional subtopologies.
This is clear in the locally compact examples, and has been shown for 
Prabir Roy's Space $\Delta$.
\end{myprob}

\begin{myprob}
\myproblem{E6}{T. Przymusi\'nsk \cite{MR83b:54007}}
If $\mathcal{F}[X]$ is the Pixley-Roy hyperspace over $X$, then is
$\mathcal{F}[X]$ strongly zero-dimensional? 
Yes, if the hyperspace is normal.
\end{myprob}

\begin{myprob}
\myproblem{E7}{K. Kunen \cite{MR81m:54044}}
Is there a locally compact, extremally disconnected space which is normal
but not paracompact? 
Yes, if there exists a weakly compact cardinal.
\end{myprob}

\begin{myprob}
\myproblem{E8}{S. Watson}
Is there a locally compact, normal, non-collectionwise normal space?
Yes, if $\text{\myaxiom{MA}}(\omega_1)$ or in models of \visl, 
\myaxiom{CH} or $\neg\text{\myaxiom{CH}}$.

\mynote{Notes}
If $\kappa$ is supercompact and $\kappa$ Cohen or random reals are added
to a model of \myaxiom{ZFC}, then the answer is negative in the resulting 
model (Z. Balogh \cite{MR91c:54030}).
It is not yet known whether a negative answer is equiconsistent with 
\myaxiom{ZFC}.
\end{myprob}

\begin{myprob}
\myproblem{E9}{S. Watson}
Is there a perfectly normal, collectionwise Hausdorff space which is not
collectionwise normal?
Yes, in some models. 
\end{myprob}

\begin{myprob}
\myproblem{E10}{E. van~Douwen}
Characterize internally the class $T_3 \vdash T_4$ of regular spaces $X$
such that every regular continuous image of $X$ is normal.

\mynote{Notes}
The class of spaces ACRIN (all continuous regular images normal)
has been studied in \cite{MR89k:54054,MR90b:54021}.
Note that $\omega_1$ and $\omega \times (\omega_1 + 1)$ are examples but
their direct sum is not.
\end{myprob}

\begin{myprob}
\myproblem{E11}{F.D.~Tall}
Levy collapse a supercompact cardinal to $\omega_2$. 
Are first countable (locally countable) $\aleph_1$-collectionwise normal
space collectionwise normal?
\end{myprob}

\begin{myprob}
\myproblem{E12}{Y. Hattori and H. Ohta, attributed to S. Nadler 
\cite{MR95k:54052}}
Must a totally disconnected separable metric space having UMP
be zero-dimensional?

\mynote{Notes}
A metric space $(X, d)$ has the \emph{unique midpoint property} (UMP)
if for every pair of distinct points $x$ and $y$ of $X$, there exists 
exactly one point $p$ such that $d(x,p) = d(y,p)$.
A positive answer to BB1 also answers this question positively.
\end{myprob}

\begin{myprob}
\myproblem{E13}{A.V.~Arhangel$'$\kern-.1667em ski\u{\i} \cite{MR89b:54004}}
Is there in \myaxiom{ZFC} a non-discrete extremally disconnected 
topological group?

\mynote{Notes}
This is an old problem; see \cite{MR36:2122}.
\end{myprob}

\section*{F. Continua theory}

\begin{myprob}
\myproblem{F1}{C. Hagopian, attributed to Bing \cite{MR22:1869}}
Is there a homogeneous tree-like continuum that contains an arc? 

\mynote{Solution}
No, (C. Hagopian \cite{MR84d:54059}).
\end{myprob}

\begin{myprob}
\myproblem{F2}{W.T.~Ingram \cite{TP1.329}}
Is there an atriodic tree-like continuum which cannot be embedded in the
plane? 

\mynote{Solution}
Yes. 
L. Oversteegen and E. Tymchatyn \cite{MR84h:54030} have given two 
examples. 
One of them consists of taking the atriodic tree-like continuum $X$ 
either of Bellamy or of Oversteegen and Rogers and adjoining to it two 
arcs at endpoints of two composants of $X$.
It is still of interest to determine if the continuum $X$ itself is an 
example, for if it were not, then it would be a solution to the 
fixed-point problem for nonseparating plane continua.
\end{myprob}

\begin{myprob}
\myproblem{F3}{W.T.~Ingram \cite{TP1.329}}
What characterizes the tree-like continua which can be embedded in the plane? 
\end{myprob}

\begin{myprob}
\myproblem{F4}{W.T.~Ingram \cite{TP1.329}}
What characterizes the tree-like continua which are in class $W$? 

A continuum $X$ is said to be in \emph{class $W$} if each continuous 
surjection from a continuum onto $X$ is weakly confluent.

\mynote{Solution}
J. Grispolakis and E. Tymchatyn have shown that a continuum $X$ is in
class $W$ if and only if it has the covering property, i.e., for any
Whitney map $\mu$ for $C(X)$ and any $t \in (0, \mu(X))$, no proper
subcontinuum of $\mu^{-1}(t)$ covers $X$. 
They have also shown that a planar tree-like continuum is in class $W$ if
and only if it is atriodic. 
\end{myprob}

\begin{myprob}
\myproblem{F5}{G.R.~Gordh and L. Lum \cite{MR58:7568}} 
Let $M$ be a continuum containing a fixed point $p$. 
Are the following conditions equivalent?
\begin{myenumerate}
\item
Each subcontinuum of $M$ which is irreducible from $p$ to some other
point is a monotone retract of $M$.
\item
Each subcontinuum of $M$ which contains $p$ is a monotone retract of $M$.
\end{myenumerate}
\end{myprob}

\begin{myprob}
\myproblem{F6}{J.T.~Rogers \cite{MR80k:54065}}
Suppose $M$ and $N$ are solenoids of pseudo-arcs that decompose to the
same solenoid.
Are $M$ and $N$ homeomorphic?

\mynote{Solution}
W. Lewis provided a positive answer to this question and completed the
classification of homogeneous, circle-like continua. 
\end{myprob}

\begin{myprob}
\myproblem{F7}{C.J.~Rhee \cite{MR84f:54010}}
Does admissibility of a metric continuum imply property c?
\end{myprob}

\begin{myprob}
\myproblem{F16}{J.T.~Rogers}
Is each nondegenerate, homogeneous, nonseparating plane continuum a
pseudo-arc?

\mynote{Notes}
Problems F16--F28 are discussed by Rogers in \cite{MR85c:54055}.
\end{myprob}

\begin{myprob}
\myproblem{F17}{J.T.~Rogers}
Is each Type $2$ curve a bundle over the Menger universal curve with 
Cantor sets as the fibers?

Here, a curve is a one-dimensional continuum and a curve is Type $2$ if it
is aposyndetic but not locally connected, and homogeneous.

\mynote{Solution}
No. J. Prajs's example (see Problem F20) is an aposyndetic, non-locally 
connected, homogeneous curve that is not the total space of a Cantor set 
bundle over the Menger curve. 
\end{myprob}

\begin{myprob}
\myproblem{F18}{J.T.~Rogers}
Is each Type $2$ curve an inverse limit of universal curves and maps? 
universal curves and fibrations? universal curves and covering maps?
\end{myprob}

\begin{myprob}
\myproblem{F19}{J.T.~Rogers}
Does each Type $2$ curve contain an arc?
\end{myprob}

\begin{myprob}
\myproblem{F20}{J.T.~Rogers}
Does each Type $2$ curve retract onto a solenoid?

\mynote{Solution}
No.  
J. Prajs \cite{MR2003f:54077} constructed a homogeneous, arcwise connected, 
non-locally connected curve.  
Such a curve must be aposyndetic. 
Since it is arcwise connected, it cannot be mapped onto a solenoid, let 
alone retracted onto one.  
\end{myprob}

\begin{myprob}
\myproblem{F21}{J.T.~Rogers}
Does each indecomposable cyclic homogeneous curve that is not a solenoid
admit a continuous decomposition into tree-like curves so that the
resulting quotient space is a solenoid?

\mynote{Solution}
E. Duda, P. Krupski and J.T.~Rogers have some partial results on this 
problem.
Krupski and Rogers \cite{MR92c:54038}
showed that the answer is yes for finitely cyclic homogeneous curves. 
Duda, Krupski, and Rogers \cite{MR92m:54061}  
show that a $k$-junctioned, homogeneous curve must be a pseudo-arc, a 
solenoid, or a solenoid of pseudo-arcs.   
See also \cite{MR90h:54038}.
\end{myprob}

\begin{myprob}
\myproblem{F22}{J.T.~Rogers}
Is every acyclic homogeneous curve tree-like?
In other words, does trivial cohomology imply trivial shape for
homogeneous curves?

\mynote{Solution}
Yes (Rogers \cite{MR88k:57016}).
\end{myprob}

\begin{myprob}
\myproblem{F23}{J.T.~Rogers}
Is every tree-like, homogeneous curve hereditarily indecomposable?
It is a pseudo-arc? 
It is weakly chainable? 
Does it have span zero?
Does it have the fixed-point property?

\mynote{Solution}
J. Prajs answered an old question of Bing by proving that each 
tree-like homogeneous continuum is hereditarily indecomposable. 
See \cite{MR90f:54054}
\end{myprob}

\begin{myprob}
\myproblem{F24}{F.B.~Jones}
Is each tree-like, homogeneous curve hereditarily equivalent?
\end{myprob}

\begin{myprob}
\myproblem{F25}{J.T.~Rogers}
Is each decomposable, homogeneous continuum of dimension greater than one
aposyndetic?
\end{myprob}

\begin{myprob}
\myproblem{F26}{J.T.~Rogers}
Is each indecomposable, nondegenerate, homogeneous continuum
one-dimensional?
\end{myprob}

\begin{myprob}
\myproblem{F27}{J.T.~Rogers}
Must the elements of Jones's aposyndetic decomposition be hereditarily
indecomposable?
\end{myprob}

\begin{myprob}
\myproblem{F28}{J.T.~Rogers}
Let $X$ be a homogeneous curve, and let $H(X)$ be its homeomorphism group.
Let $\mathcal{G}$ be a partition of $X$ into proper, nondegenerate
continua so that $H(X)$ respects $\mathcal{G}$ (this means that either
$h(G_1) = G_2$ or $h(G_1) \cap G_2 = \emptyset$, for all $G_1$ and $G_2$
in $\mathcal{G}$ and all $h$ in $H(X)$).
Are the members of $\mathcal{G}$ hereditarily indecomposable?
\end{myprob}

\begin{myprob}
\myproblem{F29}{P. Nyikos}
Is there a (preferably first countable, or, better yet, perfectly
normal) locally connected continuum without a base of open subsets with
locally connected closures?
Yes to the general question if \myaxiom{CH}.
\end{myprob}

\begin{myprob}
\myproblem{F30}{J. Grispolakis \cite{MR89g:54027}}
If $Y$ is an LC$'$ continuum with no local separating points does
$(Y, y_0)$ have the avoidable arcs property for some $y_0 \in Y$?
\end{myprob}

\begin{myprob}
\myproblem{F30}{E. Tymchatyn, attributed to D. Bellamy \cite{MR81b:54032}}
Let $S_n$ be a solenoid and let $K_n$ be the Knaster indecomposable
continuum obtained by identifying in the topological group $S_n$ each
point with its inverse.
Do there exist in $K_2$ two components without endpoints which are not
homeomorphic?

\mynote{Solution}
C. Bandt \cite{MR95m:54034a} proved that all components (except for 
the one with endpoints) of $K_2$ are homeomorphic.
\end{myprob}

\begin{myprob}
\myproblem{F32}{E. Tymchatyn \cite{MR90f:54051}}
If $S_n$ has a composant that is homeomorphic to one of $S_m$, is $S_n$
homeomorphic to $S_m$?

\mynote{Solution}
R. de Man \cite{MR96g:54040} proved that any composants of any nontrivial 
solenoids are homeomorphic.
\end{myprob}

\begin{myprob}
\myproblem{F33}{B.E.~Wilder \cite{MR94i:54073}}
Which of the known results concerning aposyndetic continua can be
extended to the class of $C$-continua?
\end{myprob}

\begin{myprob}
\myproblem{F34}{S. Mac\'{\i}as \cite{MR1255821}}
In this and the following two problems, let
$\Gamma = \mathcal{S}^1 \times \mathbb{Q}$ or
$\Gamma = \mathcal{W} \times \mathbb{Q}$,
where $\mathcal{W}$ is the figure eight, and $\mathbb{Q}$ is the Hilbert
cube.
Let $\sigma\colon \tilde{\Gamma} \to \Gamma$ be the universal covering 
space.
Let $X$ be a homogeneous continuum essentially embedded in $\Gamma$ and
let $\tilde{X} = \sigma^{-1}(X)$.
Two points of $\tilde{X}$ are said to be in the same \emph{continuum
component} if there is a continuum in $\tilde{X}$ containing them.
Do the continuum components and components of $\tilde{X}$ coincide if $X$
is homogeneous?
Without homogeneity, the answer is known to be negative.
\end{myprob}

\begin{myprob}
\myproblem{F35}{S. Mac\'{\i}as \cite{MR1255821}}
Let $\tilde{K}$ be a component of $\tilde{X}$.
If $\sigma(\tilde{K}) \neq X$, is it true that $\sigma(\tilde{K})$ is
contained in a composant of $X$?
Is it equal to a composant?
\end{myprob}

\begin{myprob}
\myproblem{F35}{S. Mac\'{\i}as \cite{MR1255821}}
Suppose that $\sigma(\tilde{K}) \neq X$.
If $\tilde K_\mathcal{C}$ is a continuum component of $\tilde{X}$,
is it true that $\sigma(K_\mathcal{C})$ is equal to a composant of $X$?
\end{myprob}

\begin{myprob}
\myproblem{F37}{The Classical Plane Fixed-Point Problem}
Does every nonseparating plane continuum have the fixed-point property?

\mynote{Notes}
See Hagopian's essay in this volume.
\end{myprob}

\begin{myprob}
\myproblem{F38}{C. Hagopian, another classical problem \cite{MR96m:54060}}
Must the cone over a tree-like continuum have the fixed-point property?
\end{myprob}

\begin{myprob}
\myproblem{F39}{C. Hagopian \cite{MR96m:54060}}
Does the cone over a spiral to a triod have the fixed-point property?
\end{myprob}

\begin{myprob}
\myproblem{F40}{\cite[J.~{\L}ysko]{MR96m:54060}}
Does there exist a $2$-dimensional contractible continuum that admits a
fixed-point-free homeomorphism?
\end{myprob}

\begin{myprob}
\myproblem{F41}{\cite[Bing]{MR96m:54060}}
If $M$ is a plane continuum with the fixed-point property, does 
$M \times [0,1]$ have the fixed-point property?
\end{myprob}

\begin{myprob}
\myproblem{F42}{C. Hagopian \cite{MR96m:54060}}
If $M$ is a simply connected plane continuum, does $M \times [0,1]$ have 
the fixed-point property? 
\end{myprob}

\begin{myprob}
\myproblem{F43}{C. Seaquist \cite{MR97h:54008}}
Does there exist a continuous decomposition of the two dimensional disk
into pseudo-arcs?
\end{myprob}

\begin{myprob}
\myproblem{F44}{J.T.~Rogers}
Let $M$ be a hereditarily indecomposable continuum. 
Assume $\dim M = n > 1$. 
Let $H(M)$ be the homeomorphism group of $M$. 
Can $H(M)$ contain a nontrivial continuum? a nontrivial connected set? 

\mynote{Notes}
For each integer $n > 1$, 
Rogers exhibited an $M$ such that $H(M)$ contains no nontrivial
connected set. 

\mynote{Solution}
Yes. 
M. Re\'nska \cite{MR2003g:54082} proved that there exist rigid 
hereditarily indecomposable continua in every dimension.  
In fact there exist continuum many such continua in each dimension.  
\end{myprob}

\begin{myprob}
\myproblem{F45}{J.T.~Rogers}
Can $M$ be rigid? 
i.e., the identity map is the only element of $H(M)$? 
\end{myprob}

\begin{myprob}
\myproblem{F46}{C. Seaquist \cite{MR2002m:54039}}
Does there exist a continuous decomposition $G$ of the plane into acyclic
continua so that for every point $x$, there is an arc $A$ and an element
$g \in G$ such that $x \in A \subset g$?
\end{myprob}

\section*{G. Mappings of continua and Euclidean spaces}

\begin{myprob}
\myproblem{G1}{D. Mauldin and B. Brechner \cite{TP1.335}}
Let $K$ be a locally connected, nonseparating continuum in $E^2$, 
$K$ not a disk. 
Let $h$ be an $\text{EC}$ homeomorphism of $K$ onto itself such that $h$ 
is extendable to a homeomorphism $\tilde{h}$ of $E^2$ onto itself.
Is $h$ necessarily periodic? 
Does there exist a homeomorphism 
$g\colon E^2 \to E^2$ such that $g\tilde{h}\colon E \to E^2$ 
is $\text{EC}^+$ with nucleus $K$? 

\mynote{Notes}
$h\colon E^n \to E^n$ is \emph{(uniformly) $\text{EC}^+$} if 
the set of non-negative iterates of $h$ forms a 
pointwise (uniformly) equicontinuous family.
$h$ is \emph{(uniformly) $\text{EC}$} iff the set of all iterates 
of $h$ forms a pointwise (uniformly) equicontinuous family.
\end{myprob}

\begin{myprob}
\myproblem{G2}{D. Mauldin and B. Brechner \cite{TP1.335}}
Let $h$ be an orientation preserving, $\text{EC}^+$ homeomorphism of
$E^2$ onto itself. 
If the nucleus of $h$ is unbounded, can $h$ be imbedded in a flow? 
\end{myprob}

\begin{myprob}
\myproblem{G3}{D. Mauldin and B. Brechner \cite{TP1.335}}
Characterize the $\text{EC}$ and $\text{EC}^+$ homeomorphisms of 
$R^\infty$.
\end{myprob}

\begin{myprob}
\myproblem{G4}{D. Mauldin and B. Brechner \cite{TP1.335}}
Characterize the nuclei of the $\text{EC}^+$ homeomorphisms of $E^n$
and characterize the action of such homeomorphisms on its nucleus. 
\end{myprob}

\begin{myprob}
\myproblem{G5}{D. Mauldin and B. Brechner \cite{TP1.335}}
Let $h$ be an orientation preserving $\text{EC}^+$ homeomorphism of
$E^n$ onto itself whose nucleus $M$ is bounded. 
If $n$ is $4$ or $5$, is it true that 
$\tilde{h}\colon {E^n}/M \to {E^n}/M$ is a topological standard 
contraction? 
\end{myprob}

\begin{myprob}
\myproblem{G6}{A. Petrus \cite{MR56:9498}}
Let $X$ be a continuum and let $\mu\colon C(X) \to [0, \infty)$ be a 
Whitney map. 
If $\mu^{-1}(t_0)$ is decomposable, must $\mu^{-1}(t)$ be decomposable
for all $t \in [t_0, \mu(X)]$? 
\end{myprob}

\begin{myprob}
\myproblem{G7}{A. Petrus \cite{MR56:9498}}
Let $\mu\colon C(X) \to [0, \infty)$ be a Whitney map. 
Characterize those continua which satisfy, for all $t \in [0, \mu(X)]$:
\begin{myenumerate}
\item 
if $\mathcal A$ is a subcontinuum of $\mu^{-1}(t)$ and 
$\sigma\mathcal{A} := \bigcup\{A : A \in \mathcal{A}\} = X$, then 
$\mathcal{A} = \mu^{-1}(t)$;
\item
if $\mathcal A$ is a subcontinuum of $\mu^{-1}(t)$, then
$A \in \mathcal{A}$ for all $A \in \mu^{-1}(t)$ such that 
$A \subset \sigma\mathcal{A}$.
\end{myenumerate}

\mynote{Notes}
See also \cite{MR80j:54029}.
\end{myprob}

\begin{myprob}
\myproblem{G8}{E.E.~Grace \cite{MR80m:54022}}
Is there a monotonely refinable map (i.e., a map that can be
$\epsilon$-approximated by a monotone $\epsilon$-map, for each positive
$\epsilon$) from a regular curve of finite order onto a topologically
different regular curve of finite order?
\end{myprob}

\begin{myprob}
\myproblem{G9}{C. Hagopian \cite{MR80k:54064}}
Is every continuous image of every $\lambda$-connected plane continuum
$\lambda$-connected?

\mynote{Notes}
Outside of the plane the answer is no since Hagopian showed that the
product of two nondegenerate, hereditarily indecomposable continua
is $\lambda$-connected. 
\end{myprob}

\begin{myprob}
\myproblem{G10}{R. Heath and P. Fletcher}
Is there a Euclidean non-Galois homogeneous continuum?

\mynote{Solution}
W. Kuperberg observed that the product of two or more copies of the Menger 
universal curve is homogeneous, but not Galois. 
An equivalent result is true for a product of pseudo-arcs. 
\end{myprob}

\begin{myprob}
\myproblem{G11}{E. Lane \cite{MR82g:54020}}
What is a necessary and sufficient condition in order for a space to
satisfy the $C$ insertion property for (nusc, nlsc)?

\mynote{Solution}
Here \emph{nusc} and \emph{nlsc} are the classes of normal lower and 
upper semicontinuous functions.
\end{myprob}

\begin{myprob}
\myproblem{G12}{E. van~Douwen}
Let $\mathbb{H}$ denote the half-line, $[0,+\infty)$.
Is every continuum of weight $\leq \omega_1$ a continuous image of
$\mathbb{H}^* = \beta\mathbb{H} \setminus \mathbb{H}$?
Yes, for metrizable continua
(J.M.~Aarts and P. van Emde Boas \cite{MR33:4904}).

\mynote{Solution}
Yes, (A. Dow and K.P.~Hart \cite{MR2001g:54037}).
\end{myprob}

\begin{myprob}
\myproblem{G13}{C.J.~Rhee \cite{MR84f:54010}}
For each fiber map $\alpha\colon X \to C(X)$, does there exist a 
continuous fiber map $\beta:X \to C(X)$ such that $\beta(x) \subset (x)$ 
for each $x \in X$?
\end{myprob}

\begin{myprob}
\myproblem{G13}{J. Mayer \cite{MR85f:54066}}
Are there uncountably many inequivalent embeddings of the pseudo-arc in 
the plane with the same prime end structure?
\end{myprob}

\begin{myprob}
\myproblem{G14}{R.G.~Gibson \cite{MR84d:26004}}
Give necessary and sufficient conditions for the extension $I^2 \to I$ of 
a connectivity function $I \to I$ to be a connectivity function.
In particular, is it necessary for the function $I \to I$ to have the
CIVP?

\mynote{Notes}
A function $f\colon I \to I$ has the \emph{Cantor intermediate value 
property} (CIVP) if for each Cantor set $K$ in $(f(x),f(y))$ there exists 
a Cantor set $C$ in $(x,y)$ for which $f(C) \subset K$. 
\end{myprob}

\begin{myprob}
\myproblem{G14}{J. Mayer \cite{MR85f:54066}}
Are there countably many inequivalent embeddings in the plane of every
indecomposable chainable continuum (with the same prime end structure)?
\end{myprob}

\begin{myprob}
\myproblem{G15}{J. Keesling}
Can the maps in \cite[Theorem 5.2]{MR84d:54060} be made monotone or 
cell-like?
\end{myprob}

\begin{myprob}
\myproblem{G15}{J.T.~Rogers \cite[Problem 467]{MR1078656}}
Can Jones's aposyndetic decomposition raise dimension? lower dimension?

\mynote{Solution}
No and yes.  
J.T.~Rogers \cite{MR1992870} proved that if $X$ is a homogeneous, 
decomposable continuum that is not aposyndetic, then the dimension of its 
aposyndetic decomposition is one.  
Hence Jones' aposyndetic decomposition can never raise dimension; in fact, 
it must lower the dimension of every such continuum of dimension greater 
than one. 
\end{myprob}

\begin{myprob}
\myproblem{G16}{J.T.~Rogers \cite{MR85c:54055}}
A homeomorphism is primitively stable if its restriction to some nonempty
open set is the identity.
Does each homogeneous continuum admit a nontrivial primitively stable
homeomorphism? 
\end{myprob}

\begin{myprob}
\myproblem{G17}{J.T.~Rogers \cite{MR85c:54055}}
Is each homogeneous continuum bihomogeneous?
That is, given points $x$ and $y$ in $X$, does there exist a homeomorphism
$h$ of $X$ onto itself such that $h(x) = y$ and $h(y) = x$?

\mynote{Solution}
Around 1921, B. Knaster asked whether every homogeneous space is 
bihomogeneous.
C.~Kuratowski \cite{Kuratowski} described a non-locally compact 
homogeneous subset of the plane which is not bihomogeneous.
D.~van~Dantzig \cite{dantzig} asked whether homogeneity implies 
bihomogeneity for continua.
H.~Cook found a locally compact, homogeneous, metric space which is not 
bihomogeneous \cite{MR88e:54045}.
K.~Kuperberg \cite{MR90m:54043} solved this long-standing problem by 
constructing a locally connected, homogeneous, $7$-dimensional continuum 
which is not bihomogeneous.
P.~Minc \cite{MR96b:54051} constructed infinite dimensional, non-locally 
connected, homogeneous continua which are not bihomogeneous.
K.~Kawamura \cite{MR96i:54028} showed that for each $n \geq 2$ there are 
$n$-dimensional homogeneous continua which are not bihomogeneous.
\end{myprob}

\begin{myprob}
\myproblem{G18}{B. Brechner \cite{MR87i:54082}}
Let $h$ be a regular homeomorphism of $B^3$ onto itself, which is the
identity on the boundary.
Must $h$ be the identity?

\mynote{Notes}
Recall that $h$ is regular iff the family of all iterates of $h$ forms an
equicontinuous family of homeomorphisms.
\end{myprob}

\begin{myprob}
\myproblem{G19}{B. Brechner \cite{MR87i:54082}}
Is every regular, orientation preserving homeomorphism of $S^3$ either
periodic, or a rotation, or a combination of rotations on two solid tori
whose union is $S^3$.
\end{myprob}

\begin{myprob}
\myproblem{G20}{J. Grispolakis \cite{MR89h:54008}}
Let $f\colon X \to Y$ be a weakly confluent mapping from a compact 
connected PL $n$-manifold $X$ onto a PL $m$-manifold $Y$ with $n,m \geq 3$.
Is $f$ homotopic to a light open mapping of $X$ onto $Y$?
\end{myprob}

\begin{myprob}
\myproblem{G21}{J. Grispolakis \cite{MR89h:54008}}
Let $f\colon M \to Y$ be a mapping from a compact connected PL 
$n$-manifold, $n \geq 3$, into an ANR $Y$ such that every simple closed 
curve can be approximated by a spiral in $Y$.
If $\Pi(Y)$ has property (Tor) relative to $f_\#\Pi(M)$, is $f$ homotopic
to a weakly confluent mapping of $M$ onto $Y$?
\end{myprob}

\begin{myprob}
\myproblem{G22}{J. Grispolakis \cite{MR89h:54008}}
Characterize all weakly confluent images of the $3$-cube.
\end{myprob}

\begin{myprob}
\myproblem{G23}{E. Tymchatyn \cite{MR90f:54051}}
Is each homeomorphism $h\colon C \to \overline{C}$ of composants of 
solenoids homotopic to a linear homeomorphism
$\overline{h}\colon \overline{C} \to \overline{C}$
(i.e., $\overline{h}(x) = ax + b$ for each $x$)?

\mynote{Notes}
A positive solution would imply positive solutions to the
problems F30, F31, F32.
\end{myprob}

\begin{myprob}
\myproblem{G24}{J. Kennedy \cite{MR92b:54071}}
If $X$ is a continuum and $x \in X$, let $G_x$ denote the set of all
points of $X$ to which $x$ can be taken by a homeomorphism of $X$. 
It is known that $G_x$ is a Borel set, even if connectedness of $X$ is
dropped. 
If $\alpha$ is a countable ordinal, does there exist continuum $X$ for
which some $G_x$ is a Borel set in class $F_\alpha$ but not in
$F\gamma$ for $\gamma < \alpha$?
\end{myprob}

\begin{myprob}
\myproblem{G25}{J. Kennedy, \cite[attributed to M.~Barge]{MR92b:54071}}
Does there exist a weakly homogeneous planar continuum $X$ with the
property that each homeomorphism it admits possesses a dense set of
periodic points?
\end{myprob}

\begin{myprob}
\myproblem{G26}{J. Kennedy \cite{MR92b:54071}}
Does there exist a homogeneous continuum $X$ with the property that each 
of its homeomorphisms, except the identity, is transitive?
\end{myprob}

\begin{myprob}
\myproblem{G27}{J. Kennedy \cite{MR92b:54071}}
Does there exist a homogeneous continuum $X$ that admits a transitive
homeomorphism and that has the property that each of its homeomorphisms
admits a dense set of periodic points?
\end{myprob}

\begin{myprob}
\myproblem{G28}{J. Kennedy \cite{MR92b:54071}}
Does there exist a homogeneous continuum $X$ with the property that for
each non-identity homeomorphism $h$ of $X$, there is some nonempty proper
open set $U$ with $h(U) \subset U$?
\end{myprob}

\begin{myprob}
\myproblem{G28}{E. E, Grace \cite{MR93f:54046}}
If $X$ is a $\theta_n$-continuum and $f\colon X \to Y$ is proximately
refinable, must $Y$ be a $\theta_{2n}$-continuum?
\end{myprob}

\begin{myprob}
\myproblem{G29}{H. Pawlak and R. Pawlak \cite{MR95m:54012}}
A function is called a \emph{Darboux} function if it takes connected 
sets to connected sets.
If $X$ is connected and locally connected space, under what additional
assumption does there exist a connected Alexandroff compactification 
$X^*$ such that a theorem analogous to the following theorem holds?
\emph{Let $X$ be a continuum having an extension $X^*$ with a one-point 
remainder ${x_0}$ such that $X^*$ has an exploding point with respect to 
$x_0$.
Then there is a closed Darboux function $f\colon X^* \to [0, 1]$
which is discontinuous at $x_0$.}
In general, what kinds of hypotheses on a space $X^*$
(weaker than compactness) allow one to prove a theorem analogous
to this one?
\end{myprob}

\begin{myprob}
\myproblem{G30}{H. Pawlak and R. Pawlak \cite{MR95m:54012}}
Do there exist, for a nondegenerate locally connected continuum $X$ and 
any homeomorphism $h\colon X \to X$, spaces $X_1, X_2$ ``close to 
compactness'' such that $X$ is a subspace of $X_1$ and $X_2$ and 
there exists a d-extension $h^*\colon X_1 \to X_2$ of the function $h$ 
such that $h^*$ is a discontinuous and closed Darboux function?
\end{myprob}

\begin{myprob}
\myproblem{G31}{D. Garity \cite{MR98d:54063}}
If a homogeneous compact metric space is locally $n$-connected for all
$n$, is the space necessarily $2$-homogeneous?
\end{myprob}

\begin{myprob}
\myproblem{G32}{J. Haywood \cite{MR99j:54034}}
If $f\colon G \to G'$ is a universal function, is it possible that $G$ is 
a graph but not a tree?
Show that if $G'$ is a graph, then it is a tree.
\end{myprob}

\begin{myprob}
\myproblem{G33}{J. Charatonik and W. Charatonik \cite{MR99k:54023}}
Let $f\colon X \to Y$ be a surjective mapping between continua.
Under what conditions about $f$ and about $Y$ the mapping $f$ is universal?
In particular, is $f$ universal if $f$ satisfies some conditions related 
to confluence and $Y$ is
a dendrite
a dendroid?
a $\lambda$-dendroid?
a tree-like continuum having the fixed point property?
\end{myprob}

\begin{myprob}
\myproblem{G34}{J. Charatonik and W. Charatonik \cite{MR99k:54023}}
Does there exist an arcwise connected, unicoherent and one-dimensional
continuum $X$ and a confluent mapping from $X$ onto a locally connected
continuum $Y$ which is not weakly arc-preserving?

\mynote{Notes}
A mapping $f\colon X \to Y$ between continua is said to be 
\emph{arc-preserving} provided that it is surjective and for each arc 
$A \subset X$ its image $f(A)$ is either an arc or a point;
it is \emph{weakly arc-preserving} provided that there is an arcwise
connected subcontinuum $X'$ of $X$ such that the restriction 
$f|X'\colon X' \to Y$ is arc-preserving.
\end{myprob}

\begin{myprob}
\myproblem{G35}{J. Charatonik and W. Charatonik \cite{MR99k:54023}}
For what continua $X$ and $Y$
is each confluent mapping $f\colon X \to Y$ weakly arc-preserving?
For what continua $X$ and $Y$
is each weakly arc-preserving mapping $f\colon X \to Y$ weakly confluent?
\end{myprob}

\begin{myprob}
\myproblem{G36}{J. Charatonik and W. Charatonik \cite{MR99k:54023}}
Is every weakly arc-preserving mapping from a continuum onto a
dendroid universal?
\end{myprob}

\begin{myprob}
\myproblem{G37}{J. Charatonik and W. Charatonik \cite{MR99k:54023}}
Is any confluent mapping from a continuum (from a dendroid)
onto a dendroid universal?

\mynote{Notes}
A (metric) continuum $X$ is said to have the \emph{property of Kelley}
provided that for each point $x \in X$, for each subcontinuum $K$ of $X$
containing $X$, and for each sequence of points $x_n$ converging to $x$,
there exists a sequence of subcontinua $K_n$ of $X$ containing
$x_n$ and converging to $K$.

Let $K$ be a subcontinuum of a continuum $X$.
A continuum $M \subset K$ is called a \emph{maximal limit continuum in}
$K$ provided that there is a sequence of subcontinua $M_n$ of $X$
converging to $M$ such that for each converging sequence of
subcontinua $M'_n$ of $X$ with $M_n \subset M'_n$ for each
$n \in \mathbb{N}$ and $\lim M'_n = M' \subset K$,
we have $M = M'$.

A continuum is said to be \emph{semi-Kelley} provided that, for each
subcontinuum $K$ of $X$ and for every two maximal limit continua
$M_1$ and $M_2$ in $K$, either $M_1\subset M_2$ or $M_2\subset M_1$.

A mapping $f\colon X \to Y$ between continua is said to be 
\emph{semi-confluent} provided that, for each subcontinuum $Q$ of $Y$ 
and for every two components $C_1$ and $C_2$ of the inverse image 
$f^{-1}(Q)$, either $f(C_1) \subset f(C_2)$ or $f(C_2) \subset f(C_1)$.
\end{myprob}

\begin{myprob}
\myproblem{G38}{J. Charatonik and W. Charatonik \cite{MR2001f:54037}}
What classes of mappings preserve the property of being semi-Kelley?
In particular, is the property preserved under
monotone mappings?
open mappings?
\end{myprob}

\begin{myprob}
\myproblem{G39}{J. Charatonik and W. Charatonik \cite{MR2001f:54037}}
Is it true that if a continuum $Y$ has the property of Kelley and $X$ is 
an arbitrary continuum, then the uniform limit of semi-confluent mappings 
from $X$ onto $Y$ is semi-confluent?
\end{myprob}

\begin{myprob}
\myproblem{G40}{F. Jordan \cite{MR2002k:54014}}
Characterize the continua which are the almost continuous images of the
reals.
\end{myprob}

\section*{H. Homogeneity and mapping of general spaces}

\begin{myprob}
\myproblem{H1}{P. Nyikos}
Is there any reasonably large class of spaces $X$, $Y$ for which 
$\operatorname{ind} X \leq \operatorname{ind} Y + n$ 
when $f\colon X \to Y$ is a perfect mapping and 
$\operatorname{ind} f^{-1}(y) \leq n$ for all $y \in Y$? 
Does it even hold for all metric spaces?
Does it hold if $\operatorname{ind} Y = 0$? 
\end{myprob}

\begin{myprob}
\myproblem{H2}{L. Janos \cite{MR81e:54038}}
Let $(X,d)$ be a compact metric space of finite dimension and 
$f\colon X \to X$ an isometry of $X$ onto itself.
Does there exist a topological embedding $i\colon X \to E^m$ of $X$ into 
some Euclidean space $E^m$ such that $f$ is transformed into Euclidean 
motion?
This would mean that there exists a linear mapping $L\colon E^m \to E^m$ 
such that $L \circ i = i \circ f$.
\end{myprob}

\begin{myprob}
\myproblem{H3}{C.E.~Aull \cite{MR80m:54044}}
Are $\gamma$-spaces, quasi-metrizable spaces, or spaces with $\sigma$-Q
bases preserved under compact open maps?
What about spaces with orthobases?
\end{myprob}

\begin{myprob}
\myproblem{H4}{C.E.~Aull \cite{MR80m:54044}}
Are $\theta$-spaces or spaces with a $\delta\theta$-base preserved under
perfect mappings?

\mynote{Solution}
Yes. 
D. Burke \cite{MR85d:54009,MR85k:54014} proved that both classes are 
closed under perfect mappings.
\end{myprob}

\begin{myprob}
\myproblem{H5}{E. van~Douwen}
Does there exist a homogeneous zero-dimensional separable metrizable
space which
cannot be given the structure of a topological group or, more strongly,
has the fixed-point property for autohomeomorphisms?

\mynote{Solution}
Yes to the first (E. van~Douwen).
\end{myprob}

\begin{myprob}
\myproblem{H6}{E. van~Douwen}
Does there exist an infinite homogeneous compact zero-dimensional space
which has the fixed-point property for autohomeomorphisms?
\end{myprob}

\begin{myprob}
\myproblem{H7}{E. van~Douwen}
Does there exist a rigid zero-dimensional separable metrizable space
which is absolutely Borel, or at least analytic?
\end{myprob}

\begin{myprob}
\myproblem{H8}{B.J.~Ball and S. Yokura \cite{MR85b:54040}}
Let $X$ be the one-point compactification of a discrete space of
cardinality $\kappa$.
If $\kappa < \aleph_\omega$, there is a subset $F$ of $C(X)$ with 
$|F| \leq \kappa$ such that every element of $C(X)$ is the composition of
an element of $F$ followed by a map of $\mathbb{R}$ into $\mathbb{R}$.
Can the restriction $\kappa < \aleph_\omega$ be dropped?
\end{myprob}

\begin{myprob}
\myproblem{H9}{E. van~Douwen}
Characterize the spaces $X$ such that the projection map 
$\pi_1\colon X^2 \to X$ preserves Borel sets.
\end{myprob}

\begin{myprob}
\myproblem{H9}{E. van~Douwen}
For a linearly ordered set $L$ define an equivalence relation 
$T_L = \{ (x,y) \in L \times L : \text{there is an order-preserving
bijection of $L$ taking $x$ to $y$}\}$.
\begin{myenumerate}
\item
Does $\mathbb{R}$ have a subset $L$ such that $T_L$ has only one
equivalence class, but $(L, \leq)$ is not isomorphic to $(L, \geq)$?
\item
Does $\mathbb{R}$ have a subset $L$ such that $T_L$ has exactly two
equivalence classes, both dense [this much is possible] but of different
cardinalities?
\end{myenumerate}

\mynote{Solution}
Yes to the first part, (J. Baumgartner \cite{Baumgartner}). 
\end{myprob}

\begin{myprob}
\myproblem{H10}{E. van~Douwen}
One can show that a compact zero-dimensional space $X$ is the continuous
image of a compact orderable space if $X$ has a clopen family 
$\mathcal{S}$ which is $T_0$-point-separating (i.e., if $x \neq y$ then
there is $S \in \mathcal{S}$ such that $|S \cap \{x,y\}| = 1$)
and of rank $1$ (i.e., two members are either disjoint or comparable).
Is the converse false?

\mynote{Solution}
No, the converse is also true (S. Purisch \cite{MR98c:54021}).
L. Heindorff \cite{MR97j:06016} also answered this question in the 
context of Boolean interval algebras.
\end{myprob}

\begin{myprob}
\myproblem{H10}{E. van~Douwen}
Does every compact space without isolated points admit an irreducible map
onto a continuum?
Does $\omega^*$?
Yes to the second part, if \myaxiom{CH}.
\end{myprob}

\begin{myprob}
\myproblem{H11}{E. van~Douwen}
Is every compact $P^\prime$-space (i.e., nonempty $G_\delta$'s have
non\-empty interiors) an irreducible continuous image of a compact
zero-dimensional $P^\prime$-space (preferably of the same weight)?
\end{myprob}

\begin{myprob}
\myproblem{H12}{E. van~Douwen}
Let $\kappa > \omega$,
let $U(\kappa)$ denote the space of uniform ultrafilters on $\kappa$ and
let $A(U(\kappa))$ be the group of autohomeomorphisms of $U(\kappa)$.
Is every member of $A(U(\kappa))$ induced by a permutation of $\kappa$?
Is $A(U(\kappa))$ simple?
Is there, for every $h \in A(U(\kappa))$, a nonempty proper clopen subset
$V$ of $U(\kappa)$ with $h^{\to}V = V$?

\mynote{Notes}
Yes to the first part would imply yes to the other two parts.
Also, if yes to the third part, then $U(\kappa)$ and $\omega^*$ are not
homeomorphic.
\end{myprob}

\begin{myprob}
\myproblem{H13}{H. Kato \cite{MR89f:54030}}
Do refinable maps preserve countable dimension?
\end{myprob}

\begin{myprob}
\myproblem{H14}{A. Koyama \cite{MR89f:54030}}
Do $c$-refinable maps between nontrivial spaces preserve Property $C$?
\end{myprob}

\begin{myprob}
\myproblem{H15}{C.R.~Borges \cite{MR89f:54059}}
Let $\theta$ be a family of gages for a set $X$,
$\theta^{**}$ the gage for $X$ generated by $\theta$.
If $f\colon X \to X$ is $(\theta, \xi)$-expansive for some $\xi > 0$, 
is $f$ also $(\theta^{**}, \xi)$-expansive?
\end{myprob}

\begin{myprob}
\myproblem{H16}{C.R.~Borges \cite{MR89f:54059}}
Let ($X, U)$ be a sequentially compact (or countably compact or
pseudocompact) uniform space and $\theta$ a subgage for $U$.
If $f\colon X \to X$ is a continuous (w.r.t.\ the uniform topology)
$(\theta, \xi_0)$-expansive map for some $\xi_0 > 0$, is $f(X) = X$?
\end{myprob}

\begin{myprob}
\myproblem{H17}{T. Wilson \cite{MR94a:54102}}
Let $A$ be a compact metric space and let
$g\colon A \to A$ be a continuous surjection.
The sequence $S = \{x_n\}^\infty_{n=0} \subset A \times [0,1]$
is a \emph{generating sequence for $g$} if:
$A$ is the derived set of $S$;
the function $T_0$ defined by $T_0 x_n = x_{n+1}$ is continuous on $S$; 
and
$T_0$ has a continuous extension to $\operatorname{cl}(S)$ such that 
$T \upharpoonright A = g$.
[We are identifying $A \times \{0\}$ with $A$.]
When do generating sequences exist?
\end{myprob}

\begin{myprob}
\myproblem{H18}{T. Wilson \cite{MR94a:54102}}
Suppose $A$ is countable and $S$ is a generating sequence for
$g\colon A \to A$.
Let $A^\alpha$ denote the $\alpha^\text{th}$ derived set of $A$,
and let $\alpha_0$ be the least ordinal such that $A^{\alpha_0}$ is finite.
If $p$ is a fixed point of $g$, is $p \in A^{\alpha_0}$?
More generally, is $A^\alpha \subset g(A^\alpha)$?
\end{myprob}

\begin{myprob}
\myproblem{H29}{A.V.~Arhangel$'$\kern-.1667em ski\u{\i}, W. Just, and H. Wicke 
\cite{MR96k:54011}}
Is there a tri-quotient (compact) mapping which is not strongly blended?
Which is not blended?
\end{myprob}

\begin{myprob}
\myproblem{H30}{A.V.~Arhangel$'$\kern-.1667em ski\u{\i}, W. Just, and H. Wicke
\cite{MR96k:54011}}
Find topological properties other than submaximality and the $I$-space
property that are inherited by subspaces and are preserved by open
mappings and by closed mappings, but are not preserved in general by
pseudo-open mappings.
\end{myprob}

\begin{myprob}
\myproblem{H31}{A.V.~Arhangel$'$\kern-.1667em ski\u{\i} \cite{MR89b:54004}}
Let $X$ be an infinite homogeneous compactum.
Is there a nontrivial convergent sequence in $X$? 
What if we assume $X$ to be $2$-homogeneous? countable dense homogeneous?
\end{myprob}

\begin{myprob}
\myproblem{H32}{A.V.~Arhangel$'$\kern-.1667em ski\u{\i} \cite{MR89b:54004}}
Is there a homogeneous compactum of cellularity greater than $2^{\omega}$?
One that is $2$-homogeneous?

\mynote{Notes}
Negative answers would imply negative ones to the respective parts of the
following problem.
\end{myprob}

\begin{myprob}
\myproblem{H33}{A.V.~Arhangel$'$\kern-.1667em ski\u{\i} \cite{MR89b:54004}}
Can every compactum be represented as a continuous image of a homogeneous
compactum?
Of a $2$-homogeneous compactum?
\end{myprob}

\begin{myprob}
\myproblem{H34}{A.V.~Arhangel$'$\kern-.1667em ski\u{\i} \cite{MR89b:54004}}
Is every first countable compactum the continuous image of a first
countable homogeneous compactum?
Yes, if \myaxiom{CH}.
\end{myprob}

\begin{myprob}
\myproblem{H35}{A.V.~Arhangel$'$\kern-.1667em ski\u{\i} \cite{MR89b:54004}}
Is every separable space [resp.\ separable compactum] the continuous image 
of a countable dense homogeneous space [resp.\ compactum]?
\end{myprob}

\begin{myprob}
\myproblem{H36}{A.V.~Arhangel$'$\kern-.1667em ski\u{\i} \cite{MR89b:54004}}
If $Y$ is a zero-dimensional compactum, is there a compactum $X$ such that
$X \times Y$ is homogeneous? $2$-homogeneous?
\end{myprob}

\begin{myprob}
\myproblem{H37}{A.V.~Arhangel$'$\kern-.1667em ski\u{\i} \cite{MR89b:54004}}
If $Y$ is a Tychonoff space, is there a Tychonoff space $X$ such that
$X \times Y$ is $2$-homogeneous?
\end{myprob}

\begin{myprob}
\myproblem{H38}{A.V.~Arhangel$'$\kern-.1667em ski\u{\i} \cite{MR89b:54004}}
Let $Y$ be a compactum.
Is there a homogeneous compactum $X$ which contains an $l$-embedded
topological copy of $Y$? A $t$-embedded topological copy?
\end{myprob}

\begin{myprob}
\myproblem{H39}{D. Garity \cite{MR98d:54063}}
Is there a compact metric space of dimension less than $(n+2)$ that is
homogeneous, locally $n$-connected, and not $2$-homogeneous?
\end{myprob}

\begin{myprob}
\myproblem{H40}{D. Garity \cite{MR98d:54063}}
If a homogeneous compact metric space is locally $n$-connected for all $n$,
is the space necessarily $2$-homogeneous?
\end{myprob}

\section*{I. Infinite-dimensional topology}

\begin{myprob}
\myproblem{I1}{J. West \cite{MR58:24276}}
Let $G$ be a compact, connected Lie group acting on itself by left 
translation. 
Is ${2^G}/G$ a Hilbert cube? 
\end{myprob}

\begin{myprob}
\myproblem{I2}{J. West \cite{MR58:24276}}
Give conditions ensuring that, if $G$ is a compact Lie group acting on a
Peano continuum $X$, the induced G action on the Hilbert cube $2^X$ is
conjugate to some standard, such as the induced translative action on
$2^G$. 
\end{myprob}

\begin{myprob}
\myproblem{I3}{J. West \cite{MR58:24276}}
In general, given a compact Lie group, give conditions on $G$ actions on
manifolds, ANRs, Peano continua, or any other class of spaces which
ensure that the induced $G$ actions on hyperspaces are conjugate. 
\end{myprob}

\begin{myprob}
\myproblem{I4}{J. West \cite{MR58:24276}}
Let $G$ be a compact Lie group acting on a Peano continuum $X$ and
consider the injection of $X \to 2^X$ as the singletons. 
Then $G$ acts on $2^X$ and we can iterate the procedure, obtaining a
direct sequence $X \to 2^X \to 2^{2^X} \to \cdots$.
If we give $X$ a $G$-invariant convex metric then the inclusions are
isometries, and, moreover, the Hausdorff metric is both $G$-invariant and
convex.
Using the expansion homotopies $A \mapsto N_t(A)$, we see that $X$ is a
$Z$-set in $2^X$. 
If we now take the direct limit, we obtain a space which is homeomorphic
to separable Hilbert space equipped with the bounded-weak topology and has
an induced $G$ action on it. 
Identify this action directly in terms of $\ell_2$. 
\end{myprob}

\begin{myprob}
\myproblem{I5}{J. West \cite{MR58:24276}}
If, in the situation of Problem I4, we take the metric direct limit, we
have a separable metric space with a $G$ action on it.
Characterize this space and/or its completion in terms of more familiar
objects. in particular, are they homeomorphic to any well-known vector
spaces?
Once the above is done, characterize the induced $G$ action. 
\end{myprob}

\begin{myprob}
\myproblem{I6}{H. Hastings \cite{MR91g:55015}}
Is every (weak) shape equivalence of compact metric spaces a strong shape
equivalence?
\end{myprob}

\begin{myprob}
\myproblem{I7}{M. Jani \cite{MR85d:55026}}
Is there a cell-like shape fibration $p\colon E \to B$ from a compactum 
$E$ onto the dyadic solenoid $B$, which is not a shape equivalence?
\end{myprob}

\begin{myprob}
\myproblem{I8}{J.T.~Rogers \cite{MR85c:54055}}
Is any nondegenerate, homogeneous contractible continuum homeomorphic to the
Hilbert cube?
\end{myprob}

\begin{myprob}
\myproblem{I9}{H. Gladdines \cite{MR95c:54012}}
Let $L(\mathbb{R}^2)$ denote the collection of Peano continua in 
$\mathbb{R}^2$.
Is $L(\mathbb{R}^2)$ homeomorphic to the product of infinitely many
circles?
\end{myprob}

\section*{J. Group actions}

\begin{myprob}
\myproblem{J1}{R. Wong \cite{MR80k:57072}}
Every finite group $G$ can act on the Hilbert cube, $Q$, semi-freely with
unique fixed point, which we term based-free.
Let $G$ act on itself by left translation and extend this in the natural
way to the cone $C(G)$.
Let $Q_G$ (which is homeomorphic to $Q$) be the product of countably
infinitely many copies of $C(G)$.
The diagonal action $\sigma$ is based-free $G$-action on $Q$, and any
other based-free $G$-action on $Q_G$ is called standard if it is
topologically conjugate to $\sigma$. 
Does there exist a non-standard based-free $G$-action on $Q_G$?
\end{myprob}

\begin{myprob}
\myproblem{J2}{W. Lewis}
Does every zero-dimensional compact group act effectively on the
pseudo-arc?
\end{myprob}

\begin{myprob}
\myproblem{J3}{W. Lewis}
If a compact group acts effectively on a chainable (tree-like) continuum,
must it act effectively on the pseudo-arc?
\end{myprob}

\begin{myprob}
\myproblem{J4}{W. Lewis}
Under what conditions does a space $X$ with a continuous decomposition
into pseudo-arcs admit an effective $p$-adic Cantor group action which is 
an extension of an action on individual pseudo-arcs of the decomposition?
\end{myprob}

\begin{myprob}
\myproblem{J5}{Z. Balogh, J. Mashburn, and P. Nyikos \cite{MR91k:54039}}
Will a space $X$ freely acted upon by a finite group of autohomeomorphisms
necessarily have a countable closed migrant cover if it is subparacompact?
What if $X$ is a Moore space or a $\sigma$-space?

\mynote{Solution}
(Peter von Rosenberg \cite{rosenberg})
Yes, in the case of semi-stratifiable spaces, hence yes in the case of
Moore spaces or $\sigma$-spaces.
\end{myprob}

\begin{myprob}
\myproblem{J6}{Z. Balogh, J. Mashburn, and P. Nyikos \cite{MR91k:54039}}
If a paracompact space $X$ with finite Ind is acted upon freely by a
finite group of autohomeomorphisms, must $X$ have a finite open or
closed migrant cover?
What if $\dim X$ is finite?
If the answer to either one is affirmative, what is the optimal bound on
the size of the cover?

\mynote{Solution}
(Peter von Rosenberg \cite{rosenberg})
If $\dim X$ is finite, then $X$ does have finite open migrant covers and
hence finite migrant closed covers.
\end{myprob}

\section*{K. Connectedness}

\begin{myprob}
\myproblem{K1}{J.A.~Guthrie, H.E.~Stone, and M.L.~Wage \cite{TP2.1.349}}
What is the greatest separation which may be enjoyed by a maximally 
connected space?
In particular, is there a regular or semi-regular Hausdorff maximally
connected space?
\end{myprob}

\begin{myprob}
\myproblem{K2}{J.A.~Guthrie, H.E.~Stone, and M.L.~Wage \cite{TP2.1.349}}
For which $\kappa$ does there exist a maximally connected Hausdorff space
of cardinal $\kappa$?
In particular, is there a countable one?
\end{myprob}

\begin{myprob}
\myproblem{K3}{P. Nyikos}
Does there exist a weakly $\sigma$-discrete, connected, normal space?

\mynote{Notes}
A space is \emph{weakly $\sigma$-discrete} if it is the union of a 
sequence $X_n$ of discrete subsets so that $\bigcup_{i < n} X_i$ is 
closed for each $n$.
\end{myprob}

\begin{myprob}
\myproblem{K4}{P. Zenor}
Does \manotch\ imply that there is no locally connected, rim-compact 
$L$-space?

\mynote{Solution}
Solved in the affirmative by G. Gruenhage.
\end{myprob}

\begin{myprob}
\myproblem{K5}{P. Collins}
Is a locally compact, $\sigma$-compact connected and locally connected
space always the union of a countable sequence of compact, connected,
locally connected subsets such that $C_i \subset \operatorname{int}C_{i+1}$ 
for all $i$?
\end{myprob}

\begin{myprob}
\myproblem{K6}{P. Nyikos, attributed to M.E.~Rudin}
Does \manotch\ imply every compact, perfectly normal, 
locally connected space is metrizable?
\end{myprob}

\begin{myprob}
\myproblem{K7}{E. van~Douwen}
Is there a connected (completely) regular space without disjoint dense
subsets?
There are Hausdorff examples.
\end{myprob}

\begin{myprob}
\myproblem{K8}{P.A.~Cairns \cite{MR96k:54057}}
Is there any space of transfinite cohesion?
\end{myprob}

\section*{L. Topological algebra}

\begin{myprob}
\myproblem{L1}{M. Henriksen \cite{MR82i:54035}}
Find a necessary and sufficient condition on a realcompact, rim-compact
space $X$ in order that $C^{\#}(X)$ will determine a compactification (and
hence the Freudenthal compactification) of $X$.
To do so, it will probably be necessary to characterize the zero-sets of
elements of $C^{\#}(X)$.
\end{myprob}

\begin{myprob}
\myproblem{L2}{D.L.~Grant \cite{MR82i:22001}}
If every finite power of a group is minimal (or totally minimal, or a
$B(A)$ group), must arbitrary powers of the group have the same property?
\end{myprob}

\begin{myprob}
\myproblem{L3}{R.A.~McCoy \cite{MR82j:54021}}
Let $X$ be a completely regular $k$-space.
If $C(X)$ with the compact-open topology is a $k$-space, must $X$ be 
hemicompact?
This would imply that $C(X)$ is completely metrizable.
\end{myprob}

\begin{myprob}
\myproblem{L4}{E. van~Douwen}
Must every locally compact Hausdorff topological group contain a dyadic
neighborhood of the identity? 

\mynote{Solution}
van~Douwen asked this question in 1986 but it had been answered long 
before by B. Pasynkov and M. Choban who proved (independently) that any 
compact $G_\delta$ subset of any topological group (not necessarily 
locally compact) is dyadic. Pasynkov never published a proof and Choban's 
proof appeared in a conference proceedings (in Russian) that were hardly 
available.
See \cite{MR89i:22005} by V. Uspenskij for a proof of this theorem.
A strengthening of the Choban-Pasynkov theorem is in \cite{MR91a:54064}.
\end{myprob}

\begin{myprob}
\myproblem{L5}{E. van~Douwen}
A quasi-group is a set $G$ with three binary operations
$\cdot$, $/$ and $\setminus$ 
such that $a/b$ and $b \setminus a$ are the unique solutions to 
$x \cdot b = a$ and $b \cdot x = a$ for all $a, b \in G$.
A topological quasi-group is a quasi-group with a topology with respect
to which these operations are jointly continuous.
\begin{myenumerate}
\item
Is there a (preferably compact) zero-dimensional topological quasi-group
whose underlying set cannot be that of a topological group?
The quasi-group of Cayley numbers of value $1$ ($S^7$) is a well-known
connected example.
\item
Is there a quasi-group which is also a (preferably compact) space such
that the $\cdot$ is jointly continuous, $/$ and $\setminus$ are separately
continuous, but are not jointly continuous?
\item
Is there a quasi-group which is a (preferably compact) space as in (2)
but whose underlying set cannot be that of a semigroup?
\end{myenumerate}
\end{myprob}

\begin{myprob}
\myproblem{L6}{E. van~Douwen}
If $\{G_i : i \in I\}$ is a collection of topological groups, the
coproduct-topology on the weak product 
$$\textstyle
\sum_{i \in I} G_i = \{x \in \prod_{i \in I} G_i : x_i = e_i\ 
\text{for all but finitely many $i$}\}$$
is defined to be the finest group topology such that the relative
topology on each finite subproduct is the product topology.
Is it possible to have families $\{G_i : i \in I\}$ and
$\{H_i : i \in I\}$ of (preferably abelian) topological groups such 
that $g_i$ and $H_i$ have the same underlying space and the same 
underlying identity for all $i \in I$, yet the coproduct topologies on the 
respective weak products are unequal? nonhomeomorphic?
\end{myprob}

\begin{myprob}
\myproblem{L7}{A.V.~Arhangel$'$\kern-.1667em ski\u{\i}}
Let $F(X)$ denote the free topological group on the space $X$. 
If $\dim \beta F(X) = 0$, does $\dim F(x) = 0$ follow?
\end{myprob}

\begin{myprob}
\myproblem{L8}{D. Shakhmatov \cite{MR94d:54010}}
Let $G$ be a countably compact (Hausdorff) topological group.
Is then $t(G \times G) = t(G)$ (here $t$ denotes tightness)?
What if $t(G)$ is countable?

\mynote{Notes}
(D. Shakhmatov)
If countable compactness is dropped, there are counterexamples under
various set-theoretic hypotheses as shown by: 
V. Malykhin under \myaxiom{CH};
Malykhin and Shakhmatov in the model obtained by adding one Cohen real to
a model of \manotch, and;
Shakhmatov in the model obtained by first adding $\omega_2$ Cohen reals 
to a model of \myaxiom{GCH} then using the Martin-Solovay poset for 
obtaining \manotch.
In the last case, examples were found of dense pseudocompact subgroups 
$G$ of $2^{\omega_1}$ for which $G^n$ is hereditarily separable and 
Fr\'echet-Urysohn but $G^{n+l}$ has uncountable tightness.
\end{myprob}

\begin{myprob}
\myproblem{L9}{D. Shakhmatov \cite{MR94d:54010}}
Let $G$ be a countably compact Fr\'echet-Urysohn topological group.
Is $G \times G$ Fr\'echet-Urysohn? 
Is $G^n$ Fr\'echet-Urysohn for all $n$?

\mynote{Notes}
A counterexample could not be $\alpha_3$ since the product of a
Fr\'echet-Urysohn $\alpha_3$ space and a countably compact
Fr\'echet-Urysohn space is Fr\'echet-Urysohn. 
See \cite{MR81i:54017}.

\mynote{Solution}
No. 
Using \myaxiom{CH}, A. Shibakov constructed a countable Fr\'echet-Urysohn 
group whose square is not Fr\'echet-Urysohn \cite{MR2000a:54062}.
\end{myprob}

\begin{myprob}
\myproblem{L10}{D. Shakhmatov \cite{MR94d:54010}}
Let $G$ be a topological group so that $G^n$ is Fr\'echet-Urysohn for
every natural number $n$.
Is $G^\omega$ Fr\'echet-Urysohn?
What if one assumes also that $G$ is countably compact?

\mynote{Notes}
A counterexample could not be $\alpha_3$
(T. Nogura \cite[Corollary~3.8]{MR87a:54036}).
\end{myprob}

\begin{myprob}
\myproblem{L11}{D. Shakhmatov \cite{MR94d:54010}}
Is every countably compact sequential topological group Fr\'echet-Urysohn?

\mynote{Notes}
An affirmative answer to L11 would imply affirmative answers for L9
and the second part of L10: the product of a countably compact
sequential space and a sequential space is sequential.
See also the background references for A22.
\end{myprob}

\begin{myprob}
\myproblem{L12}{D. Shakhmatov \cite{MR94d:54010}}
Is there a (countable) Fr\'echet-Urysohn group which is an $\alpha_3$-space
without being an $\alpha_2$-space?

\mynote{Solution}
Using \myaxiom{CH}, A. Shibakov constructed a Fr\'echet-Urysohn group that
satisfies the $\alpha_3$-property but not the $\alpha_2$-property
\cite{MR2000a:54062}.
\end{myprob}

\begin{myprob}
\myproblem{L13}{D. Shakhmatov \cite{MR94d:54010}}
Is it consistent with \myaxiom{ZFC} to have a Fr\'echet-Urysohn 
$\alpha_{1.5}$-group which is not a v-group?
\end{myprob}

\begin{myprob}
\myproblem{L14}{D. Shakhmatov \cite{MR94d:54010}}
Is there a real (requiring no additional set-theor\-et\-ic assumptions
beyond \myaxiom{ZFC}) example of a countable nonmetrizable w-group?
\end{myprob}

\begin{myprob}
\myproblem{L15}{D. Shakhmatov \cite{MR94d:54010}}
Is there a real example of a Fr\'echet-Urysohn topological group that is
not an $\alpha_3$-space?
\end{myprob}

\begin{myprob}
\myproblem{L16}{D. Shakhmatov \cite{MR94d:54010}}
Do the convergence properties $\alpha^i$ ($i = 0, 1, \dots, \infty$)
coincide for Fr\'echet-Urysohn topological groups?

\mynote{Solution}
See the results by A. Shibakov described in L12 and L18.
\end{myprob}

\begin{myprob}
\myproblem{L17}{D. Shakhmatov \cite{MR94d:54010}}
Is every Fr\'echet-Urysohn group an $\alpha^\infty$-space? 
\end{myprob}

\begin{myprob}
\myproblem{L18}{D. Shakhmatov \cite{MR94d:54010}}
Do some new implications between $\alpha_i$-properties,
$i \in \{1, 1.5, 2, 3, 4\}$, and $\alpha^k$-properties,
$k \in \omega \cup \{\infty\}$,
appear in Fr\'echet-Urysohn groups belonging to one of the following classes:
(1) countably compact groups,
(2) pseudocompact groups,
(3) precompact groups (= subgroups of compact groups), and
(4) groups complete in their two-sided uniformity?

\mynote{Solution}
A. Shibakov gave a simple proof that $\alpha_1$ and $\alpha_{1.5}$ are
equivalent for Fr\'echet-Urysohn groups \cite{MR2000a:54062}.
See L12.
\end{myprob}

\begin{myprob}
\myproblem{L19}{D. Shakhmatov \cite{MR94d:54010}}
Is every Fr\'echet-Urysohn group having a base of open neighborhoods
of its neutral element consisting of subgroups a w-space?
\end{myprob}

\begin{myprob}
\myproblem{L20}{M. Tka\v{c}enko \cite{MR94c:22002}}
Is every c.c.c.\ topological group $\mathbb{R}$-factorizable?
What if it is separable?

\mynote{Notes}
Recall that a group $G$ is said to be \emph{$\mathbb{R}$-factorizable} 
if for any continuous real-valued function $f$ on $G$ there exist a 
continuous homomorphism $\pi$ of $G$ onto a group $H$ of countable weight 
and a continuous function $h$ on $H$ such that $f = h \circ \pi$.
\end{myprob}

\begin{myprob}
\myproblem{L21}{M. Tka\v{c}enko \cite{MR94c:22002}}
Let $S$ be the Sorgenfrey line and $A(S)$ the free abelian topological
group over $S$.
Is $A(S)$ $\mathbb{R}$-factorizable?
This is a very special case of L20.
\end{myprob}

\begin{myprob}
\myproblem{L22}{M. Tka\v{c}enko \cite{MR94c:22002}}
Let $g$ be a continuous real-valued function on an $\aleph_0$-bounded
group $G$.
Are there a continuous homomorphism $\pi$ of $G$ onto a group $H$ of
weight at most $2^{\aleph_0}$ and a continuous function $h$ on $H$ such
that $g = h \circ \pi$?

\mynote{Notes}
Not every $\aleph_0$-bounded group is $\mathbb{R}$-factorizable;
but Tka\v{c}enko conjectures that the above weakening of 
$\mathbb{R}$-factorizability holds for it.
\end{myprob}

\begin{myprob}
\myproblem{L23}{M. Tka\v{c}enko \cite{MR94c:22002}}
Is every subgroup of $\mathbb{Z}^\tau$ $\mathbb{R}$-factorizable?

\mynote{Notes}
Every subgroup of $\mathbb{Z}^\tau$, for each $\tau$, is
$\aleph_0$-bounded but, by a result of V. Uspenskij, is not necessarily 
c.c.c.\ \cite{MR97f:22002}.
\end{myprob}

\begin{myprob}
\myproblem{L24}{M. Tka\v{c}enko \cite{MR94c:22002}}
Must every locally finite family of open subsets of an
$\mathbb{R}$-factorizable group be countable?
Is every $\mathbb{R}$-factorizable group $G$ weakly Lindel\"of?
That is, is it true that every open cover of $G$ has a countable
subfamilly a union of which is dense in $G$?
\end{myprob}

\begin{myprob}
\myproblem{L25}{M. Tka\v{c}enko \cite{MR94c:22002}}
Does a continuous homomorphic image of an $\mathbb{R}$-factoriz\-able
group inherit the $\mathbb{R}$-factorization property?
If the homomorphism is open as well, then yes
\cite[Theorem~3.1]{MR94c:22002}. 
\end{myprob}

\begin{myprob}
\myproblem{L26}{M. Tka\v{c}enko \cite{MR94c:22002}}
Is every $\aleph_0$-bounded group a continuous image of an 
$\mathbb{R}$-factorizable group?
Yes to L26 implies no to L25.
\end{myprob}

\begin{myprob}
\myproblem{L27}{M. Tka\v{c}enko \cite{MR94c:22002}}
Is the $\mathbb{R}$-factorization property inherited by finite products?
\end{myprob}

\begin{myprob}
\myproblem{L28}{M. Tka\v{c}enko \cite{MR94c:22002}}
Is the product of an $\mathbb{R}$-factorizable group with a compact
group $\mathbb{R}$-factorizable?

\mynote{Notes}
This is a special case of L27.
An affirmative answer to the first part of Problem L24 would imply 
an affirmative answer to Problem L28.
\end{myprob}

\begin{myprob}
\myproblem{L29}{M. Tka\v{c}enko \cite{MR94c:22002}}
Suppose $G$ is an $\mathbb{R}$-factorizable group of countable
$o$-tightness and $K$ is a compact group.
Is the product $G \times K$ $\mathbb{R}$-factorizable?
What if $G$ is a $k$-group?
\end{myprob}

\begin{myprob}
\myproblem{L30}{M. Tka\v{c}enko \cite{MR94c:22002}}
Must the product of a Lindel\"of group with a totally bounded group
be $\mathbb{R}$-factorizable?

\mynote{Notes}
We may assume without loss of generality that the totally bounded factor
is second countable.
\end{myprob}

\begin{myprob}
\myproblem{L31}{M. Tka\v{c}enko \cite{MR94c:22002}}
Is it true that the closure of a $G_{\delta, \Sigma}$ set in a $k$-group
$H$ is a $G_\delta$-set?
What if $H$ is sequential or Fr\'echet-Urysohn?
\end{myprob}

\begin{myprob}
\myproblem{L32}{D. Dikranjan and D. Shakhmatov \cite{MR95a:54058}}
Which infinite groups admit a pseudocompact topology?
In other words, what special algebraic properties must pseudocompact
groups have?
\end{myprob}

\begin{myprob}
\myproblem{L33}{D. Dikranjan and D. Shakhmatov \cite{MR95a:54058}}
If $G$ is a pseudocompact abelian group, must either the torsion subgroup
$t(G) = 
\{g \in G : ng = 0\ \text{for some}\ n \in \mathbb{N} \setminus 
\{0\}\}$
or $G/t(G)$ admit a pseudocompact group topology?
\end{myprob}

\begin{myprob}
\myproblem{L34}{D. Dikranjan and D. Shakhmatov \cite{MR95a:54058}}
If an abelian group $G$ admits a pseudocompact group topology,
must the group $G/t(G)$ admit a pseudocompact group topology?

\mynote{Notes}
(D. Dikranjan and D. Shakhmatov)
The answer to this and the preceding question is affirmative for
torsion and torsion-free groups (both trivially), for divisible groups,
for groups with $|G| = r(G)$, where $r(G)$ is the free rank of $G$,
and when $t(G)$ admits a pseudocompact topology or is bounded, 
i.e., there is some $n \in \mathbb{N} \setminus \{0\}$ such that 
$ng = 0$ for all $g \in G$.
\end{myprob}

\begin{myprob}
\myproblem{L35}{D. Dikranjan and D. Shakhmatov \cite{MR95a:54058}}
Suppose that $G$ is an abelian group,
$n \in \mathbb{N} \setminus \{0\}$ and both
$nG = \{ng : g \in G\}$ and $G/nG$
admit pseudocompact group topologies.
Must then $G$ also admit a pseudocompact topology?
\end{myprob}

\begin{myprob}
\myproblem{L36}{D. Dikranjan and D. Shakhmatov \cite{MR95a:54058}}
Let $D(G)$ denote the maximal divisible subgroup of an abelian group $G$.
If $G$ is pseudocompact, must either $D(G)$ or $G/D(G)$ admit a
pseudocompact topology?
\end{myprob}

\begin{myprob}
\myproblem{L37}{D. Dikranjan and D. Shakhmatov \cite{MR95a:54058}}
Let $G$ be an abelian group with $D(G) = \{0\}$,
i.e., a reduced abelian group. 
If $G$ admits a pseudocompact group topology,
must $G$ admit also a zero-dimensional pseudocompact group topology?
\end{myprob}

\begin{myprob}
\myproblem{L38}{D. Dikranjan and D. Shakhmatov \cite{MR95a:54058}}
Let $G$ be a non-torsion pseudocompact abelian group.
Do there exist a cardinal $\sigma$ and a subset of cardinality
$r(G)$ of $\{0,1\}^\sigma$ whose projection on every countable subproduct
is a surjection?
\end{myprob}

\begin{myprob}
\myproblem{L39}{D. Dikranjan and D. Shakhmatov \cite{MR95a:54058}}
Characterize (abelian) groups which admit a group topology which has one
of the following properties:
countably compact,
$\sigma$-compact, or
Lindel\"of.
\end{myprob}

\begin{myprob}
\myproblem{L40}{D. Dikranjan and D. Shakhmatov \cite{MR95a:54058}}
For which cardinals $\tau$ does the free abelian group with $\tau$ 
generators admit a countably compact group topology?
\end{myprob}

\begin{myprob}
\myproblem{L41}{H. Teng \cite{MR95a:46037}}
Let $X$ be a fortissimo space and $p$ the particular point of $X$.
Let $Y = X \setminus \{p\}$.
Is $C_p(Y|X)$ normal?
\end{myprob}

\begin{myprob}
\myproblem{L42}{H. Teng \cite{MR95a:46037}}
With $X$ and $Y$ as in L41, is $C_p(Y|X)$ homeomorphic to the
$\Sigma$-product of $|X|$-many real lines?
\end{myprob}

\begin{myprob}
\myproblem{L43}{J. Covington \cite{MR96m:54070}}
If $(G,t)$ is a protopological ($t$-protopological) group and $A$ and $B$
are connected (compact) subsets containing the identity, is $AB$ connected
(compact)?
\end{myprob}

\begin{myprob}
\myproblem{L44}{V. Bergelson, N. Hindman, and R. McCutcheon 
\cite{MR2001a:20114}}
In a group, if $A$ and $B$ are both right syndetic, does it follow that
$AA^{-1} \cap BB^{-1}$ necessarily contains more than the identity?
\end{myprob}

\begin{myprob}
\myproblem{L45}{V. Bergelson, N. Hindman, and R. McCutcheon
\cite{MR2001a:20114}}
If $m_l(B) > 0$ in a left amenable semigroup, and $A$ is infinite, does
$BB^{-1} \cap AA^{-1}$ necessarily contain elements different from the
identity?
\end{myprob}

\begin{myprob}
\myproblem{L46}{A. Giarlotta, V. Pata, and P. Ursino \cite{MR2003f:28035}}
Are $\mathcal{S}_\infty$ and $\mathcal{A}$ comparable as groups?
That is, does there exist an emdedding of one of them into the other one?

\mynote{Notes}
$\mathcal{A}$ is the group of measure-preserving bijections of $[0,1)$.
$\mathcal{S}_\infty$ is the group of permutations of $\mathbb{N}$.
It is known that if $\text{\myaxiom{MA}}(k)$ is assumed, and $A$ is a 
Boolean algebra with infinitely many atoms such that $|A| = k$, then 
$\mathcal{S}_\infty$ can be isomorphically embedded in 
$\operatorname{Aut}(A)$.
The authors note that above question is probably very difficult, yet the 
following weaker version of it seems to be very interesting as well. 
\end{myprob}

\begin{myprob}
\myproblem{L47}{A. Giarlotta, V. Pata, and P. Ursino \cite{MR2003f:28035}}
Are $\mathcal{S}_\infty$ and $\mathcal{A}$ comparable as subgroups of
$\text{Aut}(\mathcal{P}(\mathbb{N})/\text{fin})$?

\mynote{Notes}
(A. Giarlotta, V. Pata, and P. Ursino)
This question makes sense, since both groups can be isomorphically 
embedded into $\mathcal{P}(\mathbb{N})/\text{fin}$, as is proved in the 
paper \cite{MR2003g:06018}.
\end{myprob}

\begin{myprob}
\myproblem{L48}{A. Giarlotta, V. Pata, and P. Ursino \cite{MR2003f:28035}}
Is there a formula that gives the order of a (particular) element 
$\gamma$ in $(\mathcal{S},\circ)$ in terms of the parameters of the shifts 
of which $\gamma$
is the composition?

\mynote{Notes}
A partial answer for the composition of two rational shifts has been found
by the authors \cite{MR2003f:28035}. 
\end{myprob}

\section*{M. Manifolds}

\begin{myprob}
\myproblem{M1}{W. Kuperberg \cite{TP2.1.355}}
Is it true that the orientable closed surfaces of positive genus are the
only closed surfaces embeddable in the products of two one-dimensional
spaces?
\end{myprob}

\begin{myprob}
\myproblem{M2}{W. Kuperberg \cite{TP2.1.355}}
Suppose that $T$ is a torus surface contained in the product of 
$X \times Y$ of two one-dimensional spaces $X$ and $Y$.
Do there exist two simple closed curves $A \subset X$ and $B \subset Y$
such that $T = A \times B$?
In other words, if $\pi_1$ and $\pi_2$ are the projections, is
$T = \pi_1(T) \times \pi_2(T)$?
Here $=$ always denotes set-theoretic equality.
\end{myprob}

\begin{myprob}
\myproblem{M3}{M.E.~Rudin \cite{MR80j:54014}}
Is there a complex analytic, perfectly normal, nonmetrizable manifold?
No, if \manotch.
\end{myprob}

\begin{myprob}
\myproblem{M4}{P. Nyikos \cite{MR85i:54004,MR86f:54054}}
Is every normal manifold collectionwise normal?

\mynote{Notes}
Yes, if $\mathfrak{c}$\myaxiom{MEA} (P. Nyikos \cite{MR85i:54004}).
It is not known whether the consistency of yes requires the consistency of
an inaccessible cardinal. 
See M6. 
\end{myprob}

\begin{myprob}
\myproblem{M5}{W. Haver \cite{MR83m:57027}}
Do there exist version of the standard engulfing theorems in which
the engulfing isotopy depends continuously on the given open sets and
embeddings?
\end{myprob}

\begin{myprob}
\myproblem{M6}{F.D.~Tall \cite{MR85h:54036}}
Can one prove the consistency of ``normal manifolds are collectionwise
normal'' without assuming any large cardinal axioms?
\end{myprob}

\begin{myprob}
\myproblem{M7}{P. Latiolais \cite{MR87k:57017}}
Does there exist a pair of finite $2$-dimensional CW-complexes which are
homotopy equivalent but not simple homotopy equivalent?

\mynote{Solution}
Such examples do exist.
There were provided independently by M.~Lustig \cite{MR92a:57025}
and W.~Metzler \cite{MR91c:57028}.
There is a version of Lustig's examples in \cite[\S~VII]{MR95g:57006}.
\end{myprob}

\begin{myprob}
\myproblem{M8}{P. Latiolais \cite{MR87k:57017}}
Does there exist a finite $2$-dimensional CW-complex $K$ whose fundamental
group is finite but not abelian, which is not simple homotopy equivalent
to every $n$-dimensional complex homotopy equivalent to $K$?

\mynote{Notes}
The examples of Metzler and Lustig for problem M8 have infinite
fundamental group, so they did not answer this question. 
\end{myprob}

\begin{myprob}
\myproblem{M9}{P. Latiolais \cite{MR87k:57017}}
Do the Whitehead torsions realized by self-equi\-val\-ences of a finite
$2$-dimensional CW-complex include all of the units of the Whitehead group
of that complex?
\end{myprob}

\begin{myprob}
\myproblem{M10}{P. Nyikos}
Is every normal, or countably paracompact, manifold collectionwise
Hausdorff?
Yes, if \visl\ or $\mathfrak{c}$\myaxiom{MEA}.
Is there a model of $\text{\myaxiom{MA}}(\omega_1)$ where the answer is 
yes?
\end{myprob}

\begin{myprob}
\myproblem{M11}{C. Good \cite{MR98d:54041}}
Is there a hereditarily normal Dowker manifold?

\mynote{Solution}
D. Gauld and P. Nyikos \cite{gauldnyikos} proved that $\diamondsuit$ 
implies that there is a hereditarily normal Dowker manifold.
\end{myprob}

\begin{myprob}
\myproblem{M12}{B. Brechner and J.S.~Lee \cite{MR98d:54047}}
Characterize those bounded domains $U$ in $E^3$ which admit a prime end
structure.
\end{myprob}

\begin{myprob}
\myproblem{M13}{B. Brechner and J.S.~Lee \cite{MR98d:54047}}
Characterize those bounded domains $U$ in $E^3$ which admit a
$C$-transformation onto the interior of some compact $3$-manifold.
\end{myprob}

\section*{N. Measure and topology}

\begin{myprob}
\myproblem{N1}{W.F.~Pfeffer \cite{MR82b:28025}}
Let $\alpha \leq \gamma$ and let $\mu$ be a diffused $\gamma$-regular
$\alpha$-measure on a $T_1$ space $X$.
Is $\mu$ moderated?
\end{myprob}

\begin{myprob}
\myproblem{N2}{W.F.~Pfeffer \cite{MR82b:28025}}
Let $\alpha > \beta$ and let $\mu$ be a $\beta$-finite Borel
$\alpha$-measure on a metacompact space $X$ containing no closed discrete
subspace of measurable cardinality.
Is $\mu$ $\beta$-moderated?
\end{myprob}

\begin{myprob}
\myproblem{N3}{J. Stepr\={a}s}
Is there a measure zero subset $X$ of $\mathbb{R}$ such that any measure
zero subset of $\mathbb{R}$ is contained in some translate of $X$?
in the union of countably many translates of $X$?

\mynote{Solution}
No. S. Todor\v{c}evi\'c, F. Galvin, and D. Fremlin have independently
given general theorems which imply that the answer is negative. 
\end{myprob}

\section*{O. Theory of retracts; extensions of continuous functions}

\begin{myprob}
\myproblem{O1}{R. Wong \cite{MR80k:57072}}
An absolute retract (AR) $M$ is said to be \emph{pointed} at a point 
$p \in M$ if there is a strong deformation retract $\{\lambda_t\}$ of $M$ 
onto $\{p\}$ such that $\lambda_{t}^{-1} = p$ for all $t < 1$.
It is known that a point $p$ in a compact AR is pointed if 
$M \setminus \{p\}$ has the homotopy type of an Eilenberg-MacLane space of
type $(\mathbb{Z}_n, 1)$ where $\mathbb{Z}_n = \mathbb{Z} / (n)$.
Can we relax the condition on $M \setminus \{p\}$?
In particular, is every point pointed in a compact AR?
\end{myprob}

\begin{myprob}
\myproblem{O2}{L.I.~Sennott}
Is the converse of the theorem below \cite[Theorem 1]{MR80k:54025}
true? 
If not, is there a counterexample where $S$ is $P$-embedded?
\begin{theorem}
Let $(X,S)$ have the property that every function from $S$ to a 
complete locally convex vector space extends to $X$.
Then there exists an order preserving extender from $\mathcal{P}^*(S)$ to
$\mathcal{P}^*(X)$.
\end{theorem}

\mynote{Notes}
$\mathcal{P}^*(X)$ is the collection of all bounded pseudometrics on $X$.
\emph{$S$ is $P$-embedded} means that every pseudometric on $S$ can be 
extended to a pseudometric on $X$.
\end{myprob}

\begin{myprob}
\myproblem{O3}{L.I.~Sennott}
From \cite[Theorem 2]{MR80k:54025} it is clear that if $S$ is $D$-embedded
in $X$, then there exists a s.l.e.\ from $C_\mu(S)$ to $C_\mu(X)$ and
from $C_p(S)$ to $C_p(X)$.
Must there exist a s.l.e.\ from $C_c(S)$ to $C_c(X)$?

\mynote{Notes}
A \emph{simultaneous linear extender} (s.l.e.) from 
$C(S,L)$ to $C(X,L)$ is a linear function $\Psi\colon C(S,L) \to C(X,L)$ 
such that $\Psi(f)|S = f$ for all $f \in C(S,L)$. 
\end{myprob}

\begin{myprob}
\myproblem{O4}{L.I.~Sennott}
Give characterizations (similar to those known for $P$- and
$M$-embedding) for the other embeddings introduced in 
\cite[\S~2]{MR80k:54025}.
\end{myprob}

\begin{myprob}
\myproblem{O5}{G. Gruenhage, G. Kozlowski, and P. Nyikos}
A compact space is an absolute retract (AR) if it is a retract of every
compact (equivalently, Tychonoff) space in which it is embedded, and a
Boolean absolute retract (BAR) if it is zero-dimensional, and a retract of
every zero-dimensional (compact) space in which it is embedded.
Is a nonmetrizable AR homeomorphic to $I^\kappa$ for some $\kappa$?

\mynote{Solution}
No, E. Shchepin. 
The cone over $I^\kappa$ is also an AR and is not homeomorphic to
$I^\kappa$ for $\kappa > \omega$. 
\end{myprob}

\begin{myprob}
\myproblem{O6}{P. Nyikos}
If X is a BAR, does there exist a BAR $Y$ such that 
$X \times Y \approx 2^\kappa$ for some $\kappa$?
Is $X^\kappa \approx 2^\kappa$ for large enough $\kappa$?
Here $\approx$ denotes homeomorphism.

\mynote{Solution}
Yes, E. Shchepin.
\end{myprob}

\begin{myprob}
\myproblem{O7}{P. Nyikos}
Does there exist an intrinsic characterization of BARs either among
compact spaces or among dyadic spaces?
This is an old question.

\mynote{Solution}
(E. Shchepin)
Among dyadic spaces of weight $\leq \aleph_1$, the BARs are
characterized by the Bockstein separation property: 
disjoint open sets are contained in disjoint cozero sets. 
However, this no longer holds for BARs of higher weight.
\end{myprob}

\begin{myprob}
\myproblem{O8}{L.I.~Sennott \cite{MR80k:54049}}
If $S$ is a closed subspace of a normal space $X$ such that $(X,S)$ has 
the $\gamma$-ZIP, must $S \times Y$ be C-embedded in $X \times Y$ for 
every metric space $Y$ of weight $\leq |S|$?
See \cite[Theorem~3.1, p.~511]{MR80k:54049}.

\mynote{Notes}
A space $X$ has the \emph{$\gamma$-zero-set interpolation property}
($\gamma$-ZIP) if whenever $d$ is a $\gamma$-separable pseudometric on 
$X$ there exists a zero set $Z$ of $X$ such that 
$S 
\subset Z 
\subset \{x \in X; \text{$d(x, x_0) = 0$ for some $x_0 \in S$}\}$.
\end{myprob}

\begin{myprob}
\myproblem{O9}{L.I.~Sennott \cite{MR80k:54049}}
Characterize metric spaces $Y$ such that if $X$ is a topological space and
$S$ is C-embedded in $X$, then $S \times Y$ is C-embedded in $X \times Y$. 
Is the space of rational numbers in this class?
\end{myprob}

\begin{myprob}
\myproblem{O10}{A. Koyama \cite{MR89f:54030}}
Let $r\colon X \to Y$ be a refinable map and let $K$ be a class of ANRs.
If $Y$ is extendable with respect to $K$, then is $X$ also extendable with
respect to $K$?
\end{myprob}

\begin{myprob}
\myproblem{O11}{H. Kato \cite{MR89f:54030}}
Do refinable maps preserve FANRs?
\end{myprob}

\begin{myprob}
\myproblem{O12}{R. Levy}
Is there a \myaxiom{ZFC} example of a metric space having a subset that is
$2$-embedded (i.e., every continuous function into a two-point discrete
space has an extension to a continuous function on the whole space) but
not $C^*$-embedded?

\mynote{Solution}
The answer is trivially, yes.
For example, the open unit interval is trivially $2$-embedded in the
closed unit interval.
This question was a bad transcription of a question of R. Levy, whose
correct statement is O13. 
\end{myprob}

\begin{myprob}
\myproblem{O13}{R. Levy \cite{MR91f:54009}}
Is there a \myaxiom{ZFC} example of a metric space having a subset that is
$2$-embedded but not $\omega$-embedded?

\mynote{Notes}
A subset $S$ of a space $X$ is \emph{$\kappa$-embedded} if every
continuous function from $S$ into a discrete space of cardinality $\kappa$
has an extension to a continuous function on $X$.
There cannot be a separable example;
it is consistent that there are nonseparable examples
\cite{MR91f:54009}.
\end{myprob}

\begin{myprob}
\myproblem{O14}{L. Friedler, M. Girou, D. Pettey, and J. Porter 
\cite{MR94m:54058}}
Let $Y$ be an $R$-closed [resp.\ $U$-closed] extension of a space $X$ and
$f$ a continuous function from $X$ to an $R$-closed [resp.\ $U$-closed]
space $Z$.
Find necessary and sufficient conditions that $f$ can be extended to a
continuous function from $Y$ to $Z$.
\end{myprob}

\begin{myprob}
\myproblem{O15}{R. Pawlak \cite{MR94k:54030}}
Characterize those spaces which possess Borel Darboux retracts.
\end{myprob}

\begin{myprob}
\myproblem{O16}{K. Yamazaki \cite{MR2001b:54022}}
Let $X$ be a space, $A$ a subspace and $\gamma$ an infinite cardinal.
Find a nice class $\mathcal{P}$ of spaces such that $A$
is $P^\gamma$ (locally finite)-embedded in $X$ if and only if
every continuous map $f$ from $A$ into any $Y \in \mathcal{P}$
can be continuously extended over $X$.
\end{myprob}

\section*{P. Products, hyperspaces, and similar constructions}

\begin{myprob}
\myproblem{P1}{S. Williams \cite{TP3:Williams}}
Are $\beta\mathbb{N} \setminus \mathbb{N}$ and 
$\beta\mathbb{R} \setminus \mathbb{R}$ coabsolute?

\mynote{Solution}
Yes, if \myaxiom{MA} (S. Williams)
\cite{MR83d:54060}.
\end{myprob}

\begin{myprob}
\myproblem{P2}{S. Williams \cite{TP3:Williams}}
Are $\beta\mathbb{N} \setminus \mathbb{N}$ and 
$(\beta\mathbb{N} \setminus \mathbb{N}) \times 
(\beta\mathbb{N} \setminus \mathbb{N})$ coabsolute?
Yes, if \myaxiom{MA}
(B. Balcar, J. Pelant and P. Simon \cite{MR82c:54003}).
\end{myprob}

\begin{myprob}
\myproblem{P3}{S. Williams \cite{TP3:Williams}}
Are $\beta\mathbb{R} \setminus \mathbb{R}$ and 
$(\beta\mathbb{R} \setminus \mathbb{R}) \times 
(\beta\mathbb{R} \setminus \mathbb{R})$ coabsolute?
\end{myprob}

\begin{myprob}
\myproblem{P4}{S. Williams \cite{TP3:Williams}}
Is there a locally compact noncompact metric space $X$ of density at most
$2^\omega$ such that $\beta X \setminus X$ fails to be coabsolute with
either $\beta\mathbb{N} \setminus \mathbb{N}$ or 
$\beta\mathbb{R} \setminus \mathbb{R}$?

\mynote{Solution}
(S. Williams \cite{MR83d:54060})
If $X$ is a locally compact noncompact metric space of density at 
most $2^\omega$ then $\beta X \setminus X$ is coabsolute with
\begin{myenumerate}
\item 
$\beta\mathbb{N} \setminus \mathbb{N}$, if $X$ has a dense discrete
subspace;
\item
$\beta\mathbb{R} \setminus \mathbb{R}$, if the set of isolated points of
$X$ has compact closure;
\item 
$\beta\mathbb{N} \setminus \mathbb{N} + \beta\mathbb{R} \setminus \mathbb{R}$, 
otherwise.
\end{myenumerate}
\end{myprob}

\begin{myprob}
\myproblem{P5}{R. Heath \cite{MR55:13364}}
Is the Pixley-Roy hyperspace of $\mathbb{R}$ homogeneous?

\mynote{Solution}
Yes (M. Wage \cite{MR89f:54022}).
\end{myprob}

\begin{myprob}
\myproblem{P6}{R. Heath}
Does there exist an uncountable, non-discrete space $X$ which is
homeomorphic to its Pixley-Roy hyperspace $\mathcal{F}[X]$?

\mynote{Solution}
Yes (P. Nyikos and E. van~Douwen).
\end{myprob}

\begin{myprob}
\myproblem{P7}{M.E.~Rudin \cite{MR57:17584a}}
Can the perfect image of a normal subspace of a $\Sigma$-product
of lines be embedded in the $\Sigma$-product?
\end{myprob}

\begin{myprob}
\myproblem{P8}{T. Przymusi\'nski}
Can every space with a point-countable base be embedded into a
$\Sigma$-product of intervals?
\end{myprob}

\begin{myprob}
\myproblem{P9}{P. Nyikos}
If a product of two spaces is homeomorphic to $2^\kappa$, must one of the
factors be homeomorphic to $2^\kappa$?
This is true for $\kappa = \omega$, of course.

\mynote{Solution}
Yes (E. Shchepin).
\end{myprob}

\begin{myprob}
\myproblem{P10}{V. Saks \cite{MR82a:54013}}
Does there exist a subset $D$ of $\beta\omega \setminus \omega$ such that 
$|D| = 2^\mathfrak{c}$ and whenever $\{x_n : n \in \omega\}$ and 
$\{y_n : n \in \omega\}$ are sequences in $\beta\omega$ and 
$x,y \in \beta\omega \setminus \omega$, $d,d^\prime \in D$, and 
$x = d - \lim x_n$ and $y = d^\prime - \lim y_n$, then $x \neq y$?

\mynote{Notes}
See \cite[Example~2.3]{MR82a:54013}.
\end{myprob}

\begin{myprob}
\myproblem{P11}{V. Saks \cite{MR82a:54013}}
Does there exist a set $D$ of weak $P$-points such that 
$|D| = 2^\mathfrak{c}$ and if $x \in \operatorname{cl}A$ for some
countable subset $A$ of $\bigcup\{F(d) : d \in D\}$, then there exists 
a countable subset $C$ of $D$ such that $x \not\in \operatorname{cl}B$
for all countable subsets $B$ of 
$\bigcup\{F(d) : d \in D \setminus C\}$?

\mynote{Notes}
Here $F(d)$ is the set of all nonisolated images of $d$ under self-maps
of $\beta\omega$ induced by self-maps of $\omega$.
See \cite[\S~4]{MR82a:54013}.
\end{myprob}

\begin{myprob}
\myproblem{P11}{E. van~Douwen}
Let $\mathbb{H}$ be the half-line $[0,\infty)$.
Does there exist a characterization of $\mathbb{H}^*$ under \myaxiom{CH}?
For example:
\begin{myenumerate}
\item
Under \myaxiom{CH}, if $L$ is a $\sigma$-compact connected LOTS with 
exactly one endpoint and $\omega \leq w(L) \leq \mathfrak{c}$,
is $L^*$ homeomorphic to $\mathbb{H}^*$?
\item
More generally, does \myaxiom{CH} imply that $\mathbb{H}^*$ is (up to
homeomorphism) the only continuum of weight $\mathfrak{c}$ that is an
$F$-space, has the property that nonempty $G_\delta$'s have nonempty
interior, and is one-dimensional, indecomposable, hereditarily
unicoherent, and atriodic?
\end{myenumerate}

\mynote{Solution}
Yes to the first (A. Dow and K.P.~Hart \cite{MR95b:54031}).
No to the second, consider the pseudo-arc $\mathbb{P}$ and a component 
of $(\omega\times\mathbb{P})^*$; 
this component has all the properties but is not homeomorphic
to $\mathbb{H}^*$ (it is hereditarily indecomposable).
\end{myprob}

\begin{myprob}
\myproblem{P12}{E. van~Douwen}
Does there exist in \myaxiom{ZFC} a space that is homeomorphic to 
$\mathbb{N}^*$, but not trivially so?
Is it at least consistent with $\neg\text{\myaxiom{CH}}$ that such a space 
exists?

\mynote{Notes}
For example, under \myaxiom{CH}, 
$(\mathbb{N} \times \mathbb{N}^*)^* \approx \mathbb{N}^*$.
More generally, if \myaxiom{CH} then $X^*$ is homeomorphic to 
$\mathbb{N}^*$ whenever $X$ is locally compact, Lindel\"of, (strongly) 
zero-dimensional, and noncompact. 
This follows from Parovi\v{c}enko's theorem, see
\cite[Theorem~1.2.6]{MR86f:54027}.
\end{myprob}

\begin{myprob}
\myproblem{P13}{E. van~Douwen}
Write $X_0 \simeq X_1$ if there are open $U_i \subset X_i$ with compact
closure in $X_i$ for $i = 0,1$ such that $X_0 \setminus U_0$ and 
$X_1 \setminus U_1$ are homeomorphic.
Then $X^*$ and $Y^*$ are homeomorphic if $X \simeq Y$ but not conversely.
Does there exist in \myaxiom{ZFC} a pair of locally compact realcompact 
(preferably separable metrizable) spaces $X,Y$ such that $X^*$ is 
homeomorphic to $Y^*$, but $X \not\simeq Y$?

\mynote{Solution}
(A. Dow and K.P.~Hart \cite{MR2000d:54031})
Under \myaxiom{OCA}, if $X$ is locally compact, $\sigma$-compact and not
compact and if $X^*$ is a continuous image of $\omega^*$ then
$X\simeq\omega$.
\end{myprob}

\begin{myprob}
\myproblem{P14}{R. Pol}
Let $H$ be the hyperspace of the Hilbert cube.
The set $\{X \in H : \text{$X$ is countable-dimensional}\}$ is PCA 
(the projection of a co-analytic set) but not analytic; is it true that
this set is not co-analytic?
\end{myprob}

\begin{myprob}
\myproblem{P15}{P. Nyikos}
Is it consistent that $\beta\omega \setminus \omega$ is the union of a
chain of nowhere dense sets?

\mynote{Notes}
This cannot happen under \myaxiom{MA} since then 
$\beta\omega \setminus \omega$ cannot be covered by $\mathfrak{c}$ or 
fewer nowhere dense sets (S. Hechler \cite{MR57:12217}).

\mynote{Solution}
Yes, in fact this is equivalent to being able to cover 
$\beta\omega \setminus \omega$ by $\leq \mathfrak{c}$ nowhere dense sets.
On the one hand, any chain of nowhere dense sets covering 
$\beta\omega \setminus \omega$ must have cofinality $\leq \mathfrak{c}$ 
because $\beta\omega \setminus \omega$ has a dense subspace of
cardinality $\mathfrak{c}$.
On the other hand, B.~Balcar, J.~Pelant, and P.~Simon
\cite[Theorem~3.5(iv)]{MR82c:54003}
give a good indication of how extensive this class of models is:
if it is impossible to cover $\beta\omega \setminus \omega$ with 
$\leq \mathfrak{c}$ nowhere dense sets, then there are more than 
$\mathfrak{c}$ selective $P_\mathfrak{c}$-points 
\cite[Theorem~3.7]{MR82c:54003}.
So any model without such points (in particular, any model in which there 
is no scale of cofinality $\mathfrak{c}$ gives an affirmative solution to 
P15.
Examples are the usual Cohen real and Random real models. 
\end{myprob}

\begin{myprob}
\myproblem{P16}{W. Lewis}
Is the space of homeomorphisms of the pseudo-arc totally disconnected?
\end{myprob}

\begin{myprob}
\myproblem{P17}{E. van~Douwen}
Let $Y$ be a Hausdorff continuous image of the compact Hausdorff space $X$.
If $\Box^\omega X$ is paracompact [resp.\ normal], is $\Box^\omega Y$
paracompact [resp.\ normal]?
If the $G_\delta$-modification of $X$ is paracompact (or normal), is the
same true for that of $Y$.
\end{myprob}

\begin{myprob}
\myproblem{P18}{P. Nyikos}
Can normality of a $\Sigma$-product depend on the choice of the base 
point?
\end{myprob}

\begin{myprob}
\myproblem{P19}{P. Nyikos}
Is there a chain of clopen subsets of $\omega^*$ of uncountable cofinality
whose union is regular open?
Yes, if $\mathfrak{p} > w_1$ or $\mathfrak{b} = \mathfrak{d}$ or in any
model obtained by adding uncountably many Cohen reals.

\mynote{Solution}
No is also consistent, by a modification of Miller's rational forcing
(A. Dow and J. Stepr\={a}ns \cite{MR90c:03043, MR94b:03088}).
\end{myprob}

\begin{myprob}
\myproblem{P20}{E. van~Douwen}
Can one find in \myaxiom{ZFC} a point $p$ of $\omega^*$ such that 
$\beta(\omega^* \setminus \{p\}) \neq \omega^*$
or, better yet, $\zeta(\omega^* \setminus \{p\}) \neq \omega^*$?
It is known that $p$ is as required if it has a local base of cardinality
$\omega_1$, and that it is consistent for there to be $p \in \omega^*$ 
with $\beta(\omega^* \setminus \{p\}) = \omega^*$.

\mynote{Notes}
van~Douwen apparently had \myaxiom{PFA} in mind when he wrote
``it is consistent for there to be $p \in \omega^*$''
in posing the problem; at that time it was not yet known that 
\myaxiom{PFA} implied every point of the remainder had the stated 
property.

\mynote{Solution}
No. 
$\beta(\omega^* \setminus \{p\}) = \omega^*$ for all
$p \in \omega^*$ is consistent with 
$\text{\myaxiom{MA}} + \mathfrak{c} = \aleph_2$
(E. van~Douwen, K. Kunen, and J. van~Mill \cite{MR90b:54015}).
This is also true in any model where at least as many Cohen reals are
added as there are reals in the ground model (V.~Malykhin).
In the Miller model, $\omega^*\setminus\{p\}$ is $C^*$-embedded for some 
but not all $p$: it is iff $p$ is not a $P$-point
(A. Dow \cite{MR98c:03107}).
\end{myprob}

\begin{myprob}
\myproblem{P21}{J.T.~Rogers}
Let $f\colon X \to Y$ be a map between inverse limit spaces.
When does there exist a map induced from commuting diagrams on the inverse
sequence that has the desired properties of $f$ (such as being a 
homeomorphism taking the point $x$ onto the point $y$)?
\end{myprob}

\begin{myprob}
\myproblem{P22}{T.J.~Peters \cite{MR84m:54027}}
Must every non-pseudocompact $G$-space have remote points?
\end{myprob}

\begin{myprob}
\myproblem{P23}{T. Isiwata \cite{MR87m:54071}}
Is there a pseudocompact $\kappa$-metric space $X$ such that $\beta X$
is not $\kappa$-metrizable?
\end{myprob}

\begin{myprob}
\myproblem{P24}{T. Isiwata, attributed to Y. Tanaka \cite{MR87m:54071}}
Is there a $\kappa$-metric space $X$ such that $\upsilon X$ is not
$\kappa$-metrizable where $|X|$ is nonmeasurable?
\end{myprob}

\begin{myprob}
\myproblem{P25}{M. Tikoo \cite{MR87g:54056}}
Characterize all Hausdorff spaces for which $\sigma X = \mu X$.

\mynote{Notes}
Each semiregular Hausdorff space $X$ can be densely embedded in a 
canonical semiregular $H$-closed space $\mu X$, called its 
\emph{Banaschewski-Shanin-Fomin extension}.
Tikoo constructed an analogue to $\mu X$ for any Hausdorff space $X$.
The \emph{Fomin extension} $\sigma X$ is a strict $H$-closed extension 
of $X$.
\end{myprob}

\begin{myprob}
\myproblem{P26}{I. Ntantu \cite{MR88c:54002}}
Let $X$ be a Tychonoff space and $K(X)$ the hyperspace of its nonempty
compact subsets.
Recall that a continuous $f\colon Z \to Y$ is called 
\emph{compact-covering} if each compact subset of $Y$ is the image of 
some compact subset of $Z$.
If $K(X)$ with the Vietoris (i.e., finite) topology is a continuous image 
of $\omega^\omega$, must $X$ be a compact-covering image of 
$\omega^\omega$?
The converse is true.
\end{myprob}

\begin{myprob}
\myproblem{P27}{V. Malykhin}
A point $x \in X$ is called a butterfly point if there exist disjoint sets
$A$ and $B$ of $X$ such that $\overline{A} \cap \overline{B} = \{x\}$.
Is it consistent that there is a non-butterfly point in
$\omega^* = \beta\omega \setminus \omega$?
Is it consistent with MA?

\mynote{Notes}
If \myaxiom{PFA}, non-butterfly points of $\omega^*$ would be exactly 
those points $p$ for which $\omega^* \setminus \{p\}$ would be normal.
Malykhin has withdrawn the claim made in \cite{MR49:7985}
that \myaxiom{MA} implies $\omega^* \setminus \{p\}$ is always non-normal.
This claim would have implied the non-existence of such points under 
\myaxiom{PFA}.
It is still an unsolved problem whether there is a model in which
$\omega^* \setminus \{p\}$ is normal for some $p \in \omega^*$.

\mynote{Solution}
(A. Beslagic and E. van~Douwen \cite{MR91h:54040})
It is not consistent with \myaxiom{MA}.
In fact, if $\mathfrak{r} = \mathfrak{c}$, then every point in $\omega^*$
is a butterfly point; more strongly, 
$\omega^* = \beta\omega \setminus \omega$ is non-normal for every 
$p \in \omega^*$.
Here $\mathfrak{r}$ denotes the reaping number, i.e., the least $\kappa$
for which there is a family $\mathcal{R}$ of $\kappa$ subsets of $\omega$
such that if $A$ is a subset of $\omega$ then $A$ does not reap
$\mathcal{R}$, i.e., there is $R \in \mathcal{R}$ such that either 
$R \subset^* A$ or $R \cap A$ is finite.
\end{myprob}

\begin{myprob}
\myproblem{P28}{V. Malykhin}
Is there a model in which $\beta(\omega^* \setminus \{p\}) = \omega^*$
for some points of $\omega^*$ but not for others?

\mynote{Notes}
See P20.
\end{myprob}

\begin{myprob}
\myproblem{P29}{L.B.~Lawrence \cite{MR89f:54024}}
Let $X = \Box^\omega \mathbb{Q}$, 
$Y = \bigtriangledown^\omega \mathbb{Q}$, 
$\sigma\colon X \to Y$ the natural quotient map.
Recall that two points of $x$ have the same image iff they disagree on at
most finitely many coordinates.
Is there a closed subset $C$ of $X$ such that $\sigma[C]$ is dense in $Y$
and $C$ contains at most one point in each fiber?
\end{myprob}

\begin{myprob}
\myproblem{P30}{L.B.~Lawrence \cite{MR89f:54024}}
Replace the rationals by the irrationals in P29.
\end{myprob}

\begin{myprob}
\myproblem{P31}{R. Levy}
Let $X$ be either $[0, 1)$ or a Euclidean space of dimension at least $2$
(this is so $\beta X \setminus X$ will be connected).
Is the set of weak $P$-points of $\beta X \setminus X$ connected?
\end{myprob}

\begin{myprob}
\myproblem{P32}{R. Levy}
Is there a realcompact space $X$ such that for some 
$p \in \beta X \setminus X$ the space $\beta X \setminus \{p\}$ is
normal?
We obtain an equivalent problem if ``Lindel\"of'' is substituted
for ``realcompact''.
\end{myprob}

\begin{myprob}
\myproblem{P33}{A. Okuyama \cite{MR94d:54033}}
Let $X$ be a paracompact Hausdorff space and $Y$ a K-analytic space.
If $X \times Y$ is normal, then is $X \times Y$ paracompact?

\mynote{Notes}
(A. Okuyama)
It seems that this question concerns the property of a mapping such as
$\operatorname{id}_X \times \varphi$, where $\varphi$ is an upper 
semicontinuous, compact-valued mapping from the space $\mathbb{P}$ of
irrationals to the power set of $Y$.
\end{myprob}

\begin{myprob}
\myproblem{P33}{T. LaBerge \cite{MR95h:54015}}
Is there a countable collection $\{A_n : n \in \omega\}$ of 
non-Lindel\"of ACRIN spaces whose topological sum is ACRIN?

\mynote{Notes}
ACRIN = all continuous regular images normal. 
\end{myprob}

\begin{myprob}
\myproblem{P34}{T. LaBerge \cite{MR95h:54015}}
Is there an ACRIN space $X$ such that $X+X$ is not ACRIN?
\end{myprob}

\begin{myprob}
\myproblem{P35}{T. LaBerge \cite{MR95h:54015}}
Are there Lindel\"of spaces $X$ and $Y$ such that $X \times Y$
is ACRIN but not Lindel\"of?
\end{myprob}

\begin{myprob}
\myproblem{P36}{C. Good \cite{MR98d:54041}}
Can the square of a perfectly normal manifold be a Dowker space?

\mynote{Notes}
No, if \manotch, because it implies 
perfectly normal manifolds are metrizable.
\end{myprob}

\begin{myprob}
\myproblem{P37}{C. Good \cite{MR98d:54041}}
Does \manotch\ imply the existence of a Dowker manifold, 
or even a locally compact Dowker space?
\end{myprob}

\begin{myprob}
\myproblem{P38}{C. Good \cite{MR98d:54041}}
If $X$ is a normal, countably paracompact space and $X^2$ is normal, does
\manotch\ imply $X^2$ is countably paracompact?
What if $X$ is also perfectly normal?
\end{myprob}

\begin{myprob}
\myproblem{P39}{C. Good \cite{MR98d:54041}}
Is there a Dowker space $X$ such that $X^2$ is Dowker?
Such that $X^n$ is Dowker for all finite $n$?
\end{myprob}

\begin{myprob}
\myproblem{P40}{C. Good \cite{MR98d:54041}}
Can the square of a monotonically normal space or of a Lindel\"of space be
Dowker?
\end{myprob}

\begin{myprob}
\myproblem{P41}{M. Bonanzinga \cite{MR98k:54037}}
Does there exist a \myaxiom{ZFC} example of two star-Lindel\"of 
topological groups $G$ and $H$ such that the product $G \times H$ is not 
star-Lindel\"of?

\mynote{Notes}
See C68.
\end{myprob}

\begin{myprob}
\myproblem{P42}{D. Mattson \cite{MR98j:54045}}
Can a nowhere rim-compact space have a compactification with 
zero-dimensional remainder?
\end{myprob}

\begin{myprob}
\myproblem{P43}{W.J.~Charatonik \cite{MR99k:54005}}
Let a mapping $f\colon X \to Y$ between continua $X$ and $Y$ be such that
the induced mapping $C(f)$ is a near-homeomorphism (in particular, $C(X)$
and $C(Y)$ are homeomorphic).
Does it imply that $2^f$ is a near-homeomorphism?
The same question, if $X = Y$.

\mynote{Notes}
For a given metric continuum $X$, the symbols $2^X$ and $C(X)$ denote the
hyperspaces of all nonempty closed subsets and of all nonempty subcontinua
of $X$, respectively.
Similarly, given a mapping $f\colon X \to Y$ between continua, the 
symbols $2^f$ and $C(f)$ denote the induced mappings.
\end{myprob}

\begin{myprob}
\myproblem{P44}{W.J.~Charatonik \cite{MR99k:54005}}
Let a mapping $f$ between continua be such that the induced mapping $2^f$
is confluent.
Does it imply that the induced mapping $C(f)$ is also confluent?
\end{myprob}

\begin{myprob}
\myproblem{P45}{E. Casta\~neda \cite{MR2001a:54048}}
Does there exist an indecomposable continuum $X$ such that $F_2(X)$ is not
unicoherent?

\mynote{Notes}
The space $F_2(X)$ is the hyperspace of two-point subsets of $X$.
\end{myprob}

\begin{myprob}
\myproblem{P46}{E. Casta\~neda \cite{MR2001a:54048}}
Does there exist an hereditarily unicoherent continuum $X$ such that 
$F_2(X)$ is not unicoherent?
\end{myprob}

\begin{myprob}
\myproblem{P47}{J.J.~Charatonik \cite{MR2001a:54048}}
Does there exist an hereditarily unicoherent, hereditarily decomposable
continuum $X$ such that $F_2(X)$ is not unicoherent. 
\end{myprob}

\begin{myprob}
\myproblem{P48}{J.J.~Charatonik and W.J.~Charatonik \cite{MR2001f:54037}}
Is it true that if a continuum $X$ has the property of Kelley,
then the Cartesian product $X \times [0,1]$ is semi-Kelley?
\end{myprob}

\begin{myprob}
\myproblem{P49} {J.J.~Charatonik and W.J.~Charatonik \cite{MR2001f:54037}}
Is it true that if a continuum $X$ is semi-Kelley, then the hyperspace
$2^X$ (resp.\ $C(X)$) is contractible?

\mynote{Notes}
See G38, G39.
\end{myprob}

\begin{myprob}
\myproblem{P50}{G. Acosta \cite{MR2002j:54006}}
Let $X$ be a fan without the property of Kelley. 
Is it true that $X$ does not have almost unique hyperspace?

\mynote{Notes}
Given a continuum $X$, consider a class $\mathcal{F}_X$ of continua $Y$
such that:
no member of $\mathcal{F}_X$ is homeomorphic to $X$;
no two distinct members of $\mathcal{F}_X$ are homeomorphic; 
the hyperspaces $C(X)$ and $C(Y)$ are homeomorphic, for each
$Y\in\mathcal{F}_X$;
if $Z$ is a continuum such that the hyperspaces $C(Z)$ and $C(X)$ are 
homeomorphic, then either $Z$ is homeomorphic to $X$ or $Z$ is 
homeomorphic to some member $Y$ of $\mathcal{F}_X$.
A continuum $X$ is said to have \emph{unique hyperspace} iff the class
$\mathcal{F}_X$ is empty.
If the class $\mathcal{F}_X$ is nonempty and finite, we say that $X$ has
\emph{almost unique hyperspace}.
\end{myprob}

\begin{myprob}
\myproblem{P51}{G. Acosta \cite{MR2002j:54006}}
Let $X$ be an indecomposable continuum such that each proper and
nondegenerate subcontinuum of $X$ is a finite graph.
Does $X$ have unique hyperspace?
\end{myprob}

\begin{myprob}
\myproblem{P52}{G. Acosta \cite{MR2002j:54006}}
For a metric compactification of the space $V = (-\infty,\infty)$
and connected and nondegenerate remainder $R$, we write $X = V \cup R$ 
and define
$R_1 = \bigcap_{n\in\mathbb{N}}\text{Cl}_X((n,\infty))$
and
$R_2 = \bigcap_{n\in\mathbb{N}}\text{Cl}_X((-\infty,-n))$.
Let us assume that $R_1 \not= R_2$.
Is there a continuum $Y$, not homeomorphic to $X$, such that the
hyperspaces $C(X)$ and $C(Y)$ are homeomorphic?
What is the cardinality of the class $\mathcal{F}_X$?
\end{myprob}

\section*{Q. Generalizations of topological spaces}

\begin{myprob}
\myproblem{Q1}{R. Price \cite[E.~\v{C}ech]{Cech}}
Does there exist a \v{C}ech function?
That is, a function 
$f\colon \mathcal{P}(\omega) \to \mathcal{P}(\omega)$
such that $f \neq \operatorname{id}$, $A \subset f(A)$ for all $A$,
$f(A \cup B) = f(A) \cup f(B)$ for all $A$, $B$, and $f$ is onto?
In other words, is there a countable closure space in which every subset
is the closure of a subset?

\mynote{Notes}
Yes is consistent (R. Price)
\cite{MR84d:04003}.
See \cite[\S~4]{MR88k:04001} for a proof of a previously unpublished
related theorem of F. Galvin.
\end{myprob}

\begin{myprob}
\myproblem{Q2}{R. Herrmann}
Characterize those topological spaces $(X, \tau)$ such that 
$\operatorname{sh} = \tau$ (resp.\ $\operatorname{u} = \tau$).
\end{myprob}

\begin{myprob}
\myproblem{Q3}{R. Herrmann}
Characterize those topological spaces $(X, \tau)$ such that
$\operatorname{rc} \times \operatorname{rc} = 
\operatorname{rc}(\tau \times \tau_z)$
(the r.c.\ structure generated by $\tau \times \tau_z$) and those such
that $\operatorname{u} \times \operatorname{u} = 
\operatorname{u}(\tau \times \tau_z)$,
$\operatorname{sh} \times \operatorname{sh} = 
\operatorname{sh}(\tau \times \tau_z)$.
\end{myprob}

\begin{myprob}
\myproblem{Q4}{R. Herrmann}
Characterize those topological spaces for which $\operatorname{sh}$
[resp.\ $\operatorname{u}$] is pseudotopological, pretopological, or
topological.
\end{myprob}

\begin{myprob}
\myproblem{Q5}{R. McKee \cite{MR91i:54029}}
Let $(X, \mu)$ be a nearness space and let $K(X)$ denote the group (under
composition) of all near-homeomorphisms from $(X, \mu)$ to itself.
When is it true that if $K(X)$ and $K(Y)$ are isomorphic, then $X$ and $Y$
are near-homeomorphic?
\end{myprob}

\section*{QQ. Comparison of topologies}

\begin{myprob}
\myproblem{QQ1}{T. LaBerge \cite{MR97a:54004}}
If $X = \bigcup_{\alpha < \kappa} X_\alpha$ has the fine topology and 
$t^{+}(X_\alpha) \leq \kappa^+$, is
$t(p,X) = \sup\{ t(p, X_\alpha) : \alpha < \kappa, p \in X_\alpha \}$ 
for each $p \in X$?
\end{myprob}

\begin{myprob}
\myproblem{QQ2}{T. LaBerge \cite{MR97a:54004}}
If $s^+(X) \leq \kappa$ or $hl^{+}(X) \leq \kappa$,
is the fine topology the only compatible topology?
\end{myprob}

\begin{myprob}
\myproblem{QQ3}{T. LaBerge \cite{MR97a:54004}}
Is it possible to have a $\kappa$-chain of Hausdorff spaces
with exactly two compatible Hausdorff topologies?
\end{myprob}

\section*{R. Dimension theory}

\begin{myprob}
\myproblem{R1}{T. Przymusi\'nski}
If $X$ is a metric space in which every subset is an $F_\sigma$, 
then is $\dim X = 0$?

\mynote{Notes}
Yes, if \visl\ (G.M.~Reed \cite{MR82d:54033}).
See S8.
\end{myprob}

\begin{myprob}
\myproblem{R2}{R. Pol \cite{MR82k:54060}}
Let $\mathcal{D}$ be an upper semicontinuous decomposition of a
compactum $X$ into countable-dimensional compacta.
Is it true that 
$\sup\{\operatorname{ind} S : S \in \mathcal{D}\} < \omega_1$?
\end{myprob}

\begin{myprob}
\myproblem{R3}{R. Pol \cite{MR82k:54060}}
Let $\alpha < \omega_1$.
What is the ordinal number $\mu(\alpha)$, defined to be the minimum 
$\operatorname{ind} X$ of all $X$ such that $X$ is a countable-dimensional
compactum containing topologically all compacta $S$ with 
$\operatorname{ind} S \leq \alpha$?
\end{myprob}

\begin{myprob}
\myproblem{R4}{L. Rubin \cite[P.S.~Alexandroff's CE-problem]{MR84f:54049}}
Does there exist a separable metric space, compact or not, which has
finite cohomological dimension and infinite topological dimension?

\mynote{Solution}
Yes, there is even a compact example
(A. Dranishnikov \cite{MR90e:55004}).
\end{myprob}

\begin{myprob}
\myproblem{R5}{J. Keesling \cite{MR84d:54060}}
If $f(X) = Y$ is a mapping between compact metric spaces such that 
$m \leq \dim f^{-1}(y) \leq n$ for all $y \in Y$, then is there a closed
set $K$ in $X$ such that $\dim K \leq n-m$ and $\dim f(K) = \dim Y$?

\mynote{Solution}
Yes
(E. Kurihara \cite{MR85f:54080}).
\end{myprob}

\begin{myprob}
\myproblem{R6}{E. van~Douwen}
For which sequences $\{k_n : n \geq 1\}$ of integers is there a 
separable metrizable space $X$ such that $\dim X^n = k_n$ for all $n$?
For example, is $\lim_n k_n / n = \sqrt{2}$ possible?
What if $X$ is also compact?
\end{myprob}

\begin{myprob}
\myproblem{R7}{T. Hoshina \cite{MR93f:54013}}
Suppose $X \times Y$ is normal $T_1$, where $Y$ is a La\v{s}nev space.
Does $\dim(X \times Y) \leq \dim X + \dim Y$
hold for the covering dimension $\dim$?
Yes, if $X$ is paracompact.
\end{myprob}

\begin{myprob}
\myproblem{R8}{T. Kimura \cite{MR94k:54069}}
Does there exist a normal (or metrizable) space $X$ having
$\operatorname{trdim}$ such that every compactification of $X$ fails 
to have $\operatorname{trind}$?
\end{myprob}

\begin{myprob}
\myproblem{R9}{V.A.~Chatyrko \cite{MR98b:54046}}
If $C$ is the Cantor set, is 
$\operatorname{trdim} X = \operatorname{trdim} (X \times C)$?
Yes, if $\operatorname{trdim} (X \times C) \geq \omega^2$.
\end{myprob}

\begin{myprob}
\myproblem{R10}{V.A.~Chatyrko \cite{MR98b:54046}}
Is it true that if $X$ is a space and $\alpha$ is a countable ordinal 
number $\geq \omega^2$, then $\operatorname{trdim} X \geq \alpha$
iff $X$ admits an essential map onto Henderson's cube $H^\alpha$?
Yes, for limit ordinals.
\end{myprob}

\begin{myprob}
\myproblem{R11}{D. Garity \cite{MR98d:54063}}
Is there a homogeneous compact metric space of dimension less than
$n+2$ that is locally $n$-connected but not $2$-homogeneous?
\end{myprob}

\begin{myprob}
\myproblem{R12}{T. Kimura \cite{MR98a:54029}}
Does there exist a S-w.i.d. (i.e., weakly infinite-dimensional
in the sense of Smirnov) space $X$ such that 
$\operatorname{trdim} X \geq w(X)^+$?
\end{myprob}

\begin{myprob}
\myproblem{R21}{M.G.~Charalambous \cite{MR99g:54031}}
Is there a perfectly normal space $Y$ with $\operatorname{ind} Y = 1$ such
that no Lindel\"of (or even strongly paracompact) extension of $Y$ has
small transfinite inductive dimension?
\end{myprob}

\begin{myprob}
\myproblem{R22}{A. Dranishnikov and T. Januszkiewicz \cite{MR2001k:20082}}
Does every discrete metric space $X$ of bounded geometry (e.g., a 
finitely generated group) have the property~$A$.

\mynote{Solution}
No.
In \cite{MR2002e:53056},
M. Gromov announced that there is a finitely generated group without 
property $A$. 
The construction is presented in \cite{MR1978492}.
See G. Yu's article \cite{MR2000j:19005} for the definition of property $A$. 
J.L.~Tu \cite{MR2002j:58038} and G. Bell \cite{bell2002}
proved that property $A$ is preserved under the graph of groups 
construction, in particular by the amalgamated product and by the HNN 
extension. 
In \cite{dranishnikov2002}, it was shown that property $A$ of a space $X$ 
is equivalent to the existence of a geometric Anti-\v{C}ech approximation 
of $X$.
\end{myprob}

\begin{myprob}
\myproblem{R23}{A. Dranishnikov and T. Januszkiewicz \cite{MR2001k:20082}}
Assume that the Higson corona of a discrete metric space $X$ is finite
dimensional. 
Does $X$ have property $A$? 
\end{myprob}

\begin{myprob}
\myproblem{R24}{A. Dranishnikov and T. Januszkiewicz \cite{MR2001k:20082}}
Does every $\operatorname{CAT}(0)$ group have property $A$? 
\end{myprob}

\section*{S. Problems closely related to set theory}

\begin{myprob}
\myproblem{S1}{\cite[Rudin and Lutzer]{MR80j:54027}}
Is every $Q$-set strong?
In other words, are its finite powers $Q$-sets?

\mynote{Notes}
No is consistent.
It is consistent that there is a $Q$-set of cardinality $\aleph_2$,
but no square of a space of cardinality $\aleph_2$ is a $Q$-set
(W. Fleissner \cite{MR85h:54006}).
\end{myprob}

\begin{myprob}
\myproblem{S2}{van~Douwen and Rudin}
In \myaxiom{ZFC}, are there two free ultrafilters on $\omega$ with no 
common finite-to-one image?
Under \myaxiom{MA} there are such ultrafilters.

\mynote{Solution}
The principle of near coherence of filters (\myaxiom{NCF}) asserts 
that any two free ultrafilters have a common finite-to-one image.
\myaxiom{NCF} is consistent relative to \myaxiom{ZFC}.
See the papers by A. Blass and S. Shelah
\cite{MR88d:03094a, MR88d:03094b, MR90m:03087}.
\end{myprob}

\begin{myprob}
\myproblem{S3}{K. Hofmann \cite{MR80k:06011}}
Let $f^k$ be the permutation on the discrete space $\mathbb{Z}$ of
integers which takes $n$ to $n + k$.
For $k \in \mathbb{Z}$ and $p \in \beta(\mathbb{Z})$, let 
$p^k = \{f^k(M) : M \in p\}$,
and $O_p = \{p^k : k \in \mathbb{Z}\}$ the orbit of $p$.
Let 
$\mathcal{O} = 
\{\overline{O_p} : p \in \beta\mathbb{Z} \setminus \mathbb{Z}\}$,
and let $\mathcal{M}$ be the set of maximal members of $\mathcal{O}$.
Is there an infinite strictly increasing sequence of members of 
$\mathcal{O}$?
How long can such be?
What can be said in general about $\mathcal{O}$ and $\mathcal{M}$?
\end{myprob}

\begin{myprob}
\myproblem{S4}{E. van~Douwen}
If $D \subset \beta\omega \setminus \omega$ is nowhere dense, is it true
in \myaxiom{ZFC} that
$\{(\beta\pi)^{\to}D : \text{$\pi$ is a permutation of $\omega$}\}$
is not all of $\beta\omega \setminus \omega$?
What if $D = \bigcap\{\overline{A} : A \in \mathcal{A}\}$ where 
$\mathcal A$ is one of 
$\{A \subset \omega : \lim |A \cap n|/n = 1\}$ or
$\{A \subset \omega : \sum_{n \in A \setminus \{0\}} 1/n = \infty\}$?

\mynote{Solution}
In \cite{MR92e:54023}, A.A.~Gryzlov constructs $2^\mathfrak{c}$ many
$0$-points, where $u$ is a \emph{$0$-point} if for every permutation 
$\pi$ of $\omega$ there is $A \in u$ with 
$\lim_n {|\pi[A] \cap n|}/n = 0$; this shows that the answer for the 
first $D$ above is negative.

Also, let $\{A_n : n \in \omega\}$ be any partition of $\omega$ into
infinite sets, and let 
$D = \bigcap_n\operatorname{cl}\bigl(\bigcup_{m\ge n}A_m^*\bigr)$.
Then $u$ is a $P$-point iff $u \notin \bigcup_\pi(\beta\pi)^\to D$, so 
for this $D$ the answer is yes iff there are $P$-points.
\end{myprob}

\begin{myprob}
\myproblem{S5}{J. Stepr\={a}ns}
If there is a non-meager subset of $\mathbb{R}$ of cardinality
$\aleph_1$, is there a Luzin set?
\end{myprob}

\begin{myprob}
\myproblem{S6}{J. Stepr\={a}ns}
If there is a measure zero subset of $\mathbb{R}$ of cardinality
$\aleph_1$, is there a Sierpi\'nski set?
\end{myprob}

\begin{myprob}
\myproblem{S7}{T. Przymusi\'nski}
A $\sigma$-set is a separable metric space in which every $F_\sigma$-set
is a $G_\delta$-set.
Does there exist a $\sigma$-set of cardinality $\aleph_1$?

\mynote{Notes}
Yes if \myaxiom{MA}; under \myaxiom{MA} there even exists a $\sigma$-set 
of cardinality $\mathfrak{c}$.

\mynote{Solution}
This is independent of \myaxiom{ZFC}.
As is well known, \myaxiom{MA} implies every subset of $\mathbb{R}$ of 
cardinality $< \mathfrak{c}$ is a $Q$-set (every subset is a $G_\delta$). 
On the other hand,
it is also consistent that every separable, uncountable metric space
contains subsets that are arbitrarily far up in the Borel hierarchy
(A. Miller \cite{MR80m:04003}).
\end{myprob}

\begin{myprob}
\myproblem{S8}{T. Przymusi\'nski}
A $Q$-set is a metrizable space in which every subset is a $G_\delta$.
Is every $Q$-set strongly zero-dimensional? linearly orderable?

\mynote{Notes}
Yes to both if \visl, because then every $Q$-set is $\sigma$-discrete
(G.M.~Reed \cite{MR82d:54033}).
Every $\sigma$-discrete normal space $X$ satisfies $\dim X = 0$ by the
countable sum theorem, and every strongly zero-dimensional metric space
is linearly orderable (H.~Herrlich \cite{MR32:426}). 
\end{myprob}

\begin{myprob}
\myproblem{S9}{R. Telg\'arsky \cite{MR82m:54039}}
Let $X$ belong to the $\sigma$-algebra generated by the analytic subsets 
of an uncountable Polish space $Y$.
Is the game $G(X,Y)$ determined?
\end{myprob}

\begin{myprob}
\myproblem{S10}{R. Telg\'arsky \cite{MR82m:54039}}
Let $X$ be a Luzin set on the real line.
Does Player II have a winning strategy in the game $G(X, \mathbb{R})$?
\end{myprob}

\begin{myprob}
\myproblem{S11}{E. van~Douwen \cite{MR84b:54050} 
\cite[Question~8.11]{MR87f:54008}}
For a space $X$ let $\mathcal{K}(X)$ denote the poset (under inclusion) of
compact subsets of $X$ and let $\operatorname{cf} \mathcal{K}(X)$
denote the cofinality of $\mathcal{K}(X)$, i.e., 
$\min \{ |\mathcal{L}| : {\mathcal{L} \subset \mathcal{K}(X)},\
{\forall K \in \mathcal{K}(X)}\, 
{\exists L \in \mathcal{L}}\,
{K \subset L} \}$.
If $X$ is separable metrizable, and analytic (or at least absolutely Borel) 
but not locally compact, is 
$\operatorname{cf} \mathcal{K}(X) = \mathfrak{d}$?

\mynote{Solution}
We quote J. Vaughan \cite{MR1078647}:
By \cite[8.10]{MR87f:54008} this question is clearly intended for $X$ that
are not $\sigma$-compact, and for them $\mathfrak{d} \leq k(X) \leq
\operatorname{cf}(\mathcal{K}(X))$.
Thus, the question reduces to: is $\operatorname{cf}(\mathcal{K}(X)) \leq
\mathfrak{d}$?
Here, $\operatorname{cf}(\mathcal{K}(X))$ denotes the smallest cardinality
of a family $\mathcal{L}$ of compact subsets of $X$ such that for every
compact set $K \subseteq X$, there exists $L \in \mathcal{L}$ with 
$K \subseteq L$.
The answer to the second question is in the affirmative, but the answer to
the first question is independent of the axioms of \myaxiom{ZFC}.
H. Becker \cite{MR90e:03062} has constructed a model in which there is an
analytic space $X \subset 2^\omega$ with
$\operatorname{cf}(\mathcal{K}(X)) > \mathfrak{d}$.
On the other hand, under \myaxiom{CH}, 
$\operatorname{cf}(\mathcal{K}(X)) = \mathfrak{d} = \omega_1$.
F.~van~Engelen \cite{MR90g:54031} proved that if $X$ is co-analytic,
then $\operatorname{cf}(\mathcal{K}(X)) \leq \mathfrak{d}$.
The same follows from Fremlin's theory \cite{MR92c:54032} of Tukey's 
ordering.
Also see \cite{MR95e:06006}.
\end{myprob}

\begin{myprob}
\myproblem{S12}{P. Nyikos}
For each cardinal $\kappa$, let $u_\kappa$ be the least cardinality of a
base of for a uniform ultrafilter on a set of cardinality $\kappa$.
Is it consistent to have $\lambda < \kappa$, yet $u_\kappa < u_\lambda$?
How about in the case $\lambda = \omega$, $\kappa = \omega_1$?
\end{myprob}

\begin{myprob}
\myproblem{S13}{E. van~Douwen}
\cite{MR87d:03130}
Let \myaxiom{LN} be the axiom that every linearly orderable space is 
normal.
Does \myaxiom{LN} imply \myaxiom{AC} in \myaxiom{ZF}?

\mynote{Notes}
Birkhoff asked whether \myaxiom{LN} depends on \myaxiom{AC} 
\cite{MR82a:06001}.
It is known that \myaxiom{LN} is equivalent to ``for every complete linear 
order $L$ there is a choice function for the collection of nonempty 
intervals of $L$''.
From this, \myaxiom{ZF} $\not\Rightarrow$ \myaxiom{LN} follows easily.
\myaxiom{AC} $\Rightarrow$ \myaxiom{LN} is well known.
\myaxiom{LN} does not imply \myaxiom{AC} in $\text{ZF}^-$, i.e., without 
foundation (E.~van~Douwen \cite{MR87d:03130}).

\mynote{Solution}
No, L. Hadad and M. Morillon \cite{MR91f:03102} proved that \myaxiom{LN} 
does not imply \myaxiom{AC} in \myaxiom{ZF}.
\end{myprob}

\begin{myprob}
\myproblem{S14}{P. Nyikos}
Call a point of $\omega^*$ a simple $P$-point if it has a totally ordered
clopen base.
\begin{myenumerate}
\item
Does the existence of a simple $P$-point imply the existence of a scale,
i.e., a cofinal well-ordered subset of $({}^{\omega}\omega, <^*)$?
\item
Is it consistent that there exist simple $P$-points $p$ and $q$ with bases
of different cofinalities?
\end{myenumerate}
The cofinality of any simple $P$-point is either $\mathfrak{b}$ or 
$\mathfrak{d}$, so there can be at most two different cofinalities, and
an affirmative answer to the first question implies a negative answer to
the second question.

\mynote{Solution}
(S. Shelah \cite{MR88e:03073})
No, to the first question; yes, to the second question.
To be precise, there are models in which there are simple $P$-points and 
scales, but there is a model in which there are both simple $P$-points
with bases of cardinality $\aleph_1$ and of cardinality $\aleph_2$, and
such a model cannot contain a scale.
\end{myprob}

\begin{myprob}
\myproblem{S15}{S. Yang \cite{MR87i:54051}}
Let $I$ be a subset of $\omega^*$.
If $|I| < 2^\mathfrak{c}$, does there exist $p \in \omega^*$ such that
$p$ is incomparable in the Rudin-Keisler order with all $q \in I$?
Yes is consistent.
\end{myprob}

\begin{myprob}
\myproblem{S16}{R. Levy \cite{MR92i:54037}}
Is it consistent that there is an Isbell-Mrowka $\Psi$ space such that 
every subset of $\aleph_1$ nonisolated points is $2$-embedded, or 
$C^*$-embedded?

\mynote{Notes}
The first question is equivalent to asking for a MAD family $M$ of subsets
of $\omega$ such that, given disjoint subfamilies $S$ and $T$ of 
cardinality $|M|$ there is $A \subset \omega$ such that $A$ almost 
contains each member of $S$ and almost misses each member of $T$.
\end{myprob}

\begin{myprob}
\myproblem{S17}{J. Stepr\={a}ns \cite{MR94j:03101}}
Does there exist a Cook set in $\mathbb{N}^3$?

\mynote{Notes}
Yes, if \myaxiom{MA}.
Here we will refer to maximal antichains of monotone paths in
$\mathcal{P}(\mathbb{N}^n) / \mathcal{B}_n$ as \emph{Cook sets} for all 
$n$, not just for $n = 2$.
\end{myprob}

\begin{myprob}
\myproblem{S18}{J. Stepr\={a}ns \cite{MR94j:03101}}
Does the existence of a Cook set in $\mathbb{N}^3$ imply the existence of 
a Cook set in $\mathbb{N}^4$?
\end{myprob}

\begin{myprob}
\myproblem{S19}{J. Stepr\={a}ns \cite{MR94j:03101}}
For each $n \in \omega \setminus \{0, 1\}$, does there exist a model of
set theory in which there is a Cook set in $\mathbb{N}^{n+1}$ but not in
$\mathbb{N}^n$?
\end{myprob}

\begin{myprob}
\myproblem{S20}{J. Stepr\={a}ns \cite{MR94j:03101}}
Call a family of monotone paths in $\mathbb{N}^k$ \emph{weakly maximal} 
if any two paths are separated and the family cannot be extended to a 
larger family with this property.
Let $\mathfrak{a}^-_k$ [resp.\ $\mathfrak{a}_k$] be the least cardinality
of a weakly maximal [resp.\ maximal, assuming one exists] 
family of monotone paths in $\mathbb{N}^k$.
Does $\mathfrak{a}^-_k$ equal $\mathfrak{a}_k$ when the latter exists?
\end{myprob}

\begin{myprob}
\myproblem{S21}{J. Stepr\={a}ns \cite{MR94j:03101}}
Recall that $\mathfrak{a}$ represents the least cardinality of
an infinite maximal almost disjoint family in $\mathcal{P}(\omega)$.
What are the relationships between the cardinal $\mathfrak{a}$,
the cardinals $\mathfrak{a}_k$, and the cardinals $\mathfrak{a}^-_k$?
\end{myprob}

\begin{myprob}
\myproblem{S22}{A. Tomita \cite{MR99k:22007}}
Let $\kappa$ be the least cardinal such that if $G$ is a free abelian
group endowed with a group topology, then $G^\kappa$ is not countably
compact.
Under $\text{\myaxiom{MA}}_\text{countable}$, $\kappa > 1$, and in 
\myaxiom{ZFC}, $\kappa \leq \omega$.
Find a better bound for $\kappa$ or determine which cardinals $\kappa$ 
may be.
In particular, is it true that $\kappa > 1$ in \myaxiom{ZFC}?
Is it consistent that $\kappa > 2$?
Is $\omega$ the best upper bound for $\kappa$?
\end{myprob}

\begin{myprob}
\myproblem{S23}{A. Tomita \cite{MR99k:22007}}
Let $\lambda$ be the least cardinal such that if $S$ is a both-sided
cancellative semigroup which is not a group, endowed with a group
topology, then $S^\lambda$ is not countably compact.
Under $\text{\myaxiom{MA}}_\text{countable}$, $\lambda > 1$, and in 
\myaxiom{ZFC}, $\lambda < 2^\mathfrak{c}$.
Find a better bound for $\lambda$ or determine which cardinals $\lambda$
may be.
Is there a relation between $\lambda$ and $\kappa$?

\mynote{Notes}
See S22 for the definition of $\kappa$.
\end{myprob}

\begin{myprob}
\myproblem{S24}{A. Tomita \cite{MR99k:22007}}
Is there (consistently) a free ultrafilter $p$ over $\omega$ such that
every $p$-compact group has a convergent sequence?
Is it consistent that for every free ultrafilter $p$ over $\omega$ there
exists a $p$-compact group without nontrivial convergent sequences?

\mynote{Notes}
Under $\text{\myaxiom{MA}}_\text{countable}$ there are $2^\mathfrak{c}$ 
many free ultrafilters $p$ such that there exists for each of them a
$p$-compact group without nontrivial convergent sequences
(A. Tomita and S. Watson).
\end{myprob}

\begin{myprob}
\myproblem{S25}{J.T.~Moore \cite{MR2002e:54026}}
Is it consistent to assume that every c.c.c.\ compact topological space
without a $\sigma$-linked base maps onto $[0,1]^{\omega_1}$? 
\end{myprob}

\section*{T. Algebraic and geometric topology}

\begin{myprob}
\myproblem{T1}{R. Stern \cite{MR80k:57038}}
Is $\theta^H_3$ finitely generated?

\mynote{Notes}
In problems T1--T4, let $\theta^H_3$ denote the abelian group obtained
from the set of oriented $3$-dimensional PL homology spheres using the
operation of connected sum, modulo those which bound acyclic PL
$4$-manifolds. 
Let $\alpha\colon \theta^H_3 \to \mathbb{Z}_2$ denote the 
Kervaire-Milnor-Rokhlin surjection.
\end{myprob}

\begin{myprob}
\myproblem{T2}{R. Stern \cite{MR80k:57038}}
Does $\theta^H_3$ contain an element of nontrivial finite order?
\end{myprob}

\begin{myprob}
\myproblem{T3}{R. Stern \cite{MR80k:57038}}
Is $\alpha$ an isomorphism?
\end{myprob}

\begin{myprob}
\myproblem{T4}{R. Stern \cite{MR80k:57038}}
Suppose a homology $3$-sphere $H^3$ admits an orientation reversing PL
homeomorphism.
Is it true that
$\alpha(H^3) = 0$?
$[H^3] = 0$ in $\theta^H_3$?
\end{myprob}

\begin{myprob}
\myproblem{T5}{J. Pak \cite{MR81i:55003}}
Let $\mathcal{J} = \{E,P,B,Y\}$ be an orientable Hurewicz fibering.
Is it true that if $E$ satisfies the $J$-condition, then $B$ and $Y$ do
also?
Is the converse question true?
\end{myprob}

\begin{myprob}
\myproblem{T6}{J. Pak \cite{MR81i:55003}}
Enlarge the class of Jiang spaces.

\mynote{Notes}
Jiang spaces are those that satisfy the Jiang condition from 
\cite{MR30:1510}.
\end{myprob}

\begin{myprob}
\myproblem{T7}{B. Clark \cite{MR84e:57005}}
Does longitudinal surgery on a knot $k$ always yield a manifold of maximal
Heegard genus among those that can be obtained by surgery on $k$?
\end{myprob}

\begin{myprob}
\myproblem{T8}{K. Perko \cite{MR85b:57011}}
Is every minimal-crossing projection of an alternating knot alternating?
\end{myprob}

\begin{myprob}
\myproblem{T9}{K. Perko \cite{MR85b:57011}}
Is the minimal crossing number additive for composition of primes?
\end{myprob}

\begin{myprob}
\myproblem{T10}{K. Perko \cite{MR85b:57011}}
Does the bridge number equal the minimal number for Wirtinger generators?

\mynote{Notes}
This has been resolved for two-bridged knots by M. Boileau.
\end{myprob}

\begin{myprob}
\myproblem{T11}{J. Pak \cite{MR88g:55005}}
Let $g\colon (M^n, x) \to (M^n,x)$ be a based homeomorphism on an 
$n$-dimensional manifold at $x \in M$.
If the induced homomorphism 
$g_*\colon \prod_k (M^n, x) \to \prod_k (M^n, x)$ 
is the identity map for all $k$, is then $g$ isotopic to the identity map?
How about if $M^n$ is an aspherical manifold?
\end{myprob}

\begin{myprob}
\myproblem{T12}{N. Lu}
In \cite{MR91f:57007} the presentations of the groups 
$\mathcal{M}_g, g \geq 3$ 
are not so simple as that of $\mathcal{M}_2$ given in \cite{MR91f:57006}.
The main reason is the extra Lantern law. 
Is there a simpler equivalent form from the Lantern law in the generators 
$L$, $N$, and $T$, or a more useful presentation of $\mathcal{M}_g$ for 
$g \geq 3$?
\end{myprob}

\begin{myprob}
\myproblem{T13}{N. Lu}
D. Johnson \cite{MR85a:57005} showed the Torelli groups
$\mathcal{M}_g$ are finitely generated for $g \geq 3$.
Is there a way to write Johnson's generators in terms of the generators
$L$, $N$, and $T$ \cite{MR91f:57006, MR91f:57007} of the surface mapping
class groups which will be useful in studying the fundamental group of
homology spheres?
\end{myprob}

\begin{myprob}
\myproblem{T14}{J. Stasheff \cite{MR95g:81164, MR96e:81225}}
The structure of a (based) loop space $\Omega X$ allows the reconstruction
of a space $BY$ of the homotopy type of $X$.
The parametrization of higher homotopies by the associahedra plays a
crucial role.
Does the joining of closed strings (= free loops) described in my talk 
lead in an analogous way to constructing from a free loop space
$Z = \mathcal{L}X$ a space of the homotopy type of $X$, perhaps with the
moduli space described in the article or some variant playing the role of
the associahedra?
\end{myprob}

\section*{U. Uniform spaces}

\begin{myprob}
\myproblem{U1}{R. Levy \cite{MR82i:54067}}
Which star-like subsets of $\mathbb{R}^2$ are $U$-embedded?

See the series of papers by R. Levy M. Rice
\cite{MR84d:54029,MR84k:26003,MR86f:54026,MR87a:54015,MR87m:54046,MR88d:54035}.
\end{myprob}

\begin{myprob}
\myproblem{U2}{S. Carlson \cite{MR88d:54023}}
If a proximity space admits a compatible complete uniformity, is it rich?
\end{myprob}

\begin{myprob}
\myproblem{U3}{C.R.~Borges \cite{MR89f:54059}}
If $(X, U)$ is topologically complete, is there a subgage $\theta$ for $U$
such that each $p \in \theta$ is a complete pseudometric?
\end{myprob}

\begin{myprob}
\myproblem{U4}{H.-P. K\"unzi \cite{MR98h:54038}}
Try to characterize those properties $P$ of quasi-uni\-form spaces
$(X, \mathcal{U})$ that fulfill the following condition:
$(X, \mathcal{U})$ has Property $P$ whenever
$(\mathcal{P}_0(X), \mathcal{U}_*)$ has Property $P$.
\end{myprob}

\section*{V. Geometric problems}

\begin{myprob}
\myproblem{V1}{M. Meyerson \cite{MR84d:52015}}
Can a square table be balanced on all hills (perhaps with negative 
heights) of compact convex support?
\end{myprob}

\begin{myprob}
\myproblem{V2}{M. Meyerson \cite{MR84d:52015}}
Can a cyclic quadrilateral table be balanced on all non-negative hills
with compact convex support?
\end{myprob}

\begin{myprob}
\myproblem{V3}{M. Meyerson \cite{MR84d:52015}}
Does every planar simple closed curve contain the vertices of a square?
\end{myprob}

\begin{myprob}
\myproblem{V6}{R. Pawlak \cite{MR97a:26016}}
This problem is motivated by the following theorem in the paper
\cite{MR97a:26016}: 
\emph{Let $A$ and $B$ be convex, non-singleton and strongly disjoint 
subsets of the plane.
Then $A$ possesses the property of a $D$-extension of a homeomorphism,
with the $u$-disc on $B$, if and only if $A$ and $B$ are closed.}

It seems interesting to ask the question whether the assumption of the
convexity of the sets $A$ and $B$ can be weakened in an essential way.
It could also be interesting to obtain a result analogous to the above 
theorem, where the domain of the transformations under consideration 
would be some metric space. 
Finally it is worthwhile to raise the question: can one construct 
appropriate Borel extensions (or measurable ones of class $\alpha$)?
\end{myprob}

\section*{W. Algebraic problems}

\begin{myprob}
\myproblem{W1}{N. Lu \cite{MR91f:57006}}
Call a group $G$ \emph{balanced} if it admits a finite set $s$ of 
generators so that any two elements of $s$ can be mapped to each other 
by some automorphism of $G$ which leaves $s$ invariant.
An example is the group $\mathcal{M}_2$ 
with $s = \{\Gamma_0, \ldots, \Gamma_5\}$, the set of six Dehn twists 
given in \cite[\S~3]{MR91f:57006}.
Characterize the balanced groups. 
\end{myprob}

\section*{X. Special constructions}

\begin{myprob}
\myproblem{X1}{G. Johnson \cite{MR85e:52002}}
If $(M, S)$ is a $G$-system, is $S$ connected?
\end{myprob}

\begin{myprob}
\myproblem{X2}{G. Johnson \cite{MR85e:52002}}
If $(M, S)$ is a $G$-system, $m$ is a set in $M$ which contains two points,
$\{s\} = S \cap m$, and $p \in m \setminus \{s\}$, is
$\{(1-t)s + tp : t \geq 0\}$ a subset of $m$?
\end{myprob}

\begin{myprob}
\myproblem{X3}{G. Johnson \cite{MR85e:52002}}
If $(M, S)$ is a $G$-system for $X$ and $\{w_i : i \geq 1\}$ is a 
convergent sequence in $X$, must $\{s_i : i \geq 1\}$ be a convergent 
sequence if $s_i$ and $w_i$ belong to the same set in $m$ for all $i$?
\end{myprob}

\section*{Y. Topological games}

\begin{myprob}
\myproblem{Y1}{I. Juh\'asz \cite{MR87m:54014}}
Is there a neutral point-picking game in \myaxiom{ZFC}?

\mynote{Notes}
Yes, if $\diamondsuit$ (A. Berner and I. Juh\'asz \cite{MR86c:54005}).
Yes, if $\text{\myaxiom{MA}}(\omega_1)$ for countable posets (Juh\'asz 
\cite{MR87m:54014}).
Yes, if \myaxiom{MA} for $\sigma$-centered posets (A. Dow and G. Gruenhage 
\cite{MR93b:54004}).
\end{myprob}

\begin{myprob}
\myproblem{Y2}{I. Juh\'asz \cite{MR87m:54014}}
Is there a space $X$ such that 
$\omega \cdot \omega < \operatorname{ow}(X) < \omega_1$?
\end{myprob}

\begin{myprob}
\myproblem{Y3}{I. Juh\'asz \cite{MR87m:54014}}
Does there exist, in \myaxiom{ZFC}, a $T_3$ space $X$ for which the games
$G^D_\omega(X)$ and/or $G^{SD}_\omega(X)$ are undecided?
\end{myprob}

\begin{myprob}
\myproblem{Y4}{I. Juh\'asz \cite{MR87m:54014}}
Is it true, in \myaxiom{ZFC}, that for every compact Hausdorff space $X$ 
and every cardinal $\kappa$ the game $G^D_\kappa$ is determined?
\end{myprob}

\begin{myprob}
\myproblem{Y5}{I. Juh\'asz \cite{MR87m:54014}}
Is there a space $X$ in \myaxiom{ZFC} such that 
$\text{II} \uparrow G^D_\alpha(X)$ for every $\alpha < \omega$, but 
$\text{II} \mathbin{\not{\uparrow_M}} G^D_\omega(X)$?
\end{myprob}

\begin{myprob}
\myproblem{Y6}{M. Scheepers \cite{MR95e:04009a}}
Let $\lambda$ be an uncountable cardinal of uncountable cofinality.
Let $\kappa$ be a cardinal such that
$\lambda^{<\lambda} < \operatorname{cf}([\kappa]^\lambda, \subset) \leq
2^\lambda$.
Does TWO have a winning remainder strategy in any of
$\text{WMEG}([\kappa]^\lambda)$, $\text{WMG}([\kappa]^\lambda)$
or $\text{VSG}[\kappa]^\lambda$?
\end{myprob}

\section*{Z. Topological dynamics, fractals and Hausdorff dimension}

\begin{myprob}
\myproblem{Z1}{P. Massopust \cite{MR89h:58087}}
What is the fractal dimension of $G = \operatorname{graph}(f)$ when
$f$ is a fractal interpolation function generated by polynomials or
general $C^0$-maps?
Is it possible to calculate the fractal dimension in this case by an
approximation scheme consisting of affine and/or polynomial maps?
\end{myprob}

\begin{myprob}
\myproblem{Z2}{P. Massopust \cite{MR89h:58087}}
What are the fractal dimensions of $A(I \times X)$ and 
$\operatorname{graph}(f^*)$, when $f^*$ is a hidden variable fractal 
interpolation function generated by affine, or even more general 
$C^0$-maps, rather than by similitudes?
Is it still true that $\dim A(I \times X) = \dim(X)$, or under what
conditions does this relation remain valid?
\end{myprob}

\begin{myprob}
\myproblem{Z3}{P. Massopust \cite{MR89h:58087}}
What is the exact Hausdorff-Besicovitch dimension for the graph of a
fractal interpolation and hidden variable interpolation function?
\end{myprob}

\begin{myprob}
\myproblem{Z4}{J. Graczyk and G. Swiatek \cite{MR99k:58156}}
Is there a complex bounds theorem for all real polynomials including the
polymodal ones?
In this case, does it help to assume that all critical values are real?
Note that in the polymodal case, it is not immediately clear what the 
statement of the theorem should be.

\mynote{Solution}
W. Shen showed how to define and prove complex bounds for all real 
analytic multimodal interval maps for which all critical points are of 
even order. 
But this restriction on the critical points can be eliminated by a recent 
joint work of S. van Strien and E. Vargas.
More precisely, Shen's proof begins with a careful analysis of the 
geometry of the postcritical sets by means of cross-ratio estimates and 
the related real Koebe principle, and then the complex bounds were 
concluded by modifying an earlier work of Lyubich and Yampolsky 
\cite{MR98m:58113}. 
The first part was only done for maps without inflection critical points 
in Shen's thesis \cite{Shen}, and can be completed for all maps by van 
Strien and Vargas's work \cite{vanstrien}.
\end{myprob}

\providecommand{\bysame}{\leavevmode\hbox to3em{\hrulefill}\thinspace}

\label{tpcontributedend}

\chapter*{Peter J.\ Nyikos: Classic Problems}
\label{tpclassic}
\begin{myfoot}
\begin{myfooter}
Peter J.\ Nyikos, 
\emph{Classic Problems},\\
Problems from Topology Proceedings, Topology Atlas, 2003, 
pp.\ 69--89.
\end{myfooter}
\end{myfoot}

\newtheorem{cproblem}{Classic Problem}
\renewcommand\thecproblem{\Roman{cproblem}}

\mypreface
The eight classic problems appeared in two articles by 
Peter~J.~Nyikos in volume 1 (1976) and volume 2 (1977) of 
\emph{Topology Proceedings}. 
Nyikos wrote two more articles, \emph{Classic problems---25 years later}, 
in volume 26 (2001--2) and volume 27 (2003) of \emph{Topology Proceedings} 
with detailed accounts of the progress on the eight problems.

In volume 1, Nyikos mentioned three problems that were omitted from this 
list because they were treated in M.E.\ Rudin's problem list:
the normal Moore space conjecture;
the $S$- versus $L$-space problem; 
the question of whether there are $P$-points in
$\beta\mathbb{N} \setminus \mathbb{N}$.
Two other problems were omitted because they had been solved at about 
that time.
\begin{myitemize}
\item
Is $\dim(X \times Y) \leq \dim X + \dim Y$ for completely regular 
spaces?
This was solved in the negative by M.\ Wage and T.\ Przymusi\'nski.
\item
(O.\ Frink) Is every Hausdorff compactification of a completely 
regular space a Wallman compactification?
A negative solution was given by V.M.\ Ul$'$\kern-.1667em janov.
\end{myitemize}

This version is an amalgamation of the four article by Nyikos. This
version contains the statements of all original classic problems and 
their related problems but most of the background material has been 
omitted, particularly for solved problems. The sections \emph{Consistency 
results} and \emph{References} are from the original 1976--1977 articles. 
The sections \emph{Twenty-five years later} are taken from two articles
\emph{Classic problems---25 years later}.

\section*{Introduction}

Of the eight classic problems, numbers II, III, and VIII have been solved
outright, with examples whose existence requires nothing more than the
usual (\myaxiom{ZFC}) axioms of set theory; numbers V and VI have been
shown \myaxiom{ZFC}-independent; numbers I and VII remain half-solved,
with consistent examples but no \myaxiom{ZFC} examples, and no consistency
results denying their existence. Finally, number IV, the well-known
$M_3$-$M_1$ problem, is completely unsolved---we do not even have
consistency results for it.

\section*{Classic Problem I}

\begin{cproblem}[Efimov's Problem]
Does every compact space contain either a nontrivial
convergent sequence or a copy of $\beta\mathbb{N}$?
\end{cproblem}

In this problem only, \emph{compact} will mean \emph{infinite compact 
Hausdorff}.

\subsection*{Equivalent problems}

Does every compact space contain 
(1)
a copy of $\omega + 1$ or 
a copy of $\beta\mathbb{N} \setminus \mathbb{N}$? 
(2)
a closed metric subspace or 
an infinite discrete $C^*$-embedded subspace? 

\subsection*{Related problems}

\subsubsection*{(1)}
Does every totally disconnected compact space contain either a copy of
$\omega + 1$ or a copy of $\beta\mathbb{N}$? 
Equivalently: Does an infinite Boolean algebra have either a countable
infinite or a complete infinite homomorphic image? 
\subsubsection*{(2)}
Does every compact space contain either a point with a countable
$\pi$-base or a copy of $\beta\mathbb{N} \setminus \mathbb{N}$?
\subsubsection*{(3)}
Does every compact hereditarily normal space contain 
a nontrivial convergent sequence? 
a point with a countable $\pi$-base? 
a point with a countable $\Delta$-base? 

A \emph{$\pi$-base at a point $x$} is a collection of open sets such 
that every neighborhood of $x$ contains one; 
a \emph{$\Delta$-base at $x$} is a $\pi$-base at $x$ such that every 
member has $x$ in its closure.

\subsection*{Consistency results}

Assuming \myaxiom{CH}, V.V.\ Fedor\v{c}uk constructed a compact space of 
cardinality $2^\mathfrak{c}$ so that every infinite closed subspace is of 
positive dimension. 
Since both $\omega + 1$ and $\beta\mathbb{N}$ are zero-dimensional, this 
space cannot contain a copy of either one. 
Assuming \visl, V.V.\ Fedor\v{c}uk constructed a space having all the 
above properties of his first space which is, in addition, hereditarily 
separable and hereditarily normal. 

\subsection*{References}

\cite{MR40:6505,MR50:14701,MR51:13966,MR53:3976,MR53:14380}

\subsection*{Twenty-five years later}
Efimov posed this problem a comparatively short time (roughly nine years) 
earlier \cite{MR40:6505}, but I deemed it worthy of being called a
classic even back then because of its remarkably fundamental nature.
The class of compact spaces is arguably the most important class of
topological spaces.
Its importance transcends general topology: functional analysts have
constructed many of their own with various analysis-relevant properties,
and rings of continuous functions on compact spaces have been studied for
well over half a century.
It is also an extremely broad and varied class of spaces, and I was amazed
when I first learned in 1974 that such a fundamental question was still
open.

Consistency results were not long in coming:
V.V.\ Fedor\v{c}uk showed in \cite{MR50:14701,MR54:13827} that there are 
counterexamples under both \myaxiom{CH} and
$\mathfrak{s} = \aleph_1$ $+$ $2^{\aleph_0} = 2^{\aleph_1}$.
These were already discussed in volumes 1 and 2.
Remarkably little progress has been made on Efimov's problem since then.
There are no examples just from \myaxiom{ZFC} and no good ideas as to how
to try to obtain some.
There are no consistency results in the opposite direction, although
\myaxiom{PFA} is a reasonable candidate for an affirmative answer.
So is at least one model of the Filter Dichotomy Axiom.

\subsubsection*{The Filter Dichotomy Axiom}
This axiom, which is \myaxiom{ZFC}-independent, says that every free
filter on $\omega$ can be sent by a finite-to-one function on $\omega$ to
either an ultrafilter or the cofinite filter.
These two kinds of free filters on $\omega$ are at opposite extremes,
just as $\omega+1$ and $\beta\omega$ are at opposite extremes among
compactifications of $\omega$.
A nontrivial convergent sequence and its limit point (in other words, the
space $\omega+1$) constitute the simplest infinite compact space, while
the Stone-\v{C}ech compactification $\beta\omega$ of $\omega$ is one of
the most complicated.
There is a real sense in which $\omega+1$ is the \emph{smallest} infinite
compact space while $\beta\omega$ is the \emph{largest} compact space with 
a countable dense subspace: every separable, infinite compact space maps
surjectively onto $\omega+1$ and is the continuous image of $\beta\omega$.

Another extreme contrast exists between the algebra of finite and cofinite
subsets of $\omega$ and the algebra $\mathcal{P}(\omega)$ of all subsets
of $\omega$, and the following problem is one translation, via Stone
duality, of the restriction of Efimov's problem to totally disconnected
spaces:

\startproblem
\begin{myprob}
\myproblem{Problem 1}{}
Does every infinite Boolean subalgebra of $\mathcal{P}(\omega)$ admit a
homomorphism onto either the finite-cofinite subalgebra or
$\mathcal{P}(\omega)$?
\end{myprob}

The restriction to subalgebras of $\mathcal{P}(\omega)$ is possible
because of the reduction of Efimov's problem to those compact spaces which
have the countable discrete space $\omega$ as a dense subspace: every
infinite space contains a copy of $\omega$, and the closure of such a copy
in a counterexample is itself a counterexample.

Despite these resemblances, the Filter Dichotomy Axiom is not enough for a
positive solution to Efimov's problem, because it is compatible with
$\mathfrak{s} = \aleph_1$ $+$ $2^{\aleph_0} = 2^{\aleph_1}$.
However, it is also compatible with $\mathfrak{s} = \aleph_2$ and so it 
holds out some hope.

There is a basic equivalence which leads naturally to problems related to
Efimov's problem and Problem 1: a compact space contains a copy of
$\beta\omega$ iff it maps onto $[0, 1]^\mathfrak{c}$.
For totally disconnected spaces we can substitute
$\{0, 1\}^\mathfrak{c}$ for $[0, 1]^\mathfrak{c}$.
One direction uses the fact that $\beta\omega$ is a projective object in
the class of compact spaces: if a compact space maps onto a space that
contains a copy of $\beta\omega$, then it also contains a copy.
The other direction features $\mathfrak{c}$ applications of the Tietze
Extension Theorem and a little categorical topology pertaining to the
universal property of a product space.
In the totally disconnected case, we use the fact that these compact
spaces are the ones of large inductive dimension zero; thus, any map 
from a closed subspace onto $\{0, 1\}$ induces a partition of the subspace 
into closed subsets which can then be enlarged to a clopen partition of 
the whole space.

Another application of Stone duality now shows that Problem 1 is
equivalent to:

\startproblem
\begin{myprob}
\myproblem{Problem 2}{}
Does every infinite Boolean algebra either contain a free subalgebra of 
cardinality $\mathfrak{c}$ or have a countably infinite homomorphic image?
\end{myprob}

This gives us some interesting related problems as soon as we deny
\myaxiom{CH} (which gives us counterexamples anyway!): just substitute
\emph{uncountable cardinality} for \emph{cardinality $\mathfrak{c}$} in
Problem 2, and ask:

\startproblem
\begin{myprob}
\myproblem{Problem 3}{}
Does every infinite compact space either contain a copy of $\omega+1$ or
admit a map onto $[0, 1]^{\aleph_1}$?
\end{myprob}

\subsubsection*{A new counterexample}
Recently, A. Dow showed that there is a counterexample to Efimov's
problem if $2^\mathfrak{s} < 2^\mathfrak{c}$ and the cofinality of the 
poset $([\mathfrak{s}]^\omega, \subset)$ is equal to $\mathfrak{s}$ 
(i.e., $\operatorname{cf}[\mathfrak{s}] = \mathfrak{s}$).
Roughly speaking, Dow's construction substitutes zero-sets for points
in Fedor\v{c}uk's \myaxiom{PH} construction \cite{MR54:13827}.
The construction can be done in \myaxiom{ZFC}, and results in an
infinite compact space with no convergent sequences.
The purpose of the second condition is to insure that the space has 
cardinality $2^\mathfrak{s}$, while the purpose of the condition 
$2^\mathfrak{s} < 2^\mathfrak{c}$ is to insure there is no copy of 
$\beta \omega$ in the space.

The axiom $\operatorname{cf}[\mathfrak{s}] = \mathfrak{s}$ is very general; 
its status is similar to that of the \emph{small} Dowker space of 
C.~Good which is discussed below in connection with Problem VII. 
That is, $\operatorname{cf}[\mathfrak{s}] = \mathfrak{s}$ unless there
is an inner model with a proper class of measurable cardinals.
That is because $\mathfrak{s}$ is of uncountable cofinality, and because 
the Covering Lemma over any model of \myaxiom{GCH} is already enough to 
insure that $\operatorname{cf}[\kappa] = \kappa$ for all cardinals 
\emph{except} cardinals of countable cofinality. 
Now the Core Model satisfies \myaxiom{GCH}, and it is known that there is 
an inner model with a proper class of measurable cardinals whenever the 
Covering Lemma over the Core Model (abbreviated 
$\text{Cov}(\text{\myaxiom{V}},\text{\myaxiom{K}})$) fails.

The following well-known argument that $\mathfrak{s}$ is not of countable
cofinality was pointed out by H. Mildenberger. 
Suppose $\kappa$ has cofinality $\omega$, and no subcollection of 
$\mathcal{P}(\omega)$ of cardinality $< \kappa$ is splitting.
Let $\mathcal{A}$ be a family of $\kappa$ subsets of $\omega$, and let
$\mathcal{A} = \bigcup\{A_n: n \in \omega\}$ with 
$|\mathcal{A}_n| \le \kappa$ for all $n$.
For each $n$, there is a set $B_n$ that is not split by any member of 
$\mathcal{A}_n$ and which satisfies $B_{n+1} \subset B_n$.
Then take an infinite pseudo-intersection of the $B_n$. 
This is a set that cannot be split by any member of $\mathcal{A}$.

A trivial modification of this argument shows that
$\operatorname{cf}(\mathfrak{s}) \ge \mathfrak{t}$. 
It is still not known whether $\mathfrak{s}$ is a regular cardinal.

The axiom that $2^\mathfrak{s} < 2^\mathfrak{c}$ is more restrictive,
but still quite general. 
For example, given regular cardinals $\kappa < \lambda$, there is an 
iterated c.c.c.\ forcing construction of a model where 
$\mathfrak s = \kappa$ and $\mathfrak{c} = \lambda$ 
\cite[5.1]{MR87f:54008.1}, 
where it is easy to see that the final model satisfies 
$2^{< \lambda} = \mathfrak{c}$ ( $< 2^\mathfrak{c}$). 
Even more simply, adding $\aleph_1$ Cohen reals to a model of 
$2^{\aleph_1} < 2^\mathfrak{c}$ results in a model where 
$\mathfrak{s} = \aleph_1$ and the other cardinals are not affected. 
Many other forcings have the same effect.

It might be worth mentioning here that Efimov's problem and 
Fedor\v{c}uk's constructions are of interest to analysts. 
M. Talagrand \cite{MR82g:46029} produced a Grothendieck space such that 
no quotient and no subspace contains $\ell_\infty$. 
A Banach space is called \emph{Grothendieck} if every weak${}^*$ 
convergent sequence in the dual space $X^*$ is also weakly convergent. 
Talagrand's example was the Banach space $C(K)$ for a compact space $K$ 
which contains neither $\omega + 1$ nor $\beta\omega$; it used 
\myaxiom{CH} for the construction.

A completely different application to analysis was done by
M. D\u{z}amonja and K. Kunen \cite{MR94h:54047}. 
They used $\diamondsuit$ to construct a compact $S$-space, with no copy 
of either $\omega + 1$ or $\beta\omega$, to give a hereditarily separable 
solution to the following problem: 
If $X$ is compact and supports a Radon measure with nonseparable measure 
algebra, then does $X$ map onto $[0, 1]^{\omega_1}$? 
They were able to make the measure algebra isomorphic to the one for 
$2^{\omega_1}$.

Piotr Koszmider has called my attention to a pair of Banach space
equivalents to $K$ having a copy of $\beta\omega$. 
One is that $C(K)$ (with the uniform topology) has $\ell_\infty$ as a 
quotient. 
The other is that $C(K)$ contains a subspace Banach-isomorphic to 
$\ell_1(\mathfrak{c})$.
We do not know of conditions on $C(K)$ equivalent to $K$ having a 
nontrivial convergent sequence; a necessary condition is that $C(K)$ 
has a complemented copy of $c_0$.

\subsubsection*{Related problems}
Of the related problems listed in volume 1, only one has been solved 
since then: \emph{Does every infinite compact hereditarily normal space 
contain a nontrivial convergent sequence?}
At the time, it was already known that $\diamondsuit$ implies a negative
answer \cite{MR53:14379}, and in 1990 it was shown that \myaxiom{PFA}
implies a positive answer \cite{MR96j:54027}.

Also, one other related problem had already been answered earlier:
B.\ Shap\-irov\-ski\u{\i} \cite{MR53:14380} proved that every compact
hereditarily normal space contains a point with a countable $\pi$-base.

\section*{Classic Problem II}

\begin{cproblem}
Is there a nonmetrizable perfectly normal, paracompact space with a point
countable base? 
\end{cproblem}

\subsection*{Related problems}

Which of the following implications holds for perfectly normal spaces
with point-countable bases: 
\subsubsection*{(1)}
normal implies collectionwise normal? 
\subsubsection*{(2)}
collectionwise normal implies paracompact? 
\subsubsection*{(3)}
paracompact implies metrizable? 
\subsubsection*{(4)}
non-Archimedean implies metrizable? 
\subsubsection*{(5)}
Lindel\"of implies metrizable? 

This last is equivalent to the question of whether every hereditarily
Lindel\"of regular space with a point-countable base is metrizable, and
also to whether it is separable. 
Moreover, it is equivalent to the question of whether every first
countable regular space which is of countable spread (in other words,
every discrete subspace is countable) is separable (F.D.\ Tall). 
Hence it is also equivalent to the question of whether every first
countable, hereditarily Lindel\"of regular space is hereditarily
separable. 

\subsection*{Consistency results}

A Souslin line, whose existence is independent of the usual axioms of set
theory, is a hereditarily Lindel\"of (hence perfectly normal) linearly
ordered (monotonically normal) space which is not metrizable. 

H.R.~Bennett: If there exists a Souslin line, there exists one with a
point-countable base. 

A.V.~Arhangel$'$\kern-.1667em ski{\u\i} and P.J.~Nyikos: There exists a hereditarily
Lindel\"of non-Archimedean space (and such a space necessarily has a 
point-countable base) which is not metrizable if, and only if, there
exists a Souslin line. 

E.~van~Douwen, F.D.~Tall, and W.~Weiss \cite{MR58:24187}: 
\myaxiom{CH} implies the existence of a hereditarily Lindel\"of space with 
a point-countable base which is not metrizable. 

J. Silver: \manotch\ implies the existence of a normal Moore (hence 
perfectly normal) space with a $\sigma$-point-finite base which is not 
metrizable, hence not collectionwise normal. 

\subsection*{References}

\cite{MR37:3534.1,MR58:24187,MR50:8459,MR35:7298,MR52:9176}

\subsection*{Twenty-five years later}

The answer is yes.
In 1988, S.\ Todor\v{c}evi\'c \cite{MR93k:54005} constructed an example 
in \myaxiom{ZFC}.
The problem was the third in a natural progression recounted by R.~Hodel 
in \cite{MR50:8459}.
C.E.~Aull had observed in 1971 that every perfectly normal space with a
$\sigma$-disjoint base is metrizable;
A.V.~Arhangel$'$\kern-.1667em ski{\u\i} had shown in 1963 that every perfectly normal,
collectionwise normal space with a $\sigma$-point-finite base is 
metrizable.
By weakening the base property but strengthening the covering property,
Hodel hoped to get another metrization theorem, at least consistently. 
But the example in \cite{MR93k:54005} shows this is not possible.
The example is actually quite simple, considering that the problem
remained open for fifteen years after Hodel publicized it.

The related problems mentioned along with Problem II have had
a varied history.
Remarkably enough, the word \emph{perfectly} adds nothing to our current 
knowledge about the first two related problems as far as consistency goes.

\subsubsection*{Related problem (1)}
The first related problem sits between the questions of whether every 
metacompact normal Moore space is metrizable and that of whether every 
first countable normal space is collectionwise normal.
No known model or axiom distinguishes between these two questions, which
revolve around large cardinal axioms.
Yes to both (and hence to the first related problem) is consistent if it
is consistent that there is a strongly compact cardinal; but there is a 
metacompact normal Moore space if the covering lemma holds over the Core
Model.
See \cite{MR2003c:54001} and \cite{MR86m:54023} for more on what this
means.
I take this opportunity to correct a misleading misprint in
\cite{MR2003c:54001}:
in the last section, all the $\diamondsuit_\kappa$'s should be
$\square_\kappa$'s.

\subsubsection*{Related problem (2)} 
I have no information on the question of whether every collectionwise
normal space with a point-countable base is paracompact, with or without
perfect normality.

\subsubsection*{Related problem (3)}
This is resolved by Todor\v{c}evi\'c example.

\subsubsection*{Related problem (4)} 
In the fourth related problem, it is \emph{point-countable base} that 
adds nothing.
It is shown in \cite{MR1999943} that if there is a non-Archimedean
perfectly normal space, there is one with a point-countable base.
This related problem is shown in \cite{MR1999943} to be equivalent to an
old problem of Maurice: does every perfectly normal LOTS have a 
$\sigma$-discrete dense subset?

\subsubsection*{Related problem (5)} 
As already remarked in volume 1, the branch space of a Souslin tree is a
consistent example for both the fourth and fifth related problems, and 
the latter is equivalent to the question of whether there is a first 
countable $L$-space.
This was shown to be independent by Szentmikl\'ossy, who showed that 
\manotch\ implies there are none.

\section*{Classic Problem III}

\begin{cproblem}
Is every screenable normal space paracompact? 
\end{cproblem}

A space is \emph{screenable} if every open cover has a $\sigma$-disjoint
refinement.

\subsection*{Equivalent problems}

Is every screenable normal space 
(1) countably paracompact? 
(2) $\theta$-refinable?
(3) countably $\theta$-refinable? 

\subsection*{Results}
K.\ Nagami \cite{MR16:1138d}: A screenable normal, countably paracompact 
space is paracompact. 
For normal spaces, the concepts of countable paracompact, countable 
metacompact, countable subparacompact and countable $\theta$-refinable 
are all equivalent.

\subsection*{Related problems}

\subsubsection*{(1)}
Is a screenable normal space collectionwise normal? 
(Note: it is collectionwise Hausdorff.)
\subsubsection*{(2)}
Is a screenable, collectionwise normal space paracompact? 
\subsubsection*{(3)}
Is a normal space with a $\sigma$-disjoint base paracompact? 
\subsubsection*{(4)}
Is a screenable normal space of nonmeasurable cardinality realcompact? 
\subsubsection*{(5)}
Is every collectionwise normal, weakly $\theta$-refinable space paracompact? 
\subsubsection*{(6)}
Is every normal weakly $\theta$-refinable space of nonmeasurable cardinality 
realcompact? 
countably paracompact?

\subsection*{References}

\cite{MR16:1138d}

\subsection*{Twenty-five years later}

The answer is no.
In 1996, Z.\ Balogh \cite{MR98g:54041} constructed a counterexample using
only \myaxiom{ZFC}.
The construction is very technical and uses elementary submodels heavily.
Prior to that, M.E.\ Rudin \cite{MR84d:54042} had constructed one, 
assuming $\diamondsuit^{++}$. 
Her example is somewhat simpler to describe than Balogh's but its
properties (especially normality) are much harder to show.

Any solution to Problem III has to be a Dowker space (a normal
space which is not countably paracompact) as shown by Nagami 
\cite{MR16:1138d}, who first posed the problem.
Both Rudin's and Balogh's spaces are collectionwise normal, providing
counterexamples for the second related problem. 
Rudin \cite{MR84a:54038} showed that if there is a normal screenable 
space that is not paracompact, there is one that is not collectionwise 
normal, answering the first related problem 
with the help of Balogh's example.
Since screenable spaces are weakly $\theta$-refin\-able, this also 
answers the sixth related problem which asked whether every normal weakly
$\theta$-refinable space is countably paracompact.
So does another, simpler example of Balogh \cite{MR96j:54026}: 
a hereditarily collectionwise normal space which is not paracompact but is 
the countable union of discrete subspaces, hence weakly $\theta$-refinable.
This (along with the original screenable examples) also answers another 
related problem, which asked whether every collectionwise normal weakly 
$\theta$-refinable space is paracompact. 
If one leaves off \emph{collectionwise} then any normal metacompact
space that fails to be collectionwise normal (such as Michael's 
subspace of Bing's Example G) is a counterexample. 

Another pair of related problems was whether a screenable (or weakly
$\theta$-refin\-able) normal space of nonmeasurable cardinality is
realcompact.
This referred to the old-fashioned definition of \emph{nonmeasurable
cardinality} that set theorists would express with \emph{smaller than 
the first uncountable measurable cardinal}.
I do not know whether either Balogh's or Rudin's screenable example is
realcompact.
In volume 1, de Caux published a non-realcompact, weakly
$\theta$-refinable Dowker space using $\clubsuit$.
More recently, C.\ Good \cite{MR95c:03127} produced an example using the
Covering Lemma over the Core Model.
Unless there are real-valued measurable cardinals, every weakly
$\theta$-refinable, normal, countably paracompact space is realcompact, so
Dowker spaces are required here too for \myaxiom{ZFC} counterexamples.

Finally, one related problem has become a classic in its own right:
\emph{Is every normal space with a $\sigma$-disjoint base paracompact?}
For this we have no consistency results whatsoever.
At various times, both Balogh and Rudin thought they had examples under 
various set-theoretic hypotheses, but withdrew their claims.

Also, Balogh's example raises the following question:
\emph{Is there a first countable normal screenable space which is not 
paracompact?}

\section*{Classic Problem IV}

\begin{cproblem}
Does every stratifiable space have a $\sigma$-closure-pre\-serv\-ing open
base? 
In other words, is every $M_3$ space $M_1$?
\end{cproblem}

\subsection*{Equivalent problems}

\subsubsection*{(1)}
Does any point in any stratifiable space have a closure-preserving local 
base of open sets? 
\subsubsection*{(2)}
Does any point in any stratifiable space have a 
$\sigma$-closure-preserving local base of open sets? 
\subsubsection*{(3)}
Does any closed set in a stratifiable space have a closure-preserving 
(or: $\sigma$-closure-preserving) neighborhood base of open sets? 

\subsection*{Related problems}

\subsubsection*{(1)}
Is the closed image of an $M_1$ space $M_1$?
\subsubsection*{(2)}
Is the perfect image of an $M_1$ space $M_1$? 
\subsubsection*{(3)}
Is the closed irreducible image of an $M_1$ space $M_1$ 
\subsubsection*{(4)}
Is every closed subspace of an $M_1$ space $M_1$ 
\subsubsection*{(5)}
Is every subspace of an $M_1$ space $M_1$? 

\subsection*{Partial results}

G.\ Gruenhage and H.\ Junnila \cite[Theorem 5.27]{MR86h:54038}: 
Every stratifiable space is $M_2$. 
An \emph{$M_2$-space} is one with a $\sigma$-closure preserving 
quasi-basis, a \emph{quasi-basis} being a collection of sets which 
includes a base for the neighborhoods of each point.

G.\ Gruenhage: Every $\sigma$-discrete stratifiable space is $M_1$.

C.R.\ Borges and D.J.\ Lutzer \cite{MR50:14681}: 
The irreducible perfect image of an $M_1$ space is $M_1$.

\subsection*{References}

\cite{MR32:6409,MR50:14681,MR24:A1707,MR56:6614,MR58:2734,MR45:5952}

\subsection*{Twenty-five years later}

Shortly after this problem was posed in \cite{MR24:A1707}, J.~Nagata 
is said to have predicted that it would still be unsolved ten years 
later.
Over forty years have elapsed, and we still do not have even consistency
results either way.
This does not mean, however, that little progress has been made on the
problem; in fact, as can be seen from \cite{MR91e:54074}, it is one of
the most extensively researched problems in general topology.
One sign of this is that all but one of the problems that I listed in
volume 1 as being related problems are, in fact, equivalent to it.
This is due to a powerful theorem known as the Heath-Junnila Theorem
\cite{MR82f:54045}: 
\emph{Every stratifiable space is the perfect retract of some $M_1$ 
space.}

Now, the class of stratifiable spaces is closed under the taking of
subspaces and of closed images.
The Heath-Junnila Theorem thus shows that the $M_3$-$M_1$ problem is 
equivalent to the question of whether every closed, or every perfect
image, or every subspace, or every closed subspace of an $M_1$ space is
$M_1$.

At the same time, the Heath-Junnila Theorem shows just how poor the 
currently known preservation properties of $M_1$ spaces are.
The best that we have is a theorem that was already known a quarter of
a century ago: the perfect irreducible image of an $M_1$ space is $M_1$;
and we still do not know whether \emph{perfect} can be improved to 
\emph{closed}. 
(This is the only one of the related problems in volume 1
which has not been shown equivalent to the $M_3$-$M_1$ problem.)

Nevertheless, a great many natural classes of stratifiable spaces have
been shown to be $M_1$.
One of the most general classes, which includes all stratifiable
sequential spaces, was established by T.\ Mizokami and N.\ Shimame
\cite{MR2001f:54032}.
They extended their class still further with the help of Y.\ Kitamura
\cite{MR1876154}.
In it, they prove that every WAP stratifiable space is $M_1$.
A space $X$ is said to be \emph{WAP} iff for every non-closed $A$ there 
is $x \in \overline{A} \setminus A$ and a subset $B$ of $A$ such that $x$ 
is the only point in the closure of $B$ which is not also in $A$. 
For example, sequential spaces are WAP and so are scattered spaces.
It seems to be unknown whether $C_k(X)$ is WAP for all Polish $X$.
In particular, it is unknown whether $C_k(\mathbb{P})$ is WAP, and it
is still an open problem whether $C_k(\mathbb{P})$ is $M_1$; 
it is known to be $M_3$ \cite{MR2002b:54026}.

Other classes of stratifiable spaces known to be $M_1$ are listed below,
and others can be found in the fine survey papers \cite{MR91e:54074}
and \cite{MR1229128}.

I hope the following analogy with some classic facts about metrizable
spaces will lead to a still better appreciation of this old problem.

The celebrated Nagata-Smirnov Theorem states that a regular space is
metrizable iff it has a $\sigma$-locally finite open base.
If \emph{open} is replaced by \emph{clopen} then we have de Groot's
characterization of the metrizable spaces of covering dimension zero.
Also basic to the theory is the Morita-Hanai-Stone Theorem, which has the
corollary that the perfect image of a metrizable space is metrizable.
Moreover, a space is metrizable iff it is a perfect image of a metrizable
space of covering dimension zero.

One might reasonably hope that the corresponding results continue to hold
if \emph{$\sigma$-locally finite} is weakened to \emph{$\sigma$-closure
preserving}.
A regular space with a $\sigma$-closure preserving base consisting of
clopen sets is known as an $M_0$-space, and one might therefore expect
that the $M_1$-spaces are precisely the perfect images of $M_0$-spaces,
while the $M_0$-spaces are the $M_1$-spaces of covering dimension zero. 
Also, since metrizable spaces are closed under the taking of subspaces and
perfect images, one might expect the $M_1$-spaces to have this property
too---and the Heath-Junnila Theorem would then tell us that the $M_1$
spaces coincide with the stratifiable spaces.

Unfortunately, this lovely theory is heavily besieged by a brutal gang of
cold facts:

In the forty-plus years of research on the $M_3$-$M_1$-space problem, 
only a few of the preservation properties of stratifiable spaces have been 
proven to be shared by $M_1$-spaces.
Besides the hole revealed by the Heath-Junnila Theorem, it is also not
known whether a stratifiable space which is the countable union of closed
$M_1$-subspaces is $M_1$ (nor, for that matter, whether every stratifiable
space is the union of countably many closed $M_1$-subspaces).

Perfect images of $M_0$ spaces are indeed all $M_1$, but the converse is
an open problem. 

The class of perfect images of $M_0$ spaces has somewhat better known 
preservation properties than that of $M_1$ spaces, being hereditary and
preserved under perfect images in addition to being countably productive.
However, it is not known whether it is preserved by closed maps; also, it
is not known whether every perfect image of an $M_0$ space that is of
covering dimension zero is an $M_0$ space.

True, we do not know of any stratifiable spaces that are not perfect
images of $M_0$-spaces, but there are quite a few intermediate classes
between these two, including strong $M_1$-spaces \cite{MR91a:54040};
closed images of $M_0$-spaces \cite{MR82h:54043};
hereditarily $M_1$-spaces;
stratifiable spaces in which every point has a closure-preserving open
neighborhood base \cite{MR84c:54050}, \cite{MR86g:54039};
$EM_3$-spaces [see below];
and, of course, the $M_1$-spaces themselves.
No two of the classes listed just now are known to coincide, and all but
perhaps the $EM_3$-spaces are also intermediate between the perfect images
of $M_0$-spaces and the $M_1$-spaces.

The best that has been done for the stratifiable spaces of covering 
dimension zero is to show that they are $EM_3$-spaces.
More generally, a space is $EM_3$ iff it is the perfect (or closed) image
of a stratifiable space of covering dimension zero \cite{MR84c:54052}.
Unfortunately, it is not known whether every $EM_3$-space is $M_1$, nor
whether every $M_1$-space is $EM_3$, nor whether every stratifiable
space is $EM_3$.

However, it is true that if the $EM_3$-spaces coincide with the 
$M_1$-spaces, then every stratifiable space is $M_1$.
This follows from the Heath-Junnila Theorem, together with the fact that
$EM_3$-spaces are preserved under perfect maps.
In fact, $EM_3$-spaces have all the nice preservation properties of
stratifiable spaces listed here and in the table on p.\ 382 of
\cite{MR91e:54074}.

On the other hand, it may be that a stratifiable space that is not $M_1$
is already at hand: the space $C_k(\mathbb{P})$ of continuous real-valued
functions defined on the irrationals, with the compact-open topology.
This space has a natural neighborhood base at the zero function
$\mathbf{0}$ formed by the open sets
$B(\mathbf{0}, K, \epsilon) 
= \{f : {\forall x \in K} (|f(x)| < \epsilon)\}$,
as $K$ ranges over the compact subsets of $\mathbb{P}$.
$C_k(\mathbb{P})$ was shown to be stratifiable by P. Gartside and
E. Reznichenko \cite{MR2002b:54026} by identifying $\mathbb{P}$ with
${}^\omega\omega$, and making use of the special case where the $K$ are of
the form
$f^\downarrow = \{g : g(x) \le f(x) \forall x \in {^\omega}\omega\}$
to obtain a $\sigma$-cushioned pairbase, the existence of which 
characterizes stratifiable spaces.
However, any subcollection of the natural base fails badly to be a 
$\sigma$-closure preserving base at $\mathbf{0}$.
Moreover, G.\ Gruenhage and Z.\ Balogh have shown that no finite union of
translates of basic open sets can be a base and
$\sigma$-closure preserving at the same time, while I have shown that no 
base at $\mathbf{0}$ formed by unions of sets
$B(\mathbf{0}, K, \epsilon)$ can be $\sigma$-closure-preserving.
Since every open set is a union of (countably many) translates of such
sets, it might seem as though we are close to showing $C(\mathbb{P})$
is a counterexample to the $M_3$-$M_1$ problem, but appearances can be 
deceiving! 

If it turns out that Problem IV has a negative solution, a
reasonable place to look for a substitute theorem is the replacement
of $M_1$ with $EM_3$.
By the earlier reasoning, this would follow if we could show that every
$M_1$-space is $EM_3$.
The best that can then be hoped for is (1) that all the intermediate
classes listed above coincide with either the class of perfect images of
$M_0$-spaces, or the class of stratifiable spaces, or the class of
$M_1$-spaces, and (2) that the perfect images of $M_0$ spaces are the 
stratifiable $\mu$-spaces. This would still give us an attractive theory,
because the $M_0$-spaces coincide with the stratifiable $\mu$-spaces
of covering dimension zero \cite{MR85h:54048}.
(A \emph{$\mu$-space} is a space which can be embedded in a countable 
product of paracompact $F_\sigma$-metrizable spaces---that is, of spaces 
that are the countable union of closed metrizable subspaces.)
However, the proliferation of intermediate classes, and the many open
problems it presents, do not allow for much optimism that things will work
out even this nicely.
For example, although every stratifiable $\mu$-space is a perfect image of
an $M_0$-space \cite{MR91e:54074}, the converse is an open problem. 

\section*{Classic Problem V}

\begin{cproblem}[A.V.\ Arhangel$'$\kern-.1667em ski\u{\i}]
Does every compact hereditarily normal (abbreviated $T_5$) space of
countable tightness contain a nontrivial convergent sequence? 
\end{cproblem}

In this classic problem and the next, \emph{space} means \emph{infinite 
Hausdorff space}. 
A space $X$ is of \emph{countable tightness} if 
$\overline{A} = 
\bigcup \{ \overline{B} : B \subset A, |B| \leq \aleph_0 \}$ 
for all $A \subset X$.

\subsection*{Related problems}
Is every separable compact $T_5$ space 
\begin{myitemize}
\item[(A)]
of countable tightness? 
\item[(B)]
of cardinal $\leq \mathfrak{c}$? 
\item[(C)]
sequentially compact? 
\item[(D)]
sequential?
\end{myitemize}

\subsection*{Equivalent Problems}

Let $P$ be a closed-hereditary property: that is, one that is true for
every closed subset of a space with the property.
The problem of whether every compact space satisfying $P$ contains a 
nontrivial convergent sequence is equivalent to that of whether every 
compactification of $\mathbb{N}$ satisfying $P$ contains a nontrivial 
convergent sequence.
The problem of whether every compact space satisfying $P$ is sequentially
compact is equivalent to that of whether every compactification of 
$\mathbb{N}$ satisfying $P$ has a point $x$ and a sequence of distinct 
points of $\mathbb{N}$ converging to $x$. 
Hence \emph{separable} is redundant in the third part of the last 
related problem.

Along the same lines, here is an implication which goes only one way:
if every separable compact space satisfying a closed-hereditary property
$P$ is sequential, then every compact space of countable tightness 
satisfying $P$ is sequential.

\subsection*{Consistency results}

Under Axiom $\Phi$ (which follows from $V = L$ and resembles 
$\diamondsuit$) V.V.~Fedor\v{c}uk has constructed a hereditarily 
separable (hence of countable tightness) compact $T_5$ space of 
cardinality $2^\mathfrak{c}$ which has no nontrivial convergent 
sequence.

If $2^{\aleph_0} < 2^{\aleph_1}$, then (F.B.\ Jones) every separable 
$T_5$ space is of countable spread.
Now B.\ Shapirovski\u{\i} and A.V.\ Arhangel$'$\kern-.1667em ski\u{\i} have shown 
independently that every compact space of countable spread is of 
countable tightness.
Thus, under $2^{\aleph_0} < 2^{\aleph_1}$, every separable compact $T_5$
space $X$ is of countable tightness.

It can also be shown, assuming $2^{\aleph_0} < 2^{\aleph_1}$, that any
compact $T_5$ space which does not contain an $S$-space is sequentially
compact, and if it has countable tightness it is then Fr\'echet-Urysohn,
hence sequential.

Under \manotch, every compact space of cardinality $< 2^\mathfrak{c}$ is 
sequentially compact (V.\ Malykhin and B.\ Shapirovski\u{\i}), so that 
yes to the second part of the last related problem implies yes to 
the third. 
On the other hand, it is not known whether every separable $T_5$ compact
space is of cardinality $< 2^\mathfrak{c}$ under \manotch.
In fact, it is a mystery what happens to any of these problems under 
\manotch.
It is not even known whether the Franklin-Rajagopalan space 
$\gamma \mathbb{N}$ (a compactification of $\mathbb{N}$ with growth 
$\omega_1 + 1$, hence not of countable tightness) can be $T_5$ under
\manotch.

\subsection*{References}

\cite{MR53:14377,MR51:13966,MR53:14379,MR44:972.1,MR49:7985.1,MR53:3976}

\subsection*{Twenty-five years later}
It is consistent that there is a positive solution to Problem V. 
See the comments below on Problem VI.

\subsubsection*{Related problems}

The biggest success story pertaining to any of the eight Classic Problems
has to do with Problem V. 
Not only is the problem itself solved, but all those listed under the 
heading of \emph{Related Problems} have also been solved. 

In Volume 2, it was explained how the axiom $2^{\aleph_0} < 2^{\aleph_1}$ 
gives a positive answer to (A) while Fedor\v{c}uk's construction under 
Axiom $\Phi$ (equivalent to $\diamondsuit$) \cite{MR53:14379} 
gives negative answers to (B), (C), and (D). 
\myaxiom{PFA} gives positive answers to all four parts 
\cite{MR93i:54003}. 
A model of \manotch\ was given in \cite{MR93i:54003} where (A) is 
answered 
negatively.

To find a still-open problem in the discussion of Problem V in volume 2, 
one has to look close to the end, where it is said, ``It is not known 
whether every separable compact $T_5$ space is of cardinal 
$< 2^\mathfrak{c}$ under \manotch.''
We do know from Jones's Lemma that $2^{|D|} \le \mathfrak{c}$ for any 
discrete subset $D$ of any separable $T_5$ space, and if we could 
substitute the Lindel\"of degree of any subspace for $|D|$ when the space 
is compact, we would be done. 
However, Szentmikl\'ossy's theorem that every compact space of countable 
spread is hereditarily Lindel\"of under \manotch\ does not generalize
to arbitrary spreads $< \mathfrak{c}$. 
We also do not know of any model of \manotch\ where (B) or (C) has a 
negative answer, so we have only halfway met the challenge in the 
continuation of the above quotation: 
``In fact, it is a mystery what happens to any of these problems
under \manotch.'' 
On the other hand, the final problem at the end of the discussion of 
Problem V in volume 2 has been solved: 
\manotch\ is compatible with some version of $\gamma \mathbb{N}$ being 
$T_5$ \cite{MR93i:54003}.

\section*{Classic Problem VI}

\begin{cproblem}[Moore and Mr\'owka]
Is every compact Hausdorff space of countable tightness sequential?
\end{cproblem}

A space $X$ is of \emph{countable tightness} if for every $A \subset X$,
$\bar{A} = \bigcup\{\bar{B}: B \subset A, |B| \leq \aleph_0\}$.
A subset $A$ of $X$ is \emph{sequentially closed} if no point of $X$ 
outside $A$ has a sequence in $A$ converging to it;
$X$ is \emph{sequential} if every sequentially closed subset is closed.

\subsection*{Related problems}

\subsubsection*{(A)}
Is there a hereditarily separable, countably compact, noncompact space?
\subsubsection*{(B)}
(B.\ Efimov)
Does a compact space of countable tightness have a dense set of points of
first countability?
\subsubsection*{(C)}
(A.\ Hajnal and I.\ Juh\'asz)
Is there a hereditarily separable compact space of cardinality 
$> \mathfrak{c}$?
\subsubsection*{(D)}
Is there a compact space of countable tightness that is not sequentially
compact?
\subsubsection*{(E)}
Is every separable, countably compact space of countable tightness compact?
What if it is locally compact?
\subsubsection*{(F)}
(S.P.\ Franklin and M.\ Rajagopalan)
Is every separable, first countable, countably compact (hence sequentially 
compact) space compact?
What if it is locally compact?

\subsection*{References}

\cite{MR41:7607,MR45:9286,MR54:8566,MR31:5184,MR55:4088,
mooremrowka,MR56:6630,MR51:4128}

\subsection*{Twenty-five years later}
Problem V is a double weakening of the more famous and older Problem VI, 
so they are best considered together.

In hindsight, Problem V may seem too specialized to be called a
\emph{classic}.
However, back in 1978 we were very much in the dark as to how well 
behaved compact spaces of countable tightness or compact $T_5$ spaces 
might be under \myaxiom{ZFC}-compatible axioms. 
Back then, we could not rule out the possibility that \myaxiom{ZFC} is 
enough to give a negative solution to Problem VI while Problem V is 
\myaxiom{ZFC}-independent. 
Also, we had no idea how long we would have to wait for a final solution 
to Problem VI even if it is \myaxiom{ZFC}-independent, and I felt that 
Problem V might give us a more attainable goal to shoot for in the 
interim.

We did have Fedor\v{c}uk's sensational 1975 construction under Axiom 
$\Phi$ (later shown equivalent to $\diamondsuit$) of an infinite compact 
$T_5$ hereditarily separable, hence countably tight, space with no 
nontrivial convergent sequences, so we knew a negative solution to both 
problems is consistent. 
But \myaxiom{PFA}, which turned out to imply a positive solution to 
Problem VI (and hence to V) had not even been formulated yet. 
The strongest general tool at our disposal in that direction was 
\manotch; and that is actually compatible with a negative solution to 
Problem VI \cite{MR89i:54034.1}. 
Even now, it is still not known whether \manotch\ is compatible with a 
negative solution to Problem V. 
Also, while we now know that a positive solution to Problem V is 
compatible with \myaxiom{CH}, the status of Problem VI under \myaxiom{CH} 
is still unsolved \cite{MR1992528.1} despite its being on the list of 
26 unsolved problems in \cite{MR80i:54005.1}. 
[The statement in volume 2 that Rajagopalan had constructed a compact 
non-sequential space of countable tightness from \myaxiom{CH} was 
incorrect.]

As it turned out, the solution to Problem V only predated the one for 
VI by a couple of months; but it could easily have been otherwise. 
The \myaxiom{PFA} solution to Problem V was the culmination of five 
months of intensive research by David Fremlin and myself beginning in 
March of 1986. 
We were working from combinatorial axioms derived from Martin's Maximum 
(\myaxiom{MM}),
which we soon narrowed down to one \cite[6.8]{MR89i:54034.1} that is 
now known to follow from \myaxiom{PFA}, and does not require large 
cardinals. 
One discovery by Fremlin led to another by myself, which in turn led to 
new discoveries by Fremlin (some of which appear in \cite{MR89i:03096}). 
This continued until, on the way to the 1986 Prague International 
Topological Symposium, I showed that this axiom implies that every 
compact $T_5$ space of countable tightness is sequential 
\cite{MR89i:54034.1}. 
In Prague, I gave a copy of my proof to Zolt\'an Balogh. 
Fremlin and I continued to work on Problem VI and our joint efforts 
resulted in a proof that every compact space of countable tightness is 
sequentially compact under \myaxiom{PFA}. 

There the matter might have rested for a long time, had not Balogh 
meanwhile looked closely at Fremlin's proof that \myaxiom{MM} implies the 
axiom we were using, and thought ``outside the box'' as Gary Gruenhage 
put it last year when calling Balogh's solution to Problem VI the first 
of ``Zoli's six greatest hits''. 
Balogh did it by mixing topology into Fremlin's proof and coming up with 
a modification that even broke new set-theoretic ground. 
His solution came right at the end of 1986 and can be found in 
\cite{MR89h:03088}; 
a simplified version of the proof, using elementary submodels, can be 
found in \cite{MR1229125}. 
Dow \cite{MR93i:03071.1} later showed that large cardinals are not 
necessary for these applications \myaxiom{PFA}.

\subsubsection*{Related problems}
 
All but the last two of these problems has been solved. 
In each of the other cases, Fedor\v{c}uk's Axiom $\Phi$ (equivalent to 
$\diamondsuit$) example \cite{MR53:14379} solves the problem one way, 
while \myaxiom{PFA} solves it the other way.
In the case of Related Problem C, \manotch\ is enough to solve it in 
the other direction, as was already explained in volume 2. 
In the case of Related Problem B, V. Malykhin showed that adding a single 
Cohen real is enough to produce a compact space $X$ of countable tightness
and $\pi$-character, in which every point of $X$ has character 
$\omega_1$ \cite{MR88g:54045.1}. 
In particular, if the ground model satisfies $\mathfrak{p} > \omega_1$
then $X$ is Fr\'echet-Urysohn.
I.~Juh\'asz \cite{MR90d:54005} showed that adding a single Cohen real 
results in a 
model where a weakening $(t)$ of $\clubsuit$ holds, and that $(t)$ is 
already enough to construct a space like Malykhin's.

The \myaxiom{PFA} solution to Related Problem A for regular spaces is
due to Baumgartner and Todor\v{c}evi\'c, who showed that there are no 
$S$-spaces compatible with \myaxiom{PFA} \cite{MR86g:03084}, 
\cite{MR90d:04001.4}. 
Clearly, every countably compact noncompact space is non-Lindel\"of and 
so a regular example for Related Problem A must be an $S$-space. 
For arbitrary (Hausdorff) spaces a slight modification of posets for the 
Moore-Mr\'owka problem \cite{MR89h:03088}, \cite{MR1229125} returns a 
negative \myaxiom{PFA} solution.

The \myaxiom{PFA} solution to Related Problem B is due to A. Dow 
\cite{MR91a:54003.1}, and the one to Related Problem D is due to Fremlin 
and myself as recounted above and in \cite{MR89i:54034.1}; the proof is 
similar to that of Statement 4 of \cite{MR93i:54003}, but also uses free 
sequences of length $\omega_1$ given by Statement D of \cite{MR93i:54003} 
to complete the \emph{centrifugal saturation}.

The status of Related Problems E and F is quite different from that of 
the others. 
There is a \myaxiom{ZFC} counterexample for the first part of Related 
Problem E \cite{MR93j:54015}, but it is not even Urysohn, let alone 
Hausdorff. 
For regular spaces, almost all of what we know is already to be found
in \cite{MR1078644.1}, including the information that almost every known 
regular counterexample for Statement E is also a counterexample for 
Statement F; 
that almost every published counterexample is also locally compact; 
and that this is one of the growing list of problems for which there are 
counterexamples if $\mathfrak{c}$ is either $\aleph_1$ or $\aleph_2$: 
there are counterexamples both if $\mathfrak{p} = \aleph_1$ and if 
$\mathfrak{b} = \mathfrak{c}$, and the well-known fact that 
$\mathfrak{p} \le \mathfrak{b}$ gives us no room for loopholes if 
$\mathfrak{c} \le \aleph_2$.

Incidentally, Related Problem F is one of my personal favorites. 
At the 1986 Prague International Toplogical Symposium I offered a prize
of 500 US Dollars for a solution, and raised it to \$1000 at the 1996 
Prague Toposym. 
Despite this, almost no progress has been made on it since 1986.

\section*{Classic Problem VII}

\begin{cproblem}
Does there exist a \emph{small} Dowker space?

More precisely, does there exist a normal space which is not countably
paracompact and is one or more of the following:
\begin{myitemize}
\item[A.]
first countable?
\item[B.] 
(hereditarily) separable?
\item[C.]
of cardinality $\aleph_1$?
\item[D.]
submetrizable?
\item[E.]
locally compact?
\end{myitemize}
\end{cproblem}

\subsection*{Related problems}

\subsubsection*{(1)}
Is there a pseudonormal space (a space such that two disjoint closed 
subsets, one of which is countable, are contained in disjoint open sets) 
which is not countably metacompact, and is one or more of the above?
\subsubsection*{(2)}
Is there a realcompact Dowker space?
\subsubsection*{(3)}
Is there a monotonically normal Dowker space?

\subsection*{Consistency results}

Assuming the existence of a Souslin line, M.E.\ Rudin 
\cite{MR54:3669, MR45:2660} constructed a hereditarily separable Dowker 
space and also one that is first countable and of cardinality $\aleph_1$, 
as well as realcompact.

Assuming $\clubsuit$, P. de Caux \cite{MR56:6629} constructed a Dowker 
space of cardinality $\aleph_1$ which is separable, locally countable, 
and weakly first countable. 
It is neither first countable nor locally compact nor realcompact, but it
is weakly $\theta$-refinable, collectionwise normal, and N-compact.

It is possible to construct a pseudonormal example with all these
properties except normality (and perhaps non-realcompactness), which is
not countably metacompact, and is collectionwise Hausdorff, by the
following axiom, obviously implied by $\clubsuit$:
\emph{To each countable limit ordinal $\lambda$ it is possible to assign 
a subset $T(\lambda)$ of $[0,\lambda]$ converging to $\lambda$, such that 
if $A$ is an uncountable subset of $\omega_1$ there exists $\lambda$ such
that $A \cap T(\lambda)$ is infinite.}
One simply uses the construction in \cite{MR56:6629}, substituting 
this assignment $T(\lambda)$ for the one given by de Caux.
 
\subsection*{References}

\cite{MR56:6629,MR55:1270,MR17:391e,MR42:5217}

\subsection*{Twenty-five years later}
The word \emph{small} is very informal and one person's list of 
properties might easily differ greatly from another's. 
Most people would probably agree that \emph{of cardinality $\le 
\mathfrak{c}$} has a greater claim to being called \emph{small} than 
submetrizability or local compactness.
Had I put it in, then the most significant advance on Problem VII in the 
last twenty-five years would arguably have been Balogh's \myaxiom{ZFC} 
example in \cite{MR96j:54026}.
[Its main competitor, as ably explained in the introduction to 
\cite{MR98m:54011}, would be a Dowker space shown in \myaxiom{ZFC} to 
be of cardinality $\aleph_{\omega+1}$.] 
As it is, the most significant is clearly C. Good's construction of a 
locally compact, locally countable (hence first countable) Dowker space 
under a higher-cardinal analogue of $\clubsuit$ that follows from 
$\text{Cov}(\text{\myaxiom{V}},\text{\myaxiom{K}})$
and hence requires very large cardinals for its negation 
\cite{MR95c:03127}. 
Good gave a general construction which also works under $\clubsuit$ to 
give an example that is, in addition, of cardinality $\aleph_1$. 
Moreover, it can be embedded in a separable example using the technique 
P.~de~Caux used at the end of his paper for his very similar example 
\cite{MR56:6629}.

Good used consequences of 
$\text{Cov}(\text{\myaxiom{V}},\text{\myaxiom{K}})$
similar to those employed by
W.~Fleissner for his solution of the bigger half of the normal Moore space 
problem \cite{MR86m:54023}. 
The smallest examples in either case have cardinality $\beth_\omega^+$. 
This is the successor of the first singular strong limit cardinal, 
$\beth_\omega$, which is the supremum of the sequence of cardinals 
$\beth_n$ where $\beth_0 = \aleph_0$ and $\beth_{n+1} = 2^{\beth_n}$.

Like de Caux's example, Good's examples are all countable unions of 
discrete subspaces. 
However, they are not submetrizable. 
On the other hand, the second example in \cite{MR55:1270} is 
submetrizable, as mentioned in volume 2 already. 

As already recounted in Volume 2, there is a construction of a Dowker 
space from \myaxiom{CH} that satisfies all but the last part of Classic 
Problem VII. 
See \cite{MR55:1270}, where a $\diamondsuit$ construction was announced 
that satisfies all five conditions simultaneously, including the 
hereditary version of condition B (call this version $B^+$, the other 
$B^-$). 
This does not seem to have ever appeared in print, but there is a 
$\diamondsuit$ construction in \cite{MR58:27475} that satisfies all the
conditions except D, submetrizabilty. 
One erroneous comment from \cite{MR55:1270} carried over to the volume 
2 discussion. 
It was claimed that the $\diamondsuit$ example is $\sigma$-countably 
compact, but there is no such thing as a $\sigma$-countably compact 
Dowker space.

We still do not have a locally compact Dowker space from \myaxiom{CH} 
alone.
On the other hand, I know of only two independence results directly 
bearing on Problem VII as stated.
One is that there is no first countable, locally compact, submetrizable 
example of cardinality $\aleph_1$ under \manotch. 
This is because of Balogh's theorem \cite{MR85b:54005} that under 
\manotch, every 
first countable, locally compact space of cardinality $\aleph_1$ either 
contains a perfect preimage of $\omega_1$ (hence cannot be submetrizable) 
or is a Moore space. 
Now, Moore spaces are countably metacompact, and normal spaces are 
countably paracompact iff they are countably metacompact.

The other independence theorem has little to do with Dowker-ness.
The \emph{hereditarily} version of Part B is consistently false because 
\myaxiom{PFA} implies that there are no $S$-spaces \cite{MR86g:03084} 
\cite{MR90d:04001.4} 
and so every hereditarily separable space is Lindel\"of and therefore
(countably) paracompact. 
In contrast, \myaxiom{PFA} actually implies the existence of first 
countable Dowker spaces, and is consistent with the existence of first 
countable, locally compact Dowker spaces \cite{MR86c:54018}: 
M. Bell's first countable example \cite{MR83e:03077} exists under 
$\mathfrak{p} = \mathfrak{c}$, which is implied by \manotch\ and hence by 
\myaxiom{PFA}; and Weiss constructed a locally compact first 
countable example assuming $\mathfrak{p} = \mathfrak{c} = \omega_2 +
\diamondsuit_\mathfrak{c}(\mathfrak{c}, \text{$\omega$-limits})$
\cite{MR83d:54036} \cite{MR86c:54018}, and this combination of axioms 
is known to be compatible with \myaxiom{PFA}. 
There are also examples of first countable Dowker spaces of cardinality 
$\aleph_1$ compatible with the Product Measure Extension Axiom 
(\myaxiom{PMEA}) \cite{MR95c:03127}.

Despite all this, we seem very far from any \myaxiom{ZFC} examples, 
except perhaps for Part D of Problem VII.
At the beginning of April, 2002, less than four months before his death,
I sent Zolt\'an Balogh an e-mail in which I asked him whether any of his 
Dowker examples were submetrizable. 
In his reply, which came the same day, he wrote,
``One of my Dowker space is almost submetrizable, and I somehow thought 
it could be made submetrizable. 
Give me a couple of weeks on that and I'll let you know.''
That was the last I ever heard from him. 
Part D of Problem VII remains unsolved as far as we know.\footnote{%
Dennis Burke found some handwritten notes by Zolt\'an Balogh dated 
April 25--May 1, 2002 in which he seems to be describing a \myaxiom{ZFC} 
example of a submetrizable Dowker space. 
It is too early to tell from the notes whether the example is correct.}

\subsubsection*{Related problems}
The answers to Related Problems (2) and (3) are yes \cite{MR2003c:54047.3} 
and no \cite{MR86c:54018}, respectively. 
As for (1), there is a \myaxiom{ZFC} example of a $2$-manifold which is 
pseudonormal but not countably metacompact in \cite{MR86f:54054.2}. 
Like all manifolds, it is locally compact and first countable 
(Parts A \& E). 
It is produced by adding half-open intervals to the open first octant in 
the square of the long line. 
A routine modification of the topology on the subspace of those points 
with ordinal coordinates, together with endpoints of the added intervals 
produces a first countable, locally compact pseudonormal space of 
cardinality $\aleph_1$ which is still not countably metacompact.
Finally, this subspace can be embedded in a separable example like 
Good's $\clubsuit$ example, still in \myaxiom{ZFC}, giving Parts 
A \& B${}^-$ \& C \& E.

I am unaware of any submetrizable (Part D) examples just from 
\myaxiom{ZFC}.
Locally compact, first countable, submetrizable ones of cardinality 
$\aleph_1$ (Parts A \& C \& D \& E) are ruled out just as they are for 
Dowker spaces. 
So too are hereditarily separable examples (Part B${}^+$).

\section*{Classic Problem VIII}

\begin{cproblem}
Is every $\gamma$-space quasi-metrizable?
\end{cproblem}

Let $X$ be a space and let $\tau$ be the collection of open subsets of $X$.
Let $g\colon \omega \times X \to T$ be a function such that for each $x$ 
and $n$, $x \in g(n,x)$.
A space $X$ is a \emph{$\gamma$-space} if it admits a $g$ such that for 
each $x$ and each $n$, there exists $m \in \omega$ such that if $y \in 
g(m,x)$, then $g(m,x) \subset g(n,x)$ and such that $\{g(n,x): n \in 
\omega\}$ is a local base at $x$.
A space $X$ is quasi-metrizable if, and only if, it is a $\gamma$-space 
with a function $g$ as above such that $m = n + 1$ for all $x$ and all $n$.

\subsection*{Equivalent problem}

Does every space with a compatible local quasi-uni\-form\-ity with 
countable base have a compatible quasi-uniformity with countable base?

\subsection*{Related problems}

\subsubsection*{(1)}
Is every paracompact (or Lindel\"of) $\gamma$-space quasi-metrizable?
\subsubsection*{(2)}
Is every $\gamma$-space with an orthobase quasi-metrizable?
\subsubsection*{(3)}
Is every linearly orderable $\gamma$-space quasi-metrizable?

\subsection*{Remarks}

These problems are probably not as well known as most of the others in
this subsection, but there are a number of reasons why the main one
deserves to be called a classic.
It is old enough, going back to Ribeiro's paper of 1943 where a theorem
which says in effect that every $\gamma$-space is quasi-metrizable is
given, but the proof is at best incomplete.
The concept of a $\gamma$-space has been \emph{discovered} independently 
by quite a few researchers over the years, and \cite{MR49:6173} lists 
five aliases and thirteen conditions equivalent to being a 
$\gamma$-space, some of them bearing little resemblance to that given 
here.
Moreover, consider the equivalent problem stated above:
if one drops \emph{quasi} in both places, one gets the classic 
metrization theorem of A.H.\ Frink, and there may be a neat general theory 
to be had if this \emph{quasi} analogue turned out to be right also. 
Not to mention the convenience of having one less kind of 
\emph{generalized metric space} to deal with.
On the other hand, a $\gamma$-space that is not quasi-metrizable would 
probably break some exciting new topological ground, as did Kofner's 
example several years ago of a quasi-metrizable space which does not 
admit a non-Archimedean quasi-metric.

\subsection*{References}

\cite{frink,MR55:9040,MR48:3014,MR49:6173,MR5:149b}

\subsection*{Twenty-five years later}
The answer is no. 
R. Fox \cite{MR83h:54035} came up with a machine which outputs a 
$\gamma$-space with each $\gamma$-space input, and which produces 
non-quasi-metrizable spaces in certain cases. 
It preserves the Hausdorff separation axiom, but not regularity.
Together with J. Kofner, Fox \cite{MR86m:54033} found a Tychonoff example 
which is quasi-developable and scattered. 
In a note added in proof to their article, they announced the 
construction of a paracompact $\gamma$-space that is not 
quasi-metrizable. 
Now, H.-P.~K\"unzi \cite{MR2003d:54001} has done us the service of 
publishing a description of the example and an outline of the proof that 
it works.

\subsubsection*{Related Problems}
The answer to the \emph{paracompact} part of (1) is yes as recounted 
above. 
For \emph{Lindel\"of} it is still open. 
We also do not have a \myaxiom{ZFC} example of a Lindel\"of 
$\gamma$-space that is not non-Archimedeanly quasi-metrizable. 
A Luzin subset of the Kofner plane \cite{MR48:3014} 
\cite[Example~1]{MR82j:54054} is a consistent example: 
see \cite[Proposition~5]{MR81m:54060.1}, which was misstated 
with the omission of ``not" before ``non-Archimedean".

Kofner also provided affirmative answers to (2) \cite{MR83c:54039} and 
(3) \cite{MR82h:54050}. 
In both cases, Kofner used the fact that every $k$-transitive 
$\gamma$-space is non-Archimedeanly quasi-metrizable, for any integer 
$k$. 
The former proof uses the fact that any space with an orthobase is 
$2$-transitive, while the latter uses the fact that every GO-space is 
$3$-transitive.
His article \cite{MR82j:54054} for \emph{Topology Proceedings} is a 
very nice survey of the state of the art at the time.

\nocite{MR84k:54026}

\providecommand{\bysame}{\leavevmode\hbox to3em{\hrulefill}\thinspace}

\label{tpclassicend}

\chapter*{New Classic Problems}
\label{tpnewclassic}
\begin{myfoot}
\begin{myfooter}
Zolt\'an T.~Balogh, Sheldon W.~Davis, Alan Dow, Gary Gruenhage, 
Peter J.~Nyikos, Mary~Ellen Rudin, Franklin D.~Tall and Stephen Watson.
\emph{New classic problems},\\
Problems from Topology Proceedings, Topology Atlas, 2003, 
pp.\ 91--102.
\end{myfooter}
\end{myfoot}

\mypreface
These problems were published in volume 15 (1990) of \emph{Topology 
Proceedings}. 
Many of these problems also appeared in the book \emph{Open problem in 
topology}, edited by J.~van~Mill and G.M.~Reed \cite{MR92c:54001.1}.
This version contains the problems from the original article with 
current notes on solutions.
The original contributors have authorized these new versions.

\section*{Introduction}

Mary Ellen Rudin and Frank Tall organized a problem session at the Spring 
Topology Conference in San Marcos, Texas in 1990 and 
invited several people to come up with their ideas for problems that 
should be the worthy successors to the S \& L problems, the box product 
problems, the normal Moore space problems, etc.\ in the sense that they 
could and should be the focus of common activity during the 
1990s as the older problems had been during the 
1970s.
They hoped that these problems would counterbalance the more centrifugal 
1980s, during which there was a tendency for each 
set-theoretic topologist to do his own thing, rather than there being many 
people working on problems generally recognized as important.
Time will tell whether the title is appropriate.

\section*{Zolt\'an T.~Balogh: A problem of Kat\v{e}tov}

Given a topological space $X$, let $\operatorname{Borel}(X)$ and 
$\operatorname{Baire}(X)$ denote the $\sigma$-algebras generated by the 
families
$\operatorname{closed}(X) = \{ F : \text{$F$ closed in $X$}\}$ and
$\operatorname{zero}(X) = \{ F : \text{$F$ is a zero set in $X$}\}$, 
respectively.
The following question is due, without the phrase ``in \myaxiom{zfc}'', 
to M.~Kat\v{e}tov \cite{MR14:27c}.

\startproblem
\begin{myprob}
\myproblem{Problem}{M.~Kat\v{e}tov \cite{MR14:27c}}
Is there, in \myaxiom{zfc}, a normal $T_1$ space $X$ such that 
$\operatorname{Borel}(X) = \operatorname{Baire}(X)$ 
but $X$ is not perfectly normal 
(i.e., $\operatorname{closed}(X) \ne \operatorname{zero}(X)$)?
What if $X$ is also locally compact? 
first countable? 
hereditarily normal?

\mynote{Notes}
There are several consistency examples given by Z.~Balogh in 
\cite{MR89k:54040}.
\myaxiom{ch} implies that there is a locally compact locally countable $X$ 
satisfying the conditions of the problem.
The existence of a first countable, hereditarily paracompact $X$ is 
consistent, too.

However, as summarized by the following theorem, 
a space giving a positive answer to the question cannot satisfy 
certain properties.

\begin{theorem}
Let $X$ be a normal $T_1$ space, and let $A$ be a closed Baire 
subset of $X$.
Then $A$ is a zero set in $X$ if one of the following conditions 
hold:
\begin{myitemize}
\item
$X$ is compact \emph{(P.R.~Halmos \cite{MR11:504d})}.
\item
$X$ is paracompact and locally compact \emph{(W.W.~Comfort \cite{MR41:3695})}.
\item
$X$ is submetacompact and locally compact \emph{(D.~Burke)}.
\item
$X$ is Lindel\"of and \v{C}ech-complete \emph{(W.W.~Comfort \cite{MR41:3695})}.
\item
$X$ is a subparacompact $P(\omega)$-space \emph{(R.W.~Hansell 
\cite{MR87d:54065})}.
\end{myitemize}
\end{theorem}
\end{myprob}

\begin{myprob}
\myproblem{Problem}{K.A.~Ross and K.~Stromberg \cite{MR33:4224}}
If $X$ is a a normal locally compact Hausdorff space and $A$ is a closed 
Baire set in $X$, is $A$ a zero set?

\mynote{Solution}
In \cite{MR89k:54040}, Balogh gave a counterexample to this
related problem of K.A.~Ross and K.~Stromberg.
The construction makes use of the technique of E.~van~Douwen and 
H.H.~Wicke \cite{MR55:6390} and W.~Weiss \cite{MR81f:54021}.
\end{myprob}

\section*{Sheldon W.~Davis: Questions}

\begin{myprob}
\myproblem{Question 1}{}
Is there a symmetrizable Dowker space?

\mynote{Notes}
If $X$ is such a space, then let $\langle F_n : n \in \omega \rangle$ be a 
decreasing sequence of closed sets with $\bigcap_n F_n = \emptyset$ which 
cannot be \emph{followed down} by open sets, then attach $x_\infty \not\in X$ 
to $X$ and extend the symmetric so that $B(x_\infty, \frac{1}{n}) = F_n$, 
and the resulting space has a point, $x_\infty$, which is not a 
$G_\delta$ set.
This answers an old question of Arhangel$'$\kern-.1667em ski\u{\i} and Michael.
\end{myprob}

\begin{myprob}
\myproblem{Question 2}{A.V.~Arhangel$'$\kern-.1667em ski\u{\i}, E.~Michael}
Is every point of a symmetrizable space a $G_\delta$ set?

\mynote{Results}
S.W.~Davis, G.~Gruenhage and P.~Nyikos \cite{MR80a:54052}:
\begin{myitemize}
\item
There is a $T_3$ zero-dimensional symmetrizable space with a closed set 
which is not a $G_\delta$ (also not countably metacompact).
\item
There is a $T_2$ symmetrizable space with a point which is not a 
$G_\delta$ (constructed as above).
\item
In the example above, the sequential order, $\sigma(X)$, is $3$.
\item
If $X$ is $T_2$ symmetrizable and $\sigma(X) \leq 2$, then each point of 
$X$ is a $G_\delta$ set.
\end{myitemize}
R.M.~Stephenson \cite{MR56:1260, MR81m:54056.1}:
\begin{myitemize}
\item 
If $X$ is $T_2$ symmetrizable and $X \in X$ is not a $G_\delta$ set, then 
$X \setminus \{x\}$ is not countably metacompact.
\item
If $X$ is a regular feebly compact space which is not separable, then $X$ 
has a point which is not a $G_\delta$ set.
\end{myitemize}
D.~Burke, S.W.~Davis \cite{MR82k:54048, MR85d:54033}:
\begin{myitemize}
\item
$\mathfrak{b} = \mathfrak{c}$ implies that every regular symmetrizable 
space with a dense conditionally compact subset is separable.
\item
$\mathfrak{b} = \mathfrak{c}$ implies that every feebly compact regular 
symmetrizable space with a dense set of points of countable character is 
first countable.
\item
Let $X$ be a $T_2$ symmetrizable space.
If $x \in X$ and $\kappa$ is a cardinal of uncountable cofinality
with $\chi(x,X) \leq \kappa$, then
$\psi(x,X) < \kappa$.
Hence, an absolute example must be nonseparable and in fact have 
character $> \mathfrak{c}$.
\end{myitemize}
Y.~Tanaka \cite{MR84b:54049}:
There is a regular symmetrizable $X$ with $\chi(X) > \mathfrak{c}$.
However, this example is perfect.
\end{myprob}

\begin{myprob}
\myproblem{Question 3}{A.V.~Arhangel$'$\kern-.1667em ski\u{\i}, M.E.~Rudin}
Is every regular Lindel\"of symmetrizable space separable?
Equivalently, is there a symmetrizable $L$-space?

\mynote{Results}
S.~Nedev \cite{MR35:7293}:
Lindel\"of symmetrizable spaces are hereditarily Lindel\"of;
No symmetrizable $L$-space can have a weakly Cauchy symmetric.

J.~Kofner \cite{MR49:3850}, S.W.~Davis \cite{MR84a:54052, MR86h:54039}:
No symmetric $L$-space can have a structure remotely resembling a weakly 
Cauchy symmetric.

I.~Juh\'asz, Z.~Nagy, Z.~Szentmikl\'ossy \cite{MR87a:54041}:
\myaxiom{CH} implies that there is a $T_2$ non-regular symmetrizable 
space which is hereditarily Lindel\"of and nonseparable.

D.~Shakhmatov \cite{MR89c:54009.1}:
There is a model which contains a regular symmetrizable $L$-space.

Z.~Balogh, D.~Burke, S.W.~Davis \cite{MR91e:54076.1}:
There is (in \myaxiom{ZFC} alone) a $T_2$ non-regular symmetrizable 
space which is hereditarily Lindel\"of and nonseparable;
There is no left separated Lindel\"of symmetrizable space of uncountable 
cardinality.
\end{myprob}

\section*{Alan Dow: Questions}

A point $p \in X$ is a \emph{remote} point of $X$ if $p$ is not in the 
closure of any nowhere dense subset of $X$.
It is known that pseudocompact spaces do not have remote points 
(T.~Terada \cite{MR81b:54026}; A.~Dow \cite{MR90a:54066})
and that not every non-pseudocompact space has a remote point 
(E.~van~Douwen and J.~van~Mill \cite{MR81b:54038}).
Every non-pseudocompact metric space has remote points 
(S.B.~Smith and J.H.~Smith \cite{MR81m:54037})
(or of countable $\pi$-weight (E.~van~Douwen \cite{MR83i:54024})), 
but the statement \emph{every non-pseudocompact space of weight 
$\aleph_1$ has remote points} is independent of \myaxiom{zfc} 
(A.~Dow \cite{MR84f:54031}; K.~Kunen, J.~van~Mill and C.F.~Mills 
\cite{MR80h:54029}).
There is a model in which not all separable non-pseudocompact spaces have 
remote points (A.~Dow \cite{MR90a:54066}).
It follows from \myaxiom{ch} that all non-pseudocompact c.c.c.\ 
spaces of weight at most $\aleph_2$ have remote points (A.~Dow).

\begin{myprob}
\myproblem{Question 1}{}
Does if follow from \myaxiom{ch} (or is it consistent with \myaxiom{ch}) 
that if a non-pseudocompact space $X$ has some nonempty open subset that 
is c.c.c.\ and has non-pseudocompact closure then $X$ has remote points?

\mynote{Notes}
See Dow's article \cite{MR90a:54066}.
In \cite{MR1078637}, Dow conjectures that \myaxiom{ch} implies that all 
non-pseudocompact c.c.c.\ spaces of weight less than $\aleph_\omega$ have 
remote points.
\end{myprob}

\begin{myprob}
\myproblem{Question 2}{}
Is there a compact nowhere c.c.c.\ space $X$ such that $\omega \times X$
has remote points?

\mynote{Notes}
This question is discussed in Dow's article \cite[Problem 2]{MR1078637}:
``Of course there may not be a reasonable answer to this question in 
\myaxiom{zfc}, but it may be possible to obtain a nice characterization 
under such assumptions as \myaxiom{ch} or \myaxiom{pfa}.
For example, I would conjecture that there is a model satisfying that if 
$X$ is compact and $\omega \times X$ has remote points then $X$ has an 
open subset with countable cellularity.
See \cite{MR84f:54031, MR88m:03073}.''

\mynote{Solution}
A.~Dow \cite{MR1971307} showed (in \myaxiom{zfc}) that there is a compact 
nowhere c.c.c.\ space $X$ such that $\omega \times X$ has remote points
\end{myprob}

\begin{myprob}
\myproblem{Question 3}{}
Is there, for every space $X$, a cardinal $\kappa$ such that 
$\kappa \times X$ has remote points?

\mynote{Notes}
This is Problem 3 of Dow's list \cite{MR1078637}.
A.~Dow and T.J.~Peters \cite{MR90a:54069.1} showed that this is true if 
there are arbitrarily large cardinals $\kappa$ such that 
$2^\kappa = \kappa^+$.
\end{myprob}

\begin{myprob}
\myproblem{Question 4}{}
Are there weak $P_{\omega_2}$-points in $U(\omega_1)$, the space of
uniform ultrafilters on $\omega_1$?

\mynote{Solution}
Yes. 
J.~Baker and K.~Kunen \cite{MR2002j:54025} proved that if $\kappa$ is 
regular, then there is a uniform ultrafilter in $U(\kappa)$ which is a 
weak $P_{\kappa^+}$-point in $U(\kappa)$ and hence a weak 
$P_{\kappa}$-point in $\beta \kappa$.
The weak $P_\kappa$-point problem for singular $\kappa$ is still open.
\end{myprob}

\begin{myprob}
\myproblem{Question 5}{}
Do there exist points $p, q \in U(\omega_1)$
such that there are embeddings $f$, $g$ of $\beta \omega_1$ with
$f(p) = g(q)$, but no embedding takes $p$ to $q$ or $q$ to $p$?

\mynote{Notes}
If so, then $\beta \omega_1$ fails to have the Frol\'ik property 
(introduced in \cite{MR92k:54048}).
\end{myprob}

\begin{myprob}
\myproblem{Question 6}{}
Does there exist a compact zero-dimensional $F$-space (or basically 
disconnected space) which cannot be embedded into an extremally
disconnected (ED) space?

\mynote{Notes}
This is Problem 9 from Dow's list \cite{MR1078637}.
E.~van~Douwen and J.~van~Mill \cite{MR81b:54038} showed that it is 
consistent that there is a strongly zero-dimensional $F$-space that cannot 
be embedded in any basically disconnected space. 
A.~Dow and J.~Vermeer \cite{MR93b:54036} proved that it is consistent 
that the $\sigma$-algebra of Borel sets of the unit interval is not the 
quotient of any complete Boolean algebra.
By Stone duality, there is a basically disconnected space of weight 
$\mathfrak{c}$ that cannot be embedded into an extremally disconnected 
space.
\end{myprob}

\section*[Gary Gruenhage: Homogeneity of $X^\infty$]{Gary Gruenhage: 
Homogeneity of $\mathbf{X}^\mathbf{\infty}$}

\begin{myprob}
\myproblem{Problem}{}
Is $X^\infty$ homogeneous for every zero-dimensional first countable
regular space $X$?
What if $X$ is compact?
What if $X$ is a zero-dimensional subspace of the real line?

\mynote{Solution}
Yes, zero-dimensional subspaces of the real line have homogeneous 
$\omega$-power (B.~Lawrence \cite{MR98k:54061}).
In general, zero-dimensional first countable spaces have homogeneous
$\omega$-power (A.~Dow and E.~Pearl \cite{MR97j:54008}).
\end{myprob}

\section*[Peter J.~Nyikos: Dichotomies in compact spaces and $T_5$ 
spaces]{Peter J.~Nyikos: Dichotomies in compact spaces and 
$\mathbf{T}_\mathbf{5}$ spaces}

\begin{myprob}
\myproblem{Problem 1}{Efimov's Problem}
Is there an infinite compact $T_2$ space which contains neither a 
nontrivial convergent sequence nor a copy of $\beta \omega$?

\mynote{Notes}
This is Classic Problem I.
\end{myprob}

\begin{myprob}
\myproblem{Problem 2}{Zero-dimensional version}
Is there an infinite Boolean algebra (BA) which has neither a countably
infinite homomorphic image nor a complete infinite homomorphic image?
\end{myprob}

\begin{myprob}
\myproblem{Problem 2$'$}{}
Is there an infinite Boolean algebra (BA) which has neither a countably
infinite homomorphic image nor an independent subset of cardinality
$\mathfrak{c}$?
\end{myprob}

\begin{myprob}
\myproblem{Problem 3}{}
Is there an infinite compact $T_2$ space which cannot be mapped onto
$[0,1]^{\omega_1}$ and in which every convergent sequence is eventually
constant?
\end{myprob}

\begin{myprob}
\myproblem{Problem 4}{Hu\v{s}ek's problem}
Does every infinite compact $T_2$ space contain either a nontrivial 
convergent $\omega$-sequence or a nontrivial convergent 
$\omega_1$-sequence.
\end{myprob}

\begin{myprob}
\myproblem{Problem 5}{I.~Juh\'asz}
Does every infinite compact $T_2$ space contain either a point of first
countability or a convergent $\omega_1$-sequence.
\end{myprob}

\begin{myprob}
\myproblem{Problem 6}{}
Does every infinite compact $T_2$ space have a closed subspace with a 
nonisolated point of character $\leq \omega_1$?
\end{myprob}

\begin{myprob}
\myproblem{Problem 7}{}
Is every infinite BA of altitude $\leq \omega_1$ of pseudo-altitude 
$\leq \omega_1$?
\end{myprob}

\begin{myprob}
\myproblem{Problem 8}{}
Is \manotch\ (or even $\mathfrak{p} > \omega_1$) 
compatible with the existence of an infinite $T_2$ compact space of 
countable tightness with no nontrivial convergent sequences?
\end{myprob}

\begin{myprob}
\myproblem{Problem 9}{}
Is there a \myaxiom{zfc} example of a separable, hereditarily normal, 
locally compact space of cardinality $\aleph_1$?
\end{myprob}

\begin{myprob}
\myproblem{Problem 9$'$}{}
Is there a locally compact, locally countable, hereditarily normal 
$S$-space in every model of $\mathfrak{q} = \omega_1$?
\end{myprob}

\begin{myprob}
\myproblem{Problem 9$''$}{}
Is there a \myaxiom{zfc} example of a separable, hereditarily normal, 
locally compact, uncountable scattered space?

\mynote{Notes}
The answer to all of these is negative.
T.~Eisworth. P.~Nyikos and S.~Shelah \cite{MR1998110} showed that there 
is a model of $2^{\aleph_0} < 2^{\aleph_1}$ in which there are no first 
countable, locally compact $S$-spaces. 
Note that $2^{\aleph_0} < 2^{\aleph_1}$ implies $\mathfrak{q} = 
\omega_1$, and that every locally compact, locally countable Hausdorff 
space is first countable. 
Thus, there is a model where the answer to Problem 9$'$ is negative.
\end{myprob}

\begin{myprob}
\myproblem{Problem 10}{}
Is there a \myaxiom{zfc} example of a separable, uncountable, scattered 
hereditarily normal space?
\end{myprob}

\begin{myprob}
\myproblem{Problem 10$'$}{}
Is there a model of $2^{\aleph_0} < 2^{\aleph_1}$ in which there are no
hereditarily normal $S$-spaces?
\end{myprob}

\begin{myprob}
\myproblem{Problem 11}{}
Is it consistent that every separable compact hereditarily normal space is 
of character $< \mathfrak{c}$?
\end{myprob}

\section*{Mary Ellen Rudin: The linearly Lindel\"of problem}

\begin{myprob}
\myproblem{Problem}{\cite[A.~Mi\v{s}\v{c}encko]{MR25:4483}, 
\cite[N.~Howes]{Howes}}
Does there exist a non-Lindel\"of normal space $X$ such that every
increasing open cover of $X$ has a countable subcover?

\mynote{Notes}
This question has remained unanswered for about 40 years.
No significant partial results are known.
This is Problem 328 in Rudin's list \cite{MR1078646}.

An open cover $\mathcal{U}$ is \emph{increasing} if $\mathcal{U}$ can be 
indexed as $\{ U_\alpha : \alpha < \kappa \}$ for some ordinal $\kappa$ 
with $\alpha < \beta < \kappa$ implying that $U_\alpha \subset U_\beta$.

An example $X$ yielding a positive answer would have to be a Dowker 
space.
If $\mathcal{V} = \{ V_\alpha : \alpha < \kappa \}$ were an increasing 
open cover of $X$ with 
$V_\alpha \setminus \bigcup_{\beta < \alpha} V_\beta$ nonempty, then 
$\kappa$ must have countable cofinality.
If $A$ is a subset of $X$ having regular uncountable cardinality, then 
$A$ has a limit point $x$ every neighbourhood of which meets $A$ in a set 
having the same cardinality.
\end{myprob}

\section*[Franklin D.~Tall: The cardinality of Lindel\"of spaces with 
points $G_\delta$]{Franklin D.~Tall: The cardinality of Lindel\"of spaces 
with points $\mathbf{G}_\mathbf{\delta}$} 

\begin{myprob}
\myproblem{Problem}{A.V.~Arhangel$'$\kern-.1667em ski\u{\i}}
What are the possible cardinalities of Lindel\"of $T_2$ spaces with points
$G_\delta$?

\mynote{Notes}
A.V.~Arhangel$'$\kern-.1667em ski\u{\i} raised the question of the cardinalities of 
Lindel\"of $T_2$ spaces with points $G_\delta$ and proved that there are 
none of cardinality greater than or equal to the first measurable 
cardinal \cite{MR87i:54001}.
S.~Shelah proved there are none of weakly compact cardinality.
I.~Juh\'asz \cite{MR82a:54002} constructed such (non-$T_2$) spaces of 
arbitrarily large cardinality with countable cofinality below the first 
measurable cardinal.
Shelah showed that it is consistent with \myaxiom{gch} that there is a 
zero-dimensional such space of size $\aleph_2$.
I.~Gorelic \cite{MR93g:03046} improved this result to get such a space of 
cardinality $2^{\aleph_1}$ consistent with \myaxiom{ch}, where 
$2^{\aleph_1}$ can be arbitrarily large.
Assuming the existence of a weakly compact cardinal, Shelah showed that it 
is consistent that $2^{\aleph_1} > \aleph_2$ and there is no such 
space of cardinality $\aleph_2$ (see \cite{MR86j:54008.1} for a good 
exposition of this result).
Shelah's results were eventually published in \cite{MR96k:03120}.

Among other results in \cite{MR96i:54016}, Tall proves:

\begin{theorem}
$\operatorname{Con}(\text{there is a supercompact cardinal\/}) 
\implies
\operatorname{Con}(2^{\aleph_1}$ is arbitrarily large and there is no 
Lindel\"of space with points $G_\delta$ of cardinality $\geq \aleph_2$ 
but $< 2^{\aleph_1})$.
\end{theorem}

\begin{theorem}
$\operatorname{Con}(\text{there is a supercompact cardinal\/})
\implies
\operatorname{Con}(\text{\myaxiom{gch}}$ $+$ there is no indestructible 
Lindel\"of space with points $G_\delta$ of cardinality $\geq \aleph_2)$.
\end{theorem}

A Lindel\"of space is \emph{indestructible} if it cannot be destroyed by 
countably closed forcing.

The problem of finding a small consistent bound for the $T_2$ case or for 
the first countable non-$T_2$ case remains open.
It is not known whether such spaces can be destructible.

C.~Morgan has withdrawn the claim added in proof to Tall's article 
\cite{MR96i:54016}.

In \cite{MR96i:54016}, Tall wrote:
``Little more is known: perhaps it is consistent (probably assuming large 
cardinals) that Lindel\"of spaces with points $G_\delta$ must have 
cardinality $\leq 2^{\aleph_0}$ or of countable cofinality.
It may also be consistent that if $T_2$ is added, the singular case can be 
dropped.
It may also be consistent---or even true---that Lindel\"of $T_2$ spaces 
with points $G_\delta$ all have cardinality $\leq 2^{\aleph_1}$.''
\end{myprob}

\section*{Stephen Watson: Basic problems in general topology}

\begin{myprob}
\myproblem{Problem 1}{\cite[Problem 163]{MR1078640.1}}
Do there exist, in \myaxiom{zfc}, more than $2^{\aleph_0}$ pairwise
$T_1$-complementary topologies on the continuum?

\mynote{Notes}
In 1936, G.~Birkhoff published ``On the combination of 
topologies'' in \emph{Fundamenta Mathematicae} \cite{birkhoff}.
In this paper, he ordered the family of all topologies on a set by letting 
$\tau_1 < \tau_2$ if and only $\tau_1 \subset \tau_2$.
He noted that the family of all topologies on a set is a lattice with a 
greatest element, the discrete topology and a smallest element, the 
indiscrete topology.
The family of all $T_1$ topologies on a set is also a lattice whose 
smallest element is the cofinite topology whose proper closed sets are 
just the finite sets.
Indeed, to study the lattice of all topologies on a set is to explore the 
fundamental interplay between general topology, set theory and finite 
combinatorics.
Recent work has revealed some essential and difficult problems in the 
study of this lattice, especially in the study of complementation, a 
phenomena in these lattices akin to in spirit to the study of Ramsey 
theory in combinatorial set theory.
We say that topologies $\tau$ and $\sigma$ are \emph{complementary} if 
and only of $\tau \wedge \sigma = 0$ and $\tau \vee \sigma = 1$.

B.A.~Anderson \cite{MR43:6860} showed by a beautiful construction that 
there is a family of $\kappa$ many mutually complementary topologies on 
$\kappa$.
J.~Stepr\={a}ns and S.~Watson \cite{MR95m:54003} showed that:
\begin{myitemize}
\item
There are $\kappa$ many mutually complementary partial orders 
(and thus $T_0$ topologies) on $\kappa$.
\item
Using the partial orders above, there are $\kappa$ many mutually 
$T_1$-comp\-le\-ment\-ary topologies on $\kappa$.
\item
There are $\kappa$ many mutually complementary equivalence relations on 
$\kappa$.
\item
The maximum size of a mutually $T_1$-complementary family of topologies on 
a set of cardinality $\kappa$ may not be greater than $\kappa$, unless 
$\omega < \kappa < 2^\mathfrak{c}$.
it is consistent that there do not exist $\aleph_2$ many mutually 
$T_1$-complementary topologies on $\omega_1$;
\item
Under \myaxiom{CH}, there are $2^{\aleph_1}$ mutually $T_1$-complementary 
topologies on $\omega_1$.
\end{myitemize}

D.~Dikranjan and A.~Policriti \cite{MR98g:06015} showed that there are 
families of two mutually complementary equivalence relations on a 
finite set (with more than three elements).

J.~Stepr\={a}ns and S.~Watson \cite{MR95m:54003} asked several problems:
\begin{myenumerate}
\item
Can one establish, in \myaxiom{zfc}, that there are $\mathfrak{c}^+$ 
many (maybe even $2^\mathfrak{c}$ many) mutually $T_1$-complementary 
topologies on $\mathfrak{c}$?
\item
Are there infinitely many mutually $T_1$-complementary (completely regular) 
Hausdorff spaces?
\item
What are the possible cardinalities of maximal families of mutually 
complementary families of partial order (or $T_0$ topologies)?
\item
What are the possible cardinalities of maximal families of mutually 
complementary families of $T_1$ topologies?
\item
What are the possible sizes of mutually $3$-complementary (mutually 
$2$-complementary) preorders (partial orders) (equivalence relations)?
\end{myenumerate}
\end{myprob}

\begin{myprob}
\myproblem{Problem 2}{\cite[Problem 168]{MR1078640.1}}
Is there a linear lower bound for the maximum number of pairwise
complementary partial orders on a finite set?

\mynote{Notes}
Specifically, does there exist $\varepsilon > 0$ such that, for any 
$n \in \mathbb{N}$, there are at least $\varepsilon \cdot n$ many 
pairwise complementary partial orders on a set of cardinality $n$?

Let $\omega(n)$ denote the maximum number of mutually complementary 
partial orders on a set of size $n$. 
J.~Brown and S.~Watson \cite{MR94h:06005} estimated the asymptotic 
behaviour as $n / \log n = O(\omega(n))$.
See also \cite{MR95m:06001, MR97d:54004}.
\end{myprob}

\begin{myprob}
\myproblem{Problem 3}{\cite[Problem 172]{MR1078640.1}}
Can every lattice with $1$ and $0$ be homomorphically embedded as a
sublattice in the lattice of topologies on some set?

\mynote{Notes}
Yes, answered by J.~Harding and A.~Pogel \cite{MR2001a:54004}.
\end{myprob}

\begin{myprob}
\myproblem{Problem 4}{}
Which lattices can be represented as the lattice of all topologies between
two topologies?
Can all finite lattices be represented in this fashion?

\mynote{Notes}
See the articles by D.~McIntyre et.\ al.\ 
\cite{MR1617095,MR99c:54006,MR2000d:06013,mcintyre} for progress on this 
problem.
\end{myprob}

\begin{myprob}
\myproblem{Problem 5}{\cite[Problem 107]{MR1078640.1}}
Are para-Lindel\"of regular spaces countably paracompact?

\mynote{Notes}
There is also Watson's Problem 108 \cite{MR1078640.1}: 
Is there a para-Lindel\"of Dowker space?
\end{myprob}

\begin{myprob}
\myproblem{Problem 6}{\cite[Problem 109]{MR1078640.1}}
Are para-Lindel\"of collectionwise normal spaces paracompact?

\mynote{Notes}
This was first asked by W.~Fleissner and G.M.~Reed \cite{MR80j:54020.1}
as Topology Proceedings Problem D26.

Z.~Balogh \cite{MR2003c:54047.2} constructed a hereditarily collectionwise 
normal, hereditarily meta-Lindel\"of, hereditarily realcompact Dowker 
space. 
This answers negatively R.~Hodel's question \cite{MR50:8400.1} (also 
Watson's Problem 110 and Topology Proceedings Problem D27):
are meta-Lindel\"of, collectionwise normal space paracompact?
Balogh listed some open questions about meta-Lindel\"of and 
para-Lindel\"of Dowker spaces at the end of his article 
\cite{MR2003c:54047.2}:
\begin{myenumerate}
\item
Is there a para-Lindel\"of, collectionwise normal Dowker space?
\item
Is there a para-Lindel\"of Dowker space?
\item
Is there a meta-Lindel\"of, collectionwise normal and first countable 
Dowker space?
\item
(D.~Burke)
Is there a meta-Lindel\"of, collectionwise normal and countably 
paracompact space which is not paracompact?
\item
Is there a first countable Dowker space in \myaxiom{zfc}?
\end{myenumerate}
\end{myprob}

\begin{myprob}
\myproblem{Problem 7}{\cite[Problem 88]{MR1078640.1}}
Does \myaxiom{zfc} imply that there is a perfectly normal locally 
compact space which is not paracompact?

\mynote{Solution}
P.~Larson and F.D.~Tall \cite{larsontall} proved that if it is consistent 
that there is a supercompact cardinal, then it is consistent that every 
locally compact, perfectly normal space is paracompact. 
\end{myprob}

\begin{myprob}
\myproblem{Problem 8}{\cite[Problem 85]{MR1078640.1}}
Are locally compact normal metacompact spaces paracompact?

\mynote{Solution}
This is known as the Arhangel$'$\kern-.1667em ski\u{\i}-Tall problem.
A.V.~Arhangel$'$\kern-.1667em ski\u{\i} \cite{MR46:4472.2} proved that perfectly 
normal, locally compact, metacompact spaces are paracompact.
F.D.~Tall asked the problem in \cite{MR51:13992}.
The answer is independent of \myaxiom{zfc}.
Yes, if \visl\ (S.~Watson \cite{MR83k:54021.2}); 
or by adding supercompact many Cohen or random reals 
(Z.~Balogh \cite{MR91c:54030.1}); 
or if $\text{\myaxiom{MA}}(\omega_1)$ for $\sigma$-centered posets
(G.~Gruenhage and P.~Koszmider \cite{MR97g:54004.2}).
G.~Gruenhage and P.~Koszmider \cite{MR97i:54010.2} showed that 
consistently the answer can be no.
\end{myprob}

\begin{myprob}
\myproblem{Problem 9}{\cite[Problem 175]{MR1078640.1}}
Is there, in \myaxiom{zfc}, a linear ordering in which every disjoint 
family of open intervals is the union of countably many discrete 
subfamilies and yet in which there is no dense set which is the union of 
countably many closed discrete sets? 
Is there such a linear ordering if and only if there is a Souslin line?

\mynote{Notes}
A compact Souslin line is such a linear ordering but there may be others.
The Urysohn metrization theorem is to the Nagata-Smirnov-Stone 
metrization theorem as the Souslin problem is to this problem.

Y.-Q.~Qiao and F.D.~Tall showed that the existence of such a linear 
ordering is equivalent to the existence of a perfectly normal 
nonmetrizable non-Archimedean space (i.e., an archvillain).
Y.-Q.~Qiao \cite{MR2001g:54007} showed that there is a model of 
\manotch\ in which there is such a space (and yet no Souslin lines).
\end{myprob}

\begin{myprob}
\myproblem{Problem 10}{\cite[Problem 176]{MR1078640.1}}
Is there a topological space (or a completely regular space) in which the 
connected sets (with more than one point) are precisely the cofinite sets?

\mynote{Notes}
This problem was motivated by an interesting paper by S.F.~Cvid 
\cite{MR80b:54045}.
Cvid asked whether the connected sets in a countable connected Hausdorff 
space could form a filter. 
That problem remains unsolved.
P.~Erd\H{o}s \cite{MR6:43a} attributes to A.H.~Stone the result that 
there are no such metrizable spaces.
In fact, if a space is such that its connected sets are precisely its 
cofinite sets then the space must be $T_1$ and every infinite subset of 
the space must contain an infinite closed discrete set (in particular, the 
space cannot contain convergent sequences).

G.~Gruenhage \cite{MR95j:54024} constructed, consistently, several 
examples of spaces whose connected sets are their cofinite sets.
Assuming \myaxiom{MA}, there are completely regular as well as 
countable examples.
Assuming \myaxiom{CH}, there is a perfectly normal example.
Watson conjectured that an example (probably even completely 
regular) exists in \myaxiom{zfc} and that this will depend on some hard 
finite combinatorics.

Furthermore, Gruenhage \cite[Questions 4.8, 4.9, 1.10]{MR95j:54024} asked:
\begin{myitemize}
\item
Is there a completely regular space $X$ in which the nondegenerate 
connected sets are precisely the (co-$< |X|$))-sets?
Or co-$< \lambda$ for some uncountable cardinal $\lambda$?
\item
Is there a paracompact Hausdorff (or regular Lindel\"of) space in which 
the nondegenerate connected sets are precisely the cofinite sets?
\item
Is there in \myaxiom{zfc} a Hausdorff (or completely regular) space in 
which the nondegenerate connected sets are precisely the cofinite sets?
\end{myitemize}
\end{myprob}

\providecommand{\bysame}{\leavevmode\hbox to3em{\hrulefill}\thinspace}

\label{tpnewclassicend}

\chapter*{Problems from M.E.~Rudin's \emph{Lecture notes in set-theoretic 
topology}}
\markboth{\normalsize\textsc{\lowercase{Problems from Rudin's Lecture notes}}}{}
\label{tprudin}
\begin{myfoot}
\begin{myfooter}
Elliott Pearl, 
\emph{Problems from M.E.~Rudin's \emph{Lecture notes in set-theoretic 
topology}},\\
Problems from Topology Proceedings, Topology Atlas, 2003, 
pp.\ 103--121.
\end{myfooter}
\end{myfoot}

\mypreface
Here are the problems from the last chapter of Mary Ellen Rudin's 
\emph{Lecture notes in set-theoretic topology} \cite{MR51:4128.1}. 
The list was first published in 1975 and it was updated for the second 
printing in 1977.
This version uses the item numbering from the 1977 list and includes the 
(solved) problems from the 1975 list that were dropped from the 1977 list.
This material is reprinted here with the permission of the American 
Mathematical Society.
Some corrections to the second printing were provided in volume 2 (1977) 
of \emph{Topology Proceedings}. 
This version includes these corrections and other information on solutions 
that appeared in subsequent volumes of \emph{Topology Proceedings}.

\section*{Introduction}

Rudin wrote: ``The following problems are unsolved so far as I know.
They are being solved almost daily, of course, for they are problems which 
people are working on.
Some are very hard, basic, long unsolved and frequently worked on problems; 
others are just things someone ran across and did not know the answer to.
The names by the problem are not those of the first person to ask the 
problem or even the person currently most actively working on the problem: 
the name implies that that person once mentioned this problem to me and 
probably can fill in anyone interested in the problem on more details and 
background.''

All spaces are assumed to be Hausdorff.
A \emph{map} is a function which is continuous and onto.

\section*{A. Cardinal function problems}

\begin{myprob}
\myproblem{A1}{S.~Mr\'owka}
If every zero set is in $\mathcal{B}(\text{clopen})$, then is every zero 
set the intersection of a countable number of clopen sets?

$\mathcal{B}(\text{clopen})$ is the $\sigma$-algebra of clopen sets. 

\mynote{Notes}
This problem is due to M.~Kat\v{e}tov.
See Z.~Balogh's contribution to \emph{New Classic Problems}.
\end{myprob}

\begin{myprob}
\myproblem{A2}{I.~Juh\'asz and A.~Hajnal}
Is there a regular space $X$ with cardinality greater than $\mathfrak{c}$ 
which is not hereditarily separable but every closed subset is separable?

\end{myprob}

\begin{myprob}
\myproblem{A3}{I.~Juh\'asz and A.~Hajnal}
If $X$ is an infinite space and the number of open sets in $X$ is 
denoted by $o(X)$, then is $o(X)^\omega = o(X)$?
Yes if \myaxiom{GCH}.

\mynote{Solution}
This problem is due to I.~Juh\'asz from the 1976 Prague conference.
The answer is yes for several classes of spaces (e.g., hereditarily 
paracompact spaces \cite{MR83d:54006}; compact Hausdorff spaces).
No is consistent (S.~Shelah \cite{MR850057}).
\end{myprob}

\begin{myprob}
\myproblem{A4}{I.~Juh\'asz and A.~Hajnal}
Does every Lindel\"of space of cardinality $\aleph_2$ contain a Lindel\"of 
subspace of cardinality $\aleph_1$?
Yes if \myaxiom{GCH}.

\mynote{Solution}
No is consistent (P.~Koszmider and F.D~Tall \cite{MR2003a:54002}).
\end{myprob}

\begin{myprob}
\myproblem{A5 first printing}{I.~Juh\'asz and A.~Hajnal}
If $X$ is hereditarily separable and compact (subset of $2^{\aleph_1}$), 
then is $|X| \leq \mathfrak{c}$?

\mynote{Solution}
No if \myaxiom{CH}.
Yes if \manotch\ (Z.~Szentmikl\'ossy).
\end{myprob}

\begin{myprob}
\myproblem{A5}{I.~Juh\'asz}
If \myaxiom{GCH} holds and $X$ is a compact space, which cardinals less 
than $|X|$ can be omitted as the cardinality of closed subsets of $X$?

\mynote{Notes}
See Juh\'asz's \emph{Handbook} article \cite{MR86j:54008.2} and the series 
of articles \cite{MR94d:54012,MR94j:54002,MR99c:06015} for results on the 
cardinality and weight spectra of compact spaces.
\end{myprob}

\begin{myprob}
\myproblem{A6 first printing}{I.~Juh\'asz and A.~Hajnal}
If $X$ is a regular space of countable spread, does $X = Y \cup Z$ where 
$Y$ is hereditarily separable and $Z$ is hereditarily Lindel\"of?

\mynote{Solution}
No if \myaxiom{CH} or if there is a Souslin line 
(J. Roitman \cite{MR58:12909}).
Yes if \myaxiom{PFA}.
\end{myprob}

\begin{myprob}
\myproblem{A6}{R. Hodel}
Does every regular, hereditarily c.c.c., $w\Delta$ space with a 
$G_\delta$-diagonal have a countable base?
See D9 below.
\end{myprob}

\begin{myprob}
\myproblem{A7}{A.V.~Arhangel$'$\kern-.1667em ski\u{\i}}
If $X$ is a regular Lindel\"of space each point of which is a $G_\delta$, 
then is $|X| \leq \mathfrak{c}$?

\mynote{Notes}
This is the Lindel\"of points $G_\delta$ problem.
See F.D.~Tall's contribution to \emph{New Classic Problems}.
\end{myprob}

\begin{myprob}
\myproblem{A8}{A.V.~Arhangel$'$\kern-.1667em ski\u{\i}}
If a hereditarily normal space $X$ has countable cellularity and 
countable tightness, is $|X| \leq \mathfrak{c}$?
No if \visl\ without hereditary normality.
\end{myprob}

\begin{myprob}
\myproblem{A9}{A.V.~Arhangel$'$\kern-.1667em ski\u{\i}}
Does each compact hereditarily normal space of countable tightness contain 
a nontrivial convergent sequence?
a point of countable character?
No if \visl\ without hereditary normality.

\mynote{Solution}
Yes is consistent.
See Classic Problem V.
\end{myprob}

\begin{myprob}
\myproblem{A10}{A.V.~Arhangel$'$\kern-.1667em ski\u{\i}}
Does every compact homogeneous space of countable tightness have 
cardinality $\leq \mathfrak{c}$?
\end{myprob}

\begin{myprob}
\myproblem{A11}{Yu.M.~Smirnov}
Does every hereditarily normal compact space contain a point with a 
countable $\Delta$-base?

A $\Delta$-\emph{base} for a point $x$ is a family $B$ of open sets such 
that every neighborhood of $x$ contains a member of $B$ having $x$ in its 
closure.
\end{myprob}

\begin{myprob}
\myproblem{A12}{V.I.~Ponomarev}
Is a compact space of countable tightness a sequential space?
No if $\diamondsuit$.

\mynote{Solution}
This is the Moore-Mr\'owka Problem.
Yes if \myaxiom{PFA} (Z.~Balogh \cite{MR89h:03088.2}).
See Classic Problem VI.
\end{myprob}

\begin{myprob}
\myproblem{A13}{}
Is the product of two Lindel\"of spaces $\mathfrak{c}$-Lindel\"of?

\mynote{Solution}
No, there are consistent counterexamples.
S.~Shelah \cite{MR96k:03120.1} gave the first example.
D.~Velleman \cite{MR85e:03121} produced examples in \visl.
I.~Gorelic \cite{MR95e:54009} gave a forcing construction of Lindel\"of 
space whose square has a closed discrete subspace of size $2^{\aleph_1}$ 
(where this cardinal can be arbitrarily large regardless of 
$\mathfrak{c}$).
\end{myprob}

\begin{myprob}
\myproblem{A14}{}
Is every separable metric space, such that every nowhere dense closed 
subset is $\sigma$-compact, $\sigma$-compact?
No if \myaxiom{CH}.
\end{myprob}

\begin{myprob}
\myproblem{A15}{S.~Purisch}
Is orderable equivalent to monotonically normal for compact, separable, 
totally disconnected spaces?
\end{myprob}

\begin{myprob}
\myproblem{A16 first printing}{R.~Telg\'arsky}
Is there a compact space $X$ with no isolated points which does not 
contain a zero-dimensional closed subset with no isolated points?
No if \visl.
\end{myprob}

\begin{myprob}
\myproblem{A16}{E.~van~Douwen}
Is every point-finite open family in a c.c.c.\ space $\sigma$-centered 
(i.e., the union of countably many centered subfamilies)?

\mynote{Solution}
No (Ortwin F\"orster).
J.~Stepr\={a}ns and S.~Watson \cite{MR89e:54061.1} described a subspace 
of the Pixley-Roy space on the irrationals that is a first countable 
c.c.c.\ space which does not have a $\sigma$-linked base.
\end{myprob}

\begin{myprob}
\myproblem{A17}{R.~Telg\'arsky}
Is every image of a scattered space under a closed map scattered?
No if \myaxiom{MA}.

\mynote{Solution}
No (V.~Kannan and M.~Rajagopalan \cite{MR57:13830}).
\end{myprob}

\begin{myprob}
\myproblem{A18 first printing}{Z.~Semadeni \cite{MR21:6571}}
Is every scattered completely regular space zero-dim\-en\-sional?

\mynote{Solution}
No (R.C.~Solomon \cite{MR56:13183}).
\end{myprob}

\begin{myprob}
\myproblem{A18}{E.~van~Douwen}
Does the Sorgenfrey line have a connected compactification?

\mynote{Solution}
No (A.~Emeryk and W.~Kulpa \cite{MR57:1422}).
\end{myprob}

\begin{myprob}
\myproblem{A19}{A.~Hajnal}
Suppose $A$ and $B$ are sets with $|A| = 2^{\aleph_1}$ and 
$|B| = 2^{\aleph_0}$.
Color $A \times B$ with two colors.
Must there be $A' \subset A$ and $B' \subset B$ such that 
$|A'| = \aleph_0$, $|B'| = \aleph_1$, and $A' \times B'$ is one color?
Yes is consistent.
\end{myprob}

\begin{myprob}
\myproblem{A20 first printing}{R.M.~Stephenson}
Is the property initially $\mathfrak{m}$-compact productive for regular 
uncountable $\mathfrak{m}$?

$X$ is \emph{initially $\mathfrak{m}$-compact} if every open cover of 
cardinality $\mathfrak{m}$ on $X$ has a finite subcover.

\mynote{Solution}
No is consistent (E.~van~Douwen \cite{MR93f:54012.1}).
See Topology Proceedings Problems C21, C22, C23 for related problems.
\end{myprob}

\begin{myprob}
\myproblem{A20}{W.~Weiss}
If $X$ is a compact scattered space such that $X^\alpha - X^{\alpha+1}$ is 
countable for all $\alpha$, what are the bounds on the order (minimal 
$\alpha$ with $X^\alpha$ finite) of $X$?

\mynote{Solution}
$\alpha = \omega_2$ is consistent (J.~Baumgartner and S.~Shelah).
If \myaxiom{CH}, then $\alpha < \omega_2$ 
(I.~Juh\'asz and W.~Weiss \cite{MR82k:54005}). 
There are \myaxiom{ZFC} examples with $\alpha < \omega_2$.
\end{myprob}

\begin{myprob}
\myproblem{A21 first printing}{K.~Morita \cite{MR51:6734}}
Is every normal space $X$ countably compactifiable?
That is, is $X$ dense in a countably compact space $S$ such that 
every countably compact subset of $X$ is closed in $S$?

\mynote{Solution}
No.
D.~Burke and E.~van~Douwen \cite{MR55:13381} constructed a normal, 
locally compact $M$-space which does not have a countable
compactification.
A.~Kato \cite{MR55:11211} showed that $\beta\mathbb{R} - \beta\mathbb{N}$ 
is an $M$-space which is not countably compactifiable.
\end{myprob}

\begin{myprob}
\myproblem{A21}{E.~van~Douwen}
Is every paracompact (or metacompact or subparacompact or hereditarily 
Lindel\"of) space a $D$-space?
\end{myprob}

\section*{B. Souslin and compactness problems}

\begin{myprob}
\myproblem{B1 first printing}{K.~Kunen}
If $X$ is c.c.c.\ and $Y$ is c.c.c.\ but $X \times Y$ is not c.c.c.\ then 
is there a Souslin line?

\mynote{Solution}
No.
R.~Laver and F.~Galvin showed that such an $X$ and $Y$ can exist under 
\myaxiom{CH}.
\end{myprob}

\begin{myprob}
\myproblem{B1}{E. van~Douwen}
Does a compact homogeneous space have a nontrivial convergent sequence?
\end{myprob}

\begin{myprob}
\myproblem{B2 first printing}{}
Is there a Souslin line if there is a normal, not countably paracompact 
space (a Dowker space) which is one (or many) of the 
following:
first countable,
separable,
cardinality $\aleph_1$,
c.c.c.,
realcompact,
monotonically normal?
Yes if \myaxiom{CH}.

\mynote{Solution}
No, $\clubsuit$ also works.
\end{myprob}

\begin{myprob}
\myproblem{B2}{}
Is there a Dowker space which has any of the following properties?
(For all except the first and last the answer is yes if \myaxiom{CH} or 
if there is a Souslin line.)
Extremally disconnected,
first countable,
separable,
cardinality $\aleph_1$,
c.c.c.,
realcompact,
monotonically normal.

\mynote{Solution}
See Classic Problem VII for a discussion of small Dowker spaces.
\end{myprob}

\begin{myprob}
\myproblem{B3}{Yu.M.~Smirnov}
Does every compact space contain either a copy of $\mathbb{N}^*$ or a 
point of countable $\pi$-character?
\end{myprob}

\begin{myprob}
\myproblem{B4}{Yu.M.~Smirnov}
Is there a c.c.c., compact space $X$ with countable $\pi$-character (or 
with $|X| \leq \mathfrak{c}$) which is not separable?
\end{myprob}

\begin{myprob}
\myproblem{B5}{E.~van~Douwen \cite{MR82m:54017}}
\begin{myenumerate}
\item
Is every compact space supercompact?
\item
Is the continuous image of a supercompact space supercompact?
\item
Is a dyadic space supercompact?
\end{myenumerate}

A space is \emph{supercompact} if it has a subbasis $S$ for the closed 
sets such that, if $T \subseteq S$ and every two members of $T$ meet, then 
$\bigcap T$ is nonempty.

\mynote{Solution}
M.~Bell \cite{MR57:13846} showed that not all compact spaces are 
supercompact.

C.~Mills and J.~van~Mill \cite{MR80m:54033} showed that the continuous 
image of a supercompact space need not be supercompact: Let $X$ be the 
subspace of $(\omega_1 + 1)^2$ comprising the diagonal and everything 
below it.
$X$ is supercompact, but the quotient space obtained by collapsing, 
for each $\alpha < \omega_1$, $\{(\alpha, \alpha). (\omega_1, \alpha)\}$ 
to a point, is not supercompact.

M.~Bell \cite{MR92c:54029} showed that not all dyadic spaces are 
supercompact. 
\end{myprob}

\begin{myprob}
\myproblem{B6 first printing}{W. Fleissner}
Is there a Baire space whose square is not Baire?
Yes under \myaxiom{MA} even for metric spaces.

\mynote{Solution}
Yes, W.~Fleissner and K.~Kunen \cite{MR80f:54009} showed that there 
are even metric examples of so-called barely Baire spaces.
\end{myprob}

\begin{myprob}
\myproblem{B6}{W.~Fleissner}
Is any product of metric Baire spaces Baire?

\mynote{Solution}
No. 
See B5 above.
\end{myprob}

\begin{myprob}
\myproblem{B7 first printing}{M.~Henriksen}
Is the set of remote points in $\beta\mathbb{R}$ dense in $\mathbb{R}^*$?

\mynote{Solution}
Yes (E.~van~Douwen \cite{MR83i:54024.1}).
\end{myprob}

\begin{myprob}
\myproblem{B7}{W.~Fleissner}
Is any box product of second countable (metric) Baire spaces Baire?

\mynote{Notes}
See \cite{MR57:4116, MR80f:54009}
\end{myprob}

\begin{myprob}
\myproblem{B8}{Z.~Frol\'ik}
Is there a $P$-point in $\mathbb{N}^*$?
Yes if \myaxiom{MA}.

\mynote{Solution}
The answer is independent of \myaxiom{ZFC}.
No is consistent (S.~Shelah).
\end{myprob}

\begin{myprob}
\myproblem{B9}{A.V.~Arhangel$'$\kern-.1667em ski\u{\i}}
Does every hereditarily separable compact space have
a point of countable character?
a nontrivial convergent sequence?
a butterfly point?
No if $\diamondsuit$ (V.V.~Fedor\v{c}uk).

\mynote{Solution}
Yes if \manotch\ (Z.~Szentmikl\'ossy).
\end{myprob}

\begin{myprob}
\myproblem{B10}{V.~Saks \cite{MR51:8321, MR80h:54003}}
Is there a product of countably compact topological groups which is not 
countably compact.

\mynote{Solution}
Yes if \myaxiom{MA} (E.~van~Douwen \cite{MR82b:22002}).
See also \cite{MR91h:22001, MR91e:54025, MR98a:54033}.
\end{myprob}

\begin{myprob}
\myproblem{B11}{A.~Hager}
If $X$ is a dense subset of a compact $Y$ and every open set containing 
$X$ is $C^*$-embedded in $Y$, then is $X$ $C^*$-embedded in $Y$?

\mynote{Solution}
No (M.~Sola).
If one lets $X = \Delta$, P.~Roy's example, and its zero-dimensional 
compactification $Y = \zeta \Delta$ then every open subspace $U$ of $\zeta 
\Delta$ containing $\Delta$ is strongly zero-dimensional hence $\beta U = 
\zeta \Delta$ and $U$ is $C^*$-embedded in $\zeta \Delta$, but $\Delta$ 
is not strongly zero-dimensional and so it is not $C^*$-embedded in 
$\zeta \Delta$.
See the review by P.~Nyikos \cite{nyikosmr}.
\end{myprob}

\begin{myprob}
\myproblem{B12 first printing}{A.V.~Arhangel$'$\kern-.1667em ski\u{\i}}
Is there an infinite homogeneous extremally disconnected space?
Yes if \myaxiom{CH}.

\mynote{Solution}
No, there is not even an infinite homogeneous compact $F$-space.
See Kunen's article on van~Douwen's problem \cite{MR1078652.1}.
\end{myprob}

\begin{myprob}
\myproblem{B12}{K.~Kunen}
Is there an extremally disconnected locally compact nonparacompact space?
Yes if there is a weakly compact cardinal.
\end{myprob}

\begin{myprob}
\myproblem{B13}{R.~Blair}
If $X$ is Lindel\"of and $Y$ is realcompact, does $X$ closed in $X \cup Y$ 
imply that $X \cup Y$ is realcompact?

\mynote{Solution}
This question is due to S.~Mr\'owka, who proved that if $Y$ is also 
closed in $X \cup Y$ then $X \cup Y$ is realcompact.
A.~Kato \cite{MR81b:54025} gave a negative solution with a decomposition 
inside the Tychonoff plank: $X = \omega \times \{\omega_1\}$, 
$Y = (\omega + 1) \times D$, where $D$ is the discrete subspace of 
isolated points of $\omega_1$. 
\end{myprob}

\begin{myprob}
\myproblem{B14 first printing}{K.~Kunen}
Can a compact space be decomposed into more than $\mathfrak{c}$ closed 
$G_\delta$ sets?

\mynote{Solution}
No (A.V.~Arhangel$'$\kern-.1667em ski\u{\i}).
R.~Pol's proof of Arhangel$'$\kern-.1667em ski\u{\i}'s theorem (every first countable 
compact Hausdorff space has cardinality at most $\mathfrak{c}$) can be 
adapted here, replacing points by $G_\delta$ sets.
\end{myprob}

\begin{myprob}
\myproblem{B14}{R.~Frankiewicz}
Is $\omega_1^*$ ever homeomorphic to $\omega^*$?

\mynote{Notes}
This old problem is discussed in \cite[Problem 242]{MR1078643.1}.
\end{myprob}

\begin{myprob}
\myproblem{B15}{C.~Bandy \cite{MR49:11484}}
Are there two normal countably compact spaces whose product is not 
countably compact? 

\mynote{Solution}
Yes if \myaxiom{MA} (E.~van~Douwen \cite{MR93f:54012.1}).
\end{myprob}

\section*{C. Separable-Lindel\"of problems}

\begin{myprob}
\myproblem{C1 first printing}{I.~Juh\'asz and A.~Hajnal}
Is there a first countable c.c.c.\ space with density at most 
$\mathfrak{c}$ and uncountable spread?
Yes if \myaxiom{CH} or if there is a Souslin line.

\mynote{Solution}
Yes.
The Sorgenfrey line is a separable example.
The square the Alexandroff's double arrow space is a compact example.
This problem was certainly misstated.
\end{myprob}

\begin{myprob}
\myproblem{C1}{I.~Juh\'asz}
If \myaxiom{GCH} holds, $X$ is Lindel\"of, and $|X| = \aleph_2$, does 
there exist a Lindel\"of $Y \subset X$ with $|Y| = \aleph_1$?

\mynote{Notes} 
This is Problem A4.
\end{myprob}

\begin{myprob}
\myproblem{C2}{I.~Juh\'asz and A.~Hajnal}
Is there a regular space with cardinality greater than $\mathfrak{c}$ 
which has countable spread?
Yes is consistent with $\diamondsuit$.
\end{myprob}

\begin{myprob}
\myproblem{C3}{I.~Juh\'asz and A.~Hajnal}
Is there a regular, hereditarily Lindel\"of space with weight greater than 
$\mathfrak{c}$?
Yes is consistent with \myaxiom{CH}.
\end{myprob}

\begin{myprob}
\myproblem{C4 first printing}{I.~Juh\'asz and A.~Hajnal}
Is there a (regular) hereditarily separable space $X$ with 
$|X| > 2^{\aleph_1}$?

\mynote{Solution}
This problem was asked originally by J.~Gerlits.
No if \myaxiom{CH} for regular $X$ because $w(X) \leq \mathfrak{c}$ and so
$|X| \leq 2^\mathfrak{c}$.
S.~Todor\v{c}evi\'c \cite{MR85d:03102.1} proved that it is consistent 
that every Hausdorff space with no uncountable discrete subspace has 
cardinality $\mathfrak{c}$.
I.~Juh\'asz and S.~Shelah \cite{MR87f:03143} showed that it is consistent 
that there are regular hereditarily separable spaces of size 
$2^\mathfrak{c}$, where $\mathfrak{c}$ is arbitrarily large and 
$2^\mathfrak{c}$ is arbitrarily larger.
\end{myprob}

\begin{myprob}
\myproblem{C4}{E.~van~Douwen}
Can every first countable compact space be embedded in a separable first 
countable compact space?
Yes if \myaxiom{CH}.

\mynote{Solution}
See the papers by E.~van~Douwen and T.~Przymusi\'nski 
\cite{MR80e:54031,MR82j:54051.1}
for relevant results.
\end{myprob}

\begin{myprob}
\myproblem{C5}{I.~Juh\'asz and A.~Hajnal}
Is there a regular space which is hereditarily separable but not 
Lindel\"of (i.e., an $S$-space), or vice versa (i.e., an $L$-space).
Yes in both cases if \myaxiom{CH} or if there is a Souslin line.

\mynote{Notes}
It is consistent that there are no $S$-spaces (S.~Todor\v{c}evi\'c).
\end{myprob}

\begin{myprob}
\myproblem{C7 first printing}{}
Is density not greater than the smallest cardinal greater than spread for 
compact spaces?
regular spaces?
regular hereditarily Lindel\"of spaces?

\mynote{Solution}
B.~Shapirovski\u{\i} \cite{MR52:4213} showed that $hd(X) \leq s(X)^+$ for 
compact spaces.
\end{myprob}

\begin{myprob}
\myproblem{C8}{I.~Juh\'asz and A.~Hajnal}
Could a compact hereditarily separable space have cardinality greater than 
$\mathfrak{c}$?
Yes if $\diamondsuit$ (V.V.~Fedor\v{c}uk).

\mynote{Solution}
No if \manotch\ (Z.~Szentmikl\'ossy).
\end{myprob}

\section*{D. Metrizability problems}

\begin{myprob}
\myproblem{D1}{F.B.~Jones}
Is there a normal nonmetrizable Moore space?
Yes if \manotch.

\mynote{Solution}
The normal Moore space conjecture is the assertion that normal Moore
spaces are metrizable.
P.~Nyikos \cite{MR81k:54044} showed, under the assumption of the product
measure extension axiom (\myaxiom{PMEA}), that any normal first countable
space (hence any normal Moore space) is metrizable.
K.~Kunen showed that \myaxiom{PMEA} was consistent relative to the
consistency of the existence of a strongly compact cardinal.
Assuming \myaxiom{CH}, W.~Fleissner \cite{MR84f:54040} constructed a
normal nonmetrizable Moore space. 
It follows from Fleissner's construction that if all normal Moore 
spaces are metrizable then there is a inner model with a measurable 
cardinal. 
So, large cardinals are necessary to prove the consistency of the normal 
Moore space conjecture.
A.~Dow, F.D.~Tall and W.~Weiss \cite{MR92b:54008a,MR92b:54008b} gave a
new proof, using iterated forcing and reflection, of the normal Moore
space conjecture under the assumption of the existence of a supercompact
cardinal.
For more information, see the surveys by Fleissner \cite{MR86m:54023.1}
and Nyikos \cite{MR2003c:54001.1}.
\end{myprob}

\begin{myprob}
\myproblem{D2}{P.S.~Alexandroff}
Is there a normal nonmetrizable image of a metric space under a compact 
open map?
Yes if \manotch.

\mynote{Solution}
This is equivalent to the metacompact normal Moore problem.
See Classic Problem II and its related problems.
\end{myprob}

\begin{myprob}
\myproblem{D3}{W.~Fleissner}
Is there a first countable, normal collectionwise Hausdorff space which is 
not collectionwise normal?
Yes if \manotch.

\mynote{Solution}
No if \myaxiom{PMEA} (P.~Nyikos).
\end{myprob}

\begin{myprob}
\myproblem{D4}{G.M.~Reed}
Is every countably paracompact Moore space normal?
No if \manotch.

\mynote{Solution}
Yes if \myaxiom{PMEA} (D.~Burke \cite{MR85h:54032.1}) or 
\myaxiom{PCEA} (W.~Fleissner \cite{MR88f:54044}).
A positive solution requires large cardinals.
\end{myprob}

\begin{myprob}
\myproblem{D5}{G.M.~Reed}
Is there a countably paracompact Moore space which is not paracompact?
Yes if \manotch.

\mynote{Solution}
See D4.
\end{myprob}

\begin{myprob}
\myproblem{D6}{B.~Wilder, P.S.~Alexandroff}
Is every perfectly normal manifold metrizable?

\mynote{Solution}
No if \myaxiom{CH} (M.E.~Rudin and P. Zenor \cite{MR52:15361}).
Yes if \manotch\ (M.E.~Rudin \cite{MR80j:54014.1}).
\end{myprob}

\begin{myprob}
\myproblem{D7 first printing}{P.~Zenor}
Is every perfectly normal manifold subparacompact?
This is equivalent to D6.
\end{myprob}

\begin{myprob}
\myproblem{D7}{E.~van~Douwen}
If $X$ is $\sigma$-compact and locally compact and $f$ is one of 
cardinality, cellularity, density, spread, $\pi$-weight or weight, is
$f(X^*)^\omega = f(X^*)$?
Yes if \myaxiom{GCH}.
\end{myprob}

\begin{myprob}
\myproblem{D8 first printing}{R.~Hodel}
Is every perfectly normal collectionwise normal space paracompact?
No if $\diamondsuit$ or if \manotch.

\mynote{Solution}
No (R.~Pol \cite{MR57:4113}).
\end{myprob}

\begin{myprob}
\myproblem{D8}{E. van~Douwen}
Do spaces like $\beta\mathbb{N}$, $\mathbb{N}^*$, $\beta\mathbb{R}$,
$\mathbb{R}^*$, \ldots admit a mean?
\end{myprob}

\begin{myprob}
\myproblem{D9 first printing}{}
Is every perfect space $\theta$-refinable?
No if $\diamondsuit$ or if \manotch.

\mynote{Solution} 
No (R.~Pol \cite{MR57:4113}).
\end{myprob}

\begin{myprob}
\myproblem{D9}{R.~Hodel}
Does every regular, $\aleph_1$-compact, $w\Delta$ space with a 
$G_\delta$-diagonal (or point-countable separating open cover) have a 
countable basis?
\end{myprob}

\begin{myprob}
\myproblem{D10 first printing}{R.~Hodel}
Is every normal space with a point-countable base metrizable?
No if \manotch.

\mynote{Solution}
Under \myaxiom{CH}, E.~van~Douwen, F.D.~Tall, and W.~Weiss 
\cite{MR58:24187.1} constructed a nonmetrizable hereditarily Lindel\"of 
space with a point-countable base.
\end{myprob}

\begin{myprob}
\myproblem{D10}{K.~Kunen}
Does the existence of $P$-points in $\mathbb{N}^*$ imply the existence of 
points which are the intersection of $\mathfrak{c}$ well ordered by 
inclusion open sets?
If \manotch\ are all points which are the intersection of 
$\mathfrak{c}$ well ordered by inclusion open sets of the same type in 
$\mathbb{N}^*$?
\end{myprob}

\begin{myprob}
\myproblem{D11 first printing}{R.~Hodel}
\begin{myenumerate}
\item 
Is every perfectly normal space with a point-countable basis metrizable?
\item 
Is every perfectly normal paracompact space with a point-countable basis 
metrizable?
\item
Is every perfectly normal collectionwise normal space with a point 
countable basis metrizable?
\end{myenumerate}

(Ponomarev)
Is every regular Lindel\"of space with a point-countable basis metrizable?

Yes to all if there is a Souslin line.

\mynote{Solution}
See D10 first printing.
\end{myprob}

\begin{myprob}
\myproblem{D11}{P.~Nyikos}
\begin{myenumerate}
\item
Is every perfectly normal space with a point-countable basis metrizable?
No if \manotch\ or if there is a Souslin line.
\item
Is every perfectly normal collectionwise normal space with a point 
countable base metrizable?
No if there is a Souslin line.
\end{myenumerate}

\mynote{Solution}
No to (1) (S.~Todor\v{c}evi\'c \cite{MR93k:54005.1}).
See Classic Problem II.
\end{myprob}

\begin{myprob}
\myproblem{D12}{P.~Nyikos}
Is there a perfectly normal non-Archimedean space which is not metrizable?
Yes if there is a Souslin line.

\mynote{Notes}
Such spaces are called \emph{archvillains}.
See Watson's contribution to \emph{New Classic Problems} or
\cite[Problem 175]{MR1078640.3}.
\end{myprob}

\begin{myprob}
\myproblem{D13 first printing}{P.~Zenor}
Is every countably compact space with a $G_\delta$ diagonal metrizable?

\mynote{Solution}
Yes (J.~Chaber \cite{MR58:24189}).
\end{myprob}

\begin{myprob}
\myproblem{D13}{E.~van~Douwen}
Is $d(\beta X) = d(X)$ if $X$ is a paracompact $p$-space?
No without paracompact.
\end{myprob}

\begin{myprob}
\myproblem{D14}{A.V.~Arhangel$'$\kern-.1667em ski\u{\i}}
If $X$ is completely regular and metacompact, is $X$ the image of a 
paracompact space under a compact open map?
\end{myprob}

\begin{myprob}
\myproblem{D15 first printing}{R.~Heath}
Is every linearly ordered space with a point-countable base 
quasi-metrizable?
No if there is a Souslin line (J.~Roitman).

\mynote{Solution}
No (G.~Gruenhage \cite{MR55:9040.1}).
\end{myprob}

\begin{myprob}
\myproblem{D15}{E.~van~Douwen}
Is there a discrete subset of $\beta\mathbb{N}$ of cardinality $\aleph_1$ 
which is not $C^*$-embedded?
\end{myprob}

\begin{myprob}
\myproblem{D16}{D.~Lutzer}
Is a weakly $\theta$-refinable, collectionwise normal space paracompact?
No if $\clubsuit$.

\mynote{Solution}
No. 
See Classic Problem III.
\end{myprob}

\begin{myprob}
\myproblem{D17}{J.C.~Smith}
Are compact (or paracompact $\Sigma$) spaces with a $\delta\theta$ base 
metrizable?

$B$ is a $\delta\theta$ \emph{base} for $X$ if 
$B = \bigcup_{n \in \omega} B_n$ and $x \in X$ and $U$ is a neighborhood 
of $x$ imply that there is an $n_x$ such that 
$\{ V \in B_{n_x} : x \in V \}$ is a finite nonempty subset of $U$.

\mynote{Solution}
These questions were asked by C.E.~Aull \cite{MR54:3660}.
J.~Chaber \cite{MR56:9493.1} gave a positive answer.
Every $\Sigma$-space is a $\beta$-space, and every $\theta$-refinable 
space with a base of countable order is a Moore space.
Chaber proved that every monotonic $\beta$-space with a $\delta\theta$ 
base has a base of countable order.
One can even replace \emph{paracompact $\Sigma$} with 
\emph{collectionwise normal $\Sigma$} because they are equivalent for 
spaces with a $\delta\theta$ base.
\end{myprob}

\begin{myprob}
\myproblem{D18}{P.~Nyikos)}
In screenable spaces do normal and collectionwise normal imply countably 
paracompact?

\mynote{Solution}
No.
See Classic Problem III.
\end{myprob}

\begin{myprob}
\myproblem{D19}{I.~Juh\'asz}
Suppose that $X$ is a hereditarily Lindel\"of space of weight 
$> \mathfrak{c}$.
Is the number of closed sets in $\{ Z \subset X : w(Z) \leq \mathfrak{c} \}$
at most $\mathfrak{c}$?
\end{myprob}

\begin{myprob}
\myproblem{D20 first printing}{R.~Hodel}
Does a regular $p$-space (or a $w\Delta$-space) have a countable base if 
it is also:
(1) $\omega_1$-compact with a point-countable separating open cover?
(2) hereditarily c.c.c.?
(3) hereditarily c.c.c.\ with a $G_\delta$ diagonal?

Yes for hereditarily c.c.c.\ with a point-countable separating open cover 
for $w\Delta$-spaces.

\mynote{Solution}
No to (1) (E.~van~Douwen \cite{MR94e:54047}).
No to (2): the Alexandroff double arrow space is a compact hereditarily 
separable and hereditarily Lindel\"of counterexample.
No to (3) if \myaxiom{CH} (I.~Juh\'asz, K.~Kunen and M.E.~Rudin 
\cite{MR55:1270.1}).
Yes to (3) if \myaxiom{PFA}.
\end{myprob}

\begin{myprob}
\myproblem{D20}{G.A.~Edgar}
Suppose that $M(X)$ is the space of all regular Borel measures on a 
compact space $X$.
(1) What is the cardinality of $M(X)$?
(2) Is the density of $M(X)$ equal to the cardinality of $X$?
(3) Is the weight of $X$ equal to the density of $C(X)$?

\mynote{Solution}
For each compact infinite $X$, $|X|^\omega \le |M(X)| \le 2^{|X|}$.
D.~Fremlin and G.~Plebanek \cite{MR1990806} showed that, under 
\myaxiom{MA}, there is a compact $X$ such that $|X| = \mathfrak{c}$ and 
there is a family of cardinality $2^\mathfrak{c}$ of mutually singular 
regular probability measures on $X$.
Also, they showed that in several models of set theory there is a compact 
$X$ such that $|M(X)| > |X|^\omega$.
Regarding (2), let $X$ be the Stone space of the measure algebra of the 
Lebesgue measure on $[0,1]$; $|X|=2^\mathfrak{c}$ while $M(X)$ is separable.
Regarding (3), for a compact $X$, the weight of $X$ and the density of 
$C(X)$ are both equal to the minimal cardinality of a family in $C(X)$ 
separating points of $X$.
\end{myprob}

\begin{myprob}
\myproblem{D21 first printing}{D.~Lutzer}
Let $C_X$ be the set of all bounded real valued continuous functions on 
$X$; let $T$ be the sup-norm topology on $C_X$; and let $T'$ be the 
topology of pointwise convergence.
If $A \subset X$, an \emph{extender from $A$ to $X$} is a function 
$e\colon C_A \to C_X$ such that $e(f)$ extends $f$ to $X$ for all $f$ in 
$C_A$; $e$ is \emph{linear} if $e(f + rg) = e(f) + re(g)$ for all real 
numbers $r$.
\begin{myenumerate}
\item
Is there a continuous in $T$ extender from $\mathbb{N}^*$ into 
$\beta\mathbb{N}$?
\item
Suppose that for every closed subset $A$ of a Moore space $X$ there is a 
continuous in $T$ linear extender from $A$ into $X$.
Is $X$ c.c.c.?
\item
Suppose that every closed subset $A$ of a separable space $X$ there is a 
continuous in $T'$ linear extender from $A$ into $X$.
Must $X$ be collectionwise Hausdorff?
\end{myenumerate}

\mynote{Solution}
The problems were answered for the most part by E.~van~Douwen, D.~Lutzer 
and T.~Przymusi\'nski \cite{MR56:16577}.
\end{myprob}

\begin{myprob}
\myproblem{D21}{Y.~Benyamini}
If a compact space $X$ carries a measure equivalent to the ordinary 
product measure on $\{-1,1\}^\lambda$ for some cardinal $\lambda$, does 
$X$ have an independent family of closed sets of cardinality $\lambda$?

\mynote{Solution}
This has been called Haydon's problem.
Yes for $\lambda=\omega$ is well-known and 
R.~Haydon \cite{MR58:23514,MR80e:46013}
showed that the answer is yes for $\lambda=\mathfrak{c}^{+}$ but no, 
under \myaxiom{CH}, for $\lambda=\omega_1$.
The statement about an independent family of closed sets in $X$
is equivalent to saying that $X$ can be continuously mapped onto
$[0,1]^\lambda$.
D.~Fremlin \cite{MR99d:28019} showed that the answer is positive for 
$\lambda=\omega_1$ under \manotch.
G.~Plebanek \cite{MR98m:28025} showed that the answer is positive for 
every $\lambda \ge \omega_2$ which is the so-called precalibre of 
measure algebras, so in particular yes for $\lambda=\mathfrak{c}$ is 
consistent.
See Plebanek's article \cite{MR1979560} for a survey on this and related 
questions.
\end{myprob}

\begin{myprob}
\myproblem{D22}{E.~van~Douwen}
For which $\lambda$ are there compact homogeneous spaces of cellularity 
$\lambda$?
$\lambda = \aleph_0$ trivially is possible and $\lambda = \aleph_1$ is 
possible if $\diamondsuit$.

\mynote{Notes}
This is van~Douwen's Problem.
See Kunen's article \cite{MR1078652.1}.
\end{myprob}

\begin{myprob}
\myproblem{D23}{E.~van~Douwen}
Is a compact space $X$ nonhomogeneous if it can be mapped continuously onto 
$\beta\mathbb{N}$?
Yes if $w(X) \leq \mathfrak{c}$.
\end{myprob}

\section*{E. Moore space problems}

\begin{myprob}
\myproblem{E1 first printing}{G.M.~Reed}
Is there a collectionwise Hausdorff Moore space that is not normal?

\mynote{Solution}
Yes (M.~Wage \cite{MR53:9157}).
\end{myprob}

\begin{myprob}
\myproblem{E1}{G.M.~Reed}
Is every $\sigma$-discrete collectionwise Hausdorff Moore space metrizable?
No if \myaxiom{MA}.
\end{myprob}

\begin{myprob}
\myproblem{E2}{G.M.~Reed}
Is every $\sigma$-discrete collectionwise Hausdorff Moore space metrizable?
No if \myaxiom{MA}.
\end{myprob}

\begin{myprob}
\myproblem{E3}{G.M.~Reed}
In \visl\ is each normal Moore space completable?

\mynote{Notes}
If \manotch\ there is a normal Moore space which cannot be embedded in 
a developable space with the Baire property.
\end{myprob}

\begin{myprob}
\myproblem{E4 first printing}{G.M.~Reed}
Does every Moore space $X$ have a point-countable separating open cover?
Yes if $|X| \leq \mathfrak{c}$.

\mynote{Solution}
No (D.~Burke \cite{MR55:1318.1}).
M.~Wage constructed a similar example.
\end{myprob}

\begin{myprob}
\myproblem{E4}{E.~van~Douwen}
Can a Moore space of weight $\leq \mathfrak{c}$ (equivalently, cardinality 
$\leq \mathfrak{c}$) be embedded in a separable Moore space if
it is locally compact?
or it has a point-countable base?
or it is metacompact? (equivalently, has a $\sigma$-point-finite basis?)
\end{myprob}

\begin{myprob}
\myproblem{E5}{J.~Green \cite{MR57:7541}}
(1)
Does every noncompact Moore space which is closed in every Moore space in 
which it is embedded have a dense subset which is conditionally compact?
That is, is every noncompact Moore-closed space e-countably compact?
(2)
Does every noncompact Moore-closed space have a noncompact, 
e-count\-ably compact subspace?

\mynote{Notes}
This problem was originally misstated.
The first problem is the closest nontrivial problem.
The second problem is the question that Green seemed most interested in.

\mynote{Solution}
(1) No (R.M.~Stephenson \cite{MR81m:54056.2}).
(2) No if $\mathfrak{b} = \mathfrak{c}$
(H.-X.~Zhou \cite{MR85e:54026}) or if $\mathfrak{a} = \mathfrak{c}$
(P.~Nyikos, A.~Berner and E.~van~Douwen).
\end{myprob}

\begin{myprob}
\myproblem{E6 first printing}{H. Cook}
If $G_1, G_2, \ldots$ is a development for a Moore space $X$ and
$G^*_{n+1}(p) \subset G^*_n(p)$ for all $n$, does every conditionally 
compact subset of $X$ have compact closure?

\mynote{Solution}
No (L.~Gibson \cite{gibson}). 
\end{myprob}

\begin{myprob}
\myproblem{E6}{R.~Telg\'arsky}
Is every normal Moore space the continuous one-to-one preimage of a 
metric space?

\mynote{Solution} 
This is equivalent to the problem: 
\emph{Is every normal Moore space submetrizable?}
No is consistent (D.~Shakhmatov, F.D~Tall, and S.~Watson \cite{stw}).
\end{myprob}

\begin{myprob}
\myproblem{E7 first printing}{G.M.~Reed}
Can every first countable space $X$ of cardinality $\leq \mathfrak{c}$ 
be embedded in a separable first countable space?

\mynote{Solution}
This is independent of \myaxiom{ZFC} (E.~van~Douwen and 
T.~Przymusi\'nski~\cite{MR82j:54051.1}).
\end{myprob}

\begin{myprob}
\myproblem{E7}{W.~Fleissner}
Is there a strongly collectionwise Hausdorff Moore space which is not 
normal?
\end{myprob}

\begin{myprob}
\myproblem{E8 first printing}{G.M.~Reed}
Can every Moore space of cardinality $\leq \mathfrak{c}$ be embedded in 
a separable first countable space?

\mynote{Solution}
This is independent of \myaxiom{ZFC} (E.~van~Douwen and 
T.~Przymusi\'nski~\cite{MR82j:54051.1}).
\end{myprob}

\begin{myprob}
\myproblem{E8}{W.~Fleissner}
Is there a regular para-Lindel\"of space which is not countably 
paracompact?
(or Moore or metacompact or \ldots?)
\end{myprob}

\begin{myprob}
\myproblem{E9 first printing}{G.M.~Reed \cite{MR49:9815}}
Is there a pseudocompact Moore space which contains a copy of every 
metric space of cardinality $\leq \mathfrak{c}$?

\mynote{Solution}
Yes (G.M.~Reed and E.~van~Douwen \cite{MR92d:54035}).
\end{myprob}

\begin{myprob}
\myproblem{E9}{R.~Blair}
Is there a para-Lindel\"of completely regular space $X$ (with $|X|$ 
Ulam-nonmeasurable) that is not realcompact?
\end{myprob}

\begin{myprob}
\myproblem{E10}{F.~Tall}
Is the product of two normal Moore spaces normal?
No if \manotch.

\mynote{Solution}
Yes if \myaxiom{PMEA} (P.~Nyikos).
\end{myprob}

\begin{myprob}
\myproblem{E11}{F.~Tall}
Is every para-Lindel\"of (countably compact, Moore) normal space 
paracompact?
No if \manotch.

\mynote{Solution}
No (C.~Navy).
See the section on Nyikos's survey of two problems.
\end{myprob}

\begin{myprob}
\myproblem{E12}{F. Tall}
Is a normal, locally compact, metacompact space paracompact?

\mynote{Notes}
This is the Arhangel$'$\kern-.1667em ski\u{\i}-Tall Problem.
The answer is independent of \myaxiom{ZFC}.
Yes if \visl.
No is consistent (G.~Gruenhage and P.~Koszmider~\cite{MR97i:54010.3}).
\end{myprob}

\section*{F. Normality of product problems}

\begin{myprob}
\myproblem{F1 first printing}{T.~Przymusi\'nski}
Is there a (first countable separable) paracompact space $X$ such that 
$X^2$ is normal but not paracompact?
Yes if \manotch.

\mynote{Solution}
Yes (T.~Przymusi\'nski \cite{MR80k:54008, MR80m:54010}).
\end{myprob}

\begin{myprob}
\myproblem{F1}{K.~Kunen}
Is there a box product of infinitely many non-discrete spaces which is 
normal but not paracompact?
\end{myprob}

\begin{myprob}
\myproblem{F2}{T.~Przymusi\'nski}
Is there a nonparacompact, separable, first countable space such that 
$X^\omega$ is perfectly normal?
Yes if \manotch.
\end{myprob}

\begin{myprob}
\myproblem{F3}{T.~Przymusi\'nski}
Is there a paracompact, separable, first countable space such that 
$X^\omega$ is normal but not paracompact?
Yes if \manotch.
\end{myprob}

\begin{myprob}
\myproblem{F4 first printing}{T.~Przymusi\'nski}
Is there a locally compact normal space $X$ and a metric space $Y$ such 
that $X \times Y$ is not normal?
Yes if there is a Souslin line.

\mynote{Solution}
Yes (E.~van~Douwen \cite{MR94e:54047}).
\end{myprob}

\begin{myprob}
\myproblem{F4}{A.V.~Arhangel$'$\kern-.1667em ski\u{\i}}
For what classes of spaces is the product of two spaces of covering 
dimension zero always of covering dimension zero?
\end{myprob}

\begin{myprob}
\myproblem{F5}{N.~Howes}
Does linearly Lindel\"of imply Lindel\"of in normal spaces?

$X$ is \emph{linearly Lindel\"of} provided every open cover 
$\{U_r\}_{r \in \mathfrak{m}}$ of $X$ indexed by ordinals with 
$U_r \subset U_s$ for all $r < s$ has a countable subcover.

\mynote{Notes}
This is the linearly Lindel\"of problem. 
See Rudin's contribution to \emph{New Classic Problems}.
\end{myprob}

\begin{myprob}
\myproblem{F6}{N.~Howes}
Is every normal, finally compact in the sense of complete accumulation 
points space Lindel\"of?

A space $X$ is \emph{finally compact in the sense of complete 
accumulation points} provided, for every uncountable regular 
cardinal $\mathfrak{m}$ and $Y \subset X$ with $|Y| = \mathfrak{m}$, 
there is a point $X$ such that every $|U \cap Y| = \mathfrak{m}$ for all 
neighborhoods $U$ of $x$.

This is equivalent to F6.
\end{myprob}

\begin{myprob}
\myproblem{F7}{M.~Starbird}
If $X$ is normal and $C$ is a closed subset of $X$ and 
$f\colon (C \times I) \cup (X \times \{0\}) \to Y$ is continuous, then 
can $f$ be extended to $X \times I$ if $Y$ is an ANR(normal)?

An \emph{ANR(normal)} is an absolute neighborhood retract in every normal 
space in which it is embedded.

\mynote{Notes}
Yes if $Y$ is either an ANR(compact Hausdorff) or a separable 
topologically complete ANR(metric).
See Starbird's papers \cite{MR51:9014,MR53:9132}.
\end{myprob}

\begin{myprob}
\myproblem{F8}{M.~Starbird}
Can $X \times Y$ be Dowker without either $X$ or $Y$ being Dowker?
\end{myprob}

\begin{myprob}
\myproblem{F9}{M.~Starbird \cite{MR52:1606}}
Let $N(X)$ be the class of all spaces whose product with $X$ is normal.
Is $N(X)$ closed under closed maps for paracompact spaces?
for paracompact $p$-spaces?

\mynote{Solution}
A.~Be\v{s}lagi\'c \cite{MR87k:54016} proved that if $X$ is a paracompact 
$p$-space, $X \times Y$ is collectionwise normal, and $f$ is a closed 
map from $Y$ onto $Z$, then $X \times Z$ is collectionwise normal.
\end{myprob}

\begin{myprob}
\myproblem{F10}{K.~Kunen}
Suppose that $T$ is compact and that $Y$ is the image of $X$ under a 
perfect map, $X$ is normal, and $X \times Y$ is normal.
Is $X \times T$ normal?
\end{myprob}

\begin{myprob}
\myproblem{F11}{A.H.~Stone}
Is the box product of $\aleph_1$ copies of $\omega + 1$ normal?
paracompact?

\mynote{Solution}
No (B.~Lawrence \cite{MR96f:54027.1}):
the box product of $\aleph_1$ copies of $\omega + 1$ is neither normal 
nor collectionwise Hausdorff. 
\end{myprob}

\begin{myprob}
\myproblem{F12}{K.~Nagami}
Does $\dim(X \times Y) \leq \dim X + \dim Y$ hold for completely regular 
spaces?

\mynote{Solution}
No.
M.~Wage \cite{MR80a:54064} first constructed counterexamples under 
\myaxiom{CH}.
T.~Przymusi\'nski \cite{MR80f:54033} modified Wage's technique to produce 
many counterexamples (in \myaxiom{ZFC} alone).
The factors can be separable, first countable and either Lindel\"of or 
locally compact.
\end{myprob}

\begin{myprob}
\myproblem{F13}{K.~Nagami \cite{MR46:6327}}
Is the image of a $\mu$-space under a perfect map always a $\mu$-space?
\end{myprob}

\begin{myprob}
\myproblem{F14}{H.H.~Corson \cite{MR21:5947}}
Is a $\Sigma$-product of metric spaces always normal?

\mynote{Solution}
Yes.
This was answered by S.P.~Gul$'$\kern-.1667em ko \cite{MR57:1395}.
M.E.~Rudin \cite{MR85h:54040} proved that the $\Sigma$-product of metric 
spaces has the shrinking property.
\end{myprob}

\section*{G. Continua theory problems}

\begin{myprob}
\myproblem{G1}{P.~Erd\H{o}s}
Is there a connected set in the plane which meets every vertical line in 
precisely two points such that every nondegenerate connected subset meets 
some vertical line in two points?
\end{myprob}

\begin{myprob}
\myproblem{G2}{R.H.~Bing}
If $P$ is the pseudo-arc and $f\colon P \to P$ is continuous and fixed on 
an open set, then is $f$ a homeomorphism?
\end{myprob}

\begin{myprob}
\myproblem{G3}{P.~Erd\H{o}s} 
Is there a widely connected complete metric space?

$X$ is \emph{widely connected} if each nondegenerate connected subset is 
dense.
\end{myprob}

\begin{myprob}
\myproblem{G4}{P.~Erd\H{o}s}
Is there a biconnected space without a dispersion point?
Yes if \myaxiom{CH} (P.~Swingle \cite{MR12:627e}).

$X$ is \emph{biconnected} if it is not the union of two nondense 
connected subsets.

\mynote{Solution} 
V.~Tzannes \cite{MR2000b:54032} constructed two examples of countable, 
biconnected spaces that are not widely connected, do not have a dispersion 
point, and are not strongly connected. The first is Hausdorff and the 
second is Urysohn and almost regular.
Using \myaxiom{MA} for countable posets, M.E.~Rudin \cite{MR96m:54031}
constructed a biconnected subset of the plane the connected subsets of 
which do not have dispersion points and are not widely connected 
either.
\end{myprob}

\begin{myprob}
\myproblem{G5}{R.H.~Bing \cite[Problem 3, p.~75]{MR21:4395}}
Let $S$ be the pseudo-arc and suppose $f\colon S \to S$ is fixed on some 
nonempty open set. 
Is $f$ the identity?

\mynote{Solution}
No (W.~Lewis \cite{MR80m:54053.1}).
\end{myprob}

\begin{myprob}
\myproblem{G6}{H.~Bell}
Is there a compact continuum $K$ of the plane which does not separate the 
plane and a fixed point free map from $K$ to $K$?

\mynote{Notes}
This is the fixed point problem for nonseparating plane continua.
See the survey by C.L.~Hagopian in this volume.
\end{myprob}

\begin{myprob}
\myproblem{G7}{K.~Borsuk}
Given $X \subset \mathbb{E}^3$ such that $X$ is locally connected and 
separates $\mathbb{E}^3$ does there exist a fixed point free map from $X$ 
into $X$?
Can \emph{locally contractible} replace \emph{locally connected}?
\end{myprob}

\begin{myprob}
\myproblem{G8}{H.~Cook}
Is there a hereditarily indecomposable continuum which contains a copy 
of every hereditarily indecomposable continuum?
\end{myprob}

\begin{myprob}
\myproblem{G9}{H.~Cook, Knaster}
Is the pseudo-arc a retract of every hereditarily indecomposable 
continuum in which it is embedded?
\end{myprob}

\begin{myprob}
\myproblem{G10}{A.~Lelek}
Is the confluent image of a chainable continuum chainable?
\end{myprob}

\begin{myprob}
\myproblem{G11}{A.~Lelek}
Does the confluent image of a continuum with span zero have span zero?
\end{myprob}

\begin{myprob}
\myproblem{G12}{H.~Cook}
Suppose that $f_1\colon X_1 \to Y_1$ is confluent and that 
$f_2\colon X_2 \to Y_2$ is confluent. 
If $f_1 \times f_2\colon X_1 \times X_2 \to Y_1 \times Y_2$ confluent?

\mynote{Solution}
No, T.~Ma\'ckowiak \cite{MR53:14390} described a confluent mapping whose 
product with the identity map on the unit interval is not even locally
confluent. 
\end{myprob}

\begin{myprob}
\myproblem{G13}{H.~Cook}
Is every continuum with span zero chainable?

\mynote{Notes}
This was asked by A.~Lelek \cite{MR46:6324.1}.
Chainable continua have span zero.
\end{myprob}

\begin{myprob}
\myproblem{G14}{}
Is there a hereditarily equivalent continuum which is not tree-like?

\mynote{Notes}
Specifically, does there exist a hereditarily equivalent continuum other 
than the arc or the pseudo-arc? 
\end{myprob}

\section*{H. Mapping problems}

\subsection*{Definitions}
Let $f\colon X \to Y$ be a map.
When $S \subset Y$, $f_S$ denotes $f \restriction f^{-1}(S)$.
$f$ is \emph{quotient} if, for all $S \subset Y$, $S$ is 
closed in $Y$ whenever $f^{-1}(S)$ is closed in $X$.
$f$ is \emph{(countably) biquotient} if, for each $y \in Y$, every 
(countable) collection of open sets covering $f^{-1}(y)$ has a finite 
subcollection whose images cover a neighborhood of $y$.
$f$ is \emph{hereditarily quotient} if $f_S$ is quotient for all 
$S \subset Y$.
$f$ is an \emph{s-map} (\emph{L-map}) if $f^{-1}(y)$ is separable 
(Lindel\"of) for all $y \in Y$.
$f$ is \emph{compact covering} if every compact subset of $Y$ is the 
image of some compact subset of $X$.

A space $X$ is of \emph{point-countable type} provided each point has a
sequence $\{U_n : n \in \omega\}$ of neighborhoods such that 
$\bigcap \{U_n : n \in \omega\} = C$ is compact and every neighborhood of 
$C$ contains $U_n$ for some $n$.

A set $\mathcal{G}$ of subsets of a space $X$ is \emph{equi-Lindel\"of} 
if every open cover $\mathcal{H}$ of $X$ has an open refinement with each 
$U \in \mathcal{G}$ intersecting at most countably many 
$V \in \mathcal{H}$.

\startproblem
\begin{myprob}
\myproblem{H1}{E.~Michael}
Is every quotient $s$-image of a metric space also a compact covering 
quotient $s$-image of a metric space?

\mynote{Solution}
This question was asked by E.~Michael and K.~ Nagami \cite{MR46:6269}.
H.~Chen \cite{MR2000d:54015} constructed a counterexample.
\end{myprob}

\begin{myprob}
\myproblem{H2}{E.~Michael}
Characterize those spaces $Y$ such that every closed map 
$f\colon X \to Y$ is countably biquotient (perhaps in terms of sequences 
of subsets of $Y$).
\end{myprob}

\begin{myprob}
\myproblem{H3 first printing}{}
If $X$ is the metrizable image of a complete metric space under a 
$k$-covering map, does $X$ have a complete metric?

\mynote{Solution}
Yes if $X$ is separable (J.P.R.~Christensen \cite{MR48:12488},
A.V.~Ostrovski{\u\i} \cite{MR53:14391}).
See also Michael's article \cite{MR57:7543}.
\end{myprob}

\begin{myprob}
\myproblem{H3}{P.~Nyikos}
If $X$ is locally connected, could every subcontinuum (compact, connected, 
nontrivial) of $X$ contain a copy of $\beta\mathbb{N}$?
\end{myprob}

\begin{myprob}
\myproblem{H4}{E.~Michael \cite{MR50:11129}}
Let $f\colon X \to Y$ be a quotient map and let $E$ be a subset of $Y$ 
such that $\{f^{-1}(y) : y \in E\}$ is equi-Lindel\"of in $X$.
Assume also that $Y$ is an \emph{$A$-space} (whenever $\{F_n : n \in 
\omega\}$ is a decreasing sequence of subsets of $Y$ with a common limit 
point, then there is an $A_n \subset F_n$ with $A_n$ closed such that 
$\bigcup \{A_n : n \in \omega\}$ is not closed).
Is $f_E$ then biquotient? 

Yes if $Y$ is a Hausdorff relatively countably bi-quasi-$k$ space
(R.C.~Olson \cite{MR51:1715}). 
\end{myprob}

\begin{myprob}
\myproblem{H5}{R.C.~Olson}
Suppose that $f\colon X \to Y$ is a quotient $L$-map, $X$ has a 
point-countable base, and $Y$ is of point-countable type.
Is $f$ then biquotient?
\end{myprob}

\begin{myprob}
\myproblem{H6}{R.C.~Olson}
Is there a quotient map $f\colon X \to Y$ with $X$ locally compact and 
first countable, $Y$ compact, each $f^{-1}(y)$ compact, and $f$ 
finite-to-one but not hereditarily quotient?
\end{myprob}

\begin{myprob}
\myproblem{H7}{R.C.~Olson \cite{MR51:1715}}
Is there a paracompact $X$ of point-countable type which does not admit a 
perfect map onto a first countable space?

\mynote{Solution}
Yes, H.~Ohta \cite{MR83f:54022} described a regular Lindel\"of space of 
point-count\-able type which does not admit a perfect map into any space 
in which every point is $G_\delta$.
\end{myprob}

\begin{myprob}
\myproblem{H8}{J.~Nagata}
Is the image of a metric space under a $q$-closed map a $\sigma$-space?
\end{myprob}

\providecommand{\bysame}{\leavevmode\hbox to3em{\hrulefill}\thinspace}

\label{tprudinend}

\chapter*{Problems from A.V.~Arhangel$'$\kern-.1667em ski\u{\i}'s \emph{Structure 
and classification of topological spaces and cardinal invariants}}
\markboth{\normalsize\textsc{\lowercase{Problems from 
Arhangel$'$\kern-.1667em ski\u{\i}'s Structure and classification}}}{}
\label{tpshura}
\begin{myfoot}
\begin{myfooter}
A.V.~Arhangel$'$\kern-.1667em ski\u{\i} and Elliott Pearl, 
\emph{Problems from A.V.~Arhangel$'$\kern-.1667em ski\u{\i}'s \emph{Structure and 
classification of topological spaces and cardinal invariants}},\\ 
Problems from Topology Proceedings, Topology Atlas, 2003, 
pp.\ 123--134.
\end{myfooter}
\end{myfoot}

\mypreface
This section contains the problems that appeared in the seminal 1978 
survey article \emph{Structure and classification of topological spaces 
and cardinal invariants} by A.V.~Arhangel$'$\kern-.1667em ski\u{\i} \cite{MR80i:54005}.
The survey article ended with twenty six problems in a list titled 
\emph{Open problems}. 
Arhangel$'$\kern-.1667em ski\u{\i} wrote: ``I give here only a few problems. 
The solution of many of them seems to me to require original ideas and
methods.''

This version has been prepared with the cooperation of 
A.V.~Arhangel$'$\kern-.1667em ski\u{\i}. 
This version also contains some questions that appeared throughout the 
four chapters of the survey article; some of these questions were 
mentioned merely within the exposition of the survey but some were stated 
explicitly as open problems.
Information on solutions to these problems appeared in volumes 11, 12, 13, 
and 14 of \emph{Topology Proceedings}.
This version includes information on solutions that have appeared since
the survey article was published in 1978.

This version was prepared from the English translation in 
\emph{Russian Math.\ Surveys}.
The problems have been rewritten with current English terminology.
In particular, 
\emph{Lindel\"of} replaces \emph{finally compact};
\emph{perfectly normal} replaces \emph{completely normal};
\emph{realcompact} replaces \emph{functionally closed};
\emph{$\alpha$-expanded} replaces \emph{$\alpha$-extendable};
\emph{cellularity} replaces \emph{Souslin number}; and
\emph{metacompact} replaces \emph{weakly paracompact}.
Compact spaces are assumed to be Hausdorff.
The sectioning and item numbering from the original survey article has 
been preserved.

\section*{Cardinal invariants in broad classes of spaces}

\begin{myprob}
\myproblem{From \S 1.2}{}
Is there, in \myaxiom{ZFC}, a regular space with density greater than its 
spread?
\end{myprob}

\begin{myprob}
\myproblem{From \S 1.3}{}
Is a Moore space $\alpha$-expanded?

\mynote{Notes}
A space $X$ is said to be \emph{$\alpha$-expanded} if there is a linear 
ordering $<$ on $X$, called an $\alpha$-left ordering, such that the set 
$\{y\in X : y\leq x\}$ is closed in $X$ for every $x$ in $X$. 
$\alpha$-expanded spaces were introduced by Arhangel$'$\kern-.1667em ski\u{\i} and 
studied in \cite{MR58:24204,MR58:12885}.
This notion is sometimes translated as \emph{$\alpha$-extended} or
\emph{$\alpha$-extendable}.
See also 
\cite{MR82e:54037a,MR82e:54037b}.
\end{myprob}

\begin{myprob}
\myproblem{\S 1.3 1}{}
Is there, in \myaxiom{ZFC}, a compact radial space $X$ for which 
$c(X) < d(X)$ (i.e., cellularity is strictly less than density)?
\end{myprob}

\begin{myprob}
\myproblem{\S 1.3 2}{}
Let $X$ be a compact radial space.
Is it true, in \myaxiom{ZFC}, that $d(X) \leq (c(X))^+$?
\end{myprob}

\begin{myprob}
\myproblem{\S 1.3 3}{}
Is every pseudoradial space of countable tightness sequential?

\mynote{Solution}
No.
I.~Juh\'asz and W.~Weiss \cite{MR88a:54006} constructed a 
zero-dimensional pseudoradial space of countable tightness which is not 
sequential.
P.~Simon and G.~Tironi \cite{MR88b:54006} constructed a pseudoradial 
Hausdorff space with countable tightness which is not sequential.
Under \myaxiom{CH}, I.~Jan\'e, P.R.~Meyer, P.~Simon, and R.G.~Wilson 
\cite{MR83c:54037} had constructed a pseudoradial Hausdorff space of 
countable tightness space which is not sequential. 
\end{myprob}

\begin{myprob}
\myproblem{\S 1.3 4}{}
Let $X$ be a regular pseudoradial space.
Is it true that $|\overline{A}| \leq 2^{|A|}$ for every $A \subset X$?
\end{myprob}

\begin{myprob}
\myproblem{\S 1.3 5}{G.I.~Chertanov}
Is there, in \myaxiom{ZFC}, a Hausdorff c.c.c.\ radial space that is 
not a Fr\'echet-Urysohn space?
\end{myprob}

\begin{myprob}
\myproblem{\S 1.3 6}{}
Let $X$ be a right-separated compact space.
Is it true that $|\overline{A}| \leq |A|^\omega$ for every $A \subset X$?
\end{myprob}

\begin{myprob}
\myproblem{From \S 1.3}{}
I do not know whether each $\alpha$-expanded compact space is 
pseudoradial.
Is it true that the product of radial (pseudoradial) compact spaces is
pseudoradial?
\end{myprob}

\begin{myprob}
\myproblem{From \S 1.4}{}
Is it possible to estimate the cardinality of a space using 
cellularity, Lindel\"of degree and pseudocharacter?
For regular spaces?
\end{myprob}

\begin{myprob}
\myproblem{From \S 1.4}{}
Is it impossible to reduce (1.4.7) to the case of a regular $X$?
\begin{equation}\tag{1.4.7}
\text{For any $X \in T_1$}, |X| \leq \exp (\psi(X)s(X)).
\end{equation}

\mynote{Solution}
No. See Problem 16 below.
\end{myprob}

\begin{myprob}
\myproblem{From \S 1.4}{}
Could it be that for each $Y \in T_1$ there is a regular space $X$ such 
that $Y$ has the same spread as $X$ and $X$ condenses onto $Y$?
\end{myprob}

\begin{myprob}
\myproblem{From \S 1.5}{}
Is there, in \myaxiom{ZFC}, an uncountable cardinal $\tau$ such that 
there are a pair of spaces $X$, $Y$ such that $X \times Y$ has a pairwise
disjoint family of open sets of cardinality $\tau$ but neither 
$X$ nor $Y$ have a pairwise disjoint family of open sets of cardinality 
$\tau$.
\end{myprob}

\begin{myprob}
\myproblem{From \S 1.5}{}
Is it consistent that $c(X)^+$ is always a precalibre of $X$?

A cardinal $\tau$ is a \emph{precalibre} of $X$ if each family $\gamma$ 
of cardinality $\tau$ of nonempty open sets of $X$ contains a subfamily 
$\gamma'$ of cardinality $\tau$ such that $\gamma'$ has the finite 
intersection property.
\end{myprob}

\begin{myprob}
\myproblem{From \S 1.5}{}
What families of cardinals can be obtained as the collection of all 
calibres of a topological space? (This question has often been mentioned 
in print.)
\end{myprob}

\begin{myprob}
\myproblem{From \S 1.5}{}
Is the $K_0$ property preserved by perfect images?

A space is $K_0$ if it has a dense subspace that is $\sigma$-discrete, 
i.e., the union of a countable family of discrete subspaces.

\mynote{Solution}
No (S.~Todor\v{c}evi\'c \cite{MR84g:03078.1}).
\end{myprob}

\begin{myprob}
\myproblem{From \S 1.6}{}
Is there, in \myaxiom{ZFC}, a regular hereditarily separable countably 
compact noncompact space?

\mynote{Solution}
The answer is no since, consistently, according to S. Todor\v{c}evi\v{c},
every regular hereditarily separable space is Lindel\"of. 
In such a model of set theory every countably compact (even every 
pseudocompact) hereditarily separable space is compact. 
This observation was made in \cite{MR97a:54003}.
\end{myprob}

\begin{myprob}
\myproblem{From \S 1.6}{}
Does $2^{\aleph_0} < 2^{\aleph_1}$ imply the existence of an $L$-space?
Could this be used to get an $L$-space with a pointwise countable basis?
\end{myprob}

\begin{myprob}
\myproblem{From \S 1.6}{}
Does $2^{\aleph_0} < 2^{\aleph_1}$ imply the existence a nonseparable 
perfectly normal compact space?
\end{myprob}

\begin{myprob}
\myproblem{From \S 1.6}{}
Is there, in \myaxiom{ZFC}, a regular space $X$ such that 
$X^\omega$ is hereditarily separable and hereditarily Lindel\"of
but $X$ has uncountable net weight?
Does \manotch\ imply that there are no such spaces?

\mynote{Solution}
K.~Ciesielski \cite{MR89d:03047} constructed a model of \manotch\ where 
such a space exists.
\end{myprob}

\begin{myprob}
\myproblem{From \S 1.7}{}
Is there a nonhomogeneous (compact) space whose square is homogeneous?
Could the product of two nonhomogeneous spaces be homogeneous?
Could a compact space $X$ be nonhomogeneous whereas $X^n$ is homogeneous 
for some $n > 1$?
E.~van~Douwen asked whether a compact space that can be mapped onto 
$\beta \omega$ (or $\beta \tau$), or more generally, onto some compact 
space $Y$ of cardinality $> 2^{\pi w(Y)}$, be homogeneous?

\mynote{Solution}
J.~van~Mill \cite{MR82h:54067} described a rigid infinite-dimensional 
compact space $X$ for which $X \times X$ is homeomorphic to the Hilbert 
cube.
\end{myprob}

\section*{The structure of compact spaces and cardinal invariants}

\begin{myprob}
\myproblem{From \S 2.2}{}
Let $X$ be an infinite homogeneous compact space.
Does $X$ have a dense sequential subspace?
If $X$ is also a group, does it have a dense sequential subspace?

\mynote{Notes}
Yes, if $X$ is an abelian group.
\end{myprob}

\begin{myprob}
\myproblem{From \S 2.2}{}
Can the condition $2^\tau = \tau^+$ be removed in the theorem below?
The conclusion is still true if the set of cardinals between $\tau$ and 
$2^\tau$ is finite.

\begin{theorem}[Shapirovski\u\i]
Let $2^\tau = \tau^+$ and let $X$ be a compact space such that 
$d(Y) \leq \tau$ for each dense subset $Y$ of $X$.
Then $\pi w(X) \leq \tau$.
\end{theorem}
\end{myprob}

\begin{myprob}
\myproblem{From \S 2.3}{}
Characterize internally the class of subspaces of sequential spaces.
\end{myprob}

\begin{myprob}
\myproblem{From \S 2.3}{}
What spaces have a compactification of countable tightness?
\end{myprob}

\begin{myprob}
\myproblem{From \S 2.3}{}
Can a space of countable character be embedded in a countably compact 
space of countable character?
\end{myprob}

\begin{myprob}
\myproblem{From \S 2.3}{}
Can each regular space be embedded in a regular countably compact space of 
the same tightness?
\end{myprob}

\section*{The maps and the structure of compact spaces}

The spaces in this section as assumed to be completely regular.

\startproblem
\begin{myprob}
\myproblem{From \S 3.1}{}
What spaces can be mapped continuously onto $D^\tau$ (or onto 
$\mathbb{I}^\tau$)? 
What spaces can be embedded in a $\Sigma$-product of closed intervals?
What can be said about the continuous images of a $\Sigma$-product of 
closed intervals?

\mynote{Notes}
B.~Shapirovski{\u\i} characterized compact preimages of $D^\tau$.
G.I.~Chertanov characterized subspaces of $\Sigma$-products.
\end{myprob}

\begin{myprob}
\myproblem{\S 3.2 1}{}
Does \manotch\ imply that every compact c.c.c.\ space of countable
tightness has cardinality $\leq \mathfrak{c}$?
\end{myprob}

\begin{myprob}
\myproblem{\S 3.2 2}{}
Does \manotch\ imply that every countably compact hereditarily separable 
space has cardinality $\leq \mathfrak{c}$?
\end{myprob}

\begin{myprob}
\myproblem{\S 3.2 3}{}
Can (3.2.13), (3.2.18)--(3.2.20) and (3.2.25) be generalized to the case
of any cardinal $\tau$? 
\begin{myitemize}
\item[(3.2.13)]
There is no \myaxiom{zfc} example of a nonseparable compact space $X$ for 
which $t(X) = \omega$ and $c(X) = \omega$.
\item[(3.2.18)]
Let $X$ be compact, $t(X) = \omega$ and $c(X) = \omega$.
Assuming \manotch, we have $d(X) = \omega$.
\item[(3.2.19)]
If \myaxiom{ma}, then every compact space $X$ such that $c(X) = \omega$ 
and $\pi w(X) < \mathfrak{c}$ is separable.
\item[(3.2.20)]
If \myaxiom{ma}, then every compact space $X$ such that 
$d(X) < \mathfrak{c}$, $t(X) < \mathfrak{c}$ and $c(x) = \omega$
is separable.
\item[(3.2.25)]
There is no \myaxiom{zfc} example of a homogeneous compact space $X$ for
which $|X| \le \mathfrak{c}$, $c(X) = \omega$ and $d(X) > \omega$.
\end{myitemize}
\end{myprob}

\begin{myprob}
\myproblem{\S 3.2 4}{}
Is there, in \myaxiom{ZFC}, a compact space $X$ for which 
$d(X) > c(X)t(X)$?
\end{myprob}

\begin{myprob}
\myproblem{\S 3.2 5}{}
Is there, in \myaxiom{ZFC}, a nonseparable compact c.c.c.\ space of 
cardinality $\leq \mathfrak{c}$?

\mynote{Solution}
Yes, see Problem 9 below.
\end{myprob}

\begin{myprob}
\myproblem{\S 3.2 6}{}
Is every compact c.c.c.\ space of weight $\aleph_1$ separable?

\mynote{Solution}
Yes if $\text{\myaxiom{MA}}(\omega_1)$; on the other hand, a Souslin 
continuum would be a counterexample.
\end{myprob}

\begin{myprob}
\myproblem{\S 3.2 7}{}
Is there, in \myaxiom{ZFC}, a nonseparable compact c.c.c.\ space of 
weight $\aleph_1$?
\end{myprob}

\begin{myprob}
\myproblem{\S 3.2 8}{}
Is there, in \myaxiom{ZFC}, a nonseparable compact c.c.c.\ space $X$ of 
cardinality $\leq 2^\mathfrak{c}$?
\end{myprob}

\begin{myprob}
\myproblem{From \S 3.2}{}
Does \manotch\ imply that every Lindel\"of c.c.c.\ $p$-space 
of countable tightness is separable?
\end{myprob}

\begin{myprob}
\myproblem{From \S 3.2}{}
Let $X$ be a sequential Lindel\"of $\Sigma$-space.
Is it then true
that $|X| \leq 2^{c(X)}$?

\mynote{Solution} 
No.
The $\sigma$-product of any number of closed intervals has c.c.c., 
is Fr\'echet (hence sequential), and can be of arbitrarily large 
cardinality.
\end{myprob}

\begin{myprob}
\myproblem{From \S 3.2}{}
Let $X$ be a Lindel\"of $\Sigma$-space.
Is it true that $t(X) = 
\sup \{ \tau : \text{there is a free sequence of length
$\tau$ in $X$}\}$?

\mynote{Solution}
No (O. Okunev).
\end{myprob}

\begin{myprob}
\myproblem{From \S 3.3}{}
Is there an infinite extremally disconnected compact space whose character 
at each point is the same?
\end{myprob}

\begin{myprob}
\myproblem{From \S 3.3}{}
Is the cellularity of every reduced extremally disconnected compact 
space countable?

\mynote{Notes}
If $\gamma \subset \mathcal{P}(X)$, $a(\gamma)$ is the smallest family of 
subsets of $X$ such that:
$\gamma \subset a(\gamma)$;
if $U \in a(\gamma)$ then $X \setminus U \in a(\gamma)$;
if $\lambda \subset a(\gamma)$ then $\bigcup \lambda \in a(\gamma)$.
$T_0(X)$ denotes the family of all clopen subsets of $X$.
The \emph{algebraic weight} of an extremally disconnected space $X$ is 
$n(X) = \min \{|\gamma| : \gamma \subset T_0(X), a(\gamma) = T_0(X)\}$.
An extremally disconnected compact space is \emph{reduced} if
$n(U) = n(X)$ for each nonempty clopen subset $U$ of $X$.
\end{myprob}

\begin{myprob}
\myproblem{From \S 3.3}{}
Is there, in \myaxiom{zfc}, a non-discrete extremally disconnected group?

\mynote{Notes}
This is a major old open problem. 
It was first formulated in 1967, in \cite{MR36:5259}.
Consistent examples of such groups were constructed by V.I.~Malykhin and
S.~Sirota.

\section*{Topological properties of mapping spaces}

In this section, maps are not assumed to be continuous.
If there are no separation restrictions indicated, the spaces must be 
regarded as completely regular.

The space \emph{$C_p(X)$} is the set of all continuous real-valued 
functions on $X$ with the topology of pointwise convergence.
A space is a \emph{$\Sigma$-space} if there is a $\sigma$-locally-finite
closed collection $\mathcal{F}$ in $X$ and a cover $\mathcal{C}$ of 
closed countably compact sets such that if $C \subset U$, where $c \in 
\mathcal{C}$ and $U$ is open, then $C \subset F \subset U$ for some 
$F \in \mathcal{F}$.
A \emph{Corson compact} space is a compact subspace of a 
$\Sigma$-product of intervals.
A \emph{Gul$'$\kern-.1667em ko compact} space is a compact space $X$ such that 
$C_p(X)$ is a Lindel\"of $\Sigma$-space; 
S.~Negrepontis introduced this terminology after S.P.~Gul$'$\kern-.1667em ko proved 
that every Gul$'$\kern-.1667em ko compact space is Corson compact.

This section of Arhangel$'$\kern-.1667em ski\u{\i}'s article was the first survey 
of $C_p$ theory. 
See 
\cite{MR89c:54031,MR89m:54022,MR90k:54022,MR1078667,MR1229122,MR92i:54022,MR99i:54001} 
for more information on $C_p$ theory.
\end{myprob}

\startproblem
\begin{myprob}
\myproblem{From \S 4.1}{}
Is there, in \myaxiom{zfc}, a Corson compact space $X$ with $c(X) < w(X)$?

\mynote{Notes}
Assuming \myaxiom{CH} there is a Corson compact space $X$ such that $c(X)$ 
is countable, and $w(X)$ is uncountable. 
Assuming \manotch, no such spaces exist.
\end{myprob}

\begin{myprob}
\myproblem{From \S 4.1}{}
Is there, in \myaxiom{zfc}, a Corson compact space without a dense 
metrizable subspace?

\mynote{Solution}
Yes.
S.~Todor\v{c}evi\'c \cite{MR84g:03078.1} constructed a Corson 
compact space containing no dense metrizable subspace.
A.~Leiderman constructed an adequate example.
\end{myprob}

\begin{myprob}
\myproblem{From \S 4.1}{}
Let $C_p(X)$ be Lindel\"of.
Is $C_p(X) \times C_p(X)$ Lindel\"of?
Is $(C_p(X))^\omega$ Lindel\"of?

\mynote{Notes}
This is a hard open problem.
\end{myprob}

\begin{myprob}
\myproblem{From \S 4.1}{}
Is there an infinite (compact) space $X$ for which $C_p(X)$ is not 
homeomorphic to $C_p(X, \mathbb{R}^\omega)$

\mynote{Solution}
Yes.
This is similar to problem 22 below.
\end{myprob}

\begin{myprob}
\myproblem{From \S 4.1}{}
If $X$ is compact is then 
$l(C_p(X, \mathbb{I}^\omega)) = l(C_p(X, \mathbb{R}^\omega))$?

\mynote{Solution}
Yes.
This follows from the fact that $X$ embeds in $C_p(C_p(X,\mathbb{I}))$ 
and general fact that if $X$ is a compact subspace of $C_p(Y)$, then 
$l(C_p(X)^\omega) \le l(Y^\omega)$.
The fact is still true (and the answer to the
original question is "yes" by the same argument) if we replace
\emph{$X$ is compact} by \emph{$X$ is $\sigma$-compact} 
(O.~Okunev \cite{MR94b:54055}).
\end{myprob}

\begin{myprob}
\myproblem{From \S 4.1}{}
If $X$ is a Corson compact space, is $X$ a Gul$'$\kern-.1667em ko compact space?

\mynote{Solution}
No.
K.~Alster and R.~Pol \cite{MR81h:54019} constructed a Corson compact 
space that was not a Talagrand compact space;
G.~Sokolov showed that their example is not a Gul$'$\kern-.1667em ko compact space.
A.~Leiderman constructed an adequate Corson non-Gul$'$\kern-.1667em ko compact space. 
The results by Sokolov and Leiderman were obtained in 1981 and published 
in \cite{MR86i:54016}.
S.~Argyros \cite[Theorem~6.58]{MR86i:46018} constructed another example.
\end{myprob}

\begin{myprob}
\myproblem{From \S 4.1}{}
If $X$ is a Gul$'$\kern-.1667em ko compact space, is $X$ a Corson compact space?

\mynote{Solution}
Yes.
This is famous result of S.P.~Gul$'$\kern-.1667em ko \cite{MR81b:54017}.
\end{myprob}

\begin{myprob}
\myproblem{From \S 4.1}{}
If $X$ is compact and $C_p(X)$ is a $K_{\sigma\delta}$ space, is then $X$ 
an Eberlein compact space?

\mynote{Notes}
A \emph{$K_{\sigma\delta}$} space is a space that can be represented as 
the intersection of a countable family of spaces each of which is the 
union of a countable family of compact spaces.
A \emph{$K$-analytic} space is the continuous image of a 
$K_{\sigma\delta}$ space.

\mynote{Solution}
No.
M. Talagrand gave an example of a compact $X$ which is not Eberlein but 
$C_p(X)$ is $K_{\sigma\delta}$.
It is unknown if there is an $X$ such that $C_p(X)$ is $K$-analytic
but not $K_{\sigma\delta}$.
\end{myprob}

\begin{myprob}
\myproblem{From \S 4.1}{}
If $X$ is a separable perfectly normal nonmetrizable compact space, is 
then $X$ a Gul$'$\kern-.1667em ko compact space?

\mynote{Solution}
No.
The double arrow space is a counterexample; this was noticed by 
V.V.~Uspenskij.
\end{myprob}

\begin{myprob}
\myproblem{From \S 4.1}{}
Let $X$ be a perfectly normal Gul$'$\kern-.1667em ko compact space.
Is $X$ metrizable?

\mynote{Solution}
Yes. (S.P.~Gul$'$\kern-.1667em ko).
\end{myprob}

\begin{myprob}
\myproblem{From \S 4.2}{}
\nobreak
\begin{myitemize}
\item
If $C_p(X)$ is $Q$-closed in $\mathbb{R}^X$, is then $t_0(X) = \omega$?
\item
Is $t_0(D^\tau) = \omega$ iff $\tau$ is not Ulam measurable?
\item
If $X$ is Lindel\"of, is then $t_0(C_p(X)) = \omega$?
\item
If $X$ is compact and $C_p(X)$ is realcompact, is then $t_0(X) = \omega$?
\item
If $C_p(X)$ is realcompact, is then $t_0(X) = \omega$?
\end{myitemize}

\mynote{Notes}
The \emph{functional tightness} of $X$, $t_0(X)$, is the smallest 
cardinal $\tau$ such that each $\tau$-continuous map is continuous.
A map is \emph{$\tau$-continuous} if its restriction to any subspace 
of cardinality $\tau$ is continuous.
$C_p(X)$ is \emph{$Q$-closed} in $\mathbb{R}^X$ means that for each 
$g \in \mathbb{R}^X \setminus C_p(X)$ we can find a $G_\delta$ set $F$ in 
$\mathbb{R}^X$ such that $g \in F$ and $F \cap C_p(X) = \emptyset$.
$D$ is the discrete space consisting of two points.

The \emph{weak functional tightness} of a $X$, $t_m(X)$, is
the smallest cardinal $\tau$ such that if $f$ is a real-valued 
function on $X$ such that for every $A \subset X$ with $|A| \le \tau$ 
there exists a continuous real-valued function $g_A$ on $X$ which
coincides with $f$ on $A$, then $f$ is continuous.
The \emph{Hewitt number} of $X$, $q(X)$, is the smallest 
cardinal $\tau$ such that for every $x \in \beta X \setminus X$ 
there exists a family $\gamma$ of open subsets of $\beta X$ such
that $x \in \bigcap \gamma \subset \beta X \backslash X$ and 
$|\gamma| \leq \tau$.

Note that $q(X) = \omega$ iff $X$ is realcompact.
Arhangel$'$\kern-.1667em ski\u{\i} proved that $t_m(X) = q(C_p(X))$ and 
$t_m(C_p(X)) \ge q(X)$.
Trivially, $t_m(X) \le t_0(X)$.

\mynote{Solution}
These questions have all been answered.
V.V.~Uspenskij \cite{MR84m:54015} proved that 
$t_m(C_p(X)) \le t_0(C_p(X)) \le q(X)$.
Also, If $\tau$ is a nonmeasurable cardinal, then 
$t_0(\mathbb{R}^\tau) = \omega$.
\end{myprob}

\begin{myprob}
\myproblem{From \S 4.4}{}
Is there a nonmetrizable countable Fr\'echet-Urysohn 
Eberlein-Grothen\-dieck space?

\mynote{Notes}
A space is an \emph{Eberlein-Grothendieck space} (or EG-space) if it 
can be embedded in $C_p(X)$ for some compact $X$.

\mynote{Solution}
Yes.
E.G.~Pytkeev \cite{MR84c:54020} constructed a countable, nonmetrizable 
subspace $S \subset C_p(K,D)$ with the Fr\'echet-Urysohn property. 
Here $K$ denotes the Cantor set and $D=\{0,1\}$.
\end{myprob}

\begin{myprob}
\myproblem{From \S 4.4}{}
Is every countable bisequential space an EG-space?

\mynote{Solution}
No. 
M.~Sakai \cite{MR96i:54010} constructed a countable bisequential space 
which is not an EG-space.
Sakai asked some questions about EG-spaces and $\kappa$-metrizable spaces:
\begin{myenumerate}
\item
Is every countable EG-space $\kappa$-metrizable?
Equivalently, is every countable subspace of $C_p(C)$ 
$\kappa$-metrizable?
\item
Is every (countable) stratifiable $\kappa$-metrizable space an EG-space?
\item
Is there a universal space for all countable stratifiable 
$\kappa$-metrizable spaces?
\end{myenumerate}
\end{myprob}

\begin{myprob}
\myproblem{From \S 4.4}{}
Is the image of an EG-space under a perfect map an EG-space?
\end{myprob}

\begin{myprob}
\myproblem{From \S 4.4}{}
Is there a Lindel\"of EG-space whose square is not Lindel\"of?

\mynote{Solution} 
Yes, there are lots of examples of Lindel\"of EG-spaces
with various behaviours of the Lindel\"of property in powers.
See \cite{MR96k:54023}.
\end{myprob}

\begin{myprob}
\myproblem{From \S 4.4}{}
Is every EG-space having a uniform basis metrizable?
Are there metacompact nonparacompact EG-spaces?
\end{myprob}

\begin{myprob}
\myproblem{From \S 4.4}{H.H.~Corson} 
If $C_p(X)$ is normal, is $(C_p(X))^2$ normal?
\end{myprob}

\begin{myprob}
\myproblem{From \S 4.4}{}
Is there a compact space $X$ of uncountable tightness for which $C_p(X)$ 
is normal?

\mynote{Solution}
N.V.~Velichko \cite{MR84c:54022} proved that if $X$ is a compact space 
and $C_p(X)$ is normal then $X$ has countable tightness.
\end{myprob}

\section*{Open problems}

\begin{myprob}
\myproblem{1}{}
Does there exist, in \myaxiom{ZFC}, a compact Hausdorff space of 
countable tightness that is not sequential?

\mynote{Solution}
This is the Moore-Mr\'owka problem.
No, under \myaxiom{PFA} (Z.~Balogh \cite{MR89h:03088.3}).
See Classic Problem VI.
\end{myprob}

\begin{myprob}
\myproblem{2}{}
It is true, in \myaxiom{ZFC}, that each nonempty sequential compact space 
is first countable at some point?

\mynote{Solution}
No.
V.I.~Malykhin \cite{MR88g:54045} showed that in the model produced by 
adding one Cohen real to a model of $\mathfrak{p} = \mathfrak{c} > 
\omega_1$, there is a Fr\'echet-Urysohn compact space without points of 
countable character.

This question was asked in \cite{MR41:7607.1}.
It was known that \myaxiom{CH} implies an affirmative answer 
(S.~Mr\'owka).

A.~Dow \cite{MR91a:54003} showed that \myaxiom{PFA} implies that every 
countably tight compact space has points of first countability.
P.~Koszmider \cite{MR99m:03099}
showed that consistently even a continuous image of a first countable 
compact space (therefore, a bisequential compact space) needn't have 
points of first countability.
\end{myprob}

\begin{myprob}
\myproblem{3}{}
Is it true, in \myaxiom{ZFC}, that if $X$ is a homogeneous sequential
compact space, then $X$ is first countable? 

\mynote{Notes}
Yes if \myaxiom{CH}; 
if $X$ is a sequential homogeneous compact space, then 
$|X| \leq \mathfrak{c}$ \cite{MR41:7607.1}.
Also, if $X$ is a homogeneous compact space, then
$2^{\chi(X)} \leq 2^{\pi(X)}$.
J.~van~Mill \cite{MR1968433} has shown that the existence of a non-first 
countable homogeneous compact space of countable $\pi$-weight is 
independent of \myaxiom{ZFC}.
\end{myprob}

\begin{myprob}
\myproblem{4}{}
Let $b\mathbb{N}$ be a Hausdorff compactification of the discrete space
$\mathbb{N}$ such that $b\mathbb{N} \setminus \mathbb{N}$ is sequential
and compact.
Is it true, in \myaxiom{ZFC}, that $b\mathbb{N}$ is sequential?

\mynote{Notes}
Equivalently, is there a sequence in $\mathbb{N}$ converging to a point of 
$b\mathbb{N} \setminus \mathbb{N}$?
The existence of such a Hausdorff compactification is equivalent to the 
problem of finding a compact space $X = \bigcup \{X_n : n \in \omega\}$,
where each $X_n$ is sequential and compact, such that $X$ is not
sequential.
P.~Simon conjectured that if a Hausdorff compactification of the discrete 
space $\mathbb{N}$ is such that there is no sequence in $\mathbb{N}$ 
converging to a point of $b\mathbb{N} \setminus \mathbb{N}$ then there is 
a continuous map from $b\mathbb{N} \setminus \mathbb{N}$ onto 
$\mathbb{I}^{\omega_1}$, the Tychonoff cube of weight $\omega_1$.
\end{myprob}

\begin{myprob}
\myproblem{5}{}
Is there a nonmetrizable homogeneous Eberlein compact space?

\mynote{Solution}
Yes. 
J.~van~Mill \cite{MR84h:54036} constructed a nonmetrizable homogeneous 
Eberlein compact space which is also hereditarily normal, first countable, 
and zero-dimensional.
\end{myprob}

\begin{myprob}
\myproblem{6}{}
Let $X$ be compact.
Is $\pi\chi(x, X) \leq t(x, X)$ for every point $x \in X$?
Yes, if \myaxiom{gch}.

\mynote{Notes}
If $X$ is compact, $h\pi\chi(X) = t(X)$ (B.~Shapirovski\u\i).
\end{myprob}

\begin{myprob}
\myproblem{7}{}
Does there exist, in \myaxiom{ZFC}, a compact Fr\'echet-Urysohn space 
whose square is not Fr\'echet-Urysohn?

\mynote{Solution}
Yes.
P.~Simon \cite{MR82a:54038} constructed a compact Fr\'echet-Urysohn space 
whose square is not Fr\'echet-Urysohn.
\end{myprob}

\begin{myprob}
\myproblem{8}{}
Does there exist, in \myaxiom{ZFC}, a regular space $X$ such that
$hl(X^n) \leq \tau$ for all $n \in \mathbb{N}^+$ and $d(X) > \tau$?

\mynote{Solution}
Yes.
I.~Juh\'asz and S.~Shelah \cite{MR87f:03143.1} showed that it is 
consistent there are regular hereditarily Lindel\"of spaces of weight 
$2^\mathfrak{c}$.
Furthermore, such models can be found in which $\mathfrak{c}$ is
arbitrarily large and $2^\mathfrak{c}$ is arbitrarily larger.
\end{myprob}

\begin{myprob}
\myproblem{9}{}
Does there exist, in \myaxiom{ZFC}, a nonseparable compact c.c.c.\ 
space of cardinality $\leq \mathfrak{c}$?

\mynote{Solution}
Yes.
S.~Todor\v{c}evi\'c and B. Veli\v{c}kovi\'c \cite{MR89a:03094} 
constructed a c.c.c.\ nonseparable compact Hausdorff space of cardinality 
$\mathfrak{c}$.
\end{myprob}

\begin{myprob}
\myproblem{10}{}
Does every infinite homogeneous compact space contain a nontrivial 
convergent sequence?

\mynote{Notes}
This question is due to W.~Rudin \cite{MR18:324d.1} from 1956.
Yes, if $X$ is also a group.
\end{myprob}

\begin{myprob}
\myproblem{11}{}
Is each regular left-separated space zero-dimensional?

\mynote{Solution}
No.
M.~Tka\v{c}enko \cite{MR83d:54033} constructed examples of completely 
regular pseudocompact connected left-separated spaces.
One example was even a topological group.
Tka\v{c}enko asked if there is a normal connected left-separated space.
Using \myaxiom{CH}, I.~Juh\'asz and N.~Yakovlev \cite{MR88f:54040} 
constructed a regular, hereditarily Lindel\"of (and hence normal), 
connected, left-separated space.
\end{myprob}

\begin{myprob}
\myproblem{12}{}
Does every completely regular space contain a dense zero-dimen\-sion\-al
subspace?

\mynote{Solution}
No. K.~Ciesielski \cite{MR87d:54052} showed that for any cardinal $\kappa$ 
if $2^\kappa=\kappa^+$ then there exists a completely regular space 
without any uncountable zero-dimensional subspace. In particular,
under \myaxiom{CH} this gives an example of a left separated $L$-space of 
type $\omega_1$ without any uncountable zero-dimensional subspace.
A related result is also proved in \cite{MR93k:54010}.
\end{myprob}

\begin{myprob}
\myproblem{13}{} 
Does \manotch\ imply that regular first countably 
hereditarily separable spaces are Lindel\"of?
That is, there are no first countable $S$-spaces.

\mynote{Solution}
No.
U.~Abraham and S.~Todor\v{c}evi\'c \cite{MR86h:03092.1} showed that it is 
consistent with \manotch\ that there is a first countable $S$-space.
\end{myprob}

\begin{myprob}
\myproblem{14}{}
Is it true, in \myaxiom{ZFC}, that there is an $S$-space iff there is 
an $L$-space?

\mynote{Solution}
No, S.~Todor\v{c}evi\'c \cite{MR90d:04001.2} showed that there is a model 
of \myaxiom{MA} in which there is an $L$-space but there are no $S$-spaces. 
\end{myprob}

\begin{myprob}
\myproblem{15}{} 
Let $X$ be a regular space with a $G_\delta$-diagonal and countable
pseudocharacter (i.e., points $G_\delta$).
Is $|X| \leq 2^\mathfrak{c}$?
Is $|X| \leq 2^{2^\mathfrak{c}}$? 
\end{myprob}

\begin{myprob}
\myproblem{16}{}
Let $X$ be a regular c.c.c.\ space with a $G_\delta$-diagonal.
Is $|X| \leq 2^\mathfrak{c}$?

\mynote{Solution}
No.
This problem was asked by J.~Ginsburg and R.G.~Woods \cite{MR57:1392}.
D.B.~Shakhmatov \cite{MR86g:54007} showed that there is no upper bound on 
the cardinality of Tychonoff c.c.c.\ spaces with a $G_\delta$-diagonal.
V.V.~Uspenskij \cite{MR86f:54009} proved that for each infinite
cardinal $\kappa$, there is a completely regular space $X$ with these
properties: $|X| = \kappa$; $X$ is c.c.c.; $X$ is $F_\sigma$-discrete and
hence has a $G_\delta$-diagonal; $X$ is Fr\'echet and hence has countable 
tightness.
\end{myprob}

\begin{myprob}
\myproblem{17}{}
Does the existence of a regular Luzin space imply that there is a
nonseparable perfectly normal compact space?

\mynote{Solution}
No (S.~Todor\v{c}evi\'c \cite{MR97j:03099.1}).

\end{myprob}

\begin{myprob}
\myproblem{18}{} 
Is there, in \myaxiom{ZFC}, a regular Lindel\"of space of countable 
pseudocharacter (i.e., points $G_\delta$) and cardinality 
$> \mathfrak{c}$?

\mynote{Notes}
This problem was first formulated in 1969 in \cite{MR40:4922}.
This is the Lindel\"of points $G_\delta$ problem.
See the contribution by F.D.~Tall to \emph{New Classic Problems}.
\end{myprob}

\begin{myprob}
\myproblem{19}{} 
Is there, in \myaxiom{ZFC}, a first countable compact space whose 
density is different from its cellularity?
\end{myprob}

\begin{myprob}
\myproblem{20}{}
Is there, in \myaxiom{ZFC}, a regular semi-stratifiable Lindel\"of space 
of uncountable net-weight?
\end{myprob}

\begin{myprob}
\myproblem{21}{}
Is there an infinite-dimensional linear topological space (over
$\mathbb{R}$) that is not homeomorphic to its square?

\mynote{Solution}
Yes. (R.~Pol \cite{MR85k:57014}, J.~van~Mill \cite{MR88k:57023}).
\end{myprob}

\begin{myprob}
\myproblem{22}{}
Is $C_p(X)$ homeomorphic to $C_p(X) \times C_p(X)$ for every infinite 
compact space $X$?

\mynote{Solution}
No.
W.~Marciszewski \cite{MR89h:46039} constructed a compact separable space 
$X$ with third derived set empty with the property that $C(X)$, the 
continuous functions on $X$ with either the weak or pointwise topology, 
is not homeomorphic to $C(X) \times C(X)$. 
This problem was also solved negatively by S.P.~Gul$'$\kern-.1667em ko 
\cite{MR91c:54023}.
\end{myprob}

\begin{myprob}
\myproblem{23}{}
Is there a nonseparable regular Lindel\"of symmetrizable space?

\mynote{Solution}
D.B.~Shakhmatov \cite{MR93m:54010} showed that it is consistent 
that there is a symmetrizable, completely regular, zero-dimensional, 
hereditarily Lindel\"of, $\alpha$-left, nonseparable space of size 
$\aleph_1$.
Furthermore, the space can be condensed onto a space with a countable basis.
\end{myprob}

\begin{myprob}
\myproblem{24}{}
Is $t(X \times X) = t(X)$ for each countably compact completely regular
space $X$?
\end{myprob}

\begin{myprob}
\myproblem{25}{}
Is there, in \myaxiom{ZFC}, a compact space $X$ for which 
$c(X \times X) > c(X)$?

\mynote{Solution}
Yes.
S.~Todor\v{c}evi\'c \cite{MR88h:54009.1} showed that cellularity is not 
productive in the class of compact topological spaces.
This problem was asked by D.~Kurepa~\cite{MR31:68}.
\end{myprob}

\begin{myprob}
\myproblem{26}{}
Does \myaxiom{CH} alone imply that there exists a compact space of 
countable tightness that is not sequential?
Yes, if $\diamondsuit$.

\mynote{Notes}
That is, can \myaxiom{CH} decide the Moore-Mr\'owka problem?
T.~Eisworth \cite{MR1992528} showed that there is a totally proper 
forcing notion that will destroy a fixed counterexample to the 
Moore-Mr\'owka problem, but it is not clear if it can be iterated safely 
without adding reals.
\end{myprob}

\begin{myprob}
\myproblem{27}{} 
Suppose that $X$ is a regular c.c.c.\ symmetrizable space.
Is $|X| \le \mathfrak{c}$?

\mynote{Notes}
This problem was formulated in \cite{MR81b:54005} around 1979 and 
Arhangel$'$\kern-.1667em ski\u{\i} asked that it should be added to this version of the 
list.
\end{myprob}

\providecommand{\bysame}{\leavevmode\hbox to3em{\hrulefill}\thinspace}

\label{tpshuraend}

\chapter*{A note on P.~Nyikos's \emph{A survey of two problems in topology}}
\label{tptwonyikos}
\begin{myfoot}
\begin{myfooter}
Elliott Pearl, 
\emph{A note on P.~Nyikos's \emph{A survey of two problems in topology}},\\
Problems from Topology Proceedings, Topology Atlas, 2003, 
pp.\ 135--138. 
\end{myfooter}
\end{myfoot}

\mypreface
In volume 3 (1978) of \emph{Topology Proceedings}, Peter J.~Nyikos wrote 
\emph{A survey of two problems in topology} \cite{MR82i:54046} about 
the $S$- and $L$-space problems and problems about para-Lindel\"of 
spaces.
In this version, the statements of the problems are extracted from 
the original article by Nyikos.
Some current information will follow.

\section*[The $S$ and $L$ problem]{The $\mathbf{S}$ and $\mathbf{L}$ problem}

Is there an $S$-space?
Is there an $L$-space?

An $S$-space is a regular, hereditarily separable, not hereditarily 
Lindel\"of space.
An $L$-space is a regular, hereditarily Lindel\"of, not hereditarily 
separable space.
(In this problem, all spaces are regular Hausdorff spaces.)

\subsection*{Related problems}

\subsubsection*{A}
Does there exist a countably compact $S$-space?
(This remains unsolved if \emph{regular} is dropped in the the definition 
of $S$-space.)
Does there exist an $L$-space in which every countable subset is closed?
\subsubsection*{B}
Is there an $S$-space of cardinality $> \mathfrak{c}$?
Is there an $L$-space of weight $> \mathfrak{c}$?
\subsubsection*{C}
Does there exist a perfectly normal, or a hereditarily normal $S$-space?
\subsubsection*{D}
Does there exist a first countable $S$-space?
\subsubsection*{E}
Does there exist a locally connected $S$ or $L$-space?
\subsubsection*{F}
Does there exist a space of countable spread which is not the union of a 
hereditarily separable and a hereditarily Lindel\"of space?
\subsubsection*{G}
Does the existence of an $S$-space in a given model of set theory imply 
the existence of an $L$-space, and conversely?
\subsubsection*{H}
Does there exist a cardinal $\alpha$ for which there exists a space with 
no discrete subspace of cardinality $\alpha$, but which is not 
$\alpha$-separable? not $\alpha$-Lindel\"of?

\section*{Para-Lindel\"of spaces}

The main problem in this area is the following: 
Is every regular para-Lindel\"of space paracompact?
(A space is \emph{para-Lindel\"of} if every open cover has a locally 
countable open refinement.)
Equivalently:
Is every regular para-Lindel\"of space normal?
[This is an observation of J.~van~Mill:
if there exists a para-Lindel\"of space $X$ which is not paracompact, then 
by Tamano's theorem, $X \times \beta X$ is not normal; and clearly, 
$X \times \beta X$ is still para-Lindel\"of.]

The subject of para-Lindel\"of spaces is a wide open field, with very 
little known about which implications hold between covering or separation 
axioms (regular or beyond), besides those that hold for topological spaces 
in general.
Consider the following properties:
regular, completely regular, normal, collectionwise normal, countably 
metacompact, countably paracompact, realcompact, (weakly) submetacompact, 
metacompact, paracompact.
It is not known whether para-Lindel\"of together with any of these 
properties implies another property if it does not already do so for all 
spaces.

We do not even know whether every para-Lindel\"of normal Moore space is 
metrizable, nor whether every para-Lindel\"of Moore space is normal 
(despite being strongly collectionwise Hausdorff \cite{MR80j:54020.2})
or metacompact.

We do not know whether, on the one hand, every normal space with a 
$\sigma$-locally countable base is metrizable, or, on the other, whether 
it is consistent that there be a normal Moore space with a 
$\sigma$-locally countable base which is not metrizable.

We do not know of a \emph{real} example of a normal space with a 
point-countable base which is not paracompact.

Worst of all, we do not know what para-Lindel\"of adds to having a 
$\sigma$-locally countable base.
For all we know, every para-Lindel\"of space with a $\sigma$-locally 
countable base may be metrizable (equivalently, paracompact); on the other 
hand, there may even be ones that are not countably metacompact, or 
completely regular.

\section*{Twenty-five years later}

Here is some current information on these topics.
The information on $S$- and $L$-spaces comes from 
J.~Roitman's survey \cite{MR87a:54043.1}.
The information on para-Lindel\"of spaces comes from 
S.~Watson's articles \cite{MR1078640.2,MR94f:54048}.

\subsection*{The $\mathbf{S}$ and $\mathbf{L}$ problem}

It was long known that a Souslin line is an $L$-space and
M.E.~Rudin constructed an $S$-space from a Souslin line.
T.~Jech had proved the consistency of the existence of a Souslin line.
Many relative constructions of $S$- and $L$-spaces turned out to be
inconsistent with \manotch.
However, it was shown that \manotch\ is consistent with the existence of
$S$-spaces (Z.~Szentmikl\'ossy).
U.~Abraham and S.~Todor\v{c}evi\'c \cite{MR86h:03092.2} showed that 
\manotch\ is consistent with the existence of first countable $S$-spaces, 
and furthermore this proof dualizes to get the consistency of \manotch\ 
with the existence of $L$-spaces.
Todor\v{c}evi\'c proved that it is consistent that there are no 
$S$-spaces (even while $L$-spaces may exist).

The remaining open problem is whether there is an $L$-space or whether it
is consistent that there are no $L$-spaces.

There are several surveys with more information on these results and
related problems on $S$- and $L$-spaces:
I.~Juh\'asz \cite{MR81j:54001.1}; J.~Roitman \cite{MR87a:54043.1};
M.E.~Rudin \cite{MR81d:54003.1}; S.~Todor\v{c}evi\'c \cite{MR90d:04001.3}.

\subsection*{Navy's examples}

The main problem about para-Lindel\"of spaces was answered by C.~Navy 
in 1981 in her thesis \cite{Navy.1}.
She constructed several examples of normal para-Lindel\"of spaces that 
failed to be paracompact.
Actually, her examples were all countably paracompact and 
not collectionwise normal.

Navy's technique was rather general.
Using Bing's space $G$, she modified an example of W.~Fleissner 
which was $\sigma$-para-Lindel\"of but not paracompact
to obtain a para-Lindel\"of normal non-collectionwise-normal space.
Using normality, it was possible to separate the countably many locally 
countable families so that one locally countable refinement was obtained.
See Fleissner's \emph{Handbook of set-theoretic topology} article 
\cite[\S~6]{MR86m:54023.2} for a description of this example.

Under \manotch, she obtained a para-Lindel\"of nonmetrizable normal Moore 
space by using the Moore plane.

\subsection*{Navy's problems}

In her thesis, Navy asked some interesting questions dealing with 
regular nonparacompact para-Lindel\"of spaces.

\subsubsection*{1}
Without assuming any extra set-theoretic axioms, can one construct such a 
space which is first countable?
\subsubsection*{2}
Is there any such space which is not countably paracompact?
\subsubsection*{3}
Is there any such space which is collectionwise normal?
D.~Palenz \cite{palenzthesis} has shown that every para-Lindel\"of, 
monotonically normal space is paracompact.
She also showed that every monotonically normal space with a 
$\sigma$-locally countable base is metrizable, an extension of 
Fedor\v{c}uk's theorem.
\subsubsection*{4}
Is there any such space which is normal as well as screenable?
Is there any such space which is normal and has a $\sigma$-disjoint base?

\subsection*{Moore spaces}

Fleissner modified Navy's \manotch\ example of a para-Lindel\"of 
nonmetrizable normal Moore space to obtain a nonmetrizable 
normal Moore space under \myaxiom{CH}, thus solving the normal Moore 
space conjecture.
Fleissner's example is para-Lindel\"of too.
Watson asked if the existence of a nonmetrizable normal Moore space 
implies the existence of a para-Lindel\"of nonmetrizable normal Moore
space. 
Watson had in mind Fleissner's example of nonmetrizable normal Moore 
space under \myaxiom{SCH}.
Watson asked whether Fleissner's \myaxiom{SCH} example could be modified 
to be para-Lindel\"of, or whether a negative result could be found which 
would really illustrate the difference between Fleissner's \myaxiom{CH} 
and \myaxiom{SCH} examples.

\subsection*{Watson's example}

In \cite{MR94f:54048}, Watson constructed spaces in which the properties 
such as collectionwise normal Hausdorff or para-Lindel\"of are built 
directly into the construction.
Watson described a technique for coding a class of zero-dimensional 
para-Lindel\"of Hausdorff spaces.
Furthermore, this technique can be used to yield non-collectionwise 
normal examples.
To compare techniques, recall that Navy's space was designed to be normal 
and $\sigma$-para-Lindel\"of; para-Lindel\"of but not directly so.

\subsection*{Watson's problems}

Watson's contribution to \emph{Open Problems in Topology}
stated some open problems about para-Lindel\"of spaces.

\subsubsection*{Problem 107}
Are para-Lindel\"of regular spaces countably paracompact?
\subsubsection*{Problem 108}
Is there a para-Lindel\"of Dowker space?
\subsubsection*{Problem 109} (Fleissner and Reed \cite{MR80j:54020.2})
Are para-Lindel\"of collectionwise normal spaces paracompact?
\subsubsection*{Problem 110}
Is it consistent that meta-Lindel\"of collectionwise normal spaces are 
paracompact? 
\subsubsection*{Problem 111}
Are para-Lindel\"of screenable normal spaces paracompact?
\subsubsection*{Problem 112}
Are para-Lindel\"of collectionwise normal spaces normal?

\subsection*{Problem 107}
This is now the main open problem on para-Lindel\"of spaces.
Navy's constructions are intrinsically countably paracompact.
Watson suggested that the most likely way to obtain a (consistent) example 
of a para-Lindel\"of space which is not countably paracompact could be to 
iterate a normal para-Lindel\"of space which is not collectionwise normal 
in an $\omega$-sequence to get a para-Lindel\"of Dowker space.

\subsection*{Problem 110}
R.~Hodel \cite{MR50:8400.2} first asked if meta-Lindel\"of collectionwise 
normal spaces are paracompact.
M.E.~Rudin's \visl\ example of a normal screenable nonparacompact space
is a consistent counterexample.
Z.~Balogh constructed two \myaxiom{ZFC} counterexamples:
a hereditarily meta-Lindel\"of, hereditarily collectionwise normal 
hereditarily realcompact Dowker space \cite{MR2003c:54047.1};
a meta-Lindel\"of, collectionwise normal, countably paracompact space 
which is not metacompact.
Balogh \cite{MR2003c:54047.1} asked if there is a para-Lindel\"of 
collectionwise normal Dowker space.

\providecommand{\bysame}{\leavevmode\hbox to3em{\hrulefill}\thinspace}

\label{tptwonyikosend}

\chapter*{A note on \emph{Open problems in infinite-dimensional topology}}
\label{tpinfdim}
\begin{myfoot}
\begin{myfooter}
Elliott Pearl, 
\emph{A note on \emph{Open problems in infinite-dimensional topology}},\\
Problems from Topology Proceedings, Topology Atlas, 2003, 
pp.\ 139--140.
\end{myfooter}
\end{myfoot} 

\section*{}
There is a well-known list of problems in infinite-dimensional topology 
with a long history.
It appeared as an appendix to T.A.~Chapman's 1975 volume, \emph{Lectures 
on Hilbert Cube Manifolds} \cite{MR54:11336}, in the CBMS series of the 
American Mathematical Society.
Ross Geoghegan edited a version of the problem list \cite{MR82a:57015} in 
\emph{Topology Proceedings} as the result of a satellite meeting of 
infinite-dimensional topologists held at the 1979 Spring Topology 
Conference in Athens, OH.
The problem list was updated and revised by James West \cite{MR1078666}
in 1990 for the book \emph{Open problems in topology}. 
The book \cite{MR92c:54001} is no longer available in print but the 
publisher has made it freely available online.
For updates to this problem list, please find the series of status 
reports that have appeared in the journal \emph{Topology and its 
Applications}
\cite{opit1,MR92i:54001,MR93k:54001,MR94i:54001,MR95m:54001,
MR1467216,MR1838331,opit8}.
This long problem list will not be reproduced here. 

For a basic introduction to infinite-dimensional topology, please see
Chapman's volume \cite{MR54:11336}, Cz.~Bessaga and A.~Pelczy\'nski's 
monograph \emph{Selected topics in infinite-dimensional
topology} \cite{MR57:17657}, or the books by J.~van~Mill 
\cite{MR90a:57025,MR2002h:57031}.
See also the surveys by 
A.N.~Dranishnikov \cite{dranishnikov} and J.~Dydak 
\cite{MR2003a:55002}.

\providecommand{\bysame}{\leavevmode\hbox to3em{\hrulefill}\thinspace}

\label{tpinfdimend}

\chapter*{W.R.~Utz: Non-uniformly continuous homeomorphisms with uniformly 
continuous iterates}
\markboth{\normalsize\textsc{\lowercase{W.R.~Utz: Uniformly continuous 
iterates}}}{}
\label{tputz}
\begin{myfoot}
\begin{myfooter}
W.R.~Utz, 
\emph{Non-uniformly continuous homeomorphisms with uniformly continuous 
iterates},\\
Problems from Topology Proceedings, Topology Atlas, 2003, 
pp.\ 141--142. 
\end{myfooter}
\end{myfoot}

\mypreface
This article is reprinted whole:
W.R.~Utz, \emph{Non-uniformly continuous homeomorphisms with uniformly 
continuous iterates}, Topology Proceedings \textbf{6}, no.~2, 
(1981) 449--450.

\section*{}

It is not difficult to find examples of self-homeomorphism of a metric 
space which are not uniformly continuous but which have some uniformly 
continuous powers.

My purpose is to raise the question of what variety of powers of a 
non-uniformly continuous homeomorphism may be uniformly continuous.
My particular interest is in self-homeomorphisms of the reals.
The following theorem gives some information.

\begin{theorem}
Corresponding to any integer $n > 1$ there exists a self-homeo\-morphism, 
$f$ of the reals such that $f, f^2, f^3, \ldots, f^{n-1}$ are not 
uniformly continuous but $f^n$ is uniformly continuous.
\end{theorem}

Clearly, for such an $f$ if follows that $f^{-1}$ is not uniformly 
continuous.
Also, it is trivial that a homeomorphism and all its positive powers may 
be uniformly continuous but the negative iterates are non-uniformly 
continuous, etc.
It will be clear from the proof of the theorem of the theorem that the 
same theorem holds for any Euclidean space.

The question posed here is to describe all subsets, $z$, of $\mathbb{Z}$ 
for which one may find a self-homeomorphism, $f$, of the real which is not 
uniformly continuous but if $j \in z$ then $f^j$ is uniformly continuous.

An answer to the question would be of interest in discrete dynamical 
systems.

\begin{proof}[Proof of the theorem]
We will take the positive reals as our model and will give an example of 
an orientation preserving homeomorphism.
It will be clear that this convenience is not vital.

Let $n > 1$ be specified.
Let $x_1 = 1$.
If the integer $s$ is of the form 
$$nk, nk-1, \ldots, nk-n+2\quad (k=1,2,3,\ldots)$$
then define $x_{s+1} - x_s = 1/s$ and define $x_{s+1} - x_s = 2$ for $s$ 
of the form $nk-n+1$.

For example, for $n=4$, the values of $x_{s+1} - x_s$ are
$$1, 2, 
\frac{1}{2}, \frac{1}{3}, \frac{1}{4}, 2, 
\frac{1}{6}, \frac{1}{7}, \frac{1}{8}, 2, 
\frac{1}{10}, \frac{1}{11}, \frac{1}{12}, 2, \ldots$$
Define $f(x_s) = x_{s+1}$, $f(0) = 0$.
Define $f$ to be linear on each interval $[x_s, x_{s+1}]$ and, also on 
$[0, x_1]$.

The homeomorphisms $f$, $f^2$, $f^3$, \ldots, $f^{n-1}$ are not uniformly
continuous because in each instance a null sequence of intervals maps into 
a intervals of length $2$.
However, $f^n$ is uniformly continuous since it is piecewise linear and 
the slope of each segment is less than or equal to $1$.
\end{proof}

\label{tputzend}

\chapter*{Beverly L.~Brechner: Questions on homeomorphism groups of 
chainable and homogeneous continua}
\markboth{\normalsize\textsc{\lowercase{B.~Brechner: Questions on 
homeomorphism groups of continua}}}{}
\label{tpbrechner}
\begin{myfoot}
\begin{myfooter}
Beverly L.~Brechner, 
\emph{Questions on homeomorphism groups of chainable and homogeneous 
continua},\\
Problems from Topology Proceedings, Topology Atlas, 2003, 
pp.\ 143--144. 
\end{myfooter}
\end{myfoot}

\mypreface
This article is reprinted whole:
Beverly L.~Brechner, \emph{Questions on homeomorphism groups of chainable 
and homogeneous continua}, Topology Proceedings \textbf{7}, no.~2 (1982) 
391--393.

\section*{}

The following theorem is likely to be of importance in the solution of the 
problems posed below.

\begin{theorem}[Effros]
Let $X$ be a homogeneous metric continuum.
Then for every $\epsilon > 0$, there exist $\delta > 0$ such that if 
$d(x,y) < \delta$, then there is a homeomorphism $h\colon X \to Y$ such 
that $d(h, \text{id}) < \epsilon$ and $h(x) = y$.
\end{theorem}

In \cite{MR32:4662}, we began a study of the topological structure, in 
particular dimension properties, of homeomorphism groups of various 
continua.
In particular, it was shown that the groups of homeomorphisms of 
locally-setwise-homogeneous continua are non-zero dimensional, and, 
in fact, contain the infinite product of non-zero dimensional subgroups.
Such continua include the Sierpi\'nski universal plane curve and the 
Menger universal curve.
The homeomorphism groups of those two continua are totally disconnected, 
and it is still an \emph{open question} to determine what the dimension 
is.
Examples $M_n$ are also constructed in \cite{MR32:4662}, with the property 
that $G(M_n)$ is topologically and algebraically the product of $n$ 
one-dimensional groups.
It is still \emph{unknown} what their dimension is, too.

Here we list some questions about the homeomorphism groups of the 
pseudo-arc and other homogeneous continua.
These questions were raised by the author at the University of Texas 
Summer 1980 Topology Conference, held in Austin, Texas.

Let $P$ be the pseudo-arc, and let $X$ be any homogeneous metric 
continuum.
Let $H(X)$ denote the group of all homeomorphisms of $X$ onto itself.
It is well known and easy to see that $H(P)$ contains no arcs:
for any such arc is a homotopy $\{h_t\}$ of $P$, and if $\{x\} \times I$ 
is the track of the homotopy such that $h_1(x) \neq x$, then 
$\bigcup \{h_t(x)\}_{t \in I}$ is a subcontinuum of $P$ which is a 
continuous image of an arc, and therefore locally connected.
But $P$ contains no nondegenerate locally connected continua.
Thus we raise the following.

\startproblem
\begin{myprob}
\myproblem{1}{}
Is $H(P)$ totally disconnected?
zero-dimensional?
infinite-dimensional?
\end{myprob}

\begin{myprob}
\myproblem{2}{}
Does $H(P)$ contain a pseudo-arc?
an infinite product of pseudo-arcs?

\mynote{Solution}
Wayne Lewis \cite{MR84e:54038} has just answered this question in the 
negative, by showing that $H(P)$ contains no nondegenerate subcontinua.
\end{myprob}

\begin{myprob}
\myproblem{3}{}
Is $H(P)$ connected?
If not, does it contain a nondegenerate component?
\end{myprob}

\begin{myprob}
\myproblem{4}{}
Let $G$ denote the subgroup of $H$ keeping every composant invariant.
Then $G$ is normal in $H$.
Is $G$ minimal normal? (See \cite{MR20:4607, MR42:6793, MR27:737}.)
What is the (non-identity) minimal normal subgroup?
Is $G$ generated by those homeomorphisms supported on small open sets? 
(See \cite{MR80m:54053}.)
\end{myprob}

\begin{myprob}
\myproblem{5}{}
Let $X$ be any homogeneous metric continuum.
Is $H(X)$ non-zero dimensional?
infinite dimensional?

\mynote{Remark}
It has recently been shown by Wayne Lewis \cite{MR85c:54065} that the 
pseudo-arc admits $p$-adic Cantor group actions, as well as period $n$ 
homeomorphisms for all $n$.
\end{myprob}

\nocite{
MR20:4607,
MR32:4662,
MR42:6793,
MR39:3471,
MR80m:54053,
MR85c:54065,
MR84e:54038,
MR27:737,
MR23:A3801}

\providecommand{\bysame}{\leavevmode\hbox to3em{\hrulefill}\thinspace}

\label{tpbrechnerend}

\chapter*{Some problems in applied knot theory and geometric topology}
\label{tpsumners}
\begin{myfoot}
\begin{myfooter}
D.W.~Sumners, J.L.~Bryant, R.C.~Lacher, R.F.~Williams and J.~Vieitez, 
\emph{Some problems in applied knot theory and geometric topology},\\
Problems from Topology Proceedings, Topology Atlas, 2003, 
pp.\ 145--152.
\end{myfooter}
\end{myfoot}
\setcounter{footnote}{0}

\mypreface
In volume 13 of \emph{Topology Proceedings}, D.W.~Sumners \cite{MR91a:57006} 
edited a collection of problems in applied knot theory and geometric topology.
The collection, with contributions by D.W.~Sumners, J.L.~Bryant, 
R.C.~Lacher, and R.F.~Williams is reproduced here with a few new notes.
Jos\'e Vieitez contributed a new essay describing current results on 
expansive diffeomorphisms on $3$-manifolds.

\section*{D.W.~Sumners: Some problems in applied knot theory and some 
problems in geometric topology}

Modern knot theory was born out of physics in the 19th century.
Gauss' considerations on inductance in circular wires gave rise to the
``Gauss Integral,'' a formula for the linking number of two simple closed
curves in $3$-space \cite{sumnersG}.
William Thompson (Lord Kelvin), upon seeing experiments performed by 
P.G.~Tait involving colliding smoke rings, conceived the ``vortex theory of
atoms,'' in which atoms were modelled as configurations of knotted vortex
rings in the aether \cite{sumnersTh}
In this context, a table of the elements was---you guessed it---a knot
table!
Tait set about constructing this knot table, and the rest is history
\cite{sumnersTa}!

Given the circumstances of its birth, it is not surprising that knot 
theory has, from time to time, been of use in science.
One can think of $3$-dimensional knot theory as the study of flexible 
graphs in $\mathbb{R}^3$, with emphasis on graph entanglement (knotting 
and linking).
A molecule can be represented by its molecular graph---atoms as vertices, 
covalent bonds as edges.
A large molecule does not usually maintain a fixed $3$-dimensional 
configuration.
It can assume a variety of configurations, driven from one to the other by 
a thermal motion, solvent effects, experimental manipulation, etc.
From an initial configuration for a molecule (or a set of molecules), knot 
theory can help identify all of the possible attainable configurations of 
that molecular system.
It is clear that the notion of topological equivalence of embeddings of 
graphs in $\mathbb{R}^2$ is physically unrealistic---one cannot stretch 
or shrink molecules at will.
Nevertheless, the topological definition of equivalence is, on the one 
hand, broad enough to generate a large body of mathematical knowledge, 
and, on the other hand, precise enough to place useful an computable 
limits on the physically possible motions and configuration changes of 
molecules.
For molecules which possess complicated molecular graphs, knot theory can 
also aid in the prediction and detection of various spatial isomers 
\cite{MR87m:57007}.
As evidence for the utility of knot theory (and other mathematics) in 
chemistry and molecular biology, see the excellent survey articles 
\cite{sumnersWa, sumnersWC} and the conference proceedings 
\cite{MR89c:00037,sumnersKR,sumnersL}.

Some of the problems posed below deal with configuration of random walks 
or self-avoiding (no self-intersection) random on the integer cubic 
lattice in $\mathbb{R}^3$.
The statistics of random walks on the lattice are used to model 
configurations of linear and circular macromolecules.
A macromolecule is a large molecule formed by concatenating large 
numbers of monomers---such as the synthetic polymer polyethylene and the 
biopolymer DNA.
Conversion of circular polymers from one topological state (say unknotted 
and unlinked) to another (say knotted and linked) can occur through the 
action of various agents., chemical or biological.
Given constraints (energetic, spatial or temporal), linear polymers can 
exhibit entanglement (knotting and linking).
Moreover, linear polymers can be converted to circular polymers in various 
cyclization reactions.
If one wants a random sample of the configuration space of a macromolecule 
in $\mathbb{R}^3$, one can model the spatial configuration of a 
macromolecule as a self-avoiding random walk in $\mathbb{R}^3$, where 
the vertices represent the positions of carbon atoms, and adjacent 
vertices are connected by straight line segments (all the same length), 
representing covalent bonds.
A discrete version of random walks in $\mathbb{R}^3$ is random walks in 
the integer cubic lattice.
One studies the statistical mechanics of large ensembles of these random 
walks in hopes of detecting physically observable quantities (such as 
phase transition) of the physical system being modelled.

The problems below are stated in an informal style, and addresses of 
relevant people are included when known, in hopes that the interested 
reader will contact them.

\section*{Problems proposed by J.L.~Bryant and R.C.~Lacher}

Consider random walks on a cubic lattice in $\mathbb{R}^3$ that start 
with $0 < y < n$, $n > 1$, and end when either $y=0$ or $y=n$.
An L-walk (R-walk) is a walk that starts with $y=1$ ($y= n - 1$).
(Think of an L-walk or R-walk as a walk that starts on one of the planes 
$y=0$ or $y=n$ and takes its first step into the region between the 
planes.)
An L-loop (R-loop) is an L-walk that ends with $y=0$ ($y=n$).
Assume step probabilities are all equal to $1/6$ (pure isotropy).
Given an L-walk L and an R-walk R, defined the offset linking number 
$\operatorname{olk}(L,R)$ as follows:
If each of L and R is a loop, complete it to a closed curve by joining its 
endpoints with an arbitrary path in its base plane, offset the lattice for 
R by the vector $(-1/2, -1/2, -1/2)$, and define the 
$\operatorname{olk}(L,R)$ to be the homological linking number of the 
resulting (disjoint) closed curves.
Otherwise, set $\operatorname{olk}(L,R) = 0$.
We say L links R if $\operatorname{olk}(L,R) \neq 0$.

\startproblem
\begin{myprob}
\myproblem{Problem 1}{}
Given an L-walk L and a family $\mathcal{R}$ of R-walks with density of 
starts $d$, what is the probability $P_\text{link}(n)$ that L will link a 
member of $\mathcal{R}$?
\end{myprob}

\begin{myprob}
\myproblem{Problem 2}{}
Compute $\lim_{n \to \infty} P_\text{link}(n)$.
\end{myprob}

\begin{myprob}
\myproblem{Problem 3}{}
Find the expected value $D_\text{link}(n)$ of the number of members of 
$\mathcal{R}$ that L links.
\end{myprob}

\begin{myprob}
\myproblem{Problem 4}{}
Compute $\lim_{n \to \infty} D_\text{link}(n)$.
\end{myprob}

\begin{myprob}
\myproblem{Problem 5}{}
Find the expected sum $W_1(n)$ of the absolute values of the offset 
linking number of L with the members of $\mathcal{R}$.
\end{myprob}

\begin{myprob}
\myproblem{Problem 6}{}
Compute $\lim_{n \to \infty} W_1(n)$.
\end{myprob}

\begin{myprob}
\myproblem{Problem 7}{}
Find the expected sum $W_2(n)$ of the squares of
the offset linking number of L with the members of $\mathcal{R}$.
$W_2(n)$ should be easier to deal with than $W_1(n)$.
\end{myprob}

\begin{myprob}
\myproblem{Problem 8}{}
Compute $\lim_{n \to \infty} W_2(n)$.

Given an L-loop that starts at $(0,1,0)$, define its reach to be its 
maximum $y$-value, its range to be its maximum $x$-values, and its 
breadth $b = \text{range}/\text{reach}$.
By analogy, define the breadth of any loop.
\end{myprob}

\begin{myprob}
\myproblem{Problem 9}{}
Compute the expected value of $b$ as a function of $n$ and its 
asymptotics.
Simulation statistics seem to indicate that $b = 1.19$. 
See \cite{MR91f:92007}.
\end{myprob}

Represent a loop by an isosceles triangle parallel to the $y$ axis having 
its base on the base plane for the loop.
Its ``breadth'' $b = \text{altitude}/{2\cdot\text{base}}$.
Analogs of $D_\text{link}(n)$ and $P_\text{link}(n)$ for these simplified
loops are
\begin{align*}
D(n)&
\textstyle
= 2 b^2 d \sum_{i=1}^{n-1} d_i \sum_{j=n-1}^{n-1} 
[1 - 2 b^2 d ({i+j+1}/{2-n})^2 d_j],\ \text{and}\\
P(n)&
\textstyle
= 1 - 1/n - \sum_{i=1}^{n-1} d_i \prod_{j=n-1}^{n-1}
[1 - 2 b^2 d ({i+j+1}/{2-n}) d_j]
\end{align*}
Asymptotics for $D(n)$ are given in \cite{MR91f:92007}.

\startproblem
\begin{myprob}
\myproblem{Problem 10}{}
Compute $\lim_{n \to \infty} P(n)$.
We conjecture 
$n \cdot P(n) \sim O(\log(n))$.
\end{myprob}

\begin{myprob}
\myproblem{Problem 11}{}
Show that
$\lim_{n \to \infty} P(n) = \lim_{n \to \infty} P_\text{link}(n)$
and that
$\lim_{n \to \infty} D(n) = \lim_{n \to \infty} D_\text{link}(n)$.
\end{myprob}

\section*{Problems proposed by D.W.~Sumners}

There exist naturally occurring enzymes (topoisomerases and recombinases) 
which, in order to mediate the vital life processes of replication, 
transcription and recombination, manipulate cellular DNA in topologically 
interesting and nontrivial ways \cite{sumnersWC, MR88c:57012}.
These enzyme actions include promoting writhing (coiling up) of DNA 
molecules, passing one strand of DNA through another via an enzyme-bridged 
break in one of the strands, and breaking a pair of strands and 
recombining to different ends.
If one regards DNA as very thin string, these enzyme activities are the 
stuff of which recent combinatorial knot theory is made!
Moreover, relatively new experimental techniques (rec A enhanced 
microscopy) \cite{sumnersKS} make possible the unambiguous resolution of 
the DNA knots and links produced by reacting circular DNA with high 
concentrations of a purified enzyme \emph{in vitro} (in the laboratory).
The experimental protocol is to manufacture (by cloning techniques) 
artificial circular DNA substrate on which a particular enzyme will act.
As experimental control variables, one has the knot type(s) of the 
substrate, and the amount of writhing (supercoiling) of the substrate 
molecules.
The product of an enzyme reaction is an enzyme-specific family of DNA 
knots and links.
The reaction products are fractionated by gel electrophoresis, in which 
the molecules migrate through a resistive medium (the gel) under the 
forcing of an electric field (electrophoresis).
Molecules which are ``alike'' group together and travel together in a band 
through the gel.
Gel electrophoresis can be used to discriminate between molecules on the 
basis of molecular weight.
Given (as in the case here) that all molecules are the same molecular 
weight, it then discriminates between molecules on the basis of average 
3-dimensional ``shape''.
Following electrophoresis, the molecules are fattened with a protein (rec 
A) coating, to enhance resolution of crossovers in an electron micrograph 
of the molecule.
In this manner, the knot (link) type of the various reaction products is 
an observable.
This new observational power makes possible the building of knot-theoretic 
models \cite{sumnersWC,sumnersWMC,MR92f:92024} for enzyme action, in which 
one wishes to extract information about enzyme mechanism from the 
DNA knots and links produced by an enzyme reaction.

\startproblem
\begin{myprob}
\myproblem{Problem 1}{}
Build new models for enzyme action.

The models now existing involve signed crossover number \cite{sumnersWC}, 
polynomial invariants \cite{sumnersWMC}, and $2$-string tangles 
\cite{MR92f:92024}.
The situation is basically this: as input to a black box (the enzyme), one 
has a family of DNA circles (of known knot type and degree of 
supercoiling). 
The output of the black box is another family of DNA knots and links.
The problem: What happened inside the box?
\end{myprob}

\begin{myprob}
\myproblem{Problem 2}{}
Explain gel electrophoresis experimental results.

Gel electrophoresis is a race for molecules---they all start together, 
and the total distance travelled by a molecule when the electric field is 
turned off is determined by its gel mobility.
At the finish of a gel run, the molecules are grouped in bands, the 
slowest band nearest the starting position, the fastest band farthest 
away. 
When relaxed (no supercoils) DNA circles (all the same molecular weight) 
run under certain gel conditions, the knotted DNA circles travel according 
to their crossover number \cite{sumnersDS}!
What is it about crossover number (an artifact of $2$-dimensional knot 
projections) that determines how fast a flexible knot moves through a 
restive medium?
The theory of gel mobility of molecules (linear or circular) is rather 
difficult to work out.
See \cite{sumnersLZ} for some results on the gel mobility of unknotted 
circular molecules under pulsed field electrophoresis.
\end{myprob}

\begin{myprob}
\myproblem{Problem 3}{}
What are the properties of a random knot (of fixed length)?

Chemists have long been interested in the synthesis of molecules with 
exotic geometry in particular, the synthesis of knotted and linked 
molecules \cite{sumnersWa}.
One can imagine such a synthesis by means of a cyclization reaction 
(random closing) of linear chain molecules \cite{sumnersFW}.
Let $N$ represent the number of repeating units of the substance, or the 
equivalent statistical length of the substance.
For example, the equivalent statistical length for polyethylene is about 
3.5 monomers, and for duplex DNA, about 5000 base pairs.
A randomly closed chain of length $N$ is a random piecewise linear 
embedding of $\mathbb{S}^1$, with all the 1-simplexes of the same length.
See \cite{sumnersR1,sumnersR2} for a discussion of the topology of the 
configuration space of such PL embeddings.
In order to make predictions about the yield of such a cyclization 
reaction, one needs answers to the following mathematical questions 
\cite{sumnersS2}:
\end{myprob}

\begin{myprob}
\myproblem{Problem 3A}{}
For random simple closed curves of length $N$ (as above), what is the 
distribution of knot types, as a function of $N$?
\end{myprob}

\begin{myprob}
\myproblem{Problem 3B}{}
What is the probability of knotting, as a function of $N$?
One can show that, for simple closed curves of length $N$ inscribed on the 
cubical lattice in $\mathbb{R}^3$, the knot probability goes to one 
exponentially rapidly with $N$ \cite{MR89i:82060}.
\end{myprob}

\section*{Problems proposed by R.F.~Williams}

\subsection*{Expansive vs.\ pseudo-Anosov}

The references here are two preprints.%
\footnote{These have since been published.}
In \cite{MR91b:58184} and \cite{MR92i:58139}, the authors 
[resp.\ K.~Hiraide and J.~Lewowicz] independently prove that the concepts 
\emph{expansive} and \emph{pseudo-Anosov} coincide for surfaces.
\begin{myitemize}
\item[A.] 
What is the situation for $3$-manifolds?
\item[B.] 
Find a good example of a $3$-manifold (such as $\mathbb{S}^3$) which 
does not support an Anosov diffeomorphism.
\item[C.] 
Prove some of the beginning lemmas of Lewowicz-Hiraide for $3$-manifolds.
\end{myitemize}

\subsection*{Dynamical systems}

The two topics of zeta functions in dynamical systems and Alexander 
polynomials in knot theory are closely related; see \cite{MR39:3497}.
For flows in $\mathbb{S}^3$, periodic orbits are knots; thus there should 
be a combination such as a $2$ variable polynomial, combining knot theory 
(e.g., the degree of the Alexander polynomial) and dynamical systems (the 
length of the orbit).
See \cite{MR86a:58084}.

Branched surfaces can support Anosov endomorphisms.
However, all that are known are shift equivalent to linear maps on the 
$2$-torus, such as that induced by the $2 \times 2$ matrix 
$\left(
\begin{smallmatrix}
3&1\\
1&1
\end{smallmatrix}
\right)$.

\startproblem
\begin{myprob}
\myproblem{Conjecture}{}%
\footnote{L.~Wen \cite{MR94b:58076} has since proven a special case of
the conjecture:
If $F\colon K \to K$ is an Anosov endomorphism of branched surface, in 
which the branch set is the union of a finite collection of simple closed
curves, then $F$ is shift equivalent to a linear endomorphism of the
$2$-torus.}
Given an Anosov endomorphism $g\colon K \to K$, there is a linear map 
$f\colon T \to T$, $T$ the $2$-torus, such that $f$ is shift equivalent 
to $g$.
\end{myprob}

$f\colon X \to X$ and $g\colon Y \to Y$ are \emph{shift equivalent}
provided that there exist maps $r\colon X \to Y$ an $s\colon Y \to X$ and 
an integer $m$ such that $rf = gr$, $sg = fs$, $sr = f^m$, and 
$rs = g^m$.

$g\colon K \to K$ is \emph{Anosov}, provided there is a sub-bundle $E$ 
of the tangent bundle $TK$, such that $dg$ leaves $E$ invariant and 
contacts vectors, and such that the map induced on $TK/E$ by $dg$ expands 
vectors.

Hassler Whitney gives an example which is dear to the heart of all 
continuum theorists that know it---both of us!
It is a carefully constructed arc $A$ in the plane and smooth function 
$F\colon A \to \mathbb{R}$ with $\text{grad} f = 0$ (both partials are $0$), 
yet $f$ has different values at $A$'s endpoints.
Contact Alec Norton, Boston University for his preprints and ideas on this 
subject.
(Don't be afraid of smooth functions on manifolds.
They have beautiful pathology and are crying out for continuum theorists 
to look at them.
And they are really and truly easy to get the hang of.)

\section*[J. Vieitez: Expansive diffeomorphisms on 
$3$-manifolds]{J. Vieitez: Expansive diffeomorphisms on 
$\mathbf{3}$-manifolds}

\mypreface
R.F.~Williams's problems on expansive diffeomorphisms on $3$-mani\-folds
have been answered.
Here is a survey by Jos\'e Vieitez of the results.

\subsection*{$\mathbf{3}$-manifolds}

The answer for $3$-manifolds is no, expansive is not equivalent with 
pseudo-Anosov.
First of all we should define pseudo-Anosov for manifolds which are not 
surfaces.
A possible (rough) definition is the following:
We say that $f\colon M\to M$ is a \emph{pseudo-Anosov} homeomorphism if 
there exist two foliations with a finite set of singularities
$\mathcal{F}^s$ (stable) and $\mathcal{F}^u$ (unstable), invariant by 
$f$, and with the same finite set of singularities, such that 
$\mathcal{F}^s$ and $\mathcal{F}^u$ are transverse except at the 
singularities.
There exists $0 < \lambda<1$ such that $f$ contracts in $\mathcal{F}^s$
by a factor less or equal than $\lambda$ and $f^{-1}$ contracts in
$\mathcal{F}^u$ by a factor less or equal than $\lambda$.
At the singularities $\mathcal{F}^s$ is not locally Euclidean (but we 
should say something more. In the $2$-dimensional case we use the notion
of separatrices).
There is also a notion of measures transverse to the foliations which 
implies density of periodic points, almost of them topologically 
hyperbolic.
In particular the non-wandering set of $f$ is all of $M$.
Assuming this definition of pseudo-Anosov we should observe that only in 
the $2$-dimensional case we have a symmetric behaviour of $\mathcal{F}^s$
and $\mathcal{F}^u$.
In the $3$-dimensional case one of the foliations should separate while 
the other cannot do it.
Otherwise they should intersect in a nontrivial line which has to 
contract to the future (being in $\mathcal{F}^s$) and to contract to the 
past (being in $\mathcal{F}^u$).
But both things cannot occur at the same time.
Moreover, in the $n$-dimensional case, $n>3$, we cannot expect any of them 
to separate.
Returning to the $3$-dimensional case, using that topologically hyperbolic 
periodic points are dense in $M$, I have proved that $f$ has to be 
conjugated to Anosov and $M$ has to be $\mathbb{T}^3$, the $3$-torus, a 
rather restrictive result \cite{MR97i:58132}.
Moreover, assuming only that $\Omega(f)=M$ and that $f$ is a 
diffeomorphism and not merely a homeomorphism, I have proved the same 
result (that $f$ has to be conjugated to an Anosov diffeomorphism and
$M=\mathbb{T}^3$) \cite{MR2003b:37049}.

On the other hand, R.~Ma\~n\'e, has proved in his paper 
\emph{Expansive diffeomorphisms} \cite{MR58:31263}, that the 
$C^1$-interior of expansive diffeomorphisms are the so called quasi-Anosov 
diffeomorphisms which should not be confused with the pseudo-Anosov 
homeomorphisms.
Quasi-Anosov diffeomorphism can be defined as diffeomorphisms $f$ such 
that $Df$, the tangent map, expands to infinity the norm of any vector 
different from zero either to the past or to the future (or in both 
directions).
But J.~Franks and C.~Robinson have given in their paper
\emph{A quasi-Anosov diffeomorphism that is not Anosov} 
\cite{MR54:11399} an example of a quasi-Anosov (hence an expansive 
diffeomorphism) such that its non-wandering set is the union of a 
codimension one repeller $R$ with a codimension one attractor $A$.
This example is defined in the amalgamated sum of two $3$-tori, a 
$3$-dimensional manifold $M$.
In this case, from the topological point of view, most of $M$ is not 
foliated by transverse stable and unstable foliations.
They exist but in almost all points of $M$ (an open and dense subset) 
both are $1$-dimensional and therefore cannot be transverse.
In $A$ the unstable foliation is $2$-dimensional and the stable foliation 
is $1$-dimensional and in $R$ the situation reverses.
In particular this implies that $f$ cannot be Anosov.
And with the definition I have given of pseudo-Anosov, $f$ cannot be
pseudo-Anosov (in particular periodic points are not dense in $M$).
With J.~Rodriguez~Hertz and R.~Ures in the paper
\emph{On manifolds supporting quasi-Anosov diffeomorphisms}
\cite{MR2002m:37040} we have studied the case of quasi-Anosov diffeomorphism 
in the $3$-dimensional case proving that it is either Anosov or admits at 
least a codimension one hyperbolic repeller and a codimension one 
hyperbolic attractor. 
In this sense, the example of Franks and Robinson is minimal.

\subsection*{Examples of $\mathbf{3}$-manifolds that do not support an 
Anosov diffeomorphism}

This follows from the result of S.~Newhouse \cite{MR43:2741} which states 
that if $f$ is Anosov defined on $M$ and one of the foliations (the stable 
or the unstable) is of codimension one then $M$ is the $3$-torus (in fact 
Newhouse proves that the non-wandering has to be the whole manifold $M$ 
and a result of J.~Franks \cite{MR42:6871} implies that $M$ has to be a 
torus).
Thus Anosov diffeomorphisms live only in tori in the $3$-dimensional 
case.
The example of Franks and Robinson cited above shows that there are other
examples of expansive diffeomorphisms in $3$-manifolds.
Moreover, this example (by the result of Ma\~n\'e cited above) is 
Axiom A with the no-cycle condition and therefore $\Omega$-stable.
Being quasi-Anosov it is in the $C^1$-interior of the expansive
diffeomorphisms and so if we perturb it we obtain another expansive
diffeomorphism \emph{not} necessarily conjugate to $f$ (otherwise $f$ 
should verify the strong transversality condition).

\subsection*{Lewowicz-Hiraide lemmas for $\mathbf{3}$-manifolds}

Most of these lemmas are valid in the $3$-dimensional case.
Moreover, I use these results.
That is for instance the case of the following results,

There exist Lyapunov functions for $f$.

In the boundary of a given ball of radius $r>0$ and center $x$ the 
distance between the local stable (unstable) manifold of $x$ and the 
local unstable (resp.\ stable) manifold of points $y$ near $x$ is
bounded away from zero \emph{uniformly} in $x\in M$.
We should better say \emph{local stable and unstable sets} instead of 
\emph{local stable and unstable manifolds} because one of the main 
difficulties in the general case is to have a good topological picture of
these sets.

Given $0<\delta<\epsilon$ there exists $r>0$ such that if 
a point $y$ is 
in the $\epsilon$-stable (unstable) set of $x$ and 
$\operatorname{dist}(x,y)<r$ then it 
is in the $\delta$-stable (resp.\ unstable) set of $x$.

The connected components of the local stable and unstable manifolds of a 
point $x$ that contain $x$ contain nontrivial continua which reach the 
boundary of a ball of radius $r>0$ with $r$ independent of $x$.

There are not Lyapunov stable points for $f$ an expansive homeomorphism 
defined in a compact manifold or even in any locally connected compact
metric continuum.

Let $\mathcal{H}(M)$ be the the space of homeomorphisms defined in a 
compact metric space $M$.
If we perturb in $\mathcal{H}(M)$ in the $C^0$-topology an expansive 
homeomorphism $f$ with constant of expansivity $\alpha>0$ we cannot assure 
that the resulting homeomorphism $g$ is expansive. But a sort of dichotomy
appears: given $\delta>0$ there exists a neighbourhood $\mathcal{N}(f)$ 
such that if $g\in \mathcal{N}(f)$ then given two points $x,y\in M$ either 
there exists $n_0\in \mathbb{Z}$ such 
$\operatorname{dist}(g^{n_0}(x),g^{n_0}(y)) \ge \alpha$ or 
$\operatorname{dist}(g^n(x),g^n(y)) \le \delta$ for all $n\in \mathbb{Z}$.
Taking $\delta<\alpha/2$ we may define an equivalence relation $\sim$ in 
$M$ such that $g$ pass to the quotient space $M/{\sim}$ as an expansive 
homeomorphism $\hat g$. 
It results that $f$ is topologically stable iff $\hat g$ is conjugated to 
$f$ for all possible $g\in \mathcal{N}(f)$.

The two last mentioned results have been proved by Lewowicz in 
\cite{MR85m:58140}.
The quotient spaces $M/{\sim}$ that are obtained have a rich topological
structure that is not well understood and is object of recent research. 

\providecommand{\bysame}{\leavevmode\hbox to3em{\hrulefill}\thinspace}

\label{tpsumnersend}

\chapter*{Problems from Chattanooga, 1996}
\label{tpchattanooga}
\begin{myfoot}
\begin{myfooter}
W.W.~Comfort, F.D.~Tall, D.J.~Lutzer, C.~Pan, G.~Gruenhage, S.~Purisch 
and P.J.~Nyikos, 
\emph{Problems from Chattanooga, 1996},\\
Problems from Topology Proceedings, Topology Atlas, 2003, 
pp.\ 153--154.
\end{myfooter} 
\end{myfoot} 

\mypreface
These problems appeared in volume 20 (1996) of \emph{Topology Proceedings}.
At the AMS Regional Meeting in Chattanooga, Tennessee, October 11--12,
1996, during the Special Session in Set-theoretic Topology, there was a
problem session at which the following problems were posed.
Some of the notes are new.

\section*{}

\startproblem
\begin{myprob}
\myproblem{1}{W.W.~Comfort, attributed to N.~Noble}
Can there be an uncountable family of noncompact Tychonoff spaces whose
product is a $k$-space?

\mynote{Notes} 
N.~Noble showed in his Ph.D.\ thesis that a co-countable subfamily 
must have pseudocompact product, hence all but countably many factors
must be countably compact.
See also \cite{MR36:5883, MR44:979, MR44:4707}.
\end{myprob}

\begin{myprob}
\myproblem{2}{F.D.~Tall, attributed to W.~Fleissner}
Is there a normal $k$-space which is not collectionwise normal?

\mynote{Notes}
Peg Daniels \cite{MR92g:54031} has shown the consistency of every normal
$k'$-space being collectionwise normal, assuming large cardinal axioms.
\end{myprob}

\begin{myprob}
\myproblem{3}{D.J.~Lutzer}
Can every perfectly normal suborderable space be embedded in a perfectly 
normal LOTS?
\end{myprob}

\begin{myprob}
\myproblem{4}{D.J.~Lutzer}
Let $X$ be a suborderable space with a $\sigma$-discrete dense subspace.
Can $X$ be embedded in a perfectly normal LOTS?
a perfectly normal LOTS with a $\sigma$-discrete dense subspace?
\end{myprob}

\begin{myprob}
\myproblem{5}{Chunliang Pan}
Dowker showed that a space $X$ is normal and countably paracompact if, and 
only if, it is possible to choose, for each USC real-valued function $g$ 
and each LSC real-valued function $h$ such that $g(x) < h(x)$ for all $x$, 
a continuous function $\Phi(g, h)$ such that $g < \Phi(g, h) < h$ 
everywhere.
Can we characterize internally those spaces $X$ for which this choice can
be done monotonically, i.e., if $g < g'$ and $h < h'$ then 
$\Phi(g, h) < \Phi(g', h')$ everywhere?

\mynote{Notes}
If $\leq$ is substituted for $<$ everywhere, then we get a condition
equivalent to perfect normality.
\end{myprob}

\begin{myprob}
\myproblem{6}{G.~Gruenhage, attributed to R.~McCoy}
Find a property $P$ such that $X$ has $P$ iff $C(X)$ with the compact-open
topology is a Baire space.
Does the Moving-Off Property (MOP) provide such a characterization?

\mynote{Notes}
If $C(X)$ is Baire in the compact-open topology, then $X$ has the MOP, 
which is the property that every collection $\mathcal{L}$ of compact sets 
that moves off the compact sets contains an infinite subcollection with a
discrete open expansion.
A family $\mathcal{L}$ is said to \emph{move off the compact sets} if 
for each compact subset $K$ of $X$ there is a member of $\mathcal{L}$ that 
is disjoint from it.
See \cite{MR98f:54018}.
\end{myprob}

\begin{myprob}
\myproblem{7}{S.~Purisch}
Can we characterize the compact spaces of diversity $2$, i.e., those
compact spaces with exactly two open subspaces up to homeomorphism?

\mynote{Notes}
See the papers by J.~Mioduszewski \cite{MR80a:54040} and
by J.~Norden, S.~Purisch and M.~Rajagopalan \cite{MR97d:54043}.
\end{myprob}

\begin{myprob}
\myproblem{8}{P.J.~Nyikos, attributed to A.~Dow and K.P.~Hart}
If a continuum is the continuous image of the Stone-\v Cech remainder
$\omega^*$, is it the continuous image of the Stone-\v Cech remainder
$\mathbb{H}^*$ of the closed half-line?

\mynote{Notes}
A.~Dow and K.P.~Hart \cite{MR2001g:54037.1} proved that every continuum 
of weight $\aleph_1$ is a continuous image of $\mathbb{H}^*$.
\end{myprob}

\providecommand{\bysame}{\leavevmode\hbox to3em{\hrulefill}\thinspace}

\label{tpchattanoogaend}

\chapter*{Problems from Oxford, 2000}
\label{tpoxford}
\begin{myfoot}
\begin{myfooter}
A.V.~Arhangel$'$\kern-.1667em ski\u\i, S.~Antonyan, K.P.~Hart, L.~Ludwig, M.~Matveev, 
J.T.~Moore, P.J.~Nyikos, S.A.~Peregudov, R.~Pol, J.T.~Rogers, M.E.~Rudin, 
and K.~Shankar, 
\emph{Problems from Oxford,~2000},\\
Problems from Topology Proceedings, Topology Atlas, 2003, 
pp.\ 155--163.
\end{myfooter}
\end{myfoot}
\setcounter{footnote}{0}

\mypreface
These problems appeared in volume 25 (2000) of \emph{Topology Proceedings}.
They were collected during the problem session of the Summer Topology 
Conference at Miami University (Oxford, OH, August 2000).
In this version, a few references have been updated.
L.~Ludwig, M.~Matveev and J.T.~Moore have made minor modifications to 
their original contributions.

\section*{Alexander Arhangel$'$\kern-.1667em ski\u\i}
 
\begin{myprob}
\myproblem{Locally Compact Linearly Lindel\"of Spaces}{}
Let $X$ be a compact Hausdorff space and $a\in X$ such that for every 
uncountable subset $A$ of $X$ of regular cardinality there exists an open 
neighbourhood $U$ of $a$ such that the cardinality of $X\setminus A$ is 
the same as the cardinality of $A$. 
Is then $X$ first countable at $a$? 

\mynote{Comments} 
The above question is obviously equivalent to the following one: 
Is every linearly Lindel\"of locally compact Hausdorff space Lindel\"of?%
\footnote{\mypreface 
K.~Kunen \cite{MR2003d:54040} constructed a Hausdorff, locally compact,
linearly Lindel\"of space which is not Lindel\"of.
P.~Nyikos \cite{nyikosll} proved that it is consistent (relative to the
existence of large cardinals) that every locally compact linearly
Lindel\"of normal space is Lindel\"of.}
Certain results in the direction of the above problem were obtained in 
\cite{MR99c:54034} (where the question was formulated for the first 
time). 
Recall that a space $Y$ is said to be \emph{linearly Lindel\"of}
if every uncountable set of regular cardinality has a point of complete 
accumulation in $Y$. 
There are linearly Lindel\"of Tychonoff spaces that are not Lindel\"of. 
However, it is still unknown if there exists a normal linearly Lindel\"of 
space that is not Lindel\"of.
\end{myprob}

\begin{myprob}
\myproblem{First Countable Linearly Lindel\"of Spaces}{}
Let $X$ be a first countable linearly Lindel\"of Tychonoff (regular) 
space. 
Is then $X$ Lindel\"of?

\mynote{Comments}
The answer to the above question is yes under \myaxiom{CH} (and even 
under some weaker assumptions). 
This follows from the following result of Arhangel$'$\kern-.1667em ski\u{\i} 
and Buzyakova, proved in \myaxiom{ZFC} (see \cite{MR99j:54006}): 
the cardinality of every first countable linearly Lindel\"of Tychonoff 
space does not exceed $\mathfrak{c}$. 
One can find other results, related to the problem, in 
\cite{MR99j:54006}, where the question was formulated for the first time.
\end{myprob}

\begin{myprob}
\myproblem{Discretely Lindel\"of Spaces}{}
Is every discretely Lindel\"of Tychonoff (regular) space Lindel\"of?

\mynote{Comments}
A space $X$ is called \emph{discretely Lindel\"of} if for every discrete 
subspace $A$ of $X$ the closure of $A$ is Lindel\"of.
Discretely Lindel\"of spaces were called \emph{strongly discretely 
Lindel\"of} in \cite{MR99j:54006}. 
Every discretely Lindel\"of regular space is linearly Lindel\"of 
\cite{MR99j:54006}. 
However, it is not even known if every discretely Lindel\"of locally 
compact Hausdorff space is Lindel\"of. 
I believe, the answer to the last question should be positive, at 
least, consistently. 
\end{myprob}

\section*{Sergey Antonyan}

\begin{myprob}
\myproblem{Problem 1}{}
Let $X$ be a paracompact (metrizable, if necessary) space, $Y$ be a 
completely regular, Hausdorff space and $f\colon X\to Y$ be a continuous, 
open, 
surjective map with connected, second countable fibers. 
Furthermore, assume that $Y$ has an open cover $\{U_\alpha\}$ 
satisfying the following conditions: 
\begin{myenumerate}
\item 
each $f^{-1}(U_\alpha)$ is dense in $X$, 
\item
there exists a subset $S_\alpha\subset f^{-1}(U_\alpha)$, closed in 
$f^{-1}(U_\alpha)$, such that $f(S_\alpha)=U_\alpha$, and the 
restriction $f|_{S_\alpha}$ is a perfect and open map.
\end{myenumerate}
Is then $Y$ paracompact (normal)? 
What if $f|_{S_\alpha}$ is a homeomorphism ?

\mynote{Comments}
Some important problems in the theory of topological transformation 
groups can be reduced to this purely general-topological problem. 
Namely, let $G$ be a separable Lie group and $Z$ be a proper (in the 
sense of R.~Palais \cite{MR23:A3802}) $G$-space. 
By making use the results of Palais \cite{MR23:A3802}, it can be shown 
that if $G$ is connected then the orbit map $f\colon Z\to Z/G$ with $Z$ 
paracompact (metrizable), satisfies to the conditions of Problem~1 with 
$X=Z$ and $Y=Z/G$.

A positive answer to Problem~1 will provide a solution of the 
H\'ajek-Abels conjecture on paracompactness of the orbit space of a 
paracompact proper $G$-space (see \cite{MR42:6806} and 
\cite{MR51:11460}). 
As it was shown by H.~Abels, this conjecture will imply, for instance, 
the parallelizability of dispersive dynamical systems on arbitrary 
paracompact phase spaces (see \cite{MR42:6806} and \cite{MR51:11460}), 
a generalization of a classical Antosiewicz-Dugundji-Nemytski theorem.
 
On the other hand, as it is shown in the \cite{MR2003m:22025}, the 
paracompactness of the orbit space will imply the existence of a 
consistent $G$-invariant metric on each metrizable proper $G$-space 
$X$ (for $X$ second countable the result was established first by Palais 
\cite{MR23:A3802}). 
This fact can have important applications in equivariant theory of retracts.

Finally, we recall that a $G$-space $X$ is called \emph{proper} (see 
R.~Palais \cite{MR23:A3802}) if:
\begin{myenumerate}
\item
$G$ is a locally compact Hausdorff topological group,
\item
$X$ is completely regular Hausdorff space and
\item
every point of $X$ has a neighborhood $V$ such that for every point of 
$X$ there is a neighborhood $U$ with the property that the set 
$\langle U,V\rangle = \{g\in G : gU\cap V\not= \emptyset\}$ 
has compact closure in $G$. 
\end{myenumerate}
\end{myprob}

\section*{K.P.~Hart}

\begin{myprob}
\myproblem{Problem 1}{}
This problem is due to Alan Dow: 
Is it consistent that all extremally disconnected continuous images
of $\omega^*$ are separable?

\mynote{Comments}
Indeed, most of the \myaxiom{ZFC} results on continuous images of 
$\omega^*$ are quite general: \emph{all} separable compact spaces, 
\emph{all} compact spaces of weight $\aleph_1$ (or less) and 
\emph{all} perfectly normal compact spaces are continuous images of 
$\omega^*$.
What distinguishes the separable spaces from the rest is that for these 
spaces the proof is nearly trivial: enumerate a dense subset of $X$ with
infinite repetitions and take the \v{C}ech-Stone extension of this 
enumeration;
the repetitions ensure that $\omega^*$ gets mapped onto $X$.
The other types of spaces are, in some sense, small so that a recursive 
construction of an onto map is possible; but note that none is extremally
disconnected.

It seems likely that \myaxiom{OCA} or \myaxiom{PFA} will settle this 
problem positively because both axioms tend to dictate that many maps 
must be nearly trivial, 
in the sense that the map has a lifting to all of $\beta\omega$ and so 
the map on $\omega^*$ decides where $\omega$ must go.

In our case one would expect the following to hold under \myaxiom{OCA} or 
\myaxiom{PFA}:
if $f\colon \omega^* \to 2^\mathfrak{c}$ is continuous and 
$X=f[\omega^*]$ is extremally disconnected then any extension of $f$ to 
all of $\beta\omega$ will have to map all but finitely many elements of 
$\omega$ into $X$.
\end{myprob}

\begin{myprob}
\myproblem{Problem 2}{}
This problem is due to Eric van~Douwen:
Give a topological characterization of $\mathbb{H}^*$ under \myaxiom{CH}.
Here $\mathbb{H}=[0,\infty)$ the half-line.

\mynote{Comments}
This topological characterization should be in the spirit of 
Parovi\-\v{c}enko's characterization of $\omega^*$.
Remember that $\omega^*$ is, under \myaxiom{CH}, the only compact space of 
weight $\mathfrak{c}$ that is zero-dimensional and without isolated 
points, and which is also an $F$-space in which nonempty $G_\delta$-sets 
have nonempty interior.

The characterization of $\mathbb{H}^*$ should replace `zero-dimensional 
and without isolated points' by some, preferably finite, list of 
topological properties.
\end{myprob}

\begin{myprob}
\myproblem{Problem 3}{}
This problem is due to Stevo Todor\v{c}evi\'c:
Does \myaxiom{OCA} imply that there is no Borel lifting for the measure 
algebra?

\mynote{Comments}
The measure algebra is defined as $M=\operatorname{Bor}/\mathcal{N}$, 
where $\operatorname{Bor}$ denotes the $\sigma$-algebra of Borel sets 
and $\mathcal{N}$ is the ideal of measure-zero sets.
A lifting is a homomorphism $l\colon M \to \operatorname{Bor}$ such that 
$q\circ l$ is the identity on $M$, where 
$q\colon \operatorname{Bor} \to M$ is the quotient homomorphism.
\myaxiom{CH} implies such a lifting exists and Shelah has shown, by his 
oracle-cc method, that it is consistent that no lifting exists.

A recent metatheorem of Ilijas Farah states that if there is no reason 
(in \myaxiom{ZFC}) for two quotients of $\mathcal{P}(\omega)$ to be 
isomorphic then \myaxiom{OCA} + \myaxiom{MA} implies that they are not 
isomorphic.
A positive answer to the present question, possibly using \myaxiom{MA}, 
would reinforce the idea that \myaxiom{OCA} and \myaxiom{MA} together 
generate a quite complete theory of the reals---with the right theorems.
\end{myprob}

\section*{Lew Ludwig}

A space $X$ is called $\beta$-\emph{normal} if for any two disjoint 
closed subsets $A$ and $B$ of $X$, there exists open subsets $U$ and $V$ 
of $X$ such that $A\cap U$ is dense in $A$, $B\cap V$ is dense in $B$, and 
the closure of $U$ in $X$ and the closure of $V$ in $X$ have empty 
intersection \cite{MR2002h:54019}.

In \cite{MR1875595}, we demonstrated the existence of a $\beta$-normal 
non-normal space assuming the existence of a normal, right-separated in 
type-$\omega_1$, S-space.

\startproblem
\begin{myprob}
\myproblem{Question 1}{}
Does there exist in \myaxiom{ZFC} a $\beta$-normal, non-normal space?

\mynote{Solution}
E.~Murtinov\'a \cite{MR2003c:54045} constructed an example.
\end{myprob}

\begin{myprob}
\myproblem{Question 2}{}
Does there exist a space $X$ in \myaxiom{ZFC} which is normal, right 
separated of type $\kappa$ with $hd(X)< \kappa$?
\end{myprob}

\begin{myprob}
\myproblem{Question 3}{}
Does there exist a space $Y$ in \myaxiom{ZFC} which is not normal, has 
scattered height $2$, and the set of nonisolated points is of 
cardinality $\lambda$, such that any two disjoint closed sets, one of 
which has size less than $\lambda$, can be separated?
\end{myprob}

\begin{myprob}
\myproblem{Question 4}{}
Does there exist an uncountable $\lambda$ that gives a positive answer 
to Questions 2 and 3 simultaneously? 
\end{myprob}

\begin{myprob}
\myproblem{Question 5}{}
Does there exist a hereditarily normal ($\alpha$-normal, $\beta$-normal), 
extremally disconnected Dowker space?
\end{myprob}

\begin{myprob}
\myproblem{Question 6}{A.V.~Arhangel$'$\kern-.1667em ski{\u\i}}
If every power $X^{\kappa}$ of a $T_1$ topological space is $\alpha$- 
($\beta$)-normal, is $X$ compact?
\end{myprob}

\section*{Mikhail Matveev} 

\begin{myprob}
\myproblem{Problem 1}{}
Is there a Tychonoff space without a minimal (Tychonoff) pseudocompact 
extension? 

A space $pX$ is called a \emph{pseudocompact extension} of a space $X$ if 
$pX$ is pseudocompact and contains $X$ as a dense subspace; $pX$ is a 
\emph{minimal pseudocompact extension} of $X$ if no proper subspace of 
$pX$ is a pseudocompact extension of $X$.
\end{myprob}

\begin{myprob}
\myproblem{Problem 2}{}
Is every topological vector space B-homogeneous?

A space is called \emph{basically homogeneous} (B-homogeneous for 
short) if it has a base every element of which can be mapped onto every 
other by an autohomeomorphism of the entire space.

\mynote{Comments}
Stanislav Shkarin has given a partial positive answer for Problem~2 in 
\cite{MR2001g:54042}: yes for locally convex TVS.
\end{myprob}

\begin{myprob}
\myproblem{Problem 3}{}
Is every Hausdorff monotonically compact space metrizable?

A space $X$ is \emph{monotonically compact} (\emph{monotonically 
Lindel\"of}) if there is a mapping that assigns to every open cover 
$\mathcal{U}$ of $X$ a finite (resp.\ countable) open refinement 
$R(\mathcal{U})$ so that $R(\mathcal{U}_1)$ refines $R(\mathcal{U}_2)$
as soon as $\mathcal{U}_1$ refines $\mathcal{U}_2$.

\mynote{Comments}
The author and Jerry Vaughan have shown (unpublished) that some well-known 
examples of nonmetrizable compacta are not monotonically compact. 
\end{myprob}

\begin{myprob}
\myproblem{Problem 4}{}
Is every countable monotonically Lindel\"of space metrizable? 
\end{myprob}

\begin{myprob}
\myproblem{Problem 5}{}
Is there a Hausdorff (regular, Tychonoff, normal) inversely compact space 
which is not compact?

Let $F$ be a family of subsets of a set $X$. 
A \emph{partial inversement} of $F$ is a family $\{p(A): A\in F\}$ such 
that for every $A\in F$, $p(A)$ is either $A$ or $X\setminus A$.
A space $X$ is called \emph{inversely compact} (\emph{inversely 
Lindel\"of}) if every open cover of $X$ has a partial inversement which 
contains a finite (resp.\ countable) subcover of $X$.
In other words, a space is inversely compact if every independent family 
of closed sets has nonempty intersection.
\end{myprob}

\begin{myprob}
\myproblem{Problem 6}{}
Is the discrete sum of two inversely Lindel\"of spaces 
inversely Lindel\"of? 

\mynote{References}
\cite{MR96b:54036, MR98g:54050}
\end{myprob}

\begin{myprob}
\myproblem{Problem 7}{}
Is the inequality $2^{|K|} \leq (\chi(X))^{\operatorname{St-l}(X)}$
true for every closed discrete subset $K$ in a normal space $X$? 

For a topological space $X$, $\operatorname{St-l}(X)$ denotes the 
minimum of such cardinals $\tau$ that for every open cover $U$ of $X$ 
there is a subset $A$ of $X$ with $|A|\leq\tau$ and 
$\operatorname{St}(A,U)=X$.

\mynote{Reference}
See the author's preprint \emph{Some Questions} \cite{matveev}. 
\end{myprob}

\section*{Justin Tatch Moore}

\begin{myprob}
\myproblem{Problem 1}{J.T.~Moore}
If $2^{\aleph_0} < 2^{\aleph_1}$ is there a nontrivial automorphism of
$\mathbb{N}^*$?
\end{myprob}

\begin{myprob}
\myproblem{Problem 2}{}
Does \myaxiom{OCA} imply that all automorphisms of $\mathbb{N}^*$ are 
trivial?

\mynote{Comments} 
If \myaxiom{OCA} is supplemented with \myaxiom{MA} then yes
\cite{MR94a:03080}.
It is also known that \myaxiom{OCA} implies that every automorphism of 
$\mathbb{N}^*$ is somewhere trivial \cite{MR94a:03080}.
This question was mentioned explicitly in I.~Farah's Ph.D.\ thesis 
but probably is older.
\end{myprob}

\begin{myprob}
\myproblem{Problem 3}{S.~Todor\v{c}evi\'c}
Assume $\text{\myaxiom{MA}}_{\aleph_1}$.
After forcing with an arbitrary measure algebra, must every nonmetrizable
compact space contain an uncountable discrete set in its square?

\mynote{Comments} 
If the measure algebra is trivial then yes, since any counterexample 
would have to have a square which is compact and contains an $S$-space.
See also \cite{MR93k:03050} and \cite{MR97j:03099}.
\end{myprob}

\begin{myprob}
\myproblem{Problem 4}{J.T.~Moore}
Assume $\text{\myaxiom{MA}}_{\aleph_1}$.
After adding one random real, does $\text{\myaxiom{MA}}_{\aleph_1}$ hold 
for all partial orders which have a c.c.c.\ product with every c.c.c.\ 
partial order?
\end{myprob}

\begin{myprob}
\myproblem{Problem 5}{S.~Todor\v{c}evi\'c}
If the countable chain condition is productive, must 
$\text{\myaxiom{MA}}_{\aleph_1}$ hold?

\mynote{Comments} 
While the answer seems to be no, one can still prove partial positive 
results.
For instance $\mathfrak{b}$ and $\operatorname{cf}(\mathfrak{c})$ are 
both greater than $\aleph_1$ if the countable chain condition is 
productive.
The most interesting open question of this sort is probably:
``If the countable chain condition is productive must 
$2^{\aleph_0} = 2^{\aleph_1}$?''
See, e.g., \cite{MR87b:04003.1, MR88h:54009, MR90d:04001.1}.
\end{myprob}

\begin{myprob}
\myproblem{Problem 6}{M.E.~Rudin}
After adding $\aleph_2$ Laver reals, is there a Lindel\"of space which 
has a non-Lindel\"of product with the irrationals?

\mynote{Comments} 
A solution in either direction would be interesting.
A negative solution would solve Michael's problem.
A positive solution would give a fundamentally new construction of a 
Michael space since the argument could not involve Baire category
(or the assumption $\mathfrak{b}={\aleph_1}$).
See \cite{MR99j:54008} for a survey of this problem.
M.E.~Rudin also has an unpublished note surveying Michael's problem 
(which I regrettably was unaware of when I wrote \cite{MR99j:54008}).
\end{myprob}

\begin{myprob}
\myproblem{Problem 7}{}
Is it consistent that for every $c\colon [\omega_1]^2 \to 2$ there is a 
pair of uncountable sets $A,B \subseteq \omega_1$ such that $c$ is 
constant on $A \otimes B = 
\{\{\alpha,\beta\}:\alpha \in A, \beta \in B, \alpha < \beta\}$?

\mynote{Comments} 
I believe this problem is due to F.~Galvin. 
If the answer is yes then it is consistent that every regular space 
is hereditarily separable iff it is hereditarily Lindel\"of iff it does 
not contain an uncountable discrete set (by results of Woodin, if this 
statement is consistent at all, then it is consistent with 
$\text{\myaxiom{MA}}_{\aleph_1}$). 
See \cite{MR85d:03102}.
\end{myprob}

\section*{Peter Nyikos} 

\begin{myprob}
\myproblem{Problem}{}
Is every separable normal manifold $\omega_1$-compact?

\mynote{Comments}
A space is said to be $\omega_1$-\emph{compact} if every closed 
discrete subspace is countable.

An equivalent problem is: In a separable normal manifold, is every closed 
discrete subspace a $G_\delta$---a countable intersection of open sets? 
This problem is unsolved even for hereditarily normal manifolds.

Jones' Lemma shows that the answer to this problem is yes if 
$2^{\aleph_0} < 2^{\aleph_1}$. 
Some theorems on manifolds would be improved upon if we could even show 
that \manotch\ or even \myaxiom{PFA} is consistent with an affirmative 
answer. 
A yes answer would follow from one to Mary Ellen Rudin's favorite problems 
about manifolds: is every normal manifold collectionwise Hausdorff? 

There are plenty of separable normal manifolds from \myaxiom{ZFC} 
alone, but all the ones I know of are $\omega_1$-compact.

\mynote{References}
\cite{MR85i:54004.1, MR86f:54054.1}
\end{myprob}

\section*{Stanislav Peregudov}

\begin{myprob}
\myproblem{Problem 1}{}
Is every Lindel\"of regular space that has a weakly uniform base first 
countable? 
No, in the class of Lindel\"of Hausdorff spaces.
\end{myprob}

\begin{myprob}
\myproblem{Problem 2}{}
Is there an $L$-space with a weakly uniform base? 
No, under \manotch\ \cite{MR2002m:03096}.
\end{myprob}

\begin{myprob}
\myproblem{Problem 3}{}
Is every pseudocompact space that has a weakly uniform base compact? 
The space is \v{C}ech complete first countable \cite{MR2000h:54033}.
\end{myprob}

\section*{Roman Pol}

\begin{myprob}
\myproblem{Problem}{}
Let $f\colon X \to Y$ be a light map (i.e., with zero-dimensional 
fibers), where $X$ is a metrizable continuum with $\dim X > 2$.
Does there exist a nontrivial continuum $C$ in $X$ such that $f$ 
restricted to $C$ is injective ?

\mynote{Comments} 
For some related information, see \cite{MR98e:54014}.
In particular, the answer is positive for $Y$ finite-dimensional. 
\end{myprob}

\section*{James T. Rogers, Jr.}

Let $M$ be a hereditarily indecomposable continuum. 
Assume $\dim M=n>1$.
Let $H(M)$ be the homeomorphism group of $M$.

\startproblem
\begin{myprob}
\myproblem{Question 1}{}
Can $H(M)$ contain a nontrivial continuum? 
A nontrivial connected set?

For each integer $n>1$, Rogers has exhibited an $M$ such that $H(M)$ 
contains no nontrivial connected set.
\end{myprob}

\begin{myprob}
\myproblem{Question 2}{}
Can $M$ be rigid? 
i.e., the identity map is the only element of $H(M)$?
\end{myprob}

\section*{Mary Ellen Rudin}

\begin{myprob}
\myproblem{The Linearly Lindel\"of Problem}{}
Is there a Hausdorff, normal, linearly Lindel\"of not Lindel\"of space? 

\mynote{Comments} 
A space is \emph{linearly Lindel\"of} if every increasing open cover 
has a countable subcover. 
There are several regular but not normal examples known. 
It is normality that is hard to achieve.

\mynote{References}
\cite{MR99j:54006, MR25:4483.1}
\end{myprob}

\begin{myprob}
\myproblem{The Point Countable Base Problem}{}
We are given a $T_1$ space $X$ such that each point $x$ of $X$ has a 
countable open basis $\mathcal{B}(x)$ having the property that for all 
open $U$ with $x \in U$ there is an open $V$ with $x \in V \subset U$ 
such that $y \in V$ implies there is $B \in \mathcal{B}(y)$ with 
$x \in B \subset U$. 
Does this imply that the space must have a point countable base (i.e., 
that each point is in only countably many members of the base)?

\mynote{Comments} 
It is true if the density of the space is $\leq \omega_1$.

\mynote{References}
\cite{MR1078650.1, MR1229128.1}
\end{myprob}

\begin{myprob}
\myproblem{Max Burke's Spaces}{}
Max Burke would like a nice characterization of the class of spaces which 
are the continuous image of an arbitrary product of compact, linearly 
ordered spaces. 

\mynote{Comments}
Burke has recently proved that if $X$ is in this class then every 
separately continuous function $f\colon X \times Y \to \mathbb{R}$ (where 
$Y$ is any space) is Borel. 
Also, if $I\colon C(X) \to \mathbb{R}$ is integration with respect to 
finite Borel measure on $X$, then $I$ is Borel measurable when $C(X)$ has 
the topology of pointwise convergence.
See \cite{MR2003m:54045}.
\end{myprob}

\section*{Krishnan Shankar} 

\begin{myprob}
\myproblem{Problem}{}
Is the Berger space an $S^3$ bundle over $S^4$?

\mynote{Comments}
$S^3$ bundles over $S^4$ have been important ever since 
J.~Milnor~\cite{MR18:498d} showed that the total spaces of such bundles 
with Euler class $\pm 1$ are homeomorphic to the standard sphere $S^7$ 
but not always diffeomorphic to it.
In 1974, D.~Gromoll and W.~Meyer \cite{MR51:11347} constructed a metric 
of non-negative curvature on one of these exotic spheres (a generator in 
the group of homotopy $7$-spheres) and it remained the only exotic sphere 
known to admit non-negative curvature. 
More recently, K.~Grove and W.~Ziller \cite{MR2001i:53047} constructed 
metrics of non-negative curvature on all the total spaces of $S^3$ 
bundles over $S^4$ thus showing that all the Milnor $7$-spheres admit 
non-negative curvature. 
They also asked whether the homogeneous Berger space admits the structure 
of such a bundle. 
They showed that it cannot be a principal $S^3$ bundle over $S^4$.

The Berger space is described as the quotient $Sp(2)/Sp(1) = SO(5)/SO(3)$ 
where the embedding of $SO(3)$ is nonstandard. 
If we think of the space of $\mathbb{R}^5$ as the space of symmetric, 
$3\times 3$, traceless matrices, then $SO(3)$ acts on this space by 
conjugation which gives a representation of $SO(3)$ into $SO(5)$. 
The Berger space is described as the resulting quotient space.

In a recent paper, Kitchloo and Shankar \cite{MR2001m:55044} showed that 
the Berger space $SO(5)/SO(3)$ is PL-homeomorphic to an $S^3$ bundle over 
$S^4$. 
They were unable to settle the diffeomorphism question since this 
requires computing the Eells-Kuiper invariant. 
To do this one needs to exhibit the Berger space as a spin coboundary. 
It remains open whether the Berger space is diffeomorphic to an $S^3$ 
bundle over $S^4$.
\end{myprob}

\providecommand{\bysame}{\leavevmode\hbox to3em{\hrulefill}\thinspace}

\label{tpoxfordend}

\chapter*{Wayne Lewis: Continuum theory problems}
\label{tplewis}
\begin{myfoot}
\begin{myfooter}
Wayne Lewis, 
\emph{Continuum theory problems},\\
Problems from Topology Proceedings, Topology Atlas, 2003, 
pp.\ 165--182.
\end{myfooter}
\end{myfoot}

\mypreface
Wayne Lewis collected this list of problems in volume 8 (1983) of 
\emph{Topology Proceedings} \cite{MR86a:54038} and updated it in 
volume 9 \cite{TP9.2.375}. 
This version includes a few new notes.

\section*{Introduction}

The problems listed below have come from a number of sources. 
Some were posed at the Texas Topology Symposium, 1980 in Austin, 
some at the American Mathematical Society meeting in Baton Rouge in 1982, 
some at the Topology Conference in Houston in 1983, 
some at discussions at the University of Florida in 1982, 
and some at the International Congress of Mathematicians in Warsaw in 1983. 
Some are classical, while others are more recent or primarily of technical 
interest. 
Preliminary versions of subsets of this list have been circulated, 
and an attempt has been made to verify the accuracy of the statements of the 
questions, comments, and references given. 
In many cases, variations on a given question have been asked by many people 
on diverse occasions. 
Thus the version presented here should not be considered definitive. 
Any errors or additions which are brought to my attention will be noted at a 
later date. 

The division of the questions into categories is only intended as a 
rough guide, and many questions could properly be placed in more than one 
category. 
A number of these questions have appeared in the University of Houston 
Problem Book (UHPB), a good reference for further problems. 
Assistance in compiling earlier versions of subsets of this list was provided 
by Bellamy, Brechner, Heath, and Mayer.

\section*{Chainable continua}

\begin{myprob}
\myproblem{1}{Brechner, Lewis, Toledo}
Can a chainable continuum admit two non-con\-jugate homeomorphisms of 
period $n$ with the same fixed-point set? 

\mynote{Notes}
Earlier (Brechner): Are every two period $n$ homeomorphisms of the 
pseudo-arc conjugate? 
Lewis has since shown that the pseudo-arc admits homeomorphisms of every 
period, and Toledo has shown that it admits such homeomorphisms with 
nondegenerate fixed-point sets. 

\mynote{Update}
The question should be rephrased to require the sets of fundamental
periods of points under the two homeomorphisms to be identical. 
Toledo has shown that for any sequence of positive integers 
$1 \leq n_0 < n_1 < n_2 < \cdots < n_k$ where $n_i$ divides $n_j$ for 
$i < j$, there is a period $n_k$ homeomorphism of the pseudo-arc
with points of each of the fundamental periods $n_i$
\end{myprob}

\begin{myprob}
\myproblem{2}{Brechner}
Classify, up to conjugacy, the periodic homeomorphisms of the pseudo-arc.
\end{myprob}

\begin{myprob}
\myproblem{3}{Anderson}
Does every Cantor group act effectively on the pseudo-arc? 

\mynote{Notes}
Lewis has shown that every inverse limit of finite solvable groups acts 
effectively on the pseudo-arc.
\end{myprob}

\begin{myprob}
\myproblem{4}{Nadler}
Does the pseudo-arc have the complete invariance property? 

\mynote{Notes}
A continuum $X$ has the \emph{complete} invariance property if every 
nonempty closed subset of $X$ is the fixed-point set of some continuous 
self-map of $X$. 
Martin and Nadler have shown that every two-point set is a fixed-point 
set for some continuous self-map of the pseudo-arc. 
Cornette has shown that every subcontinuum of the pseudo-arc is a retract. 
Toledo has shown that every subcontinuum is the fixed-point set of a 
periodic homeomorphism. 
Lewis has shown that there are proper subsets of the pseudo-arc with 
nonempty interior which are the fixed-point sets of homeomorphisms.
\end{myprob}

\begin{myprob}
\myproblem{5}{Brechner and Lewis}
Do there exist stable homeomorphisms of the pseudo-arc which are extendable 
(or essentially extendable) to the plane? 
How many, up to conjugacy? 

\mynote{Notes}
This is a rewording of a question earlier posed by Brechner. 
Lewis has shown that there are non-identity stable homeomorphisms of the 
pseudo-arc.
\end{myprob}

\begin{myprob}
\myproblem{6}{Brechner}
Let $M$ be a particular embedding of the pseudo-arc in the plane, and let $G$ 
be the group of extendable homeomorphisms of $M$. 
Does $G$ characterize the embedding? 
\end{myprob}

\begin{myprob}
\myproblem{7}{Lewis}
Are the periodic (resp.\ almost periodic, or pointwise-periodic) 
homeomorphisms dense in the group of homeomorphisms of the pseudo-arc? 

\mynote{Notes}
The conjecture is that answer is no. 
For each $n \geq 2$, the period $n$ homeomorphisms do act transitively 
on the pseudo-arc.
\end{myprob}

\begin{myprob}
\myproblem{8}{Brechner}
Does each periodic homeomorphism $h$ of the pseudo-arc have a square root 
(i.e., a homeomorphism $g$ such that $g^2 = h$)? 

\mynote{Notes}
It is known that some periodic homeomorphisms have an infinite sequence of 
$p_i$-roots, for any sequence $\{p_i\}$ of positive integers.
\end{myprob}

\begin{myprob}
\myproblem{9}{Toledo}
Can a pointwise-periodic, regular homeomorphism on a chainable 
(indecomposable), or tree-like (indecomposable) continuum, or the 
pseudo-arc, always be induced by square commuting diagrams on inverse 
systems of finite graphs? 

\mynote{Notes}
Fugate has shown that such homeomorphisms on chainable continua cannot 
always be induced by square commuting diagrams on inverse systems of 
arcs. Toledo has shown that periodic homeomorphisms of the pseudo-arc can 
always be induced by square commuting diagrams on finite graphs (not 
necessarily trees).
\end{myprob}

\begin{myprob}
\myproblem{10}{Toledo} 
Can a homeomorphism of a chainable continuum always be induced by square 
commuting diagrams on inverse systems of finite graphs? 

\mynote{Notes}
See remark after question 9.
\end{myprob}

\begin{myprob}
\myproblem{12}{Duda} 
Characterize chainability and/or circular chainability without using span. 

\mynote{Notes}
Oversteegen and Tymchatyn have a technical partial characterization, 
but a complete, useful, and satisfying characterization remains to be 
developed.
\end{myprob}

\begin{myprob}
\myproblem{12}{Duda} 
What additional conditions make the following statement true? 
If $X$ has only chainable proper subcontinua and (?) then $X$ is either 
chainable or circularly chainable. 

\mynote{Notes}
Ingram's examples show that an additional condition is needed. 
If $X$ is decomposable, no additional condition is needed. 
If $X$ is hereditarily indecomposable, then either homogeneity 
(or the existence of a $G_\delta$ orbit under the action of its homeomorphism 
group) or weak chainability is a sufficient condition. 
Hereditary indecomposability alone is insufficient. 
Also (Fugate, UHPB 106): 
If $M$ is tree-like and every proper subcontinuum of $M$ is chainable, 
is $M$ almost chainable?
\end{myprob}

\begin{myprob}
\myproblem{13}{Fugate, UHPB 104} 
If $X$ is circularly chainable and $f\colon X \to Y$ is open, then is $Y$ 
either chainable or circularly chainable? 

\mynote{Notes}
Yes if $Y$ is decomposable.

\mynote{Update}
Krupski has shown that if $X$ is a solenoid, then $Y$ is either a point, 
a solenoid, or a Knaster continuum, i.e., an inverse limit of arcs with 
open bonding maps. 
\end{myprob}

\begin{myprob}
\myproblem{14}{Duda} 
Can the following theorem be improved---say by dropping ``hereditarily 
decomposable''?

\begin{theorem}[Duda and Kell]
Let $f\colon X \to Y$ be a finite-to-one open mapping of an hereditarily 
decomposable chainable continuum onto a $T_2$ space. 
Then $X = \bigcup_{j=1}^n K_j$ where each $K_j$ is a continuum and 
$f |_{\operatorname{Int}K_j}$ is a homeomorphism. 
\end{theorem}
\end{myprob}

\begin{myprob}
\myproblem{15}{Cook and Fugate, UHPB 105} 
Suppose $M$ is an atriodic one-dimensional continuum and $G$ is an upper 
semi-continuous decomposition of $M$ such that $M/G$ and every element of $G$ 
are chainable. 
Is $M$ chainable? 

\mynote{Notes}
Michel Smith has shown that if ``one-dimensional'' is removed and ``$M/G$ 
is an arc'' is added to the hypothesis, then the answer is yes.

It follows from a result of Sher that even if $M$ contains a triod, if $M/G$ 
and every element of $G$ are tree-like, then $M$ is tree-like.
If $M$ is hereditarily indecomposable and $G$ is continuous then the 
answer is yes.
\end{myprob}

\begin{myprob}
\myproblem{16}{Mohler}
Is every weakly chainable, atriodic, tree-like continuum chainable? 

\mynote{Notes}
A positive answer would imply that the classification of homogeneous plane 
continua is complete.
\end{myprob}

\section*{Decompositions}

\begin{myprob}
\myproblem{17}{Rogers} 
Suppose $G$ is a continuous decomposition of $E^2$ into nonseparating 
continua. 
Must some element of $G$ be hereditarily indecomposable? 
What if all of the decomposition elements are homeomorphic? 
Must some element have span zero? 
be chainable? 

\mynote{Notes}
This is a revision of a question by Mayer. 
Possibly related to this, Oversteegen and Mohler have recently shown that 
there exists an irreducible continuum $X$ and an open, monotone map 
$f\colon X \to [0,1]$ such that each nondegenerate subcontinuum of $X$ 
contains an arc, and so no nondegenerate $f^{-1}(t)$ is hereditarily 
indecomposable. 
Oversteegen and Tymchatyn have shown that there must exist an $f^{-1}(t)$ 
which contains arbitrarily small indecomposable subcontinua.
\end{myprob}

\begin{myprob}
\myproblem{18}{Krasinkiewicz, UHPB 158}
Let $X$ be a nondegenerate continuum such that there exists a continuous 
decomposition of the plane into elements homeomorphic to $X$. 
Must $X$ be the pseudo-arc?
\end{myprob}

\begin{myprob}
\myproblem{19}{Mayer} 
How many inequivalent embeddings of the pseudo-arc are to be found in the 
Lewis-Walsh decomposition of $E^2$ into pseudo-arcs? 
\end{myprob}

\begin{myprob}
\myproblem{20}{Ingram} 
Does there exist a tree-like, non-chainable continuum $M$ such that the plane 
contains uncountably many disjoint copies of $M$? 
Is there a continuous collection of copies of $M$ filling up the plane? 

\mynote{Notes}
W.T.~Ingram has constructed an uncountable collection of disjoint, 
nonhomeomorphic, tree-like, non-chainable continua in the plane.
\end{myprob}

\begin{myprob}
\myproblem{21}{Lewis} 
Is there a continuous decomposition of $E^2$ into Ingram continua (not 
necessarily all homeomorphic)? 
\end{myprob}

\begin{myprob}
\myproblem{22}{Lewis} 
If $M$ is an hereditarily equivalent or homogeneous, nonseparating plane 
continuum, does there exist a continuous collection of continua, each 
homeomorphic to $M$, filling up the plane? 
Does the plane contain a (homogeneous) continuous circle of copies of $M$, 
as in the Jones Decomposition Theorem? 
\end{myprob}

\begin{myprob}
\myproblem{23}{Lewis} 
If $X$ and $Y$ are one-dimensional continua with continuous decompositions $G$ 
and $H$, respectively, into pseudo-arcs such that $X/G$ and $Y/H$ are 
homeomorphic, then are $X$ and $Y$ homeomorphic? 

\mynote{Notes}
It follows from arguments of Lewis that if every element of $G$ and $H$ is a 
terminal continuum in $X$ and $Y$ respectively then $X$ and $Y$ are 
homeomorphic.
\end{myprob}

\begin{myprob}
\myproblem{24}{Burgess}
Is there a continuous decompositions $G$ of $E^3$ into pseudo-arcs such 
that $E^3/G \approx E^3$ and the pre-image of each one-dimensional 
continuum is one-dimensional?
If so, is the pre-image of a homogeneous curve under such a decomposition
itself homogeneous?
Can this process produce any new homogeneous curves?

\mynote{Notes}
It is known that for every one-dimensional continuum $M$ there exists a 
one-dimensional continuum $\hat{M}$ with a continuous decomposition $G$ into
pseudo-arcs such that $\hat{M}/G \approx M$.
If $M$ is homogeneous, then $\hat{M}$ can be constructed to be homogeneous.
This method can produce new homogeneous continua.
\end{myprob}

\section*{Fixed points}

\begin{myprob}
\myproblem{25}{Bellamy} 
Allowing singletons as degenerate indecomposable continua, are the following 
statements true? 
\begin{myenumerate}
\item
Suppose $X$ is a tree-like continuum and $f\colon X \to X$ is continuous. 
Then there is an indecomposable subcontinuum $W$ of $X$ such that 
$f(W) \subseteq W$. 
\item
The same with hereditarily unicoherent replacing tree-like in the 
hypothesis. 
\end{myenumerate}

\mynote{Notes}
Bellamy has constructed a tree-like indecomposable continuum without the 
fixed-point property. 
Manka has shown that every $\lambda$-dendroid (hereditarily decomposable, 
hereditarily unicoherent continuum) has the fixed-point property. 
Cook has shown that $\lambda$-dendroids are tree-like. 

\mynote{Solution}
Ma\'ckowiak has described such a hereditarily unicoherent continuum $X$ 
and map $f$ so that the statement is false. 
\end{myprob}

\begin{myprob}
\myproblem{26}{Bellamy} 
Suppose $X$ is a tree-like continuum and every indecomposable subcontinuum has 
the fixed-point property. 
Does $X$ have the fixed-point property? 
\end{myprob}

\begin{myprob}
\myproblem{27}{Bellamy} 
Suppose $X$ is a tree-like continuum and $f\colon X \to X$ is a function 
homotopic to the identity on $X$. 
Must $f$ have a fixed-point? 
\end{myprob}

\begin{myprob}
\myproblem{28}{Bellamy} 
Suppose $X$ is a tree-like continuum. 
Does there exist $\epsilon > 0$ such that every self-map of $X$ within 
$\epsilon$ of the identity has a fixed-point? 
\end{myprob}

\begin{myprob}
\myproblem{29}{Knaster} 
Does every hereditarily indecomposable tree-like continuum have the fixed 
point property? 

\mynote{Solution}
No (P. Minc \cite{MR2000k:54029}).
\end{myprob}

\begin{myprob}
\myproblem{30}{Cook} 
Does every hereditarily equivalent continuum have the fixed-point 
property? 

\mynote{Notes}
A continuum is \emph{hereditarily equivalent} if it is homeomorphic with 
each of its nondegenerate subcontinua. 
Cook has shown that every nondegenerate hereditarily equivalent continuum 
other than the arc or pseudo-arc is hereditarily indecomposable and tree-like.
\end{myprob}

\begin{myprob}
\myproblem{31}{Bellamy} 
Suppose $X$ is triod-like (or $K$-like for some fixed tree $K$). 
Must $X$ have the fixed-point property? 

\mynote{Notes}
Marsh has shown that an inverse limit of fans $\{F_i\}$---where each
bonding map preserves ramification points and is except for one branch, a
homeomorphism of each branch of $F_{i+1}$ onto a branch of $F_i$---has 
the fixed-point property. 
\end{myprob}

\begin{myprob}
\myproblem{32}{Bellamy} 
Does every inverse limit of real projective planes with homotopically 
essential bonding maps have the fixed-point property? 
for homeomorphisms? 
\end{myprob}

\begin{myprob}
\myproblem{33}{}
Suppose $X$ is a nonseparating plane continuum with each arc component 
dense. 
Is $X$ an almost continuous retract of a disc? 

\mynote{Notes}
If $X \subseteq D$, a function $f\colon D \to X$ is 
\emph{almost (quasi-) continuous} if every neighborhood in $D \times D$ 
(in $D \times X$) of the graph of $f$ contains the graph of a continuous 
function with domain $D$.
Akis has shown that the disc with a spiral about its boundary is neither 
an almost continuous nor quasi-continuous retract of a disc.
\end{myprob}

\begin{myprob}
\myproblem{34}{Bellamy} 
Suppose $f$ is a self-map of a tree-like continuum which commutes with 
some homeomorphism of period greater than one, or with every member of 
some nondegenerate compact group of homeomorphisms. 
Must $f$ have a fixed-point? 

\mynote{Notes}
Fugate has shown that if a compact group acts on a tree-like continuum, 
then all the homeomorphisms in the group have a common fixed-point. 
\end{myprob}

\begin{myprob}
\myproblem{35}{Edwards} 
Does every self-map (homeomorphism) of a tree-like continuum have a periodic 
point? 
\end{myprob}

\begin{myprob}
\myproblem{36}{Bellamy} 
Does every weakly chainable tree-like continuum have the fixed-point 
property? 
What about tree-like continua which are the continuous image of circle-like 
continua? 
\end{myprob}

\begin{myprob}
\myproblem{37}{Rosenholtz} 
Suppose $f$ is a map from a nonseparating plane continuum $M$ to itself which 
is differentiable (i.e., $f$ can be extended to a neighborhood of $M$ 
with 
partial derivatives existing). 
Must $f$ have a fixed-point? 
\end{myprob}

\begin{myprob}
\myproblem{38}{Sternbach, Scottish Book 107} 
Does every nonseparating plane continuum have the fixed-point property? 
\end{myprob}

\begin{myprob}
\myproblem{39}{Bellamy} 
Do each two commuting functions on a simple triod have a common incidence 
point? 
\end{myprob}

\begin{myprob}
\myproblem{40}{Manka} 
Let $C$ be the composant with an endpoint in the simplest Knaster 
indecomposable continuum. 
Does $C$ have the fixed-point property? 

\mynote{Notes}
Also: If $f\colon C \to C$ is continuous with noncompact image, is $f$ 
onto? 
An affirmative answer gives an affirmative answer to the previous question. 
\end{myprob}

\begin{myprob}
\myproblem{41}{Oversteegen and Rogers} 
Does the cone over $X$ have the fixed-point property, where $X$ is the 
tree-like continuum without the fixed-point property constructed by 
Oversteegen and Rogers? 
\end{myprob}

\begin{myprob}
\myproblem{42}{{\L}ysko} 
Does there exist a continuum $X$ with the fixed-point property such that 
$X \times P$ ($P$ = pseudo-arc) does not have the fixed-point property? 
\end{myprob}

\begin{myprob}
\myproblem{43}{Gordh} 
If $X$ is an irreducible continuum and each tranch has the fixed-point 
property, must $X$ have the fixed-point property? 

\mynote{Notes}
If $X$ is an irreducible continuum such that each indecomposable 
subcontinuum of $X$ is nowhere dense, then there exists a finest monotone 
map $f\colon X \to [0,1]$. 
Point-inverses under $f$ are nowhere dense subcontinua of $X$ and are called 
the \emph{tranches} of $X$.
\end{myprob}

\begin{myprob}
\myproblem{44}{Bell} 
Is there a map $f\colon K \to K$, where $K$ is a continuum in 
$\mathbb{R}^2$ and $K$ is minimal with respect to $f(K) \subset K$, such 
that $\operatorname{Index}(f,K) = 0$? 

\mynote{Notes}
If $g\colon A \to \mathbb{R}^{n+1}$ is a fixed-point free map where $A$ 
is an $n$-sphere in $\mathbb{R}^{n+l}$, then $\operatorname{Index}(g,A)$ is 
the degree of $h(z) = \frac{g(z) - z}{||g(z) - z||}$. 
If $K$ is a point-like continuum in $\mathbb{R}^n$ and $f$ is a 
fixed-point free map $f\colon \operatorname{Bd}K \to \mathbb{R}^n$ then 
$f$ has an extension to a map $F\colon \mathbb{R}^n \to \mathbb{R}^n$ 
that is fixed-point free on $\mathbb{R}^n - K$. 
$\operatorname{Index}(f, \operatorname{Bd}X) = \operatorname{Index}(F,B)$, 
where $B$ is any $n-1$ sphere in $\mathbb{R}^n$ that surrounds $K$. 
\end{myprob}

\begin{myprob}
\myproblem{45}{Bell} 
Let $B$ be a point-like continuum in $\mathbb{R}^n$, $n > 2$, 
$f\colon \operatorname{Bd}(B) \to B$, and 
$\operatorname{Index}(f, \operatorname{Bd}B) = 0$. 
Must there be a continuum $K \subset \operatorname{Bd}B$ such that $K = f(K)$? 

\mynote{Notes}
The answer is no if there is a fixed-point free map on a point-like 
continuum $X$, 
where $\operatorname{Bd}K$ contains no invariant subcontinua. 
\end{myprob}

\begin{myprob}
\myproblem{46}{Minc} 
Is there a planar continuum $X$ and $f\colon X \to X$ such that $f$ 
induces the zero homomorphism on the first \v{C}ech cohomology group and 
$f$ is fixed-point free? 
\end{myprob}

\section*{Higher-dimensional problems}

\begin{myprob}
\myproblem{47}{Ancel} 
If $f\colon S^2 \to \mathbb{R}^3$ is continuous and $U$ is the unbounded 
component of $\mathbb{R}^3 - f(S^2)$, is 
$f\colon S^2 \to \mathbb{R}^3 - U$ homotopically trivial in 
$\mathbb{R}^n - U$?

\mynote{Notes}
The analogous result is true one dimension lower and false one dimension 
higher. 
\end{myprob}

\begin{myprob}
\myproblem{48}{Ancel} 
If $X$ is a cellular subset of $\mathbb{R}^3$ is $\pi_2(X) = 0$? 
\end{myprob}

\begin{myprob}
\myproblem{49}{Burgess} 
Is a $2$-sphere $S$ in $S^3$ tame if it is homogeneously embedded?

$S$ is \emph{homogeneously embedded} in $S^3$ if for each $p$,$q$ in $S$ 
there is a homeomorphism $h\colon (S^3, S, p) \to (S^3, S, q))$.
\end{myprob}

\begin{myprob}
\myproblem{50}{Burgess} 
Is a $2$-sphere $S$ in $S^3$ tame if every homeomorphism of $S$ onto 
itself can be extended to a homeomorphism of $S^3$ onto itself? 
\end{myprob}

\begin{myprob}
\myproblem{51}{Bing} 
If $S$ is a toroidal simple closed curve in $E^3$ (i.e., an intersection of 
nested solid tori with small meridional cross-sections) such that over each 
arc $A$ in $S$ a singular fin can be raised, with no singularities on $A$,
must $S$ be tame? 

A \emph{fin} is a disc which contains $A$ as an arc on its boundary and 
is otherwise disjoint from $S$. 
It follows from a result of Burgess and Cannon that $S$ is tame if the fin 
can always be chosen to be non-singular.
\end{myprob}

\begin{myprob}
\myproblem{52}{Bing} 
Is a simple closed curve $S$ in $E^3$ tame if it is isotopically homogeneous 
(i.e., for each $p$, $q$ in $S$ there is an ambient isotopy of $E^3$, 
leaving $S$ invariant at each stage, with the $0$-th level of the isotopy the 
identity and the last level a homeomorphism taking $p$ to $q$)? 

\mynote{Notes}
Compare this with this dissertation and related work of Shilepsky. 
Shi\-lepsky has conjectured that the answer is yes. 
Shilepsky and Bothe have independently constructed wild simple closed
curves in $E^3$ which are homogeneously embedded in $E^3$ but not
isotopically homogeneous. 
\end{myprob}

\begin{myprob}
\myproblem{53}{J. Heath, Jack Rogers}
If $r\colon X \to Y$ is refinable and $X$ is an ANR, must $Y$ be an ANR? 

\mynote{Notes}
A map $r\colon X \to Y$ is \emph{refinable} if for each $\epsilon > 0$ 
there is an $\epsilon$-refinement, i.e., an $\epsilon$-map 
$g\colon X \to Y$ such that $\operatorname{dist}(g(x), r(x)) < \epsilon$ 
for each $x \in X$. 
Heath and Kozlowski have shown: If $X$ is finite dimensional, then $Y$ 
must be an ANR if either:
each $r^{-1}(y)$ is locally connected;
each $r^{-1}(y)$ is nearly $1$-movable; 
each $r^{-1}(y)$ is approximately $1$-connected;
$Y$ is $\operatorname{LC}^1$ at each point, or;
there is a monotone $\epsilon$-refinement of $r$ for each $\epsilon > 0$.
\end{myprob}

\begin{myprob}
\myproblem{54}{J. Heath, Kozlowski} 
If $r\colon S^3 \to Y$ is refinable, is $Y$ an ANR?
\end{myprob}

\begin{myprob}
\myproblem{55}{J. Heath, Kozlowski}
If $r\colon S^n \to S^n/A$ is refinable and $n > 3$, must $A$ be 
cellular? 

\mynote{Notes}
The answer is yes if $n \leq 3$.
\end{myprob}

\begin{myprob}
\myproblem{56}{Edwards} 
If $f\colon S^3 \to S^2$ is a continuous surjection must there exist 
$\Sigma^2$ (an embedded copy of $S^2$ in $S^3$) such that $f |_{\Sigma^2}$ 
is a surjection?

\mynote{Notes}
The analogous question for a map $f\colon S^2 \to S_1$ has an affirmative 
solution. 

\mynote{Solution}
Bestvina and Walsh have shown that the answer is no. 
\end{myprob}

\begin{myprob}
\myproblem{57}{Boxer} 
Do ARI maps preserve property K? 

\mynote{Notes}
A continuous surjection of compacta $f\colon X \to Y$ is 
\emph{approximately right invertible} (ARI) if there is a null sequence 
$\{\epsilon_n\}$ of positive numbers and a sequence of maps 
$g_n\colon Y \to X$ such that $d(fg_n, \operatorname{id}_Y) < \epsilon_n$ 
for each $n$ ($d$ = sup-metric). 
A continuum $X$ has \emph{property K} if for each $\epsilon > 0$ there 
exists $\delta > 0$ such that for each $p \in X$ and each $A \in C(X)$ with 
$p \in A$, if $q \in X$ and $\operatorname{dist}(p,q) < \delta$, then there 
exists $B \in C(X)$ with $q \in B$ and $H(A,B) < \epsilon$. 
($H$ = Hausdorff metric. $C(X)$ = hyperspace of subcontinua of $X$.) 
If $\{g_n\}$ is equicontinuous, the question has a positive answer. 
This does not represent new knowledge unless the next question has a negative 
answer for a continuum $X$ with property K. 
The above question is a special case of a question in Nadler's book. 
\end{myprob}

\begin{myprob}
\myproblem{58}{Boxer} 
If $f\colon X \to Y$ is an ARI map with an equicontinuous sequence as in 
the above comments, is $f$ an $r$-map? 
\end{myprob}

\section*{Homeomorphism groups}

\begin{myprob}
\myproblem{59}{Duda} 
Let $G(P)$ is the group of homeomorphisms of the pseudo-arc $P$.
Is the map $h\colon G(P) \to \mathbb{R}$ defined by 
$H(g) = \operatorname{dist}(g, \operatorname{id})$, a surjection 
onto $[0, \operatorname{diam}P]$, or does the image at least contain a 
neighborhood (relative to $[0, \operatorname{diam}P]$) of $0$?
\end{myprob}

\begin{myprob}
\myproblem{60}{Brechner} 
Is the homeomorphism group of the pseudo-arc totally disconnected? 

\mynote{Notes}
Brechner and Anderson have proven an analogous result for the Menger universal 
curve. 
The homeomorphism group of the pseudo-arc contains no nondegenerate 
subcontinua, by a result of Lewis. 
\end{myprob}

\begin{myprob}
\myproblem{61}{Lewis} 
Is the homeomorphism group of every hereditarily indecomposable continuum 
totally disconnected? 
\end{myprob}

\begin{myprob}
\myproblem{62}{Lewis} 
Must the homeomorphism group of a homogeneous continuum either contain an arc 
or be totally disconnected? 
\end{myprob}

\begin{myprob}
\myproblem{63}{Brechner} 
If a homogeneous continuum $X$ has a homeomorphism group which contains an arc 
(or admits nontrivial isotopies), must $X$ admit a nontrivial flow? 
\end{myprob}

\begin{myprob}
\myproblem{64}{Brechner} 
Is the homeomorphism group of the pseudo-arc infinite dimensional? 
\end{myprob}

\begin{myprob}
\myproblem{65}{Lewis} 
Is the homeomorphism group of every nondegenerate homogeneous continuum 
infinite dimensional? 

\mynote{Notes}
Keesling has shown that if the homeomorphism group $G(X)$ of a compact metric 
space $X$ contains an arc, then $G(X)$ is infinite dimensional. 
\end{myprob}

\begin{myprob}
\myproblem{66}{Lewis} 
Is every connected subset of the space of continuous maps of the pseudo-arc 
into itself which contains a homeomorphism degenerate? 

\mynote{Notes}
The analogous result for the Menger universal curve is true. 
\end{myprob}

\begin{myprob}
\myproblem{67}{Lewis} 
Is there a natural measure which can be put on the space $M(P)$ of self-maps of 
the pseudo-arc? 
If so, what is the measure of the subspace $\hat{H}(P)$ of maps which are 
homeomorphisms onto their image? 
Is it the same as the measure of $M(P)$? 

\mynote{Notes}
$\hat{H}(P)$ is a dense $G_\delta$ in $M(P)$. 

Kallman has shown that there is no standard Borel structure on
$H(P)$---the full autohomeomorphism group of the pseudo-arc $P$---with 
respect to which $H(P)$ is a Borel group and which admits a
$\sigma$-finite Borel measure which is quasi-invariant under left
translations. 
This seems to imply a negative answer to this question. 
\end{myprob}

\begin{myprob}
\myproblem{68}{Lewis} 
Does the pseudo-circle have uncountably many orbits under the action of its 
homeomorphism group? 
What about other non-chainable continua all of whose nondegenerate proper 
subcontinua are pseudo-arcs? 

\mynote{Notes}
It can be shown that each orbit of such a continuum is dense, and that no 
such continuum has a $G_\delta$ orbit. 

\mynote{Solution}
Lewis proved that no such continuum has a $G_\delta$-orbit under the
action of its homeomorphism group. 
Kennedy and Rogers observed that a version of Effros' theorem implies a
positive answer to both questions. 
\end{myprob}

\begin{myprob}
\myproblem{69}{Wechsler} 
If $X$ and $Y$ are homogeneous continua with isomorphic and homeomorphic 
homeomorphism groups, are $X$ and $Y$ homeomorphic? 

\mynote{Notes}
Whittaker has shown that compact manifolds with or without boundary are 
homeomorphic if and only if their homeomorphism groups are isomorphic. 
Rubin has shown that if $X$ and $Y$ are locally compact and strongly locally 
homogeneous, then they are homeomorphic if and only if their homomorphism 
groups are isomorphic. 
Sharma has shown that there are (nonmetric) locally compact Galois spaces 
$X$ and $Y$ with isomorphic homeomorphism groups such that $X$ and $Y$ are not 
homeomorphic. 
van~Mill has constructed non-locally compact, connected subsets of the 
$2$-sphere which are strongly locally homogeneous and have algebraically 
(but not topologically) isomorphic homeomorphism groups, but which are not 
themselves homeomorphic.
\end{myprob}

\begin{myprob}
\myproblem{70}{Ancel} 
Let $G$ be the space of homeomorphisms of $S^n$ and $X$ the embeddings of 
$S^{n-1}$ into $S^n$.
Is there some condition analogous to $1$-ULC which will detect when the 
orbit in $G$ of a given embedding is a $G_\delta$? 
Is there also a way to distinguish non-$G_\delta$ orbits? 
Are there, in some reasonable sense, more non-$G_\delta$ orbits than 
$G_\delta$ orbits? 
\end{myprob}

\begin{myprob}
\myproblem{71}{Brechner} 
If $G$ is the collection of homeomorphisms of the pseudo-arc $P$ which 
leave every composant invariant, is $G$ dense in the full homeomorphism 
group of $P$? of first category? 
\end{myprob}

\begin{myprob}
\myproblem{72}{Brechner}
Do minimal normal subgroups of the groups of homeomorphisms characterize 
chainable continua? 

\mynote{Notes}
Also (Brechner): Find a nice characterization of normal subgroups of the 
homeomorphism group of the pseudo-arc. 
\end{myprob}

\begin{myprob}
\myproblem{73}{Jones} 
What is the structure of the collection of homeomorphisms leaving a given point 
of the pseudo-arc fixed? 
\end{myprob}

\section*{Homogeneity}

\begin{myprob}
\myproblem{74}{Jones} 
Is every homogeneous, hereditarily indecomposable, nondegenerate continuum a 
pseudo-arc? 

\mynote{Notes}
Rogers has shown that it must be tree-like. 
\end{myprob}

\begin{myprob}
\myproblem{75}{Jones} 
Is each nondegenerate, homogeneous, nonseparating plane continuum a 
pseudo-arc? 

\mynote{Notes}
Jones and Hagopian have shown that it must be hereditarily indecomposable. 
Rogers has shown that it must be tree-like. 
\end{myprob}

\begin{myprob}
\myproblem{76}{Fitzpatrick, UHPB 88} 
Is every homogeneous continuum bihomogeneous? 

\mynote{Notes}
$X$ is \emph{bihomogeneous} if for every $x_0, x_1 \in X$ there exists a 
homeomorphism $h\colon X \to X$ with $h(x_i) = x_{1-i}$.

\mynote{Solution}
K. Kuperberg \cite{MR90m:54043.1} constructed a homogeneous continuum 
which is not bihomogeneous.
\end{myprob}

\begin{myprob}
\myproblem{77}{Jones} 
What effect does hereditary equivalence have on homogeneity in continua? 
\end{myprob}

\begin{myprob}
\myproblem{78}{Hagopian} 
If a homogeneous continuum $X$ contains an arc must it contain a solenoid 
or a simple closed curve? 
What if $X$ contains no simple triod? 

\mynote{Notes}
Ma\'ckowiak and Tymchatyn have shown that the answer is yes if $X$ is
atriodic. 
Connor presented a candidate for a counterexample. 
\end{myprob}

\begin{myprob}
\myproblem{79}{Rogers} 
Is each acyclic, homogeneous, one-dimensional continuum tree-like? 
hereditarily indecomposable? 

\mynote{Solution}
Rogers proved that acyclic, homogeneous curves are tree-like.
Krupski and Prajs showed that tree-like homogeneous continua are 
hereditarily indecomposable.
\end{myprob}

\begin{myprob}
\myproblem{80}{K. Kuperberg} 
Does there exist a homogeneous, arcwise-connected continuum which is not 
locally connected? 

\mynote{Solution}
Prajs \cite{MR2003f:54077.1} constructed a a homogeneous arcwise-connected 
curve which is not locally connected.
\end{myprob}

\begin{myprob}
\myproblem{81}{Minc} 
Is the simple closed curve the only nondegenerate, homogeneous, hereditarily 
decomposable continuum? 
\end{myprob}

\begin{myprob}
\myproblem{82}{Gordh} 
Is every hereditarily unicoherent, homogeneous, $T_2$ continuum indecomposable? 

\mynote{Notes}
Jones has shown that the metric version of this question has an affirmative 
answer. 
\end{myprob}

\begin{myprob}
\myproblem{83}{Ungar} 
Is every finite-dimensional, homeotopically homogeneous continuum a manifold? 

\mynote{Notes}
$X$ is homeotopically homogeneous if for each $x,y \in X$ there is a 
homeomorphism $h\colon (X,x) \to (X,y)$ and an isotopy connecting $h$ to 
the identity. 
\end{myprob}

\begin{myprob}
\myproblem{84}{Burgess} 
Is every $n$-homogeneous continuum ($n + 1$)-homogeneous for each 
$n \geq 2$? 

\mynote{Notes}
A continuum $X$ is \emph{$n$-homogeneous} if for each pair of collections 
$\{x_i : 1 \leq i \leq n\}$ and $\{y_i : 1 \leq i \leq n\}$ of $n$ 
distinct points of $X$ there is a homeomorphism $h\colon X \to X$ with 
$\{h(x_i) : 1 \leq i \leq n\} = \{y_i : 1 \leq i \leq n\}$.
Ungar has shown that if $X$ is $n$-homogeneous and $X \not= S^1$ then $h$ 
can also be chosen such that $h(x_i) = y_i$ for each $1 \leq i \leq n$. 
Kennedy has shown that if $X$ is $n$-homogeneous ($n \geq 2$) and admits a 
non-identity stable homeomorphism then $X$ is $m$-homogeneous for each 
positive integer $m$ (and in fact countable dense homogeneous and 
representable). 
Ungar has shown that if $X$ is $n$-homogeneous ($n \geq 2$) then $X$ is 
locally connected. 
\end{myprob}

\begin{myprob}
\myproblem{85}{Kennedy} 
Does every nondegenerate homogeneous continuum admit a non-identity stable 
homeomorphism? 
\end{myprob}

\begin{myprob}
\myproblem{86}{Bing} 
Is every homogeneous tree-like continuum hereditarily indecomposable? 

\mynote{Notes}
Jones and Hagopian have shown that in the plane the answer is yes. 
Jones has shown that such a continuum must be indecomposable. 
Hagopian has shown that it cannot contain an arc. 
Each of the following variants has been asked by various persons at various 
times. 
Is each such nondegenerate continuum 
a pseudo-arc? 
weakly chainable? 
hereditarily equivalent? 
of span zero? 
a continuum with the fixed-point property? 

Krupski has shown that if $X$ is a homogeneous continuum which contains a
local endpoint, then either $X$ is hereditarily indecomposable or $X$
admits a continuous decomposition into mutually homeomorphic,
nondegenerate, homogeneous, hereditarily indecomposable subcontinua with
decomposition space a homogeneous continuum with no local
endpoints. 
\end{myprob}

\begin{myprob}
\myproblem{87}{Rogers} 
Does every indecomposable, homogeneous continuum have dimension at most one? 
\end{myprob}

\begin{myprob}
\myproblem{88}{Rogers} 
Is each aposyndetic, non-locally-connected, one-dimensional, homogeneous 
continuum an inverse limit of Menger curves and continuous maps? 
Menger curves and fibrations? 
Menger curves and covering maps? 
Is each a Cantor set bundle over the Menger curve? 
\end{myprob}

\begin{myprob}
\myproblem{89}{Minc} 
Can each aposyndetic, non-locally connected, one-dimensional homogeneous 
continuum be mapped onto a solenoid? 

\mynote{Notes}
Rogers: Can such a continuum be retracted onto a nontrivial solenoid? 
Does each such continuum contain an arc?
\end{myprob}

\begin{myprob}
\myproblem{90}{Rogers} 
Is each pointed-$1$-movable, aposyndetic, homogeneous 
one-dim\-en\-sion\-al continuum locally connected? 
\end{myprob}

\begin{myprob}
\myproblem{91}{Rogers} 
Must each cyclic, indecomposable, homogeneous, one-dimensional continuum 
either be a solenoid or admit a continuous decomposition into tree-like, 
homogeneous continua with quotient space a solenoid? 
\end{myprob}

\begin{myprob}
\myproblem{92}{Rogers} 
Is every decomposable, homogeneous continuum of dimension greater than one 
aposyndetic? 
\end{myprob}

\begin{myprob}
\myproblem{93}{Rogers} 
Can the Jones Decomposition Theorem be strengthened to give decomposition 
elements which are hereditarily indecomposable? 
Can such a decomposition raise dimension? 
lower dimension? 

\mynote{Solution}
J. Rogers \cite{MR1992870.1} proved that if $X$ is a homogeneous,
decomposable continuum that is not aposyndetic and has dimension greater
than one, then the dimension of its aposyndetic decomposition is one.
\end{myprob}

\begin{myprob}
\myproblem{94}{}
Let $X$ be a nondegenerate, homogeneous, contractible continuum. 
Is $X$ an AR? 
Is $X$ homeomorphic to the Hilbert cube? 
\end{myprob}

\begin{myprob}
\myproblem{95}{Patkowska} 
What are the homogeneous Peano continua in $E^3$?
\end{myprob}

\begin{myprob}
\myproblem{96}{Patkowska} 
Does there exist a $2$-homogeneous continuum $X = X_1 \times X_2$ where $X_1$ 
and $X_2$ are nondegenerate, which is not either a manifold or an infinite 
product of manifolds? 
\end{myprob}

\begin{myprob}
\myproblem{97}{Bellamy} 
Is the following statement false? 
Statement: Suppose $X$ is a homogeneous compact connected $T_2$ space. 
Then for every open cover $U$ of $X$ there is an open cover $V$ of $X$ such 
that whenever $x$ and $y$ belong to the same element of $V$ there is a 
homeomorphism $h\colon X \to X$ such that $h(x) = y$ and such that for 
every $p \in X$, $p$ and $h(p)$ belong to the same element of $U$. 
\end{myprob}

\begin{myprob}
\myproblem{98}{Bellamy} 
If $X$ is an arcwise connected homogeneous continuum other than a simple 
closed curve, must each pair of points be the vertices of a 
$\theta$-curve in $X$? 

\mynote{Notes}
Bellamy and Lum have shown that each pair of points of $X$ must lie on a 
simple closed curve.
\end{myprob}

\begin{myprob}
\myproblem{99}{Bellamy} 
Does each finite subset of a nondegenerate arcwise connected homogeneous 
continuum lie on a simple closed curve? 
\end{myprob}

\begin{myprob}
\myproblem{100}{Bellamy} 
Does each nondegenerate arcwise connected homogeneous continuum other than the 
simple closed curve contain simple closed curves of arbitrarily small diameter? 
\end{myprob}

\begin{myprob}
\myproblem{101}{Wilson} 
Does there exist a uniquely arcwise connected homogeneous compact $T_2$ 
continuum, with an arc being defined either as a homeomorph of $[0,1]$ or as a 
compact $T_2$ continuum with exactly two nonseparating points? 

\mynote{Notes}
By a result of Bellamy and Lum, such a continuum cannot be metric. 
\end{myprob}

\begin{myprob}
\myproblem{102}{Lewis}
Does there exist a homogeneous one-dimensional continuum with no nondegenerate 
chainable subcontinua? 

\mynote{Notes}
If there exists a nondegenerate, homogeneous, hereditarily indecomposable 
continuum other than the pseudo-arc, the answer is yes.
\end{myprob}

\begin{myprob}
\myproblem{103}{Bennett} 
Is each open subset of a countable dense homogeneous continuum itself countable 
dense homogeneous? 

\mynote{Notes}
$M$ is countable dense homogeneous if for each two countable dense subsets 
$S$ and $T$ of $M$ there is a homeomorphism $h\colon M \to M$ with $h(S) = T$.
\end{myprob}

\begin{myprob}
\myproblem{104}{Fearnley} 
Is every continuum a continuous image of a homogeneous continuum? 
In particular, is the spiral around a triod such an image? 

\end{myprob}

\begin{myprob}
\myproblem{105}{J. Charatonik} 
Is the Sierpi\'nski curve homogeneous with respect to open surjections? 
\end{myprob}

\section*{Hyperspaces}

In each of the following, $X$ is a metric continuum, and $C(X)$ (resp.\ 
$2^X$) is the hyperspace of subcontinua (resp.\ closed subsets) of $X$ 
with the Hausdorff metric. 

\startproblem
\begin{myprob}
\myproblem{106}{Rogers} 
If $\dim X > 1$, is $\dim C(X) = \infty$? 
What if $X$ is indecomposable? 

\mynote{Notes}
Rogers raised the question and conjectured at the USL Mathematics 
Conference in 1971 that the answer is yes. 
The answer is known to be yes if any of the following are added to the 
hypothesis: 
$X$ is locally connected; 
$X$ contains the product of two nondegenerate continua; 
$\dim x > 2$;
$X$ is hereditarily indecomposable. 
\end{myprob}

\begin{myprob}
\myproblem{107}{Rogers} 
If $\dim X = 1$ and $X$ is planar and atriodic, is $\dim C(X) = 2$? 
Is $C(X)$ embeddable in $\mathbb{R}^3$? 

\mynote{Notes}
The answer is yes if $X$ is either hereditarily indecomposable or locally 
connected.
\end{myprob}

\begin{myprob}
\myproblem{108}{Rogers}
If $\dim X = I$ and $X$ is hereditarily decomposable and atriodic, is 
$\dim C(X) = 2$? 
\end{myprob}

\begin{myprob}
\myproblem{109}{Rogers} 
If $X$ is tree-like, does $C(X)$ have the fixed-point property? 
\end{myprob}

\begin{myprob}
\myproblem{110}{Nadler} 
When does $2^X$ have the fixed-point property?
\end{myprob}

\begin{myprob}
\myproblem{111}{Dilks} 
Is $C(X)$ or $2^X$ locally contractible at the point $X$? 

\mynote{Solution}
No. 
H. Kato \cite{MR91b:54014} constructed a chainable continuum $X$ such that 
$C(X)$ and $2^X$ are not locally contractible at $X$; and a dendroid $Y$ 
such that $C(Y)$ is locally contractible at $Y$ but $2^Y$ is not locally 
contractible at $Y$. 
A. Illanes \cite{MR91g:54014} constructed a continuum $X$ such that $2^X$, 
as well as $C(X)$, is not locally contractible at any of its points.
\end{myprob}

\begin{myprob}
\myproblem{112}{Rogers} 
Are any of the following Whitney properties: 
$\delta$-connected,
weakly chainable, or 
pointed-one-movable? 

\mynote{Notes}
Krasinkiewicz and Nadler have asked which of the following are Whitney 
properties: acyclic, ANR, AR, contractibility, Hilbert cube, homogeneity, 
$\lambda$-connected, $Sh(X) < Sh(Y)$, and weakly chainable. 
W. Charatonik has recently shown that homogeneity is not a Whitney property. 
\end{myprob}

\begin{myprob}
\myproblem{113}{Dilks and Rogers} 
Let $X$ be finite-dimensional and have the cone = hyperspace property. 
Must $X$ have property K? 
belong to class W? 
be Whitney stable? 
\end{myprob}

\section*{Inverse limits}

\begin{myprob}
\myproblem{114}{Young}
Is there for each $k \geq 1$ an atriodic tree-like continuum which is level 
$(k + 1)$ but not level $k$ 
(equivalently: Burgess' $(k + I)$-junctioned but not $k$-junctioned). 
What about the equivalent question for $(k + I)$-branched but not 
$k$-branched?
Find a useful way to characterize level $n$. 

\mynote{Notes}
A tree-like continuum $M$ is \emph{level $n$} if for every 
$\epsilon > 0$ there exists an $\epsilon$-map of $M$ onto a tree with $n$ 
points of order greater than two.
\end{myprob}

\begin{myprob}
\myproblem{115}{Young} 
Is there a continuum which is $4$-od like, not $T$-like, and every 
nondegenerate proper subcontinuum of which is an arc? 
\end{myprob}

\begin{myprob}
\myproblem{116}{}
Under what conditions is the inverse limit of dendroids a dendroid? 

\mynote{Notes}
A dendroid is an arcwise connected, hereditarily unicoherent continuum. 
\end{myprob}

\begin{myprob}
\myproblem{117}{Bellamy} 
Define $f_a\colon [0,1] \to [0,1]$ by $f(t) = at(1 - t)$ for $0 \leq a \leq 4$. 
Is there a relationship between the existence of periodic points of $f_a$ of 
various periods and the topological nature of the inverse limit continuum 
obtained by using $f_a$ as each one-step bonding map? 
In particular, is the inverse limit continuum indecomposable if and only 
if $f_a$ has a point of period $3$? 
\end{myprob}

\section*{Mapping properties}

\begin{myprob}
\myproblem{118}{W. Kuperberg, UHPB 31} 
Is it true that the pseudo-arc is not pseudo-contractible? 

\mynote{Notes}
A continuum $X$ is \emph{pseudo-contractible} if there exists a 
continuum $Y$, points $a,b \in Y$ and a map $h\colon X \times Y \to X$ 
such that $h_a\colon x \times \{a\} \to X$ is a homeomorphism and 
$h_b\colon x \times \{b\} \to X$ is a constant map. 
Also (W. Kuperberg, UHPB 29): 
Does there exist a one-dimensional continuum which is pseudo-contractible 
but not contractible? 
\end{myprob}

\begin{myprob}
\myproblem{119}{Ma\'ckowiak} 
Does there exist a chainable continuum $X$ such that if $H$ and $X$ are 
subcontinua of $X$ then the only maps between $H$ and $K$ are the identity or 
constants? 

\mynote{Notes}
Ma\'ckowiak has constructed a chainable continuum which admits only the 
identity or constants as self maps.

\mynote{Solution}
Ma\'ckowiak has constructed a nondegenerate chainable continuum
with the desired property. 
\end{myprob}

\begin{myprob}
\myproblem{120}{Lewis} 
Is every subcontinuum of a weakly chainable, atriodic, tree-like continuum 
weakly chainable? 
\end{myprob}

\begin{myprob}
\myproblem{121}{Lewis} 
If $P$ is the pseudo-arc and $X$ is a nondegenerate continuum, is 
$P \times X$ Galois if and only if $X$ is isotopy Galois? 

\mynote{Notes}
$X$ is \emph{Galois} if for each $x \in X$ and open $U$ containing $x$ 
there exists a homeomorphism $h\colon X \to X$ with $h(x) \not= x$ and 
$h(z) = z$ for each $z \in U$. 
If in addition $h$ can be chosen isotopic to the identity, each level of 
the isotopy satisfying $h(z) = z$ for each $z \not\in U$, then $X$ is 
\emph{isotopy Galois}. 
The parallel question for the Menger curve has a positive answer. 
\end{myprob}

\begin{myprob}
\myproblem{122}{Lewis} 
If $h$ is a homeomorphism of $\Pi_{\alpha \in A} P_\alpha$ where each 
$P_\alpha$ is a pseudo-arc, is $h$ necessarily of the form 
$h = \Pi_{\alpha \in A} h_{s(\alpha)}$, where $s$ is a permutation of $A$ 
and $h_s(\alpha)$ is a homeomorphism of $P_\alpha$ onto $P_{s(\alpha)}$? 

\mynote{Notes}
Bellamy and Lysko have given a positive answer when $A$ contains at most two 
elements. 
Cauty has shown the parallel question has a positive answer for any product of 
one-dimensional continua each open subset of which contains a simple closed 
curve (e.g., Menger curves or Sierpi\'nski curves). 

\mynote{Solution}
Bellamy provided a positive answer if $A$ is finite, and Bellamy and
Kennedy provided a positive answer for arbitrary $A$. 
\end{myprob}

\begin{myprob}
\myproblem{123}{Eberhart} 
If $X$ is a locally compact, metric space with every proper subcontinuum of 
$X$ hereditarily indecomposable, and $f$ is a local homeomorphism on $X$, 
is $f$ a homeomorphism on proper subcontinua of $X$?
\end{myprob}

\begin{myprob}
\myproblem{124}{Bellamy}
Conjecture: Let $X$ be a nondegenerate metric continuum, $p \in X$. 
Then there exist mappings $H\colon C \to C(X)$, ($C$ = Cantor set, 
$C(X)$ = hyperspace of subcontinua of $X$) and $h\colon C \to X$ such 
that $H$ and $h$ are embeddings and for each $x \in C$, $H(x)$ is 
irreducible from $p$ to $h(x)$ and if $x,y \in C$, $x < y$ (in ordering as 
a subset of [0,1]), then $H(x) \subsetneq H(Y)$. 
\end{myprob}

\begin{myprob}
\myproblem{125}{Minc} 
Suppose $X$ is a plane continuum such that for each $x,y \in X$ there is a 
weakly chainable subcontinuum of $X$ containing both $x$ and $y$. 
Is $X$ weakly chainable? 

\mynote{Notes}
Special case: 
Suppose $X$ is arcwise connected. 
The answer may be no if $X$ is non-planar.
\end{myprob}

\begin{myprob}
\myproblem{126}{Young} 
Suppose that $f$ is a light map of a tree $T_1$ onto a tree $T$ with the 
following property: 
Given light maps $g$, $h$ from the unit interval $I$ onto $T_1$. 
There exist maps $a, a: I \twoheadrightarrow I$ such that 
$fg\alpha = fh\beta$. 
Does $f$ factor through an arc? 
What if all maps are piecewise linear? 
\end{myprob}

\begin{myprob}
\myproblem{127}{Oversteegen}
Suppose $X$ is a weakly chainable, tree-like continuum. 
Do there exist inverse sequences 
$\varprojlim (I_n, g_n) \approx P$ ($P$ = pseudo-arc, $I$ = unit interval), 
$\varprojlim (T_n, f_n) \approx X$ (each $T_n$ a tree), and maps 
$h_n\colon I_n \to T_n$ such that 
$h = \varprojlim h_n = P \twoheadrightarrow X$?

\mynote{Notes}
Mioduszewski has shown that the answer is yes if $X$ is arc-like.
\end{myprob}

\begin{myprob}
\myproblem{128}{Oversteegen} 
Suppose $X$ is a continuum such that for each $x \in X$ there exists a 
neighborhood $U_x$ of $x$ such that $U_x \approx (0,1) \times A$ 
($A$ = compact, zero-dimensional set). 
Is $X$ not tree-like? 
\end{myprob}

\begin{myprob}
\myproblem{129}{Bellamy} 
Suppose $X$ is a non-pointed-one-movable continuum. 
Is there a non-pointed-one-movable continuum $K(X)$ which is either 
circle-like or figure-eight-like onto which $X$ can be mapped? 
\end{myprob}

\begin{myprob}
\myproblem{130}{Krasinkiewicz} 
Is there a finite-to-one map of an hereditarily indecomposable continuum onto 
an hereditarily decomposable continuum? 
\end{myprob}

\begin{myprob}
\myproblem{131}{Bellamy} 
For countable non-limit ordinals $\alpha$, what are the continuous images of 
$C(\alpha)$, the cone over a? 
For $\alpha \geq \omega_2 + 1$, what are the continuous pre-images of 
$C(\alpha)$? 

\mynote{Notes}
Katsuura has characterized the continuous images of the harmonic fan.
\end{myprob}

\begin{myprob}
\myproblem{132}{Bellamy} 
Is every continuous image of the cone over the Cantor set g-contractible? 

\mynote{Notes}
A continuum is \emph{g-contractible} if and only if it admits 
null-homotopic self surjection.
\end{myprob}

\begin{myprob}
\myproblem{133}{Bellamy} 
If an hereditarily indecomposable continuum admits an essential map onto a 
circle, does it admit map onto a pseudo-circle? 
\end{myprob}

\begin{myprob}
\myproblem{134}{Bellamy} 
Does every finite dimensional, hereditarily indecomposable continuum embed 
into a finite product of pseudo-arcs? 
\end{myprob}

\begin{myprob}
\myproblem{135}{Bellamy} 
Does every one dimensional hereditarily indecomposable continuum embed in 
a product of three (or maybe even two) pseudo-arcs? 
\end{myprob}

\begin{myprob}
\myproblem{136}{Bellamy} 
Does every tree-like hereditarily indecomposable continuum embed into a 
product of two (or three) pseudo-arcs? 
Does every planar hereditarily indecomposable continuum embed in a product 
of two pseudo-arcs? 
\end{myprob}

\begin{myprob}
\myproblem{137}{Bellamy} 
Is the pseudo-circle-a retract of every one-dimensional hereditarily 
indecomposable continuum containing it? 
\end{myprob}

\section*{Problems in the plane}

\begin{myprob}
\myproblem{138}{Lewis} 
Does every hereditarily indecomposable plane continuum have 
$\mathfrak{c} = 2^{\omega_0}$ distinct embeddings in $E^2$. 
Does each such continuum have, for each integer $n > 1$, an embedding with 
exactly $n$ accessible composants? 
Does every such continuum have an embedding with no two accessible points 
in the same composant? 
\end{myprob}

\begin{myprob}
\myproblem{139}{Burgess} 
Which continua in $E^2$ have the property that all of their embeddings in 
$E^2$ are equivalent? 
\end{myprob}

\begin{myprob}
\myproblem{140}{Nadler and Quinn} 
If $p$ is a point of the chainable continuum $M$, is there an embedding of 
$M$ in $E^2$ which makes $p$ accessible? 
\end{myprob}

\begin{myprob}
\myproblem{141}{Mayer}
Are there uncountably many inequivalent embeddings of every chainable 
indecomposable continuum in $E^2$? 
\end{myprob}

\begin{myprob}
\myproblem{142}{Mayer} 
Can every chainable indecomposable continuum be embedded in $E^2$ 
non-principally (i.e., without a simple dense canal)? 

\mynote{Notes}
This is known for such continua with at least one endpoint. 
\end{myprob}

\begin{myprob}
\myproblem{143}{Brechner and Mayer} 
Does there exist a nonseparating plane continuum such that every embedding 
of it in $E^2$ has a simple dense canal? 
\end{myprob}

\begin{myprob}
\myproblem{144}{Ancel} 
Is every embedding of a Peano continuum in $\mathbb{R}^2$ micro-unknotted? 
Is the standard inclusion $S^3 \to S^4$ micro-unknotted? 

\mynote{Notes}
Suppose $M$ and $N$ are compact, metric spaces, $G$ is the homeomorphism 
group of $N$, and $X$ is the space of embeddings of $M$ in $N$. 
An embedding $e\colon M \to P$ is \emph{micro-unknotted} if for each 
$\epsilon > 0$ there exists $\delta > 0$ such that if $h \in G$ and 
$\operatorname{dist}_X(e, h\circ e) < \delta$, then there exists 
$h' \in G$ with $\operatorname{dist}_G (1_N, h') < \epsilon$ 
and $h' \circ e = h \circ e$. 
$e\colon M \to N$ is micro-unknotted iff 
acts micro-transitively on the orbit $G \circ e$ iff 
$G \circ e$ is $G_\delta$ in $X$ (Effros' theorem). 
\end{myprob}

\begin{myprob}
\myproblem{145}{Jones} 
What characterizes dendroids that are embeddable in $E^2$? 
What characterizes dendroids that are contractible? 
\end{myprob}

\begin{myprob}
\myproblem{146}{Ancel} 
Is there a recognizable family of nonseparating plane continua such that every 
nonseparating plane continuum is a retract of a member of this family? 
\end{myprob}

\begin{myprob}
\myproblem{147}{Bellamy} 
When is the inverse image $S$ of an indecomposable plane continuum $X$ under 
a complex power map ($f(z) = zn$ for some $n$) itself an indecomposable 
continuum? 
In particular, if $0$ lies in an inaccessible composant of $X$, is $S$ 
indecomposable? 
\end{myprob}

\begin{myprob}
\myproblem{148}{}
Suppose $M$ is a nondegenerate connected subset of $E^2$, such that the 
complement of each point in $M$ is connected but the complement of each pair 
of points in $M$ is disconnected. 
Can $E^2 - M$ be arcwise connected? 
\end{myprob}

\section*[Set function $T$]{Set function $\mathbf{T}$}

Let $S$ be a compact Hausdorff space, and let $A$ be a subset of $S$. 
$T(A)$ is the set of points which have no continuum neighborhood missing $A$. 
$K(A)$ is the intersection of all continuum neighborhoods of $A$. 
The following problems are unsolved for compact Hausdorff continua, with the 
possible exception of number 157. 
Except for number 158, they are unsolved for compact metric continua. 
The phrase `$T$ is continuous for $S$' means that $T$ is continuous 
considered as a function from the hyperspace of closed subsets of $S$ to 
itself; similarly for $K$. 
`S is $T$-additive' means that for closed sets $A, B \subseteq S$, 
$T(A \cup B) = T(A) \cup T(B)$. 
All questions in this section were posed by Bellamy unless indicated 
otherwise. 

\startproblem
\begin{myprob}
\myproblem{149}{}
If $T$ is continuous for $S$, is $K$ also continuous for $S$?
\end{myprob}

\begin{myprob}
\myproblem{150}{}
If $T$ is continuous for $S$ and $S$ is decomposable, is it true that for each 
$p \in S$, $\operatorname{Int}(T(p)) = \emptyset$? 

\end{myprob}

\begin{myprob}
\myproblem{151}{}
If $T$ is continuous for $S$, is $S$ $T$-additive? 

\mynote{Notes}
Bellamy has offered a prize for the solution of this question---one bushel of 
extra fancy Stayman Winesap apples, delivered in season.
\end{myprob}

\begin{myprob}
\myproblem{152}{}
If $S/T$ denotes the finest decomposition space of $S$ which shrinks each 
$T(p)$ to a point, is $S/T$ locally connected? 

\mynote{Notes}
This is not difficult to show if $S$ is also $T$ additive.
\end{myprob}

\begin{myprob}
\myproblem{153}{Jones} 
If $X$ is indecomposable and $W$ is a subcontinuum of $X \times X$ with 
nonempty interior, is $T(W) = X \times X$? 
\end{myprob}

\begin{myprob}
\myproblem{154}{Cook} 
If $X$ is atriodic (or contains no uncountable collection of pairwise disjoint 
triods) and $X$ has no continuum cut point, does this imply that there is a 
continuum $W \subset X$ such that $\operatorname{Int}(W) \not= \emptyset$ 
and $T(W) \not= X$? 
\end{myprob}

\begin{myprob}
\myproblem{155}{}
If $T$ is continuous for $S$ and $f\colon S \to Z$ is a continuous and 
monotone surjection, is $T$ continuous for $Z$ also? 
\end{myprob}

\begin{myprob}
\myproblem{156}{}
If $X$ is one-dimensional and homogeneous is $T$ continuous for $S$? 
\end{myprob}

\begin{myprob}
\myproblem{157}{}
Call a continuum $S$ \emph{strictly point $T$ asymmetric} if for $p \not= q$ 
and $p \in T(q)$ we have $q \not\in T(p)$. 
In dendroids, does this property imply smoothness? 
\end{myprob}

\begin{myprob}
\myproblem{158}{H. Davis and Doyle} 
If $S$ is almost connected im kleinen, is $S$ connected im kleinen at some 
point? 

\mynote{Notes}
Almost connectedness im kleinen can be expressed in terms of the set function 
$T$ as follows: 
$S$ is \emph{almost connected im kleinen} at $p \in S$ if and only if for 
each closed $A$ for which $p \in \operatorname{Int}(T(A))$ we have 
$p \in \operatorname{Int}(A)$. 
This question is known to be true for the metric case. 
\end{myprob}

\begin{myprob}
\myproblem{159}{}
Suppose the restriction of $T$ to the hyperspace of subcontinua of $S$ is 
continuous. 
Does this imply that $T$ is continuous for $S$? 

\mynote{Notes}
This is true if $T$ is the identity on subcontinua.
\end{myprob}

\begin{myprob}
\myproblem{160}{}
Do open maps preserve $T$-additivity? 
$T$-symmetry? 

$S$ is \emph{$T$-symmetric} if and only if for all closed sets $A$ and 
$B$ in $S$, if $A \cap T(B) = \emptyset$ then $B \cap T(A) = \emptyset$.
\end{myprob}

\section*{Span}

\begin{myprob}
\myproblem{161}{Lelek, Cook, UHPB 81} 
Is each continuum of span zero chainable? 
\end{myprob}

\begin{myprob}
\myproblem{162}{Duda} 
To what extent does span zero parallel chainability?
\begin{myenumerate}
\item
Is the open image of a continuum of span zero a continuum of span zero? 
\item
(Lelek, UHPB 84) 
Is the confluent image of a chainable continuum chainable? 
\item
(Lelek, Cook, UHPB 86)
Do confluent maps of continua preserve span zero? 
\end{myenumerate}

\mynote{Notes}
Also (Lelek, UHPB 85): If $f$ is a confluent mapping of an acyclic 
(or tree-like or arc-like) continuum $X$ onto a continuum $Y$, is 
$f \times f$ confluent? 
An affirmative solution to (2) would show that the classification of 
homogeneous plane continua is complete. 
McLean has shown that the confluent image of a tree-like continuum is 
tree-like, and Rosenholtz has shown that the open image of a chainable 
continuum is chainable.

\mynote{Solution}
K. Kawamura \cite{MR89j:54037} proved that (1) has an affirmative answer.
\end{myprob}

\begin{myprob}
\myproblem{163}{Cook, UHPB 92} 
If $M$ is a continuum with positive span such that each of its proper 
subcontinua has span zero, does every nondegenerate, monotone, continuous 
image of $M$ have positive span? 
\end{myprob}

\begin{myprob}
\myproblem{164}{Cook, UHPB 173} 
Do there exist, in the plane, two simple closed curves $J$ and $K$ such that 
$X$ is in the bounded complementary domain of $J$, and the span of $K$ is 
greater than the span of $J$?
\end{myprob}

\begin{myprob}
\myproblem{165}{Bellamy}
Suppose $X$ is a homogeneous, aposyndetic continuum which contains two
disjoint subcontinua with interior.
Is $X$ mutually aposyndetic?
What if $X$ is also arcwise connected?
\end{myprob}

\begin{myprob}
\myproblem{166}{Bula}
Suppose $F\colon X \to Y$ is an open map, with each of $X$ and $Y$ 
compact metric and each $F^{-1}(y)$ infinite.
Do there exist disjoint closed subsets $F$ and $H$ of $X$ such that
$f(H) = f(H) = Y$?

\mynote{Notes}
It is known that if each point inverse is perfect and $Y$ is
finite-dim\-en\-sion\-al then there exists a continuous surjection 
$g\colon X \times Y \to [0,1]$ such that $f = \pi_Y \circ g$, where 
$\pi_Y$ is the projection of $Y \times [0,1]$ onto $Y$.
\end{myprob}

\begin{myprob}
\myproblem{167}{Lewis}
Under what conditions does there exist a wild embedding of the $k$-sphere
$S^k$ in $E^k$ which is a homogeneous embedding?

\mynote{Notes}
Compare with questions 49 and 50. 
\end{myprob}

\begin{myprob}
\myproblem{168}{Lewis} 
Does there ever exist a wild embedding of $S^k$ in $E^n$ which is 
isotopically homogeneous? 

\mynote{Notes}
Compare with question 52. 
\end{myprob}

\begin{myprob}
\myproblem{169}{Lewis} 
Does there exist a nondegenerate continuum $K$ which can be embedded in
$E^n$, $n \geq 3$, such that every embedding of $K$ in $E^n$ is a
homogeneous embedding? 
\end{myprob}

\begin{myprob}
\myproblem{170}{Minc} 
Does there exist an hereditarily indecomposable continuum which is
homogeneous with respect to continuous surjections but not homogeneous
with respect to homeomorphisms?

\mynote{Notes}
The pseudo-circle and pseudo-solenoids are known not to have this 
property. 
\end{myprob}

\begin{myprob}
\myproblem{171}{Bellamy} 
Does there exist an hereditarily indecomposable nonmetric continuum with
only one composant? 

\mynote{Notes}
D. Bellamy and Smith have independently constructed indecomposable,
nonmetric continua with only one or two composants. 
Smith has constructed an hereditarily indecomposable, nonmetric
continuum with only two composants. 
\end{myprob}

\begin{myprob}
\myproblem{172}{Van Nall} 
Is it true that an atriodic continuum in class W is hereditarily in
class W if and only if each $C$-set in it is in class W? 
\end{myprob}

\providecommand{\bysame}{\leavevmode\hbox to3em{\hrulefill}\thinspace}

\label{tplewisend}

\chapter*{Janusz R.\ Prajs: Problems in continuum theory}
\label{tpprajs}
\begin{myfoot}
\begin{myfooter}
Janusz R.\ Prajs, 
\emph{Problems in continuum theory},\\
Problems from Topology Proceedings, Topology Atlas, 2003, 
pp.\ 183--189. 
\end{myfooter}
\end{myfoot} 

\mypreface
The material in this section is taken from the article \emph{Several 
old and new problems in continuum theory} \cite{MR2003g:54079} 
by J.J.~Charatonik and J.R.~Prajs and from the website 
\emph{Open problems in continuum theory} edited by 
J.R.~Prajs \cite{prajs}.
Three short essays were contributed by J.J.~Charatonik and 
C.L.~Hagopian and are included below.

\section*{Introduction}

Properties of continua (i.e., compact connected Hausdorff spaces) have
been concentrating much attention since the very beginning of topology
studies. Now, when foundations of general topology are already
established, a great number of natural questions about continua remain
open. Many of them are easy to formulate and understand even for
beginners. Nevertheless, they turned out to be difficult and they are
still a great challenge and inspiration to current research. Below we
present a sample of these questions. For other collections of continuum
theory problems see historically the first such set \cite{MR96f:54042},
and also \cite{MR1078655.1}, \cite{MR86a:54038.2} and \cite{TP9.2.375.1}.

The presented questions are divided into two parts. 
First, we list some old and well known problems that should be reminded 
whenever important questions in topology are discussed. 
Second, we recall twelve newer questions that are connected with authors' 
recent research. 
All problems presented below concern metric spaces only. 
All mappings are assumed to be continuous.

\section*{Classic problems}

\begin{myprob}
\myproblem{Fixed point problem for nonseparating plane continua}{}
Does every nonseparating plane (tree-like) continuum have the fixed-point 
property?

\mynote{Notes}
A space $X$ is said to have the \emph{fixed-point property} provided 
that for every continuous function $f\colon X \to X$ there exists a 
point $p$ in $X$ such that $f(p)=p$.
For more information see the survey paper \cite{MR92i:54033.1} by 
Charles L.~Hagopian.

See also the short survey below about this problem by C.L.~Hagopian.
\end{myprob}

\begin{myprob}
\myproblem{Hereditary equivalence}{}
Assume that a nondegenerate continuum $X$ is hom\-eomorphic to each of 
its proper nondegenerate subcontinua. 
Must then $X$ be either an arc or a pseudo-arc? 

\mynote{Notes}
Such continua $X$ are named \emph{hereditarily equivalent}. 
As early as 1921 S. Mazurkiewicz posed a question as to whether every 
hereditarily equivalent continuum is an arc, \cite{mazurkiewicz19212}. 
In 1948 E.E.~Moise constructed the pseudo-arc which is hereditarily 
equivalent and hereditarily indecomposable, \cite{MR10:56i2}, and thus 
answered Mazurkiewicz's question in the negative. 
Later G.W.~Henderson showed that a hereditarily equivalent decomposable 
continuum is an arc, \cite{MR22:9949}. 
H.~Cook proved that a hereditarily equivalent continuum is tree-like, 
\cite{MR42:1072.1}. 
Compare \cite[Section~2, p.~307]{MR1078656.1}. 
\end{myprob}

\begin{myprob}
\myproblem{Homogeneous tree-like continua}{}
Is each nondegenerate homogeneous tree-like (planar, weakly chainable) 
continuum a pseudo-arc? 

\mynote{Notes}
Research directed to classify homogeneous continua was initiated by the 
question of B.~Knaster and K.~Kuratowski in 1920, \cite{knastekuratowksi}, 
whether the simple closed curve is the only homogeneous nondegenerate 
plane continuum. 
A continuum $X$ is said to be \emph{homogeneous} provided that for 
every two points $x$ and $y$ of $X$ there exists a homeomorphism 
$h\colon X \to X$ such that $h(x)= y$. 

A \emph{weakly chainable} continuum is meant a continuous image of 
the pseudo-arc. 
J.T.~Rogers, Jr., proved in \cite{MR84b:54072.1} that a hereditarily 
indecomposable homogeneous continuum is tree-like. 
Answering an old question of R.H.~Bing, the second named author showed 
(the proof is presented in the joint paper \cite{MR90f:54054.1}) that 
tree-like homogeneous continua are hereditarily indecomposable. 
A positive answer to any of these questions would finally classify, 
after eight decades of study, all nondegenerate homogeneous plane continua 
as: the circle, the pseudo-arc and the circle of pseudo-arcs. 
For more detailed information on classifications of homogeneous continua, 
see \cite[Chapter~8]{MR2002d:54017}, \cite{MR93d:54049} and 
\cite{MR94c:54061}. 
For the definition of the pseudo-arc and for more information about this 
continuum see \cite{MR2000f:54029}.
\end{myprob}

\begin{myprob}
\myproblem{Homogeneous indecomposable continua}{}
Is each nondegenerate homogeneous indecomposable (cell-like) continuum 
one-dimensional? 

\mynote{Notes}
This question was asked by James T.~Rogers, Jr. 
The pseudo-arc, solenoids and solenoids of pseudo-arcs are the only known 
nondegenerate homogeneous indecomposable continua, and all they are 
one-dimensional. 
If the answer to any of these questions is yes, then an essential progress 
in the study of the structure of homogeneous higher dimensional continua
would be obtained, namely the completely regular decompositions described
in \cite{MR17:180e}, \cite{MR84e:54040} and 
\cite[Theorem~7.1, p.~18]{MR87a:54048} would be trivial (in particular 
such continua would be aposyndetic and they would contain no proper 
nondegenerate terminal subcontinua). 
On the other hand an example of a higher dimensional homogeneous 
indecomposable continuum would be of a great importance in this area.

In the nonmetric case the answer is negative (J.~van~Mill, 
\cite{MR91k:54063}).
\end{myprob}

\begin{myprob}
\myproblem{Confluent image of arc-like continua}{}
Is a confluent image of an arc-like continuum (of a pseudo-arc) 
necessarily arc-like? 

\mynote{Notes}
It is known that a positive answer to this question implies that every
nondegenerate, planar, homogeneous, tree-like continuum is a pseudo-arc.
This question was raised by A.~Lelek in 
\cite[Problem~4, p.~94]{MR46:6324}.
\end{myprob}

\begin{myprob}
\myproblem{Property of Kelley}{}
Assume that a continuum $X$ has the property of Kelley. 
Does the product $X \times [0,1]$ necessarily have this property? 

\mynote{Notes}
A continuum $X$ is said to have the \emph{property of Kelley} provided
that for each point $x \in X$, for each sequence of points $x_n \in X$
converging to $x$ and for each continuum $K$ such that $x \in K \subset X$
there exists a sequence of continua $K_n \subset X$ such that $x_n \in
K_n$ and $\lim K_n = K$. The property is a one of the most extensively
studied and useful in continuum theory. All hereditarily indecomposable,
all (openly) homogeneous continua, all locally connected continua and all
absolute retracts for hereditarily unicoherent continua have this property
(see \cite[pp.~167--175, 277--279 and 405--406]{MR99m:54006.1}; 
\cite{MR86a:54039} and \cite[Corollary~3.7]{ccpapk}).

The recalled problem arose from the original question of S.B.~Nadler, 
Jr., \cite[16.37, p.~558]{MR58:18330}, whether the property of Kelley of 
a continuum $X$ implies the property of the hyperspace $C(X)$ of all
nonempty subcontinua of $X$ with the Hausdorff metric. 
In \cite[Corollary~3.3, p.~1147]{MR91j:54013}, H.~Kato proved that 
Nadler's question is equivalent to the considered problem. 
Since Kato's variant of the problem is more intuitive for 
non-specialists, we have chosen it here.
\end{myprob}

\begin{myprob}
\myproblem{Dendroids and small retractions onto dendrites}{}
Let $X$ be a dendroid. 
Do there exist, for each $\varepsilon > 0$, a tree $T \subset X$ and a 
retraction $r\colon X \to T$ with $d(x, r(x)) < \varepsilon$ for each 
point $x \in X$? 

\mynote{Notes}
See the short essay below about this problem by Janusz J. Charatonik.
\end{myprob}

\begin{myprob}
\myproblem{Span $\mathbf{0}$ vs. arc-like}{}
Let $X$ be a continuum with span $0$. 
Must $X$ be arc-like?

\mynote{Notes}
For any two maps $f,g\colon Z \to Y$, where $Y$ is a metric space, define
$m(f,g)=\inf\{d(f(z),g(z)) \mid z\in Z\}$. 
For any continuum $X$ the number
$$
\sigma(X) = 
\sup\{ m(f,g) \mid f,g\colon Z\to X,\ \text{$Z$is a continuum, and}\ 
f(Z)\subset g(Z)\}
$$
is called the \emph{span} of $X$. 
Note that $\sigma(X)=0$ is a topological property of a continuum $X$. 
The concept of the span of a continuum is due to A.~Lelek. 
The above question was posed by A.~Lelek in \cite{MR46:6324}.
\end{myprob}

\begin{myprob}
\myproblem{Homogeneous $\mathbf{n}$-dimensional ANRs}{}
Let $X$ be a homogeneous, $n$-dimen\-sional continuum. 
If $X$ is an absolute neighborhood retract (ANR), must $X$ be an 
$n$-manifold? 

\mynote{Notes}
This question is due to R.H.~Bing and K.~Borsuk.
A positive answer to this question was given by Bing and Borsuk for 
$n < 3$.
\end{myprob}

\section*{Some new questions}

The next three problems below are related to each other. 
They deal with a more general question: 
Given continua $X$ and $Y$, does there exist a continuous surjection of 
$X$ onto $Y$?
 
Among initial famous results in this area there is the construction of a 
continuous surjection of $[0,1]$ onto $[0,1]^2$ by G.~Peano and its 
generalization, the Hahn-Mazurkiewicz theorem saying that each locally 
connected continuum is a continuous image of $[0,1]$. 

In this area we study invariants and inverse invariants of continuity for
continua (sometimes called generalized continuous invariants). The study
of generalized continuous invariants (e.g., local connectedness, uniform
pathwise connectedness, various types of so called ``indices of local
disconnectivity'', see e.g., \cite{prajsinv}, \cite{MR1347229},
\cite{MR2003b:54017}, \cite{MR91g:54047}, and compare also $\sigma$-local
connectedness in \cite{MR2000e:54023}), did not allow yet to exclude the
existence of continuous surjections questioned in the next three problems

\startproblem
\begin{myprob}
\myproblem{Mappings onto hyperspaces of subcontinua}{}
Does there exist a continuum $X$ admitting no continuous surjection onto 
its hyperspace $C(X)$ of all nonempty subcontinua?

\mynote{Notes}
Originally, a related problem was considered by S.B.~Nadler, Jr. in
\cite[Question~4.6, p.~243]{MR58:18330}. No tools are known to prove
non-existence of a continuous surjection from any continuum $X$ onto
$C(X)$. On the other hand, no natural tools promising to construct such
mappings for all continua are developed either. A (possible) continuum $X$
with no such mapping must be non-locally connected, and each of its open
subsets must have countably many components only, see a remark in
\cite[Question~4.6, p.~243]{MR58:18330}.
\end{myprob}

\begin{myprob}
\myproblem{Mappings between hyperspaces of subcontinua}{}
Assume that there exists a continuous surjection $f\colon X \to Y$ 
between continua $X$ and $Y$. Does there exist a continuous surjection 
$g\colon C(X) \to C(Y)$ between their hyperspaces $C(X)$ and $C(Y)$?

\mynote{Notes}
If the mapping $f$ is weakly confluent, then the induced mapping 
$A \mapsto f(A)$ between $C(X)$ and $C(Y)$ is surjective, 
\cite[Theorem~0.49.1, p.~24]{MR58:18330}. 
However, there are pairs of continua $X$ and $Y$ admitting a continuous 
surjection $f$ and such that there is no weakly confluent mapping from $X$ 
onto $Y$.
\end{myprob}

\begin{myprob}
\myproblem{Mappings between Cartesian squares}{}
Does there exist a pair of continua $X$ and $Y$ with a continuous 
surjection $f\colon X^2 \to Y^2$ that admits no continuous surjection 
from $X$ onto $Y$?

\mynote{Notes}
An example of such a pair for locally compact, noncompact, connected 
spaces was found by M.~Morayne (an oral communication).
\end{myprob}

\begin{myprob}
\myproblem{Tree-likeness of absolute retracts}{}
Is every absolute retract $X$ for the class of all hereditarily unicoherent 
continua a tree-like continuum?

\mynote{Notes}
In the recent paper \cite{ccpar} an extensive study of absolute retracts
for hereditarily unicoherent continua was presented. 
This problem and the next seem to be the most important among those that 
arose from this research.

Such a continuum $X$ has the property of Kelley, and each of its arc 
components is dense in $X$ (in particular $X$ is approximated from within 
by trees). 
Proofs of these properties, together with many other ones, are presented 
in \cite{ccpar}.
\end{myprob}

\begin{myprob}
\myproblem{Absolute retracts and inverse limits}{}
Does there exist an absolute retract $X$ for tree-like continua such that 
$X$ cannot be represented as an inverse limit of trees with confluent 
bonding mappings? 

\mynote{Notes}
The arc-like continuum having exactly three end points as constructed in
\cite[1.10, p.~7, and Figure~1.10, p.~8]{MR93m:54002.1} is our candidate 
for such a continuum $X$. It is proved in \cite[Theorem~3.6]{ccpar} that 
the inverse limit of trees with confluent bonding mappings is an absolute 
retract for hereditarily unicoherent continua.

\mynote{Solution}
Recently, W.J.~Charatonik and J.R.~Prajs found examples of absolute 
retracts for hereditarily unicoherent continua that cannot be represented 
as the inverse limit of trees with confluent bonding mappings. These 
examples are dendroids and thus they are tree-like. 
Thus the above question is answered in the positive.
\end{myprob}

\begin{myprob}
\myproblem{Continuous decomposition of a $\mathbf{3}$-book}{}
Let $T$ be a simple triod. Does there exist a continuous decomposition of 
the product $T \times [0,1]$ into pseudo-arcs?

\mynote{Notes}
For motivation of studying continuous decompositions into pseudo-arcs see
the introduction of \cite{MR2000b:54043}. In \cite{MR58:2750} and in the 
recent papers \cite{MR2000b:54043} and \cite{MR2000k:54009} it was shown 
that the plane and each locally connected continuum in a $2$-manifold with 
no local separating point, as well as the Menger curve, admit a continuous 
decomposition into pseudo-arcs (compare also \cite{MR97h:54008.1} and 
\cite{MR96a:54006}). Among Peano continua local separating point is the 
only known true obstacle to construct such a decomposition, 
\cite[Proposition~15, p.~34]{MR2000b:54043}. The methods developed in
the above quoted papers cannot be directly extended to the $3$-book case.
\end{myprob}

\begin{myprob}
\myproblem{Homogeneous Peano continua in the $\mathbf{3}$-space}{}
Does there exist a homogeneous locally connected $2$-dimensional 
continuum in the Euclidean $3$-space that is neither a surface nor the 
Pontryagin sphere? 

\mynote{Notes}
We can define the Pontryagin sphere as the quotient space of the standard
Sierpi\'nski universal plane curve $S$ in $[0,1] \times [0,1]$. Namely we
identify each pair of points belonging to the boundary of one component of
$\mathbb{R}^2 \setminus S$ having either $x$-coordinates or 
$y$-coordinates equal. 
The Pontryagin sphere can also be seen as the quotient space of the 
disjoint union of two Pontryagin discs $\mathbb{D}^2$ (see 
\cite[Section~3, pp.~608--609]{MR93f:57024}) with each pair of the 
corresponding points in the boundary $\partial \mathbb{D}^2$ identified.

S.~Mazurkiewicz had shown that the only nondegenerate locally connected 
homogeneous plane continuum is the simple closed curve, 
\cite{mazurkiewicz1924}. 
Locally connected $1$-dimensional homogeneous continua are characterized 
as the simple closed curve and the Menger universal curve (see e.g., 
\cite[12.2, p.~29]{MR87a:54048}). 
Therefore, a negative answer to this question would provide a complete 
classification of locally connected homogeneous continua in $3$-space. 
A continuum in question could not contain a $2$-cell, see 
\cite{MR90f:54055}, and it would not be an ANR, see 
\cite[Theorem~16.10, p.~194]{MR35:7306}.
\end{myprob}

\begin{myprob}
\myproblem{Continuous decomposition of the plane}{}
Let $X$ be a nondegenerate continuum such that the plane admits a 
continuous decomposition into topological copies of $X$.
Must then $X$ be hereditarily indecomposable? 
Must $X$ be the pseudo-arc? 

\mynote{Notes}
The existence of a continuous decomposition of the plane into pseudo-arcs 
was announced by R.D.~Anderson in 1950. 
The first known proof of this fact, given by W.~Lewis and J.~Walsh, 
appeared in 1978, \cite{MR58:2750}. 
\end{myprob}

\begin{myprob}
\myproblem{Simply connected, homogeneous continua in $\mathbf{R^3}$}{}
Let $X$ be a simply connected, nondegenerate, homogeneous continuum in 
the $3$-space $\mathbb{R}^3$. 
Must $X$ be homeomorphic to the unit sphere $S^2$?

\mynote{Notes}
A continuum $X$ is called \emph{simply connected} provided that $X$ 
is arcwise connected and every map from the unit circle $S^1$ into $X$ is 
null-homotopic.
If $X$ either is an ANR, or topologically contains a $2$-dimensional 
disk, then the answer is positive.
\end{myprob}

\begin{myprob}
\myproblem{Local connectedness of simply connected homogeneous continua}{}
Let $X$ be a simply connected, homogeneous continuum. 
Must $X$ be locally connected?

\mynote{Notes}
This question is related to a question by K.~Kuperberg whether an arcwise
connected, homogeneous continuum must be locally connected. 
This last question was recently answered in the negative by J.R.~Prajs.
\end{myprob}

\begin{myprob}
\myproblem{Disks in simply connected homogeneous continua}{}
Let $X$ be a homogeneous, simply connected (locally connected) 
nondegenerate continuum.
Must $X$ contain a $2$-dimensional disk?

\mynote{Notes}
This question is due to Panagiotis Papazoglou.
\end{myprob}

\begin{myprob}
\myproblem{Path connectedness of homogeneous continua}{}
Let $X$ be an arcwise connected, homogeneous continuum. 
Must $X$ be uniformly path connected? 
(Equivalently, is $X$ a continuous image of the Cantor fan?)

\mynote{Notes}
A continuum $X$ is called \emph{uniformly path connected} provided that
there is a compact collection $P$ of paths in $X$ such that each pair of
points $x$, $y$ in $X$ is connected by some member of $P$. 
The \emph{Cantor fan} is defined as the cone over the Cantor set. 
It is known that a homogeneous arcwise connected continuum need not be 
locally connected \cite{MR2003f:54077.2}. 
The strongest result in the direction of this question has been obtained 
by D.P.~Bellamy, \cite{MR88h:54048}. 
See also \cite{MR86m:54047} and \cite{MR82j:54071}.
\end{myprob}

\providecommand{\bysame}{\leavevmode\hbox to3em{\hrulefill}\thinspace}

\label{tpprajsend}

\chapter*{Charles L.~Hagopian: The plane fixed-point problem}
\label{tphagopian}
\begin{myfoot}
\begin{myfooter}
Charles L.~Hagopian, 
\emph{The plane fixed-point problem},\\ 
Problems from Topology Proceedings, Topology Atlas, 2003, 
pp.\ 191--193. 
\end{myfooter}
\end{myfoot}

Does every nonseparating plane continuum have the fixed-point property? 
This is the plane fixed-point problem. 
It has been called the most interesting outstanding problem in plane 
topology \cite{MR38:5201}. 
A positive answer would provide a natural generalization to the 
$2$-dimensional version of the Brouwer fixed-point theorem.

A space $S$ has \emph{the fixed-point property} if for every map 
(continuous function) $f$ of $S$ into $S$ there exists a point $x$ of $S$ 
such that $f(x) = x$.
A \emph{continuum} is a nondegenerate compact connected metric space.
A continuum in the plane that has only one complementary domain is a 
\emph{nonseparating plane continuum}.
Every nonseparating plane continuum is the intersection of a nested 
sequence of topological disks.

To summarize related results, suppose $\mathcal{C}$ is a nonseparating 
plane continuum and $f$ is a fixed-point-free map of $\mathcal{C}$ into 
$\mathcal{C}$.
Ayres \cite{ayres} in 1930 proved $\mathcal{C}$ is not locally 
connected if $f$ is a homeomorphism. 
In 1932 Borsuk \cite{borsuk1935} proved $\mathcal{C}$ cannot be locally 
connected (even if $f$ is not a homeomorphism). 
He accomplished this by showing that every locally connected 
nonseparating plane continuum is a retract of a disk. 
Stallings and Borsuk \cite{MR22:8485} pointed out that the plane 
fixed-point problem would be solved if it could be shown that every 
nonseparating plane continuum is an almost continuous retract of a disk. 
This approach was eliminated by Akis in \cite{MR86b:54044}.

Hamilton \cite{MR1501958} in 1938 proved the boundary of $\mathcal{C}$ is 
not hereditarily decomposable if $f$ is a homeomorphism. 
Bell \cite{MR35:4888}, Sieklucki \cite{MR39:2139}, and Iliadis 
\cite{MR44:4726} in 1967--1970 independently proved the boundary of 
$\mathcal{C}$ contains an indecomposable continuum that is left invariant 
by $f$. 
The methods used to establish this theorem led to (but did not answer) 
the following questions.
Can the plane fixed-point problem be solved by digging a simple dense
canal in a disk? 
Can $f^2$ be fixed-point free?

In 1971 Hagopian \cite{MR42:8469} proved $\mathcal{C}$ is not arcwise 
connected. 
Hagopian \cite{MR97a:54047} in 1996 improved this theorem by showing that 
an arcwise connected plane continuum has the fixed-point property if and 
only if its fundamental group is trivial.

It is not known if the fixed-point-free map $f$ can be a homeomorphism.
Bell \cite{MR58:18386} in 1978 proved $f$ cannot be a homeomorphism that 
is extendable to the plane.
Akis \cite{MR2001k:54068} and Bell \cite{bellprivate} proved $f$ is not a 
map that has an analytic extension to the plane.
In 1988 Hagopian \cite{MR89d:54022} proved $f$ cannot send each 
arc-component of $\mathcal{C}$ into itself. 
Hence $f$ is not a deformation.
Must $f$ permute every arc-component of $\mathcal{C}$?

In 1951 Hamilton \cite{MR12:627f} proved $\mathcal{C}$ is not chainable. 
We do not know if $\mathcal{C}$ can be triod-like 
\cite{MR85b:54052,MR87m:54100}. 
More generally, can $\mathcal{C}$ be tree-like \cite[p.~653]{MR13:265a}? 
Bellamy \cite{MR81h:54039} in 1979 defined a nonplanar tree-like continuum 
that admits a fixed-point-free map (also see \cite{MR82j:54075,MR83b:54044}
and \cite{MR94a:54108,MR96h:54029,MR99e:54024,MR2000k:54029.1}).
Using this example and an inverse limit technique of 
Fugate and Mohler \cite{MR80k:54062},
Bellamy \cite[p.~12]{MR81h:54039} defined a second tree-like continuum $M$ 
that admits a fixed-point-free homeomorphism. 
It is not known if $M$ can be embedded in the plane. 
Note that such an embedding would solve the plane fixed-point problem. 
Every proper subcontinuum of Bellamy's continuum $M$ is an arc. 
This motivates another open question. 
Must a nonseparating plane continuum with only arcs for proper subcontinua 
have the fixed-point property?

In 1990 Minc \cite{MR90d:54067} proved $\mathcal{C}$ is not weakly 
chainable (a continuous image of a chainable continuum). 
Minc \cite{MR99e:54024} in 1999 defined a weakly chainable tree-like 
continuum that does not have the fixed-point property.

Kuratowski \cite{kuratowski1927} defined a continuum $K$ to be of 
\emph{type $\lambda$} if $K$ is irreducible and every indecomposable 
continuum in $K$ is a continuum of condensation. 
Every continuum $K$ of type $\lambda$ admits a unique monotone upper 
semi-continuous decomposition to an arc with the property that each 
element of the decomposition has void interior relative to $K$ 
\cite[Th.~3, p.~216]{MR41:4467}. 
The elements of this decomposition are called \emph{tranches}.

Can $\mathcal{C}$ be a continuum of type $\lambda$ with the property that 
each of its tranches has the fixed-point property?
In answer to a question of Gordh \cite[Prob.~43, p.~371]{MR86a:54038.1}, 
Hagopian \cite{MR1990801} defined a nonplanar continuum 
$\mathcal{M}$ of type $\lambda$ such that each tranche of $\mathcal{M}$
has the fixed-point property and $\mathcal{M}$ does not.
Recently Hagopian and Ma\'nka \cite{hm} defined a planar continuum with 
these properties.

A fundamental exposition on the plane fixed-point problem is given in 
\cite[pp.~66 and 145]{MR92k:00014} (also see 
\cite{bingscottish}, \cite{MR80k:54064.1}, and \cite{MR92i:54033}).

\providecommand{\bysame}{\leavevmode\hbox to3em{\hrulefill}\thinspace}

\label{tphagopianend}

\chapter*{Janusz J.~Charatonik: On an old problem of Knaster}
\label{tpcharatonik1}
\begin{myfoot}
\begin{myfooter}
Janusz J.~Charatonik, 
\emph{On an old problem of Knaster},\\
Problems from Topology Proceedings, Topology Atlas, 2003, 
pp.\ 195--196. 
\end{myfooter}
\end{myfoot} 

When the definition of dendroids began to be formulated, in 1958/1959 
and in the early 1960s at the Wroc{\l}aw Higher Topology Seminar of the 
Polish Academy of Sciences (conducted by Bronis{\l}aw Knaster), 
Knaster saw this class of arcwise connected curves as ones that can be 
retracted onto their subdendrites or even onto their subtrees under small 
retractions, i.e., retractions that move points a little. 
Later, the contemporary definition of a \emph{dendroid} as an arcwise 
connected and hereditarily unicoherent continuum was formulated and 
commonly accepted because it is much more convenient to work with. 
But the problem if the two concepts coincide is still open. 

\begin{question}
Let $X$ be a dendroid. Do there exist, for each $\varepsilon > 0$, a tree 
(a dendrite) $T \subset X$ and a retraction $r\colon X \to T$ with 
$d(x, r(x)) < \varepsilon$ for each point $x \in X$?
\end{question}

Some partial positive answers can be found in 
\cite[Theorem~2, p.~261]{MR45:5965} for smooth dendroids and in 
\cite[Theorem~1, p.~120]{MR45:5963} for fans.
See also \cite{MR45:5964}.

Recall that if the assumption on the mapping of being a retraction onto a 
tree $T$ contained in $X$ is omitted, then the answer to the question is 
affirmative, since each dendroid, being a tree-like continuum, admits for 
each $\varepsilon > 0$ an $\varepsilon$-mapping onto a tree, see 
\cite{MR42:1072}.

The property of having ``small'' retractions onto trees is related to the 
following concept of an approximative absolute retract.
A compact metric space $X$ is called an \emph{approximative absolute 
retract} (abbr.\ AAR) if, whenever $X$ is embedded into another metric 
space $Y$, then for every $\varepsilon >0$ there exists a mapping 
$f_{\varepsilon}\colon Y \to X$ such that 
$d(x,f_{\varepsilon}(x)) < \varepsilon$ for each $x \in X$. 
Since trees are absolute retracts, it is clear that any compact space 
that admits ``small'' retractions onto trees must be an AAR.

The two following questions are closely related to Knaster's question 
discussed here. They are formulated at the end of \cite{cpaanr}.

\begin{question}
Is every dendroid an AAR?
\end{question}

\begin{question}
Is each dendroid the inverse limit of an inverse sequence of (nested) 
trees with retractions as bonding mappings? 
\end{question}

More information on dendroids and some open problems related to them is 
in~\cite{MR96g:54044}.

\providecommand{\bysame}{\leavevmode\hbox to3em{\hrulefill}\thinspace}

\label{tpcharatonik1end}

\chapter*{Janusz J.~Charatonik: Means on arc-like continua}
\label{tpcharatonik2}
\begin{myfoot}
\begin{myfooter}
Janusz J.~Charatonik, 
\emph{Means on arc-like continua},\\
Problems from Topology Proceedings, Topology Atlas, 2003, 
pp.\ 197--200.
\end{myfooter}
\end{myfoot} 

A \emph{mean} on a topological space $X$ is defined as a mapping 
$\mu\colon X \times X \to X$ such that $\mu (x,y) = \mu (y,x)$ and 
$\mu (x,x) = x$ for every $x, y \in X$ (in other words, it is a symmetric, 
idempotent, continuous binary operation on $X$). 
In \cite[p.~285]{MR58:18330.1} an approach to this concept is presented 
from the standpoint of the theory of hyperspaces (a mean on a continuum 
$X$ can be defined as a retraction of the hyperspace $F_2(X)$ onto 
$F_1(X)$, see also \cite[Section~76, p.~371]{MR99m:54006}; compare also 
\cite[Section~5, p.~18]{MR98j:54061} and
\cite[Section~6, p.~496]{MR2000m:54035}).

A natural problem that is related to this concept is what spaces, in 
particular what metric continua, admit a mean? 
No characterization is known yet. 

It is easy to give an example of a mean on the closed unit interval 
$[0,1]$ (e.g., the arithmetic mean $\mu (x,y) = \tfrac{x+y}{2}$). 
Means on $[0,1]$, even in a more general setting, were studied by 
A.N.~Kolmogoroff who described a structural form of these mappings in 
\cite{Kolmogoroff}. 
Functional equations of the type 
\begin{equation}
f(\mu(x,y)) = \mu(f(x), f(y))\tag{$\star$}
\end{equation}
with a given mean $\mu$ on $[0,1]$ and unknown mapping 
$f\colon [0,1] \to [0,1]$ have been studied extensively, see 
\cite{MR34:8020}. 
Inversely, a question about the existence of a mean on $[0,1]$ for a 
given mapping $f$ such that ($\star$) holds for all $x, y \in [0,1]$ is 
also discussed in some papers. 
E.g., in \cite{MR85c:39005} it is shown that equation ($\star$) has no 
solutions $\mu$ for the tent map $f$ (see \cite{MR2000i:39020} for an 
extension) and it is asked if a surjection $f$ on $[0,1]$ satisfying 
($\star$) for some mean $\mu$ must necessarily be monotone. 

A study on basic properties of means defined on arbitrary spaces started
with the habilitation thesis of G.~Aumann \cite{0008.05601, 0012.25205}, 
and it was developed in \cite{MR6:277g}, where it is shown that the 
circle, or even $k$-dimensional sphere for each $k \ge 1$ does not admit 
any mean, while each \emph{dendrite} (i.e., a locally
connected metric continuum containing no simple closed curve) does. 
An outline of a quite different proof that the circle does not admit any 
mean is given in \cite[(0.71.1), p.~50]{MR58:18330.1}.
These fundamental results have been generalized later in several ways.

Given a mapping $f\colon X \to Y$, a mapping $h\colon Y \to X$ is called 
a \emph{right inverse of} $f$ provided that 
$f \circ h = \operatorname{id} \restriction Y$. 
If, for a given $f$, there exists a right inverse of $f$, then $f$ is 
called an \emph{r-mapping}. 
Each r-mapping is surjective. 
Let $f\colon X \to Y \subset X$ be a retraction (i.e., such that 
$f \restriction Y = \operatorname{id} \restriction Y$; then $Y$ is called 
a \emph{retract} of $X$). 
Then $h = f \restriction Y$ is a right inverse of $f$, so each retraction 
is an r-mapping. 
It is known that if a space $X$ admits a mean and $f\colon X \to Y$ is an 
r-mapping, then $Y$ also admits a mean, \cite{MR51:14000}. 
In particular, each retract of $X$ admits a mean, \cite{MR37:3529}. 

A continuum $X$ is said to be \emph{unicoherent} provided that for each 
decomposition of $X$ into two subcontinua, their intersection is
connected. 
It is known that if a locally connected metric continuum admits a mean, 
then it is unicoherent; if, in addition, it is $1$-dimensional, then it 
is a dendrite, see \cite{MR37:3529} (compare also 
\cite[Theorem~5.31, p.~22]{MR98j:54061}). 
Local connectedness is essential in this result because the dyadic
solenoid is $1$-dimensional, unicoherent, and admits a mean, see 
\cite[76.6, p.~374]{MR99m:54006} (also \cite[5.47, p.~24]{MR98j:54061}; 
it admits an open and monotone mean, 
\cite[Example~5]{illanessimon}). 
For further progress see 
\cite{MR39:6274, MR41:9230, MR51:14000, MR57:10686, MR84c:54013}. 

In an early period of studies on means, the majority of results was 
related to locally connected spaces. 
One of the first examples of non-locally connected continua that admit no 
mean was the $\sin (1/x)$-curve, \cite{MR41:6168} (for an extension of 
this result see \cite{MR95m:54027}).
This curve is acyclic (in the sense that all its homology groups are
trivial). 
All known examples of locally connected continua that do not admit any 
mean are cyclic. 
So, a question arises if cyclicity is the only obstruction which does not 
let a locally connected continuum to admit a mean, \cite{MR39:6274}.

A (metric) continuum $X$ is said to be \emph{arc-like} provided that for 
each $\varepsilon > 0$ it has an $\varepsilon$-chain cover; or, 
equivalently, if it is the inverse limit of an inverse sequence of arcs 
with surjective bonding mappings. 

Let an inverse sequence $\{X_n, f_n : n \in \mathbb{N}\}$ be given 
each coordinate space $X_n$ of which admits a mean 
$\mu_n\colon X_n \times X_n \to X_n$ such that for each 
$n \in \mathbb{N}$ the functional equation 
$f_n(\mu_{n+1}(x,y)) = \mu_n(f_n(x), f_n(y))$ is satisfied for all 
$x, y \in X_{n+1}$. 
Then the inverse limit space 
$X = \varprojlim \{X_n, f_n: n \in \mathbb{N}\}$ 
admits a mean $\mu\colon X \times X \to X$ defined by 
$\mu(\{x_n\}, \{y_n\}) = \{\mu_n(x_n,y_n)\}$. 
Some special results concerning this concept are in \cite{MR85c:39005} and 
\cite{MR94b:39021}. 
As an answer in the negative to a question whether every mean on an 
arc-like continuum is an inverse limit mean, \cite{MR85c:39005}, a 
suitable example showing that inverse limit means are not preserved under 
homeomorphisms has been constructed in \cite{MR90j:39008}. 

In connection with the main result of \cite{MR41:6168} that the 
$\sin (1/x)$-curve does not admit any mean, P.~Bacon asked the following. 

\begin{question}[{\cite[p.~13]{MR41:6168}}]
Is the arc the only arc-like continuum that admits a mean?
Is the arc the only continuum containing an open dense half-line that 
admits a mean? 
\end{question}

After more than thirty years, the questions still remain unanswered. 
However, a sequence of important partial answers has been obtained. 

The above mentioned result of Bacon (that the $\sin (1/x)$-curve does not 
admit any mean) has been essentially extended in \cite{MR98g:54079}, where 
some criteria are obtained for the existence as well as for the 
non-existence of means on continua (the non-existence criterium is also 
presented in \cite[Section~76, p.~374--376]{MR99m:54006}). 
A further generalization was obtained in \cite{MR97i:54045}. 
It runs as follows. 

Two points $a$ and $b$ of an arc-like continuum are called 
\emph{opposite end points} of the continuum provided that for each 
$\varepsilon > 0$ there is an $\varepsilon$-chain cover of the continuum 
such that only the first link of the chain contains $a$ and only the last 
link of the chain contains $b$. 
Let a continuum $X$ contain an arc-like continuum $A$ with opposite end 
points $a$ and $b$ of $A$. 
A sequence $\{A_n: n \in \mathbb{N}\}$ of subcontinua of $X$ is called a 
\emph{folding sequence with respect to the point $a$} provided that for 
each $n \in \mathbb{N}$ there are two subcontinua $P_n$ and $Q_n$ of $A_n$ 
such that 
$A_n = P_n \cup Q_n$, 
$\operatorname{Lim} {(P_n \cap Q_n)} = \{a\}$, 
and
$\operatorname{Lim} P_n = \operatorname{Lim} Q_n = A$. 

\begin{theorem}[{\cite[p.~99]{MR97i:54045}}] 
Let a hereditarily unicoherent continuum $X$ contain an arc-like 
subcontinuum $A$ with opposite end points $a$ and $b$ of $A$. 
If there exist folding sequences $\{A_n\}$ and $\{B_n\}$ with respect to 
$a$ and $b$ correspondingly, then $X$ admits no mean. 
\end{theorem}

The concept of a folding sequence is a generalization of the concept of 
type $N$ \cite[p.~837]{MR80e:54045} which in turn generalizes the concept 
of a zigzag \cite[p.~78]{MR82i:54065} and is related to the notion of a 
bend set \cite[p.~548]{MR80a:54025}. 
These concepts were exploited to obtain some criteria for 
noncontractibility and nonselectibility of \emph{dendroids} (i.e., 
hereditarily unicoherent and arcwise connected continua) as well as for 
non-existence of means on these curves. 
For details see \cite[p.~23--32]{MR98j:54061} and 
\cite[p.~496-498]{MR2000m:54035}. 

The above theorem does not apply to hereditarily indecomposable continua, 
because it assumes the existence of decomposable subcontinua. 
The non-existence of means on the pseudo-arc (and on each hereditarily 
indecomposable circle-like continuum) follows from the following result 
that is shown also in \cite{MR97i:54045}. 

\begin{theorem}[{\cite[p.~102]{MR97i:54045}}]
If a hereditarily indecomposable contains a pseudo-arc, then it admits no 
mean. 
\end{theorem}

Another famous arc-like continuum is the simplest indecomposable 
continuum $D$ \cite[Fig.~4, p.~205]{MR41:4467.1} also called the 
buckethandle continuum or the Brouwer-Janiszewski-Knaster continuum. 
It can be defined as the inverse limit of arcs with tent bonding 
mappings. 
$D$ has exactly one end point, each of its proper subcontinua in an arc, 
and it again is an example to which the above theorem (on folding 
sequences of arcs) does not apply. 
Answering my question \cite{MR94d:54043}, A.~Illanes has shown that $D$ 
does not admit any mean \cite{illanes}. 
Similarly constructed indecomposable continua with $k$ end points (where 
$k \ge 2$; for $k = 3$ see \cite[p.~142]{MR23:A2857} and 
\cite[1.10, p.~7]{MR93m:54002}) also do not admit any mean, 
\cite[Corollary~3.15]{charatonik}. 
Recently, D.P.~Bellamy \cite{bellamy} presented an outline of a proof 
that each Knaster-type continuum (i.e., the inverse limit of arcs with 
open bonding mappings) different from an arc admits no mean. 

\providecommand{\bysame}{\leavevmode\hbox to3em{\hrulefill}\thinspace}

\label{tpcharatonik2end}

\chapter*{James T.~Rogers, Jr.: Classification of homogeneous continua}
\label{tprogers}
\begin{myfoot}
\begin{myfooter}
James T.~Rogers, Jr., 
\emph{Classification of homogeneous continua},\\ 
Problems from Topology Proceedings, Topology Atlas, 2003, 
pp.\ 201--216.
\end{myfooter}
\end{myfoot}

\newtheorem*{adttheorem}{Aposyndetic Decomposition Theorem}
\newtheorem*{tdttheorem}{Terminal Decomposition Theorem}

\mypreface
In volume 8 (1983) of \emph{Topology Proceedings}, J.T.~Rogers 
\cite{MR85c:54055.1} proposed a complete classification of homogeneous 
curves and a strategy to prove that all homogeneous continua of dimension 
$n > 1$ are aposyndetic.
That survey was updated six years later for the Proceedings of the 
Symposium on General Topology and Applications (Oxford, 1989) in 
\cite{MR94c:54061.1}.
This version contains a summary of both surveys and some new information 
provided by Rogers. 
This version was edited by Elliott Pearl with the approval of 
J.T.~Rogers; Rogers is the first person narrator here.

\section*{Introduction}

\emph{Fundamenta Mathematica} was the first journal devoted to set theory.
Lebesgue, among others, applauded the effort but worried that a dearth of 
publishable work might doom the enterprise \cite{bk}.
Perhaps it was to avoid this calamity and to prime the pump that the 
editors included a list of questions at the end of each volume.

The first question in the first volume in 1920 was answered almost 
immediately, but the second was a dilly.
Knaster and Kuratowski \cite{knastekuratowksi2} asked if each homogeneous, 
plane continuum must be a simple closed curve.
Mazurkiewicz \cite{mazurkiewicz19242} proved in 1924 that the answer is 
yes provided the continuum is locally connected.
This was the only significant progress on the problem for over a quarter 
century, even though the problem did not sit on the back burner.

In 1948, R.H.~Bing \cite{MR10:261a} proved that the pseudo-arc is 
homogeneous.
This remarkable result initiated a spate of activity on the problem.
In some sense, this period of intense activity was concluded in 1961 by 
another paper \cite{MR22:1869.1} of Bing, in which he showed that the answer 
to the question is yes provided the continuum contains an arc.
This could be called the classical or planar period in the study of 
homogeneous continua.
Although homogeneous continua in general were also investigated, the 
predominant results continued to be spawned by the original question of 
Knaster and Kuratowski.

Later in the decade two additional and important results were obtained, 
results that concluded the classical period.
In 1968, L.~Fearnley \cite{MR39:3460} and Rogers \cite{MR41:1018} 
independently proved that the pseudo-circle is not homogeneous, and in 
1969, at the Auburn Topology Conference, F.B.~Jones \cite{MR52:11862} 
announced that indecomposable, homogeneous, plane continua must be 
hereditarily indecomposable.
The pseudo-circle, defined by Bing \cite{MR13:265b} almost 20 years 
earlier, had emerged as the leading candidate for a new homogeneous 
continuum.
The fact that it is not homogeneous suggested that new homogeneous 
continua in the plane would be hard to come by.

The proofs of these two results told the tale on the state of the art at 
that time.
Jones never wrote up his proof---he told me once that it would have 
been so complicated that he feared no one would read it.
In the same vein, I felt that the ideas in the proof of the nonhomogeneity 
of the pseudo-circle should extend to some other separating plane 
continua, but the details were formidable, and I was never tempted more 
than briefly to attack them.
The reader should recall that, in those days, to prove the pseudo-circle 
nonhomogeneous, certain points $x$ and $y$ were precisely described, and 
it was shown that no homeomorphism of the continuum could move the point 
$x$ to the point $y$.

Clearly, new techniques were needed if the study of the homogeneous 
continua were to remain a viable field.
The most important such technique was already available, although we 
didn't know it.
In 1965, E.G.~Effros \cite{MR30:5175} proved an important result about 
Polish transformation groups.
When applied to the homeomorphism group of a homogeneous continuum, it 
yields a powerful and effective tool.

G.~Ungar \cite{MR52:6684} was the first to apply the Effros result to 
continua; with it, he showed that $2$-homogeneity implies local 
connectivity. 
It is significant that this is a nonplanar result (C.E.~Burgess, one of 
the pioneers in the study of homogeneous continua, had already shown the 
result in the plane and had raised the question in general 
\cite{MR20:1961}.)

In 1975, then, the study of homogeneous continua entered its current 
state---the modern or nonplanar period---a second period of intense 
activity, marked by extensive use of the Effros result and punctuated by 
the introduction of other new techniques as well.

\subsection*{Definitions and goals of this paper}

The goals of this paper are to summarize the present state of knowledge 
of homogeneous continua, to present a possible classification of all 
homogeneous continua, to ask some questions whose answers are important 
in obtaining further progress, and to mention some of the new techniques 
currently being used in the investigation of these continua.

The classification scheme rests on the cornerstone of Jones' Aposyndetic 
Decomposition Theorem.
We present the scheme first for plane continua, then for curves, and 
finally, for continua of dimension greater than one.

A \emph{continuum} is a compact, connected nonvoid metric space.
A \emph{curve} is a one-dimensional continuum.

A space $X$ is \emph{homogeneous} if, for each pair of points $x$ and $y$ 
of $X$, there exists a homeomorphism $f$ of $X$ such that $f(x) = y$.

A continuum $X$ is \emph{decomposable} if it is the union of two of its 
proper subcontinua; otherwise $X$ is \emph{indecomposable}.
A continuum is \emph{hereditarily indecomposable} if it does not contain a 
decomposable continuum.

A \emph{pseudo-arc} is a chainable, hereditarily indecomposable continuum.

\section*{Aposyndetic decompositions}

The notion of an aposyndetic continuum is crucial to the investigation of 
homogeneous continua.
Aposyndesis is a weak form of local connectivity and the following 
implications hold and none are reversible:
locally connected $\Rightarrow$ aposyndetic $\Rightarrow$ decomposable.

A continuum $X$ is \emph{aposyndetic at $x$ with respect to $y$} if $X$ 
contains an open set $G$ and a subcontinuum $H$ such that 
$x \in G \subset H \subset X \setminus \{y\}$.
A continuum is said to be \emph{aposyndetic} if it is aposyndetic at each 
point with respect to any other point.

For each $x$ in $X$, let 
$$
L_x 
= \{x\} \cup \{z : \text{X is not aposyndetic at $z$ with respect to 
$x$}\}.
$$
$L_x$ is always a subcontinuum of $X$.
If $X$ is indecomposable, then $L_x = X$ for all $x$.
If $X$ is decomposable, then $L_x$ is a proper subcontinuum of $X$ for 
some $x$.
Jones has used the $L_x$'s to fashion an important decomposition theorem 
for homogeneous decomposable continua.

\begin{adttheorem}
Let $X$ be a homogeneous continuum such that $X$ is decomposable but not 
aposyndetic.
If $G = \{L_x : x \in X\}$, then
\begin{myenumerate}
\item
the collection $G$ is a monotone, continuous decomposition of $X$,
\item
the elements of the decomposition are mutually homeomorphic homogeneous 
continua,
\item
the quotient space is a homogeneous continuum, and
\item
the quotient space is an aposyndetic continuum.
\end{myenumerate}
\end{adttheorem}

Rogers has added the following improvements to this theorem.
\begin{myenumerate}
\setcounter{myenum}{4}
\item
\emph{The elements of the decomposition are cell-like, indecomposable
continua of the same dimension as $X$.}
\item
\emph{The quotient space is a curve.}
\end{myenumerate}

In case $X$ is planar, the quotient space is homeomorphic to the circle 
$\mathbb{S}^1$.

\section*{Jones' classification of homogeneous plane continua}

In 1949, Jones \cite{MR10:468d} proved that an aposyndetic, homogeneous 
plane continuum is either a point or a simple closed curve.

In 1951, Jones \cite{MR13:573a} made the first use of decompositions of 
homogeneous continua by showing that a nonseparating homogeneous plane 
continuum must 
be indecomposable.
In 1954, he divided homogeneous plane continua into three types:
\newcounter{adt}
\begin{list}{(Type \Alph{adt})}{\usecounter{adt}}
\item
nonseparating (hence indecomposable);
\item
separating and decomposable;
\item
separating and indecomposable.
\end{list}
Furthermore, he showed \cite{MR17:180e.1} that each Type B continuum is a 
circle of Type A continua.
Rogers \cite{MR83b:54045} proved that the set of Type C continua is empty.
C.L.~Hagopian \cite{MR54:8586,MR85m:54031} proved that Type A continua 
are hereditarily indecomposable.

There are, at present, four known homogeneous plane continua:
the point, the pseudo-arc, the circle and the circle of pseudo-arcs.

An affirmative answer to the following question of Jones would imply that 
these four are the only homogeneous plane continua.

\begin{question}
Is each nondegenerate homogeneous nonseparating plane continuum a 
pseudo-arc?
\end{question}

Oversteegen and Tymchatyn \cite{MR84h:54030.1,MR85h:54068} showed that 
each Type A continuum has span zero and is a continuous image of the 
pseudo-arc.

The classification of homogeneous plane continua is summarized in 
Figure 1.

\newlength{\mycon}
\settowidth{\mycon}{connected}
\newlength{\mydne}
\settowidth{\mydne}{\emph{\small Do not exist.}}
\newlength{\mynotapo}
\settowidth{\mynotapo}{not aposyndetic}
\newlength{\mynotloc}
\settowidth{\mynotloc}{not locally}
\newlength{\mynonsep}
\settowidth{\mynonsep}{nonseparating}
\addtolength{\mynonsep}{29pt}

\begin{figure}
\begin{tabular}{|m{\mycon}|m{\mydne}|m{\mynotapo}|m{\mydne}|m{\mynonsep}|}
\hline
\multicolumn{3}{|c|}{decomposable}&
\multicolumn{2}{|c|}{indecomposable}\\
\hline
\multicolumn{2}{|c|}{aposyndetic}&
not aposyndetic&
separating&
nonseparating\\ 
\hline
\parbox{\mycon}{%
\setlength{\baselineskip}{0.75\baselineskip}%
\raggedright%
\vspace{2pt}%
locally connected%
\vspace{2pt}%
}&
\parbox{\mynotloc}{%
\setlength{\baselineskip}{0.75\baselineskip}%
\raggedright%
\vspace{2pt}%
not locally connected%
\vspace{2pt}%
}&
&
&
\\
\cline{1-2}
\parbox{\mycon}{%
\setlength{\baselineskip}{0.75\baselineskip}%
\raggedright%
\emph{\small Must be a point or a circle.}%
}&
\vspace{-1.333\baselineskip}
\emph{\small Do not exist.}&
\parbox{\mynotapo}{%
\setlength{\baselineskip}{0.75\baselineskip}%
\raggedright%
\vspace{-45pt}
\emph{\small Must be a circle of nonseparating continua.}%
}&
\vspace{-62pt}
\emph{\small Do not exist.}&
\parbox[t]{\mynonsep}{%
\setlength{\baselineskip}{0.75\baselineskip}%
\raggedright%
\vspace{-23pt}%
\emph{\small
Are~tree-like,~hereditarily~indecomposable,
span~zero,~weakly
chainable.
Pseudo-arc is the only known example.}%
\vspace{2pt}
}\\
\hline
\end{tabular}
\caption{A classification of homogeneous plane continua}
\end{figure}

\section*{Homogeneous curves outside the plane}

Homogeneous, nonplanar curves are a more yeasty mixture.
There exists, for instance, a collection of cardinality $\mathfrak{c}$ of 
solenoids.
A \emph{solenoid} is defined as an inverse limit of circles with covering 
maps as the bonding maps.
Each solenoid is an indecomposable continuum as well as an abelian 
topological group.
Hence each solenoid is an indecomposable homogeneous continuum with 
nontrivial cohomology.

If $f\colon \mathbb{S} \to \mathbb{S}^1$ is the projection of the 
solenoid $\mathbb{S}$ onto the factor space $\mathbb{S}^1$, then $f$ is a 
morphism of topological groups with kernel a topological group $G$ whose 
underlying space is a Cantor set.
The collection $\mathcal{S} = (\mathbb{S},f,\mathbb{S}^1,G)$ is a 
principal fiber bundle.

In 1958, R.D.~Anderson \cite{MR20:2675} showed that the Menger universal 
curve (the so-called ``Swiss Cheese Space'') is homogeneous, and that the 
circle and the Menger curve are the only homogeneous locally connected 
curves.

In 1961, J.H.~Case \cite{MR24:A1714} constructed a new homogeneous curve 
as an inverse limit of Menger universal curves and double-covering maps.
Case's construction was quite complicated, and in 1982 Rogers 
\cite{MR85f:54073} provided a simpler, geometric construction of similar 
continua and then \cite{MR86c:54030} a bundle-theoretic construction of 
such spaces.
These continua are simply the total spaces of bundles induced from 
solenoid bundles by a retraction of the Menger universal curve onto a 
\emph{core} circle.

It can be proved from these constructions that there are $\mathfrak{c}$ 
such continua (one for each solenoid), that each is aposyndetic but not 
locally connected, and that none is arcwise-connected, hereditarily 
decomposable, or pointed-one-movable.

In 1983, P.~Minc and Rogers \cite{MR88i:54027} constructed even more 
homogeneous, aposyndetic curves.
The geometric idea is to spin the Menger curve around several of its holes 
at the same time.
Each finite sequence of solenoids $\mathbb{S}_1, \dots, \mathbb{S}_n$ 
determines one of these continua $M$.
If $M^\prime$ is another such continuum determined by the sequence 
$\mathbb{S}_1^\prime, \dots, \mathbb{S}_m^\prime$ and if $M^\prime$ is 
homeomorphic to $M$, then $n = m$ and $\mathbb{S}_i$ is homeomorphic to 
$\mathbb{S}_i^\prime$ for some reindexing.

In 2002, J.R.~Prajs \cite{MR2003f:54077.3} constructed a homogeneous, 
arcwise connected curve that is not locally connected. 
This important example answers an old question of K.~Kuperberg. 
We will see in the next section that this example answers two more 
questions of the original survey.

Prajs' example is constructed as an inverse limit of Menger curves and
covering maps; the maps, however, are chosen differently than those in
the examples of Rogers and of Minc and Rogers. 
In particular, he spins the Menger curve around infinitely many of its 
holes.

K.~Villarreal \cite{MR94h:54043} has shown that spinning the Menger 
curve $M$ around infinitely many of its holes in the style of Minc and 
Rogers leads to a continuum that is not homogeneous.

\section*{Rogers' classification scheme for homogeneous curves}

We propose here a classification of homogeneous curves by dividing them 
into six types. 
This classification is summarized in Figures 2 and 3.

\subsection*{Type 1. Locally connected}
The Menger universal curve and the circle are the only ones 
\cite{MR20:2675}, so this type is completely understood.

\subsection*{Type 2. Aposyndetic but not locally connected}
The examples of Case, Rogers, Minc \& Rogers, and Prajs are the only 
examples known.

\subsection*{Type 3. Decomposable but not aposyndetic}
Jones' Aposyndetic Decomposition Theorem says that each Type 3 curve 
admits a decomposition into Type 6 curves such that the quotient space is 
a Type 1 or a Type 2 curve.

\subsection*{Type 4. Indecomposable and cyclic}
A curve is \emph{cyclic} if its first \v{C}ech cohomology group with 
integral coefficients does not vanish; otherwise it is \emph{acyclic}.
A curve is cyclic if and only if it admits an essential map onto 
$\mathbb{S}^1$.
The solenoids and the solenoids of pseudo-arcs are the only known 
continua of Type 4.

\subsection*{Type 5. Acyclic but not tree-like}
A curve is \emph{tree-like} if it admits finite open covers of 
arbitrarily small mesh whose nerves are trees.
A curve is tree-like if and only if it has trivial shape.
A tree-like curve is acyclic.
Bing showed that acyclic planar curves are tree-like.
In 1987, Rogers \cite{MR88k:57016.1} proved that acyclic homogeneous 
curves are tree-like, that is, there are no Type 5 curves.

\subsection*{Type 6. Tree-like}
The pseudo-arc is the only known Type 6 curve.

\newlength{\mywid}
\setlength{\mywid}{109pt}
\begin{figure}
\begin{tabular}{|m{\mywid}|m{\mywid}|m{\mywid}|}
\hline
\vspace{2pt}%
Type~1:~locally~connected\newline&
\vspace{2pt}%
Type~2:~aposyndetic,~not locally~connected&
\vspace{2pt}%
Type 3: not aposyndetic\newline\\
\hline
\vspace{2pt}%
\parbox{\mywid}{%
\setlength{\baselineskip}{0.75\baselineskip}%
\raggedright%
\small Must be a Menger curve or a circle (Anderson, 1958).
\newline\newline%
}%
\vspace{1pt}%
&
\vspace{2pt}%
\parbox{\mywid}{%
\setlength{\baselineskip}{0.75\baselineskip}%
\raggedright%
\small Only known examples are of Case, Rogers, Minc \& Rogers, Prajs.
\newline%
}%
\vspace{1pt}%
&
\vspace{2pt}%
\parbox{\mywid}{%
\setlength{\baselineskip}{0.75\baselineskip}%
\raggedright%
\small Jones~Aposyndetic~Decomposition~applies:~decomposes~%
into~Type~6~with~quotient~space~Type~1~or~Type~2.%
}%
\vspace{1pt}%
\\
\hline
\end{tabular}
\caption{Types of decomposable homogeneous curves}
\end{figure}

\begin{figure}
\begin{tabular}{|m{\mywid}|m{\mywid}|m{\mywid}|}
\hline
\vspace{2pt}%
Type 4: cyclic\newline&
\vspace{2pt}%
Type~5:~acyclic,~not tree-like&
\vspace{2pt}%
Type 6: tree-like\newline\\
\hline
\vspace{2pt}%
\parbox{\mywid}{%
\setlength{\baselineskip}{0.75\baselineskip}%
\raggedright%
\small
E.g.,~Solenoid,~solenoid of pseudo-arcs.
Terminal~decomposition theorem~applies.
}%
\vspace{1pt}%
&
\vspace{2pt}%
\parbox{\mywid}{%
\setlength{\baselineskip}{0.75\baselineskip}%
\raggedright%
\small
Do~not~exist~(Rogers,~1987)
\newline\newline\newline%
}%
\vspace{1pt}%
&
\vspace{2pt}%
\parbox{\mywid}{%
\setlength{\baselineskip}{0.75\baselineskip}%
\raggedright%
\small
E.g.,~Pseudo-arc.
Must~be~hereditarily indecomposable.
\newline%
}%
\vspace{1pt}%
\\
\hline
\end{tabular}
\caption{Types of indecomposable homogeneous curves}
\end{figure}

\section*{Classifying homogeneous curves}

\subsection*{Classifying Type 2 curves}

There are no Type 2 curves in the plane, but there
are examples in $\mathbb{R}^3$.
All the known examples of Type 2 continua can be obtained as inverse 
limits of Menger universal curves and covering maps.

\begin{question}
Is each Type 2 curve an inverse limit of Menger universal curves and maps?
and fibrations?
and covering maps?
\end{question} 

\begin{question} 
Does each Type 2 curve contain an arc?
\end{question} 

In both surveys, we asked as Question 2 if each Type 2 curve is a bundle
over the universal curve with Cantor sets as the fibers. 
In both surveys, we asked as Question 5 if each Type 2 curve retracts onto 
a solenoid.

The example constructed by Prajs is aposyndetic, but it is not a bundle
over the universal curve. 
Since it is arcwise connected, it cannot be mapped onto a solenoid, let 
alone retracted onto one. 
Hence both Question 2 and Question 5 of the surveys have negative answers.

\subsection*{Classifying Type 3 curves; decompositions into pseudo-arcs}

Jones' theorem tells us, in a sense, not to worry about Type 3 curves 
until we know enough about Type 1, Type 2, and Type 6 curves.
There is only one known Type 6 curve, the pseudo-arc, so it is natural to 
ask if each Type 1 or Type 2 curve can be realized as a decomposition of 
a Type 3 curve into pseudo-arcs.
More generally, there is the problem, given a homogeneous curve $X$, of 
\emph{blowing up} its points into pseudo-arcs to obtain a homogeneous 
curve $\tilde{X}$.

Bing and Jones \cite{MR20:7251} solved this problem for the circle.
It follows from their construction that, to any finite, connected graph
$G$, there corresponds a curve $\tilde{G}$ and a decomposition of 
$\tilde{G}$ into pseudo-arcs with quotient space $G$.

Rogers \cite{MR57:4128} solved this problem for solenoids.
The idea in that paper is to express a curve $X$ as an inverse limit of 
graphs $(G,g)$, use Bing-Jones to blow up the graphs $G$ to graphs of 
pseudo-arcs $\tilde{G}$, and obtain $\tilde{X}$ as an inverse limit of 
$(\tilde{G}, \tilde{g})$ such that $\tilde{X}$ admits a continuous 
decomposition into pseudo-arcs with quotient space $X$.

The problem then is to show that $X$ homogeneous implies $\tilde{X}$ 
homogeneous.
In the case (such as for the solenoids) that $X$ is homogeneous by 
homeomorphisms induced by commuting diagrams of maps on the inverse 
sequence $(G,g)$, the desired homeomorphisms on $\tilde{X}$ can be 
obtained by a straightforward lifting process \cite{MR57:4128}.

But it is not known that there are always enough induced homeomorphisms on 
$X$ to do the job, and in fact, it seems unlikely that this is always so.
In the absence of induced homeomorphisms, one must fall back to 
Mioduszewski's $\epsilon$-commutative diagrams \cite{MR29:4035}, and then 
appears the sticky problem of whether the lift to 
$(\tilde{G}, \tilde{g})$ of an \emph{almost commutative} diagram 
involving $(G,g)$ is still \emph{almost commutative enough}.
Fortunately, by a careful use of the Bing-Jones paper, Wayne Lewis 
\cite{MR86e:54038} has proved that this is indeed possible, and that 
hence, for each homogeneous curve $X$, there is a homogeneous curve 
$\tilde{X}$ that admits a continuous decomposition into pseudo-arcs with 
quotient space $X$.

Incidentally, the problem of replacing a map between inverse limit spaces 
by a map induced from commuting diagrams on the inverse sequences is an 
important problem in continua theory.
One would wish the induced map to have any desirable property (such as 
being a homeomorphism taking the point $x$ to the point $y$) possessed by 
the original map.
More about this possibility and its limitations is needed.

\subsection*{Classifying Type 4 curves; decompositions for indecomposable 
continua}

\begin{question}
Does each Type 4 curve that is not a solenoid admit a continuous 
decomposition into Type 6 curves so that the resulting quotient space is 
a solenoid?
\end{question}

Hagopian \cite{MR85m:54031} has shown that the answer is yes for atriodic 
curves.

Rogers \cite{MR90a:54091} proved the Terminal Decomposition Theorem, which 
gives an analogue to the Aposyndetic Decomposition Theorem.
Intuitively, the decomposition is the one sought to answer the question 
above, but it has not been proved that the quotient space is a solenoid.
We consider this further below.

How can we get an \emph{aposyndetic decomposition} for indecomposable 
continua?
If $x$ is a point of the indecomposable continuum $X$, then $L_x = X$, and 
so the aposyndetic decomposition itself is the trivial one yielding a 
degenerate quotient space.
Something else must be tried.

A subcontinuum $Z$ of $X$ is \emph{terminal} in $X$ if each subcontinuum 
$Y$ of $X$ that meets $Z$ satisfies either $Y \subset Z$ or $Z \subset Y$.
If $X$ is the topologist's $\sin 1/x$ curve, then the limit bar is a 
terminal subcontinuum of $X$.
Each point of a continuum is a terminal subcontinuum.
All subcontinua of a hereditarily indecomposable continuum are terminal 
subcontinua.

Implicit in the proof of Jones' decomposition for a decomposable, 
homogeneous continuum $X$ is the fact that
$$
\{L_x : x \in X \} =
\{Z : \text{$Z$ is a maximal, terminal proper subcontinuum of $X$}\}.
$$
The idea for an indecomposable continuum $X$ is to decompose $X$ by 
maximal, terminal proper subcontinua.
Of course the following question arises immediately:
Must an indecomposable homogeneous curve contain a maximal, terminal 
proper subcontinuum?
The answer is no, since the pseudo-arc is homogeneous.

The answer is yes, however, for cyclic homogeneous curves, and the proof 
is quite interesting. 
Here is an outline.

If $X$ is a homogeneous, cyclic curve, then $X$ can be embedded in 
$\mathbb{S}^1 \times D$, where $D$ is a $3$-cell, so that the inclusion 
map is not homotopic to a constant map.
Let $p\colon R \times D \to \mathbb{S}^1 \times D$ be the usual covering 
space, and let $\tilde{X} = p^{-1}(X)$.

We show that each component $K$ of $\tilde{X}$ is homogeneous and 
unbounded in both the positive and negative directions.
Compactify $K$ with the two-point set $\{\pm \infty\}$.
We show that the continuum $\tilde{K} = K \cup \{\pm \infty\}$ admits an 
aposyndetic decomposition and that $L_\infty = \{\infty\}$ and 
$L_{-\infty} = \{-\infty\}$.
We push the decomposition elements $K$ downstairs and show that they fit 
together to yield the following decomposition theorem.

\begin{tdttheorem}[\cite{MR90a:54091}]
Let $X$ be a homogeneous curve such that $H^1(X) \not= 0$.
If $G$ is the collection of maximal terminal proper subcontinua of $X$, 
then
\begin{myenumerate}
\item
$G$ is a monotone, continuous, terminal decomposition of $X$,
\item
the nondegenerate decomposition elements of $G$ are mutually 
homeomorphic, hereditarily indecomposable, tree-like, terminal, 
homogeneous continua,
\item
the quotient space is a homogeneous curve, and
\item
the quotient space does not contain any proper, nondegenerate terminal 
subcontinuum.
\end{myenumerate}
\end{tdttheorem}

A decomposable, homogeneous continuum is aposyndetic if and only if it 
does not contain any proper nondegenerate terminal subcontinuum.
Thus the last conditions of both the Aposyndetic Decomposition Theorem 
and the Terminal Decomposition Theorem are saying the same thing when the 
homogeneous curve $X$ is both decomposable and cyclic.

\subsection*{Classifying Type 6 curves}

Jones \cite{MR13:573a} showed that Type 6 curves are indecomposable.
Lewis \cite{MR85h:54066} showed that a new example of a Type 6 curve must 
be infinitely-branched and infinitely-junctioned and must contain a 
proper nondegenerate subcontinuum different from a pseudo-arc.
Hagopian \cite{MR84d:54059.1} showed that no example can contain an arc.
Rogers \cite{MR84b:54072} showed that any hereditarily indecomposable 
homogeneous continuum must be a type 6 curve.
Oversteegen and Tymchatyn \cite{MR88m:54044} gave a new proof that the 
pseudo-arc is homogeneous.
Bing \cite{MR21:3818} showed that the pseudo-arc is the only chainable 
homogeneous continuum.

Krupski and Prajs \cite{MR90f:54054.2} answered an old question of Bing by 
showing that Type~6 curves are hereditarily indecomposable.

We summarize the questions that have been asked by different investigators 
in seeking further restrictions on Type 6 curves.

\begin{question}
Are Type 6 curves pseudo-arcs? weakly chainable? hereditarily equivalent? 
Do they have span zero? 
Do they have the fixed point property?
\end{question}

Oversteegen \cite{MR85m:54034,MR91g:54049} obtained some very good 
results concerning the problem of determining when a tree-like continuum 
is chainable.
He has shown, for example, that all continua with zero span that are the 
image of a chainable continuum under an induced map are chainable.

\section*{A characterization of homogeneous curves}

In broader strokes from the classification above, we choose the following 
three questions:

\begin{question}
\mbox{}
\begin{myenumerate}
\item
Is each aposyndetic, nonlocally connected, homogeneous curve an inverse 
limit of Menger universal curves and covering maps?
\item
Does each homogeneous, cyclic, indecomposable curve that is not a solenoid 
admit a continuous decomposition by tree-like, homogeneous curves so that 
resulting quotient space is a solenoid?
\item
Are tree-like, homogeneous curves pseudo-arcs?
\end{myenumerate}
\end{question}

Why do we choose these three questions?
If the answer to each of these three questions is yes, then we can 
classify homogeneous curves according to the following scheme:
Each homogeneous curve would be
\begin{myenumerate}
\item
a simple closed curve or a Menger universal curve, or
\item
an inverse limit of Type 1 curves and covering maps, or
\item 
a curve admitting a continuous decomposition into pseudo-arcs such that 
the quotient space is a curve of Type 1 or Type 2, or
\item
a pseudo-arc.
\end{myenumerate}

If we could answer these three questions affirmatively, then we would 
have completed the classification of homogeneous curves.
Of course, a negative answer to any one of these questions would mean 
that there are additional homogeneous curves and the classification must 
be refined.

\section*{Other classifications of homogeneous curves}

The pseudo-arc is the only homogeneous arc-like continuum.
The circle and the circle of pseudo-arcs are the only homogeneous, 
separating, planar, circle-like continua, and the solenoids and solenoids 
of pseudo-arcs are the only homogeneous nonplanar circle-like continua.
This is the beginning of a classification of homogeneous continua 
according to the graphs used in their inverse limit representations.

A curve is \emph{simply cyclic} if it is an inverse limit of graphs each 
of which contains only one cycle.
Continuing this sort of classification, Rogers \cite{MR90h:54038.1} has 
proved that each simply cyclic homogeneous curve that is not tree-like 
either is a solenoid or admits a decomposition into mutually homeomorphic, 
tree-like, homogeneous curves with quotient space a solenoid.
More along this line should be possible.

A curve is said to be \emph{finitely cyclic} if it is the inverse 
limit of graphs of genus $\leq k$, where $k$ is some integer.
Krupski and Rogers \cite{MR92c:54038.1} have proved that each finitely 
cyclic, homogeneous curve that is not tree-like is a solenoid or admits a 
decomposition into mutually homeomorphic, tree-like, homogeneous curves 
with quotient space a solenoid.
Since the Menger curve is homogeneous, the restriction to finitely cyclic 
curves is essential.

A curve is said to be \emph{$k$-junctioned} if it is the inverse limit of 
graphs each of which has at most $k$ branchpoints.
Duda, Krupski, and Rogers \cite{MR92m:54061.1} have proved that a 
homogeneous, $k$-junctioned curve must be a pseudo-arc, a solenoid or a 
solenoid of pseudo-arcs.

Finally, more about Type 2 and Type 4 curves should be forthcoming if one 
could detect the right embedding of such a curve into $F \times D$, where 
$D$ is a $3$-cell and $F$ is a closed, hyperbolic surface of sufficiently 
high genus.

\section*{Homogeneous continua of higher dimension}

Homogeneous continua of dimension greater than one can be divided similarly 
into six types, but in general they form a rather intractible class with 
questions arising from varied sources.

\subsection*{Type 1.~Locally connected}

Closed $n$-manifolds (for $n > 1$), countable products of locally 
connected, homogeneous, nondegenerate continua, and the Hilbert cube are 
Type 1 continua.
K.~Kuperberg has shown that certain mapping tori are Type 1 
continua.
Higher dimensional analogues of the Menger curve may be homogeneous.

\subsection*{Type 2.~Aposyndetic but not locally connected}

All nontrivial products of continua are aposyndetic.
Hence any nontrivial countable product of homogeneous, nondegenerate 
continua one of whose factors is not locally connected is a Type 2 
continuum.
If $M$ is a closed $n$-manifold (for $n > 1$) that admits a retraction 
onto a finite wedge of circles, then the bundle machines of 
\cite{MR86c:54030,MR88i:54027} 
automatically provide an $n$-dimensional Type 2 continuum.
Some have speculated that certain mapping tori are Type 2 continua.
Karen Villarreal \cite{MR93j:54024} has constructed additional 
two-dimensional homogeneous continua using fibered products.

\subsection*{Type 3.~Decomposable but not aposyndetic}

Again, the Jones' Aposyndetic Decomposition Theorem comes into play, this 
time in its full generality.

\begin{theorem}
Each decomposable, homogeneous continuum admits a continuous decomposition 
into mutually homeomorphic, cell-like, indecomposable, homogeneous 
continua such the quotient space is an aposyndetic, homogeneous continuum.
\end{theorem}

A continuum is \emph{cell-like} if it has the shape of a point.
No Type 3 continuum is known, which suggests the following question:

\begin{question}
Is each decomposable, homogeneous continuum of dimension greater than one 
aposyndetic?
\end{question}

An affirmative answer to the next question would imply an affirmative
answer to the previous question and strengthen the Decomposition Theorem
enormously.

\begin{question}
Must the elements of this aposyndetic decomposition be hereditarily 
indecomposable?
\end{question} 

The first survey asked if the aposyndetic decomposition could raise
dimension or could lower dimension. 
In other words, suppose $X$ is a decomposable, homogeneous continuum of 
dimension greater than one, and suppose $X$ is not aposyndetic. 
If $Y$ denotes the quotient space of the aposyndetic decomposition of $X$, 
can the dimension of $Y$ be greater than that of $X$? 
less than that of $X$?

Rogers \cite{MR1992870.2} answered this question recently by showing 
that the dimension of $Y$ must be one.

\subsection*{Type 4.~Indecomposable and cyclic}

Type 4 continua are the indecomposable and cyclic continua.
\emph{Cyclic} means that some (reduced) \v{C}ech cohomology group is 
nontrivial; otherwise the continuum is \emph{acyclic}.

\subsection*{Type 5.~Acyclic but not tree-like}

Type 5 continua are the acyclic but not cell-like continua.

\subsection*{Type 6.~Tree-like}

Type 6 continua are the cell-like continua.

No continuum of Type 4, 5, or 6 is known; these are really uncharted 
waters.

\begin{question}
Is each indecomposable, nondegenerate homogeneous continuum
one-dimensional?
\end{question}

Rogers \cite{MR85f:54073} showed that all hereditarily indecomposable, 
nondegenerate homogeneous continua are one-dimensional.
Hagopian \cite{MR85m:54031} showed that every indecomposable, 
nondegenerate homogeneous continuum of dimension greater than one must 
contain a triod.

Some formidable obstacles lie in the path of a complete classification of 
high-dimensional homogeneous continua.
For instance, Bing and Borsuk \cite{MR30:2475} conjectured in 1965 that 
an $n$-dimensional, homogeneous, compact ANR is an $n$-manifold, and they 
proved the conjecture true for $n=1$ or $2$.
In 1980, however, W.~Jakobsche \cite{MR81i:57007} showed that the 
validity of the Bing-Borsuk conjecture for $n=3$ implies the validity of 
the Poincar\'e Conjecture!

Consider also this baffling question from infinite-dimensional topology:

\begin{question}
Is each nondegenerate, homogeneous contractible continuum homeomorphic to 
the Hilbert cube?
\end{question}

Krupski \cite{MR90k:54049} has shown that homogeneous continua are Cantor 
manifolds.
Prajs \cite{MR89e:54069} showed that homogeneous continua in 
$\mathbb{R}^{n+1}$ that contain an $n$-cell are locally connected; this 
extends a planar theorem of Bing.

M.~Re\'nska \cite{MR2003g:54082.1} proved that there exist rigid 
hereditarily indecomposable continua in every dimension.  
In fact there exist continuum many such continua in each dimension.  

\section*{A decomposition theorem}

One of the most useful tools in studying homogeneous continua is 
decompositions.
Here we state a version of the decomposition theorems of 
\cite{MR83c:54045} for homogeneous curves.

\begin{theorem}
Let $X$ be a homogeneous curve, and let $H(X)$ be its homeomorphism group.
Let $\mathcal{G}$ be a partition of $X$ into proper, nondegenerate 
continua such that $H(X)$ respects $\mathcal{G}$ (this means that either 
$h(G_1) = G_2$ or $h(G_2) \cap G_1 = \emptyset$, for all $G_1$, $G_2$ in 
$\mathcal{G}$ and all $h$ in $H(X)$). Then
\begin{myenumerate}
\item
$\mathcal{G}$ is a continuous decomposition of $X$,
\item
there is a continuum $G$ such that each element of $\mathcal{G}$ is 
homeomorphic to $G$,
\item
$G$ is homogeneous, hereditarily indecomposable and tree-like, 
\item
the quotient image of this decomposition is a homogeneous curve.
\end{myenumerate}
\end{theorem}

Here are some applications of this decomposition theorem.

\subsection*{Application 1}

Suppose $X$ is a homogeneous curve that contains an arc.
Let $\mathcal{G}$ be the set whose elements are closures of arc 
components of $X$.
One shows that $\mathcal{G}$ is a partition of $X$ which $H(X)$, of 
course, respects.
Since no homogeneous, tree-like continuum can contain an arc 
\cite{MR84d:54059.1}, it follows that $\mathcal{G}$ contains only the one 
element $X$.
Therefore, if a homogeneous curve contains an arc, then it contains a 
dense arc component \cite{MR86m:54046}.

\subsection*{Application 2}

(W.~Lewis \cite{MR85h:54066})
Suppose $X$ is hereditarily indecomposable homogeneous curve (this implies 
$X$ is tree-lke), not a pseudo-arc, and contains a pseudo-arc.
For each point $x$ of $X$, let $P_x$ be the closure of the union of all 
pseudo-arcs containing $x$. 
The quotient space of the decomposition $\{ P_x : x \in X \}$ is a 
tree-like, homogeneous continuum containing no pseudo-arc.

\subsection*{Application 3}

Suppose $X$ is decomposable homogeneous curve.
Let $L_x$ be the set of all points $z$ of $X$ such that $X$ is not 
aposyndetic at $z$ with respect to $x$.
$\{ L_x : x \in X\}$ is the decomposition in Jones' Aposyndetic 
Decomposition Theorem \cite{MR17:180e.1}.

\subsection*{Application 4}

There does not exist, for instance, a circle of solenoids.
This means that no homogeneous curve admits a decomposition into solenoids 
such that the quotient space is a simple closed curve.

\section*{Completely regular maps}

A surjective map $f \colon X \to Y$ between metric spaces is said to be 
\emph{completely regular} if, for each $\epsilon > 0$ and point $y$ in 
$Y$, there exists a $\delta > 0$ such that $d(y,y') < \delta$ implies 
there exists a homeomorphism of $f^{-1}(y)$ onto $f^{-1}(y')$ moving no 
point as much as $\epsilon$.

Projection maps of products are completely regular, and completely 
regular maps are open.
In general, neither of the converse statements is true.

Dyer and Hamstrom introduced completely regular maps in \cite{MR19:1187e} 
with the idea of showing that spaces on which certain open maps are 
defined are locally products.
They considered, for instance, maps whose fibers are $2$-spheres.
Kim \cite{MR44:3346} has shown that their techniques, together with 
current knowledge about the homeomorphism group of a compact manifold, 
imply that each completely regular map with fibers homeomorphic to a 
compact manifold is locally trivial.

Completely regular maps arise naturally in the study of homogeneous 
continua, frequently as a consequence of using the Effros result.
Moreover, these maps are often not locally trivial.
Consider the two following theorems.

\begin{theorem}
In the decomposition theorem above, the quotient map is completely 
regular.
\end{theorem}

The second theorem is an immediate corollary of 
\cite[Theorem 9]{MR20:7251}

\begin{theorem}
If $f$ is an open, surjective map between compacta with the property that 
each point inverse is a pseudo-arc, then $f$ is completely regular.
\end{theorem}

Completely regular maps have some special properties.
We close with two of them.
The first, due to Mason and Wilson \cite{MR85c:54063}, is crucial in part 
of the proof of the Decomposition Theorem.

\begin{theorem}
If $f\colon X \to Y$ is a completely regular, monotone map between 
curves, then $f^{-1}(y)$ is a tree-like continuum for all $y$ in $Y$.
\end{theorem}

The second is due to Dyer and Hamstrom \cite{MR19:1187e}.

\begin{theorem}
Let $f\colon X \to Y$ be a completely regular mapping between compacta.
Let $f^{-1}(y)$ be homeomorphic to the compactum $M$, for all $y$ in $Y$.
Let $H(M)$ be the homeomorphism group of $M$.
Suppose 
$\dim Y \leq n+1$ and 
$H(M)$ is $LC^n$, and
$\Pi_1(H(M)) = 0$ for $0 \leq i \leq n$.
Then $X$ is homeomorphic to $Y \times M$.
\end{theorem}

The Dyer-Hamstrom result requires, for most applications, a well-behaved 
homeomorphism group $H(M)$.
If $Y$ is a Cantor set, however, then $n=-1$ and the last two conditions 
are vacuously satisfied.

An application of this is the following.
Call a compactum a \emph{Cantor set of pseudo-arcs} if it admits an open 
map into a Cantor set with pseudo-arcs as the fibers.
Then we have an alternate proof of a result of Wayne Lewis 
\cite{MR85f:54072}:
\emph{Each Cantor set of pseudo-arcs is a product of a Cantor set and a 
pseudo-arc.}

These ideas may have application again in continua theory.

Here is another application to the zero-dimensional case.
If $X$ is a homogeneous compactum, then the decomposition space $Y$ 
obtained by shrinking the components of $X$ to points is homogeneous and 
zero-dimensional, and the quotient map is completely regular (see, for 
instance, the proof of Theorem 3 of \cite{MR83a:54012}).
Hence we have the following result.

\begin{theorem}
Each homogeneous compactum $X$ is homeomorphic to $M \times Y$, where $M$ 
is one of the components of $X$, and $Y$ is a homogeneous zero-dimensional 
compactum.
In particular, $Y$ is either a finite set or a Cantor set.
\end{theorem}

Mislove and Rogers \cite[Theorem 2.4]{MR90m:54044a} or 
\cite{MR90m:54044b} have another technique that can be used to prove the 
theorem above.
Aarts and Oversteegen \cite{MR91m:54006} have generalized this theorem by 
replacing ``compactum'' by ``locally compact separable metric space'' in 
both hypothesis and conclusion.

A compact metric space $(X,d)$ is said to have the \emph{Effros property}
if, for every $\epsilon > 0$, there exists a $\delta > 0$ such that if 
$d(x,y) < \delta$, for two points $x$ and $y$ of $X$, then there is an 
$\epsilon$-homeomorphism $h$ from $X$ onto itself such that $h(x) = y$.
Zhou \cite{MR1397100} has used a decomposition technique to determine 
when a compactum with the Effros property must be homogeneous.

\section*{Hereditarily equivalent continua}

A continuum is \emph{hereditarily equivalent} if it is homeomorphic to 
each of its nondegenerate subcontinua.
In 1921, S.~Mazurkiewicz \cite{mazurkiewicz19213} asked if each 
finite-dimensional, hereditarily equivalent continuum is an arc.
In 1930, G.T.~Whyburn \cite{whyburn1930} proved that a planar, 
hereditarily equivalent continuum does not separate the plane.
Although the problem was posed as worthy of attention by Klein in 1928 and 
Wilder \cite{wilder1937} in 1937, no further progress occurred until 1948, 
when E.E.~Moise \cite{MR10:56i3} constructed a pseudo-arc.
The pseudo-arc is a hereditarily indecomposable, hereditarily equivalent 
continuum in the plane, and so the answer to Mazurkiewicz's question is no.

The arc and the pseudo-arc are the only known hereditarily equivalent, 
nondegenerate continua.
G.W.~Henderson \cite{MR22:9949.1} showed that any new example must be 
hereditarily indecomposable, and H.~Cook \cite{MR42:1072.2} showed that 
any new example must be tree-like.
Rogers \cite{MR86k:54054} observed that each continuum of dimension 
greater than one contains uncountably many topologically distinct 
subcontinua.

\begin{question}
Is every hereditarily equivalent, nondegenerate continuum chainable?
\end{question}

If the answer to this question is yes, then it is known that the arc and 
the pseudo-arc are the only such examples.

\begin{question}
Does each hereditarily equivalent continuum have span zero?
\end{question}

Oversteegen and Tymchatyn \cite{MR86a:54042} have shown that planar, 
hereditarily equivalent continua have symmetric span zero.

\begin{question}
Does each hereditarily equivalent continuum have the fixed-point property?
\end{question}

\begin{question}
Is each indecomposable, hereditarily equivalent continuum homogeneous?
\end{question}

\providecommand{\bysame}{\leavevmode\hbox to3em{\hrulefill}\thinspace}

\label{tprogersend}

\end{document}